\documentclass{article}
\usepackage{epsfig,epic,eepic,latexsym, amssymb}

\input xy
\xyoption{all}
\xyoption{2cell}
\UseAllTwocells
\xyoption{all}
\makeindex

\def \pf{PROOF:}
\def    \QED    {\hfill\hbox{\hskip 4pt
                \vrule width 5pt height 6pt depth 1.5pt}}
\def\epf{\QED\\}

\newcommand{\CAT}{ {\bf Cat}}

\newcommand{\Set}{ {\bf Set}}
\newcommand{\V}{ {\cal V} } 
\newcommand{\SMC}{SMC}
\newcommand{\SSMC}{StrSMC}
\newcommand{\A}{ {\cal A} }
\newcommand{\B}{ {\cal B} } 
\newcommand{\C}{ {\cal C} }
\newcommand{\D}{ {\cal D} }

\newcommand{\T}{ {\cal T} }

\newcommand{\bA}{ {\it A}}
\newcommand{\bB}{ {\it B}}
\newcommand{\bC}{ {\it C}}
\newcommand{\cH}{\check{H}}
\newcommand{\cotimes}{@}

\newcommand{\unc} { {\cal I} }
\newcommand{\termcat}{ {\bf 1}   }
\newcommand{\Diag}{ \bigtriangleup }
 
\newcommand{\Gra}{ {\cal H} }
\newcommand{\Grd}{ {\cal K}} 
\newcommand{\F}{ {\cal F} }
\newcommand{\Obj}{ Obj }
\newcommand{\Ver}{ { \cal O} }

\newcommand{\Ten} { Ten    }

\newcommand{\gen}{\star}

\newcommand{\ud}{\bullet}
\newcommand{\Fs}{ F^*}
\newcommand{\Gs}{ G^*}
\newcommand{\Fss}{ F^{**}}

\newcommand{\gF}{ \tilde{F}}
\newcommand{\gG}{ \tilde{G}}
\newcommand{\Fb} { F^{\bigtriangledown} }
\newcommand{\Gb} { G^{\bigtriangledown} }
\newcommand{\sigb} { \sigma^{\bigtriangledown} }

\newcommand{\ac}{ ass } 
\newcommand{\lc}{ l }
\newcommand{\rc}{ r }
\newcommand{\aci}{\overline{ass} } 
\newcommand{\lci}{\bar{l} }
\newcommand{\rci}{ \bar{r} }
\newcommand{\syc}{ s } 
\newcommand{\un}{I}

\newcommand{\rev}{rev}

\newcommand{\Res}{Rn}
\newcommand{\Ext}{En}
\newcommand{\homSMC}{\underline{Hom}}
\newcommand{\tenSMC}{\underline{Ten}}
\newcommand{\Pre}{Pre}
\newcommand{\Post}{Post}
\newcommand{\pre}{pre}
\newcommand{\post}{post}
\newcommand{\Dual}{D}
\newcommand{\ev}{ev}
\newcommand{\resev}{q}
\newcommand{\vv}{v}
\newcommand{\vs}{v^*}
\newcommand{\Eval}{Eval}

\newtheorem{theorem}{Theorem}[section]

\newtheorem{proposition}[theorem]{Proposition}
\newtheorem{lemma}[theorem]{Lemma}

\newtheorem{corollary}[theorem]{Corollary}
\newtheorem{remark}[theorem]{Remark}

\newtheorem{tag}[theorem]{}

\title{Tensor product for symmetric monoidal categories}
\author{Vincent Schmitt,\\
University of Leicester,\\
University road,\\
Leicester LE1 7RH, England.
}

\begin{document}
\maketitle                      

\begin{abstract}
We introduce a tensor product for symmetric monoidal categories
with the following properties.
Let $\SMC$ denote the 2-category with objects small symmetric monoidal 
categories, arrows {\em symmetric monoidal functors} and 2-cells
monoidal natural transformations. Our tensor product together with a
suitable unit is part of a structure on $\SMC$ that is a 2-categorical
version of the symmetric monoidal closed categories.
This structure is surprisingly simple. In particular
the arrows involved in the associativity and symmetry
laws for the tensor and in the unit cancellation laws are 2-natural
and satisfy coherence axioms which are {\em strictly commuting diagrams}.
We also show that the category quotient
of $\SMC$ by the congruence generated by 
its 2-cells admits a symmetric monoidal closed structure.
\end{abstract}

\begin{section}{Summary of results}
Thomason's famous result claims that 
symmetric monoidal categories model all connective spectra
\cite{Tho95}. The discovery of a symmetric monoidal structure
on the category of structured spectra \cite{EKMM97} suggests that 
a similar structure should exist on an adequate category with 
symmetric monoidal categories as objects.
The first aim of this work is to give a reasonable candidate
for a tensor product of symmetric monoidal categories.\\
 
We define such a tensor product for two symmetric monoidal categories by 
means of a generating graph and relations. It has the following properties.
Let $\SMC$ denote the 2-category with objects symmetric 
monoidal categories, with 1-cells {\em symmetric monoidal functors} and 
2-cells monoidal natural transformations. The tensor product yields
a {\em 2-functor} $\SMC \times \SMC \rightarrow \SMC$. This one is
part of a 2-categorical structure on $\SMC$ that is 2-categorical 
version of the symmetric monoidal closed categories.
Moreover this structure is rather simple since:\\
- its ``canonical'' arrows, i.e those 
involved for the associativity, the symmetry
and the left and right unit cancellation
laws, are {\em 2-natural};\\
- all coherence axioms for the above arrows
are {\em strictly commuting diagrams.}\\
Actually this last point was quite unexpected. 
Eventually from the above structure 
one can deduce a symmetric monoidal closed structure on the category 
$\SMC_{/\sim}$ quotient of $\SMC$ by the congruence $\sim$ generated 
by its 2-cells.\\ 

Here are now, in brief, the technical results in the 
order in which they occur in the paper.
The existence of an {\em internal hom}, a {\em tensor} and a 
{\em unit} for $\SMC$ and the fundamental properties defining and 
relating those are established first.  
From this, further properties of the above structure 
can be derived,  such as the existence of associativity, symmetry
and unit laws involving 2-natural
arrows satisfying coherence axioms. 
The proofs are rather computational for establishing 
a few key facts but become hopefully convincingly 
short and abstract for the rest of the paper.\\

Sections 2,3 and 4 are preliminaries.
A first point is that the 2-category $\SMC$ admits
internal homs in the following sense.
For any symmetric monoidal categories $\A$ and $\B$,
the mere category $\SMC(\A,\B)$ admits 
a symmetric monoidal structure
denoted $[\A,\B]$. 
This data extends 
to a 2-functor $\homSMC: \SMC^{op} \times \SMC \rightarrow \SMC$
(Proposition \ref{hom2fun}).
Moreover the classical isomorphism $\CAT(\A,\CAT(\B,C)) \cong
\CAT(\B,\CAT(\A,\C))$ induces a 2-natural isomorphism 
$\Dual: [\A,[\B,\C]] \cong [\B,[\A,\C]]$ between 2-functors
$\SMC^{op} \times \SMC^{op} \times \SMC \rightarrow \SMC$
(Proposition \ref{dualitynatABC}).
The isomorphism $\Dual$ is natural even in a stronger
sense as it is precisely stated in Lemma \ref{tbc1}
and \ref{tbc3}. 
The results above are presented in the sections from \ref{inthomSMC} 
to \ref{dualitynatsec}.
A brief section \ref{evalfunsec} treats 
the functors $\ev_a: [\A,\B] \rightarrow \B$ given by the evaluation
at some object $a$ of $\A$. 
They are the images of objects $a$ by the functor 
$\A \rightarrow [[\A,\B],\B]$ that
corresponds via $\Dual$ to the identity at $[\A,\B]$.\\

The tensor $\A \otimes \B$ of any symmetric monoidal categories
$\A$ and $\B$ is defined in section \ref{tensorABsec}.
In section \ref{P1}, two functors are defined
for any symmetric monoidal categories $\A$, $\B$ and $\C$, 
namely
$$\Ext:\SMC(\A, [\B, \C]) \rightarrow \SMC(\A \otimes \B, C)$$
and
$$\Res:\SMC(\A \otimes \B, \C) \rightarrow \SMC(\A, [\B, \C]).$$
They both admit symmetric monoidal structures yielding a monoidal 
adjunction
$$\Ext \dashv \Res : [\A \otimes \B,\C] \rightarrow [\A ,[\B,\C]]$$
with $\Res \circ \Ext =1$.  
It is also shown that the functors $\Ext_{\A,\B,\C}$ 
factor as 
$$
\xymatrix{ 
\SMC(\A,[\B,\C]) 
\ar[r]^-{\cong}
&
\SSMC(\A \otimes \B,\C) 
\ar[r] 
&
\SMC(\A \otimes \B,\C)
}
$$ 
where $\SSMC$ denotes the sub-2-category of $\SMC$
with the same objects but with {\em strict} functors as 1-cells
and monoidal natural transformations as 2-cells (Proposition
\ref{characExt} and Corollary \ref{coroadjten}). As explained 
in section \ref{tens2fun}, this 
universal characterisation of the tensor 
of symmetric monoidal categories serves 
to define its extension 
to a 2-functor $\SMC \times \SMC \rightarrow \SMC$.
Then the naturality issues for the collection of arrows $\Res$
and $\Ext$ are treated in section \ref{2natExtissue}.\\

Technical lemmas regarding the interaction between 
the internal hom, the tensor and the 
isomorphism $\Dual$ in $\SMC$ are grouped in section \ref{lasttecresults}.\\

The free symmetric monoidal category over the one 
point category, denoted $\unc$, is studied in section
\ref{free1}. By the universal property defining $\unc$,
for any symmetric monoidal category $\A$,
there is a unique strict symmetric monoidal functor 
$\vv: \unc \rightarrow [\A,\A]$
sending the generator of $\unc$ to the identity $1: \A \rightarrow
\A$. 
The symmetric monoidal functor $\vs: \A \rightarrow [\unc,\A]$ 
corresponding to $\vv$ via $\Dual$ will reveal useful 
and, as a 1-cell in $\SMC$, admits as right adjoint 
in $\SMC$ the functor $\ev_{\gen}:[\unc,\A] \rightarrow \A$
evaluation at the generator $\gen$ of $\unc$.\\

Then what can be seen as a ``lax symmetric monoidal 2-categorical
structure'' is defined on $\SMC$. The tensor and unit are known.
 
Arrows for the associativity, unit and symmetry laws
are introduced in section \ref{canarrows}:
$A'$ for associativity, 
$L'$, $R'$ for the left and right unit cancellations and $S'$ for the symmetry.
It is shown that the usual coherence axioms
for symmetric monoidal categories hold for 
$\SMC$ (see points \ref{waxiom12},
\ref{waxiom22}, \ref{waxiom3}, \ref{waxiom42} and \ref{SRpLp})
and that the canonical arrows are 2-natural
with ``lax'' inverses (points \ref{Ainverse2}, \ref{Ainverse}, 
\ref{Rinverse}).\\

Eventually last section \ref{closedquo} recaps everything for 
proving the symmetric monoidal closeness of the category
$\SMC / \sim$. This last point is Theorem \ref{mainresult}.\\ 

The following remark is in order.
A related result was found by M.Hyland and J.Power 
in their paper \cite{HyPo02} which treats in particular the 
2-category with objects symmetric monoidal categories, {\em but with 
strong functors as 1-cells}, and monoidal natural 
transformations as 2-cells. 
This was done though in a general 2-categorical setting extending
A.Kock's work on commutative monads. Hyland and Power discovered
a 2-categorical structure, namely that of {\em pseudo-monoidal closed 
2-category}, which exists on the 2-category of $T$-algebras with 
their pseudo morphisms for any {\em pseudo-commutative doctrine} $T$ 
on $\CAT$. 
Their pseudo monoidal closed structure is also a 2-categorical 
generalisation of the Eilenberg-Kelly's monoidal closed structure 
\cite{EiKe66}, with a tensor which is a {\em pseudo-functor},
with canonical morphisms which are {\em pseudo-natural equivalences} 
and coherence diagrams which commute {\em up to coherent isomorphic
2-cells}.
In particular the tensor for $T$-algebras 
in \cite{HyPo02} is defined from a collection 
of left biadjoints $\A \otimes -$ to internal 
homs $[\A,-]$. 
Actually this point does not seem to generalise
to the 2-categories of $T$-algebras with their lax morphisms. 
Eventually the tensors here and in \cite{HyPo02} are not isomorphic as 
symmetric monoidal categories due to the fact that
we consider lax morphisms of algebras whereas Hyland
and Power considered
the pseudo ones.

The author's view is that the case of symmetric monoidal categories 
together with monoidal functors should be elucidated before 
considering any generalisation in the line
of Kock's or Hyland-Power's works.\\

I wish to gratefully acknowledge Teimuraz Pirashvili 
first for suggesting the main result of this paper 
and also for fruitful conversations.
I also wish to thank John Power for his kind visit to Leicester 
in March 2007 and for giving then an extensive account
of his work with M.Hyland. 
\end{section}

\begin{section}{Preliminaries}
The purpose of this section is to recall
basic notions, introduce notations and give
references to some results used through the paper.\\ 

We shall write $\Set$ for the category of 
sets, and $\CAT$ for the 2-category of {\em small}
categories.\\ 

A {\em monoidal category} 
$(\A, \otimes, \un , \ac,\rc,\lc)$ 
consists of a category $\A$ together with:\\
- a functor $\otimes : \A \times \A \rightarrow \A$,
sometimes denoted for convenience by $\Ten$;\\
- an object $\un$ of $\A$, sometimes explicitly written $\un_{\A}$;\\ 
- natural
isomorphisms $\ac_{a,b,c}: a \otimes (b \otimes c) \rightarrow (a
\otimes b) \otimes c$, $\rc_a: a \otimes \un \rightarrow a$
and $\lc_a: \un \otimes a  \rightarrow a$
that satisfy the two coherence axioms
\ref{mcaxiom1} and \ref{mcaxiom2} below.\\
\begin{tag}\label{mcaxiom1}
The diagram in $\A$
$$
\xymatrix{
a \otimes (b \otimes (c \otimes d))
\ar[r]^{\ac}
\ar[d]_{1 \otimes \ac}
&
(a \otimes b) \otimes (c \otimes d)
\ar[r]^{\ac}
&
((a \otimes b) \otimes c) \otimes d
\\ 
a \otimes (( b \otimes c ) \otimes d)
\ar[rr]_{\ac}
&
&
(a \otimes (b \otimes c)) \otimes d
\ar[u]_{\ac \otimes 1}
}
$$
commutes for any $a$,$b$,$c$ and $d$.
\end{tag}
\begin{tag}\label{mcaxiom2} 
The diagram in $\A$  
$$\xymatrix{
a \otimes (\un \otimes b) 
\ar[rr]^{\ac}
\ar[rd]_{1 \otimes \lc} 
&
&
(a \otimes \un) \otimes b
\ar[ld]^{\rc \otimes 1}
\\
&
a \otimes b
&
}$$
commutes for any $a$ and $b$.
\end{tag}

A symmetric monoidal category $(\A,\otimes, \un, \ac,\lc,\rc,\syc)$
consists of a monoidal category $(\A, \otimes,\un,\ac,\rc,\lc)$
with a natural transformation $\syc_{a,b} : a \otimes b \cong b
\otimes a$ such that the following coherence axioms
\ref{smcaxiom3}, \ref{smcaxiom4} and \ref{smcaxiom5}
are satisfied.
\begin{tag}\label{smcaxiom3} 
The composite in $\A$
$$\xymatrix{
a \otimes b
\ar[r]^{\syc}
&
b \otimes a
\ar[r]^{\syc} 
&
a \otimes b}$$
is the identity at $a \otimes b$ for any $a$ and $b$.
\end{tag}
\begin{tag}\label{smcaxiom4}
The diagram in $\A$
$$
\xymatrix{
a \otimes (b \otimes c)
\ar[r]^{\ac}
\ar[d]_{1 \otimes \syc}
&
(a \otimes b) \otimes c
\ar[r]^{\syc}
&
c \otimes (a \otimes b)
\ar[d]^{\ac}
\\
a \otimes (c \otimes b)
\ar[r]_{\ac}
&
(a \otimes c) \otimes b
\ar[r]_{\syc \otimes 1}
&
(c \otimes a) \otimes b
}
$$
commutes for any $a$,$b$,$c$.
\end{tag}
\begin{tag}\label{smcaxiom5}
The diagram in $\A$
$$
\xymatrix{
a \otimes \un
\ar[rr]^{\syc}
\ar[rd]_{\rc}
&
&
\un \otimes a
\ar[ld]^{\lc}
\\
&
a
&
}
$$
commutes for any $a$.
\end{tag}
A symmetric monoidal category 
is said {\em strict} when its isomorphisms
$\ac$, $\lc$, $\rc$ and $\syc$ are identities.\\ 

For convenience, we adopt the view that all the monoidal 
categories considered further are by default {\em small}.\\ 

A {\em monoidal} functor between 
monoidal categories $\A$ and $\B$, 
consists of a triple
$(F,F^0,F^2)$ where:\\
- $F: \A \rightarrow \B$ is a functor;\\
- $F^0$ is an arrow $\un \rightarrow F(\un)$ in $\B$;\\
- $F^2$ is a natural transformation 
$F^2_{x,y}: Fx \otimes Fy \rightarrow F(x \otimes y)$ in $x$ and $y$;\\  
and that is subject to the axioms
\ref{mofun3}, \ref{mofun41},
\ref{mofun42} below.
\begin{tag}\label{mofun3}
For any objects $a$, $b$, $c$ of $\A$,
the diagram in $\B$
$$\xymatrix@C=4pc{Fa \otimes (Fb \otimes Fc) 
\ar[r]^{\ac_{Fa,Fb,Fc}} 
\ar[d]_{1 \otimes F^2_{b,c}}
&
(Fa \otimes Fb) \otimes Fc
\ar[d]^{ F^2_{a,b} \otimes 1 }
\\ 
Fa \otimes F(b \otimes c)
\ar[d]_{F^2_{a, b \otimes c}}
&
F(a \otimes b) \otimes Fc
\ar[d]^{F^2_{a \otimes b,c}}
\\
F(a \otimes (b \otimes c))
\ar[r]_{F(\ac_{a,b,c})}
&
F((a \otimes b) \otimes c)
}$$  
commutes.
\end{tag}
\begin{tag}\label{mofun41}
For any object $a$ in $\A$, the diagram in $\B$
$$\xymatrix{
Fa \otimes \un
\ar[r]^{\rc_{Fa}}
\ar[d]_{1 \otimes F^0}
&
Fa
\\
Fa \otimes F\un
\ar[r]_{F^2_{a,\un}}
&
F(a \otimes \un)
\ar[u]_{F\rc_a}
}$$
commutes.
\end{tag}
\begin{tag}\label{mofun42}
For any object $a$ in $\A$, the diagram in $\B$
$$\xymatrix{
\un \otimes Fa 
\ar[r]^{\lc}
\ar[d]_{F^0 \otimes 1} 
&
Fa
\\
F(\un) \otimes Fa
\ar[r]_{F^2_{\un,a}}
&
F(\un \otimes a)
\ar[u]_{F\lc}
}$$
commutes.
\end{tag}

A monoidal functor $(F,F^0,F^2)$ is {\em strict} when
the arrows $F^0$ and the $F^2$ are identities and it 
is {\em strong} when those are isomorphisms.
When the monoidal categories
$\A$ and $\B$ are symmetric, a monoidal
functor $F: \A \rightarrow \B$ is {\em symmetric} when
it satisfies the following axiom.
\begin{tag}\label{symofun5}
For any objects $a$, $b$ of $\A$, the diagram in $\B$
$$\xymatrix{
Fa \otimes Fb
\ar[r]^{F^2_{a,b}}
\ar[d]_{\syc_{Fa,Fb}}
&
F(a \otimes b)
\ar[d]^{F\syc_{a \otimes b}}
\\
Fb \otimes Fa
\ar[r]_{F^2_{b,a}}
&
F(b \otimes a)
}$$
commutes.
\end{tag}
   
Though symmetric monoidal functors
are triples we will sometimes just mention
``the symmetric monoidal $F$'' to mean
actually the triple $(F,F^0,F^2)$.\\ 

Given monoidal functors $F,G: \A \rightarrow \B$,
a natural transformation $\sigma: F \rightarrow G$ 
is {\em monoidal} when it satisfies axioms
\ref{monat6} and \ref{monat7} below.
\begin{tag}\label{monat6}
For any objects $a$, $b$ of $\A$, the diagram in $\B$
$$
\xymatrix{
Fa \otimes Fb \ar[r]^{F^2_{a,b}}
\ar[d]_{\sigma_a \otimes \sigma_b} 
&
F( a \otimes b)
\ar[d]^{\sigma_{a \otimes b}}
\\
Ga \otimes Gb 
\ar[r]_{G^2_{a,b}} 
&
G ( a \otimes b)
}
$$
commutes.
\end{tag}
\begin{tag}\label{monat7}
The diagram in $\B$
$$
\xymatrix{
\un 
\ar[r]^{F^0}
\ar[rd]_{G^0}
&
F(\un)
\ar[d]^{\sigma_{\un}}
\\
&
G(\un)
}
$$
commutes.
\end{tag}

Recall that symmetric monoidal categories, symmetric monoidal functors 
and monoidal transformations form a 2-category that we denote $\SMC$.
In particular given two symmetric monoidal functors $F: \A \rightarrow \B$
and $G: \B \rightarrow \C$,
the symmetric monoidal structure
of the composite $GF$ is the following:\\
- the arrow ${(GF)}^0$ is 
$
\xymatrix@C=2pc{
\un_{\C} 
\ar[r]^{G^0}
&
G(\un_{\B})
\ar[r]^{G(F^0)}
&
GF(\un_{\A})
}
$;\\
- for any objects $a$, $a'$ in $\A$,
the arrow ${(GF)}^2_{a,a'}$ is 
$$
\xymatrix@C=3pc{
GF(a) \otimes GF(a') 
\ar[r]^-{G^2_{Fa,Fa'}}
&
G(F(a) \otimes F(a'))
\ar[r]^{G(F^2_{a,a'})}
&
GF(a \otimes a').
}
$$
We write $\SSMC$ for the sub-2-category of $\SMC$,
with the same objects
but with 1-cells the {\em strict} symmetric monoidal functors
and monoidal transformations as 2-cells.\\

One has a forgetful 2-functor $U: \SMC \rightarrow \CAT$
sending symmetric monoidal categories, 
symmetric monoidal functors and monoidal transformations
respectively to their underlying categories, functors
and natural transformations.
There exists a notion of 
free symmetric monoidal category in the following sense:
the forgetful
2-functor $\SSMC \rightarrow \CAT$ admits a left 2-adjoint.\\

Remember from enriched category theory 
the following (\cite{Kel82} p.43).
\begin{lemma}\label{natcompfun}
For any symmetric monoidal closed complete and cocomplete
category $\V$
and any $\V$-functor $F: \A \rightarrow \B$
the collection of its components $F_{a,a'}: \A(a,a') \rightarrow
\B(Fa,Fa')$ in $\V$ is $\V$-natural in $a$ and $a'$.
\end{lemma}
As a consequence the collection of 
functors $U_{\A,\B}: \SMC(\A,\B) \rightarrow \CAT(U \A, U \B)$, 
components in $\A$ and $\B$ of the above forgetful 2-functor $U$, is   
2-natural in $\A$ and $\B$.\\ 

There is an important result regarding coherence 
for symmetrical monoidal categories, which is 
explained for instance in the second 
version of Mac Lane's handbook \cite{McLaCWM} in chapter 11. We will use 
implicitly this result, 
sometimes mentioning
the coherence for symmetric monoidal categories.
Also we will often omit to explicit the canonical isomorphisms
in symmetric monoidal categories when they are the expected 
ones.\\ 

G.M. Kelly characterised in \cite{Kel74} the 
{\em symmetric monoidal adjunctions} which are the adjunctions 
in the 2-category $\SMC$. In particular let us recall that
{\em a symmetric monoidal functor is left adjoint in $\SMC$ 
if and only if it is strong and left adjoint in $\CAT$.} 
Though monoidal adjunctions are a central 
concept in the present paper, we will use only the following 
result from Kelly.
\begin{proposition}\label{monadj}
Given any symmetric monoidal functor $(G, G^0,G^2): \B \rightarrow \A$
that admits a left adjoint $F$ in $\CAT$, 
$G$ admits a left adjoint in the 2-category $\SMC$
if and only if the following two conditions hold:\\ 
- $(1)$ For any objects $a,a'$ in $\A$, the arrows in $\B$
{\small
$$\xymatrix@C=3pc{
F( a \otimes a' ) 
\ar[r]^-{F ( \eta_{a} \otimes \eta_{a'} )}
&
F( GFa \otimes GFa')
\ar[r]^-{F ( G^2_{Fa,Fa'})   }
&
FG(Fa \otimes Fa' )
\ar[r]^-{\epsilon_{Fa \otimes Fa'}}
&
Fa \otimes Fa'
}$$
} are invertible;\\
- $(2)$ The arrow in $\B$
$$
\xymatrix{
F(\un) \ar[r]^-{F(G^0)} 
&  
FG(\un)
\ar[r]^-{\epsilon_{\un}}
&
\un}
$$ is invertible.\\
In this case, the functor $F$ admits a strong monoidal structure
with the $F^2_{a,a'}$ and $F^0$ respectively inverses of the 
above arrows $(1)$ and $(2)$, and $(F,F^0,F^2)$ is the left adjoint 
of $G$ in $\SMC$.
\end{proposition}

We will use the Yoneda machinery in the enriched context.
Let us consider a base symmetric monoidal closed category $\V$
that is complete and cocomplete.
(The base considered further in the paper is always the cartesian 
closed $\CAT$.)  
Recall that for any $\V$-functor $F: \A \rightarrow \V$ the 
Yoneda lemma defines a bijection 
between $\V$-natural transformations in the argument $a$, 
$\A(a,k) \rightarrow Fa$ and arrows $\un \rightarrow Fk$ in $\V$. 
We will use extensively the two following results from \cite{Kel82}.
\begin{lemma}
Given any $\V$-functors $F: \A^{op} \otimes \B \rightarrow \V$ and
$K: \A \rightarrow \B$, with a collection of arrows
$\phi_{a,b}: \A(a,Kb) \rightarrow_{a} F(a,b)$,
$\V$-natural in the argument $a$ for each object $b$ in $\B$,
this collection 
is also $\V$-natural in $b$
if and only if the collection of arrows 
$\un \rightarrow F(Kb,b)$
corresponding by Yoneda 
is $\V$-natural in $b$. 
\end{lemma}
\begin{lemma}[Representability with parameters]
For any $\V$-functor $F: \A^{op} \otimes \B \rightarrow \V$
with a collection of objects $K(b)$ and $\V$-natural
isomorphisms $\phi_b: \A(a,Kb) \cong_{a} F(a,b)$ in the argument $a$
for each object $b$ of $\B$,
there is a unique extension of $K$ into a $\V$-functor
$\B \rightarrow \A$ such that the collection of isomorphisms 
$\phi_b$ is also $\V$-natural in the argument $b$.
\end{lemma}

\end{section}

\begin{section}{Terminal object and product in $\SMC$}
\label{termprod}
The terminal category $\termcat$ (with one object
and one arrow) admits a unique symmetric monoidal structure,
which is the strict one.
Still denoting by $\termcat$ the latter object 
in $\SMC$, one has an isomorphism in $\CAT$,
2-natural in $\A$
\begin{tag}\label{2natisoter}
$$\SMC(\A, \termcat) \cong \termcat.$$
\end{tag}
In particular $\termcat$ is also the terminal
object of the underlying category of $\SMC$. 
For any symmetric monoidal category $\B$,
the constant functor 
$\Delta_{\un_{\B}}: \termcat \rightarrow \B$
to $\un_{\B}$
admits a strong symmetric monoidal structure
as follows. 
For the unique object $\gen$ of $\termcat$,
the arrow ${( \Delta_{\un_{\B}} )}^2_{\gen,\gen}$ 
is defined as the canonical one $r_{\un} = \l_{\un}: 
\un \otimes \un \rightarrow \un$ in $\B$
and the arrow
${( \Delta_{\un_{\B}} )}^0$ is defined as the identity 
at $\un_{\B}$.
So that for any symmetric monoidal categories
$\A$ and $\B$ the constant functor 
$\A \rightarrow \B$ to $\un_{\B}$ 
admits a symmetric strong monoidal structure
as the composite 
$\xymatrix{
\A 
\ar[r] 
& 
\termcat 
\ar[r]^{\Delta_{\un_{\B}}} 
& \B
}$
in $\SMC$.\\

The cartesian product of categories satisfies the following
property.
For any categories $\A$, $\B$ and $\C$, one has an isomorphism
of categories
\begin{tag}\label{Catproduct}
$$\CAT(\A,\B \times \C) \cong \CAT(\A, \B) \times \CAT(\A,\C)$$
\end{tag}
and for some arbitrary fixed $\B$ and $\C$ 
this collection of isomorphisms is 
2-natural in $\A$.
There is a unique 2-functor, say the {\em cartesian product}
2-functor, $- \times -: \CAT \times \CAT \rightarrow \CAT$
that makes the isomorphism \ref{Catproduct} also 
2-natural in $\B$ and $\C$. (This can be checked in an ad-hoc
way but this also results from the Yoneda Lemma in its enriched version
by considering 2-categories and 2-functors as categories and
functors enriched over the cartesian closed category $\CAT$.)\\

When the categories $\B$ and $\C$ have a symmetric monoidal structure,
their cartesian product admits the following  ``pointwise'' symmetric 
monoidal structure.\\
- The unit $\un_{\B \times \C} = (\un_{\B}, \un_{\C})$.
Equivalently one can define the constant functor
$\Delta_{\un_{\B \times \C}}: \termcat \rightarrow \B \times \C$
to $\un_{\B \times \C}$ as the composite in $\CAT$
$$\xymatrix@C=3pc{
\termcat
\ar[r]^-{\Diag}
&
\termcat \times \termcat
\ar[r]^-{ \Delta_{\un_B} \times \Delta_{\un_C} }
&
\B \times \C  
}$$
where $\Diag$ stands for the diagonal functor.
Note that the above composite
is the unique functor $F: \termcat \rightarrow \B \times \C$
such that $F \circ p_{\B} = \Delta_{\un_\B}$
and $F \circ p_{\C} = \Delta_{\un_{\C}}$, where $p_{\B}$ 
and $p_\C$ are respectively the projections $\B \times \C \rightarrow \B$
and $\B \times \C \rightarrow \C$.\\
- The tensor $\Ten: (\B \times \C) \times (\B \times \C)
 \rightarrow (\B \times \C)$ is 
the composite functor
$$
\xymatrix@C=3pc{
(\B \times \C) \times (\B \times \C)  
\ar[r]^-{\cong} 
&
(\B \times \B) \times (\C \times \C)
\ar[r]^-{ \Ten \times \Ten } 
&
\B \times \C
}$$
which is also the unique functor 
$F: (\B \times \C) \times (\B \times \C) \rightarrow \B \times \C$
such that $p_{\B} \circ F = Ten \circ (p_{\B} \times p_{\B})$
and $p_{\C} \circ F = Ten \circ (p_{\C} \times p_{\C}$).
This is to say that
for any objects $b,b'$ in $\B$ 
and $c,c'$ in $\C$,
the object $(b,c) \otimes (b',c')$ is $(b \otimes b', c \otimes c')$ and
for any arrows $f:b_1 \rightarrow b_2$, 
$f':b_1' \rightarrow b_2'$ in $\B$
and $g:c_1 \rightarrow c_2$,
$g':c_1' \rightarrow c_2'$ in $\C$,
the arrow $(f,g) \otimes (f',g'): (b_1,c_1) \otimes (b_1',c_1')
\rightarrow (b_2,c_2) \otimes (b_2',c_2')$ is 
$(f \otimes f', g \otimes g'):(b_1 \otimes b_1', c_1 \otimes c_1') \rightarrow 
(b_2 \otimes b_2', c_2 \otimes c_2')$.\\
- For any objects $b,b',b''$ in $\B$
and $c,c',c''$ in $\C$, the arrow
$\ac_{(b,c), (b',c'), (b'',c'')}$ is 
$$\xymatrix{
(b,c) \otimes ( (b',c') \otimes (b'',c''))
\ar@{=}[d] 
\\
( (b \otimes (b' \otimes b''), c \otimes (c' \otimes c'') )
\ar[d]^{( \ac_{b,b',b''} , \ac_{c,c',c''}  ) }
\\
( ( b \otimes b' ) \otimes b'', (c \otimes c') \otimes c'' )
\ar@{=}[d]  
\\
( (b,c) \otimes (b',c') ) \otimes (b'', c'').
}
$$
- For any objects $b$ in $\B$ and $c$ in $\C$,
the arrow 
$\rc_{(b,c)}$ is defined as
$$\xymatrix{
(b,c) \otimes (\un_{\B}, \un_{\C})
\ar@{=}[r]
&
(b \otimes \un_{\B}, c \otimes \un_{\C})
\ar[r]^-{(\rc_{\B}, \rc_{\C})}
&
(b,c)
}$$
and 
$\lc_{(b,c)}$ is
$$\xymatrix{
(\un_{\B}, \un_{\C}) \otimes (b,c)
\ar@{=}[r]
&
(\un_{\B} \otimes b, \un_{\C} \otimes c)
\ar[r]^-{(\lc_{\B}, \lc_{\C})}
&
(b,c)
}.$$
- For any objects $b,b'$ in $\B$ and $c,c'$ in $\C$,
the arrow 
$\syc_{(b,c), (b',c')}$ is defined as
$$\xymatrix@C=3pc{ (b,c) \otimes (b',c')
\ar@{=}[r]
&
(b \otimes b', c \otimes c')
\ar[r]^-{( \syc_{b,b'}, \syc_{c,c'} )}
& 
(b' \otimes b, c' \otimes c)
\ar@{=}[r]
&
(b',c') \otimes (b,c)
.}
$$
This symmetric monoidal structure, which we still write 
$\B \times \C$, is also 
the unique one on the mere category $\B \times \C$ 
that makes the two projections 
$p_{\B}: \B \times \C \rightarrow \B$
and
$p_{\C}: \B \times \C \rightarrow \C$
into symmetric strict monoidal functors.\\

Given any symmetric monoidal functors
$(F,F^0,F^2): \A \rightarrow \B$ and 
$(G,G^0,G^2): \A \rightarrow \C$,
writing $<F,G>$ for the corresponding 
functor $\A \rightarrow \B \times \C$ via 
the isomorphism \ref{Catproduct}, this one admits 
a symmetric monoidal structure as follows:\\
- Given objects $a$,$a'$ in $\A$, 
${<F,G>}^2_{a,a'}$ is the arrow of $\B \times \C$
$$\xymatrix{
<F,G>(a) \otimes <F,G>(a')
\ar@{=}[d] 
\\
(F(a),G(a)) \otimes (F(a'),G(a'))
\ar@{=}[d]
\\
(F(a) \otimes F(a'),G(a) \otimes G(a'))
\ar[d]^{(F^2_{a,a'},G^2_{a,a'})}
\\
(F(a \otimes a'), G(a \otimes a'))
\ar@{=}[d]
\\
<F,G>(a \otimes a');
}
$$  
- ${<F,G>}^0$ is the arrow
$$
\xymatrix{
\un 
\ar@{=}[r] 
&
(\un_{\B}, \un_{\C})
\ar[r]^-{(F^0, G^0)}
&
(F(\un_{\A}), G(\un_{\A}))
\ar@{=}[r]
&
<F,G>(\un_{\A}). 
}
$$  
This symmetric monoidal structure is the unique one on 
the functor $<F,G>$ such that the two equalities
$$p_{\B} \circ (<F,G>, {<F,G>}^0, {<F,G>}^2) = (F,F^0,F^2)$$ 
and 
$$p_{\C} \circ (<F,G>, {<F,G>}^0, {<F,G>}^2) = (G,G^0,G^2)$$ 
hold in $\SMC$. When talking 
of $<F,G>$ as a symmetric monoidal functor, we will always
consider this structure. Note that it is strict 
if and only both $(F,F^0,F^2)$ and $(G,G^0,G^2)$ are.\\
 
Given a pair of monoidal natural transformations between 
symmetric functors
$\sigma: F \rightarrow F' : \A \rightarrow \B$
and
$\tau: G \rightarrow G': \A \rightarrow \C$,
the natural transformation 
$<\sigma,\tau>: <F,G> \rightarrow <F',G'>: \A \rightarrow \B \times \C$ 
that corresponds to this pair by \ref{Catproduct} is also monoidal.\\

The remarks above  are enough to show that for any symmetric monoidal
categories $\B$ and $\C$, the isomorphism \ref{Catproduct} 
induces a 2-natural isomorphism between 2-functors
\begin{tag}\label{2natisopro}
$$\SMC(-,\B \times \C) \cong \SMC(-,\B) \times \SMC(-,\C):
\SMC \rightarrow \CAT.$$
\end{tag}

This isomorphism in $\CAT$ is also 2-natural 
in $\B$ and $\C$ for a unique 2-functor 
$- \times -: \SMC \times \SMC \rightarrow \SMC$.
This 2-functor acts on 1-cells as follows.
Given any symmetric monoidal functors $(F,F^0,F^2): \B \rightarrow \B'$
and any $(G,G^0,G^2): \C \rightarrow \C'$, the monoidal 
structure on $F \times G: \B \times \C \rightarrow \B' \times \C'$,
is the following:\\ 
- The arrow ${(F \times G)}^0$ is
{\small
$$\xymatrix{
\un_{\B \times \C }
\ar@{=}[r]
&
(\un_{\B},\un_{\C})
\ar[r]^-{(F^0,G^0)}
&
(F(\un_{\B'}), G(\un_{\C'}))
\ar@{=}[r]
&
(F \times G )(\un_{\B '},\un_{\C '})
\ar@{=}[r]
&
(F \times G) (\un_{\B ' \times \C'})
};
$$
}
- For any objects $b$,$b'$ in $\B$ and $c$,$c'$ in $\C$,
the arrow ${(F \times G)}^2_{(b,c),(b',c')}$ is 
{\small
$$
\xymatrix{
(F \times G) (b,c) \otimes (F \times G)(b',c')
\ar@{=}[d]
\\
(Fb \otimes Fb', Gc \otimes Gc')
\ar[d]^{(F^2_{b,b'}, G^2_{c,c'})}
\\
( F( b \otimes b') , G(c \otimes c') )
\ar@{=}[d]
\\
(F \times G )( (b,c) \otimes (b',c') ).
}
$$
}
Note that for any symmetric monoidal category $\A$
the diagonal functor $\A \rightarrow \A \times \A$
admits a strict symmetric monoidal structure by the above results.\\

\end{section}

\begin{section}{Tensor as a symmetric monoidal functor}
For any symmetric monoidal category $\A$, 
its tensor $\Ten: \A \times \A \rightarrow \A$
admits the following symmetric monoidal 
structure.\\
- For any objects $a$,$b$,$c$ and $d$ in $\A$,
the arrow  $\Ten^2_{(a,b),(c,d)}$ is the canonical 
arrow commuting the $b$ and $c$
$$\xymatrix{
\Ten(a,b) \otimes \Ten(c,d)
\ar@{=}[d]\\  
(a \otimes b) \otimes (c \otimes d)
\ar[d]\\ 
(a \otimes c) \otimes (b \otimes d)
\ar@{=}[d]\\
\Ten ( a \otimes c, b \otimes d)
\ar@{=}[d]\\
\Ten( (a,b) \otimes (c,d)).
}$$
- The arrow 
$Ten^0: \un_{\A} \rightarrow \un_{\A} \otimes \un_{\A}$
is also the canonical one.\\

The collection of arrows $Ten^2_{(a,b),(c,d)}$ is natural 
in $(a,b)$ and $(c,d)$ and 
Axioms \ref{mofun3}, \ref{mofun41}, \ref{mofun42},
and \ref{symofun5} for the functor $\Ten$ all amount to commutations 
of canonical diagrams in $\A$.\\

Observe then that the natural isomorphisms
part of the monoidal structure of $\A$, namely
$\ac$, $\rc$, $\lc$ and $\syc$ are {\em monoidal},
precisely they are the respective following 2-cells in $\SMC$
$$
\xymatrix{
\A \times (\A \times \A)
\ar[rr]^{\cong}
\ar[d]_{1 \times \Ten}
& &
(\A \times \A) \times \A
\ar[d]^{\Ten \times 1}
\\
\A \times \A
\ar[rd]_{\Ten}
\ar@{=>}[rru]^{\ac}
& & 
\A \times \A
\ar[ld]^{\Ten}
\\
& 
\A
& 
}
$$
\begin{center}
$
\xymatrix{
\A \times \termcat
\ar[d]_{1 \times \Delta_{\un}}
&
\A
\ar@{=}[d]
\ar[l]_-{\cong}
\\
\A \times \A
\ar[r]_{\Ten}
\ar@{=>}[ru]^{\rc}
&
\A
}
$
\hspace{2cm}
$
\xymatrix{
\termcat \times \A
\ar[d]_{\Delta_{\un} \times 1}
&
\A
\ar@{=}[d]
\ar[l]_-{\cong}
\\
\A \times \A
\ar[r]_{\Ten}
\ar@{=>}[ru]^{\lc}
&
\A
}
$
\end{center}
$$
\xymatrix{
\A \times \A 
\ar[r]^{\Ten}
\ar[d]_{\rev}
&
\A
\ar@{=>}[ld]_{\syc}
\ar@{=}[d]
\\
\A \times \A
\ar[r]_{\Ten}
&
\A
}
$$
where $\rev$ stands for the automorphism of $\A \times \A$
in $\SMC$ 
such that $p_1 \circ \rev = p_2$ and $p_2 \circ \rev = p_1$
for the two projections $p_1, p_2: \A \otimes \A \rightarrow \A$.   
All this results straightforwardly again from the coherence 
in $\A$.

\end{section}

\begin{section}{Internal homs in $\SMC$}
\label{inthomSMC}
Given arbitrary symmetric monoidal categories $\A$ and $\B$,
the category $\SMC(\A,\B)$ of
symmetric monoidal functors $\A \rightarrow \B$
and monoidal transformations between them admits 
a symmetric monoidal structure denoted $[\A,\B]$
that is described below.\\

The unit $\un_{[\A,\B]}$ is given by the constant functor to
the unit $\un_{\B}$ of $\B$ with the (strong) monoidal structure as 
described previously.\\
 
For any symmetric monoidal functors $F,G: \A \rightarrow \B$, 
their tensor
$F \Box G$  is the composite in $\SMC$
$$\xymatrix{ \A \ar[r]^-{\Diag} & 
\A \times \A \ar[r]^-{F \times G} & 
\B \times \B \ar[r]^-{\Ten} & 
\B}$$ where $\Diag$ is the diagonal
functor of $\A$ and $\Ten$ the tensor in $\B$.
Therefore ${(F \Box G)}^0$ is the arrow 
$$\xymatrix{ \un_{\B} \ar@{-}[r]^-{\cong} & 
\un_{\B} \otimes \un_{\B} \ar[r]^-{F^0 \otimes G^0} & 
F(\un_{\A}) \otimes G(\un_{\A}) 
\ar@{=}[r] & F \Box G(\un_{\A}) }$$
and
${(F \Box G)}^2$ has component in any pair $(a,b)$ the arrow
$$\xymatrix{F \Box G(a)  \otimes F \Box G(b) \ar@{=}[d]\\ 
   (F(a) \otimes G(a)) 
\otimes (F(b) \otimes G(b)) \ar@{-}[d]^{\cong}\\
(F(a) \otimes F(b)) \otimes (G(a) \otimes G(b))
\ar[d]^{F^2_{a,b} \otimes G^2_{a,b}}\\
F(a \otimes b) \otimes G(a \otimes b) \ar@{=}[d]\\
F \Box G ( a \otimes b ).}$$

One defines for any two 2-cells 
$\sigma: F \rightarrow G : \A \rightarrow \B$
and $\tau: F'\rightarrow G': \A \rightarrow \B$ in $\SMC$
the 2-cell 
$\sigma \Box \tau : F \Box G \rightarrow F' \Box G'$
in $\SMC$ as
$$\xymatrix{
&&
\ar@{=>}[dd]^{\sigma \times \tau}
&
\\
\A
\ar[r]^-{\Diag}
&
\A \times \A
\ar@/^35pt/[rr]^-{F \times G}
\ar@/_35pt/[rr]_-{F' \times G'}
& &
\B \times \B
\ar[r]^-{\Ten}
& 
\B,
\\
&&&
}$$
which means that for any object $a$ of $\A$,
${(\sigma \Box \tau)}_a = \sigma_a \otimes \tau_a$.\\

That the assignments $\Box$ define a bifunctor 
$\SMC(\A,\B) \times \SMC(\A,\B) \rightarrow \SMC(\A,\B)$
results from the fact that the cartesian product defines a 2-functor 
$\SMC \times \SMC \rightarrow \SMC$.\\

For any symmetric monoidal $F,G,H: \A \rightarrow \B$,
the 2-cell $\ac_{F,G,H}$ in $\SMC$ is
$$\xymatrix{
& \A 
\ar[d]^{\Diag}
\\
&
 \A \times \A
\ar[ld]_{1 \times \Diag}
\ar[rd]^{\Diag \times 1}
\ar@{}[d]|{=}
\\
 \A \times (\A \times \A)
\ar[d]_{F \times (G \times H)}
\ar[rr]^{\cong}
&
\ar@{}[d]|{=}
&
(\A \times \A) \times \A
\ar[d]^{(F \times G) \times H}
\\
\B \times (\B \times \B)
\ar[d]_{ 1 \times \Ten}
\ar[rr]^{\cong}
& 
&
(\B \times \B) \times \B
\ar[d]^{ \Ten \times 1}
\\
\B \times \B 
\ar[dr]_{ \Ten }
\ar@{=>}[rru]_{\ac_{\B}}
 &  & 
\B \times \B
\ar[dl]^{ \Ten}
\\
& \B.
}$$

For any symmetric monoidal $F: \A \rightarrow \B$,
the 2-cell $\rc_F$ in $\SMC$ is
$$
\xymatrix{
\A
\ar@{}[rd]|{=}
\ar[r]^F
\ar[d]_{\Diag}
&
\B
\ar@{=}[r]
&
\B
\\
\A \times \A
\ar[r]_{ F \times !}
&
\B \times \termcat
\ar[u]|{\cong}
\ar@{=>}[ru]_{\rc}
\ar[r]_{ 1 \times \Delta_{\un_{\B}} }
&
\B \times \B
\ar[u]_{\Ten}
}
$$
and
$\lc_F$ is 
$$
\xymatrix{
\A
\ar@{}[rd]|{=}
\ar[r]^F
\ar[d]_{\Diag}
&
\B
\ar@{=}[r]
& 
\B
\\
\A \times \A
\ar[r]_{ ! \times F}
&
\termcat \times \B
\ar[u]|{\cong}
\ar@{=>}[ru]_{\lc}
\ar[r]_{ \Delta_{\un_{\B}} \times 1 }
&
\B \times \B.
\ar[u]_{\Ten}
}
$$

For any symmetric monoidal $F,G: \A \rightarrow \B$,
the 2-cell $\syc_{F,G}: F \Box G \rightarrow G \Box F$ 
in $\SMC$ is
$$\xymatrix{
\A
\ar@{}[rd]|{=}
\ar[r]^{\Diag}
\ar@{=}[d]
& 
\A \times \A
 \ar@{}[rd]|{=}
\ar[r]^{F \times G}
\ar[d]|{\rev}
&
\B \times \B
 \ar[d]|{\rev}
 \ar[r]^{\Ten}
&
\B 
 \ar@{=}[d]
 \ar@{=>}[ld]_{\syc}
\\
\A
 \ar[r]_{\Diag}
& 
\A \times \A
 \ar[r]_{G \times F}
&
\B \times \B
 \ar[r]_{\Ten}
&
\B.
}$$

The previous definitions amount to say that the 
canonical isomorphisms $\ac$, $\rc$, $\lc$, $\syc$
for $[\A,\B]$ are ``pointwise''.
That is, for any object $a$ in $\A$, the isomorphisms:\\
- $\ac_{F,G,H}: 
F \Box ( G \Box H ) \cong  (F \Box  G) \Box H$\\
- $\rc_F: F \Box \un \cong  F$\\ 
- $\lc_F: \un \Box F \cong F$\\
and\\
- $\syc_{F,G}: F \Box G  \cong  G \Box F$\\
are respectively defined at any object $a$ by:\\ 
- $\ac_{F(a),G(a),H(a)}: Fa \otimes (Ga \otimes Ha)
\rightarrow (Fa \otimes Ga) \otimes Ha$\\
- $\rc_{F(a)}:  F(a) \otimes \un \rightarrow F(a)$\\
- $\lc_{F(a)}:  \un \otimes F(a) \rightarrow F(a)$\\
and\\ 
- $\syc_{F(a),G(a)} : F(a) \otimes G(a) \rightarrow G(a) \otimes F(a).$\\

From this observation and the fact that $\B$ is
symmetric monoidal, it follows easily that 
the collections $\ac_{F,G,H}$, 
$\rc_F$, $\lc_F$ and $\syc_{F,G}$ are natural
respectively in $F$,$G$,$H$, in $F$, in $F$ and in $F$,$G$
and that these together
satisfy
Axioms \ref{mcaxiom1}, \ref{mcaxiom2},
\ref{smcaxiom3}, \ref{smcaxiom4} and \ref{smcaxiom5}
with $\Box$ as tensor and unit $\un_{[\A,\B]}$ as specified.

\end{section}

\begin{section}{The isomorphism $\Dual: [\A,[\B,\C]] \rightarrow
[\B,[\A,\C]]$ }
\label{dualitysec}
Given any symmetric monoidal categories
$\A$, $\B$ and $\C$,
we are going to build up from the classical isomorphism 
\begin{tag}\label{classicdual}
$$\CAT(\A, \CAT(\B,\C)) \cong \CAT(\B, \CAT(\A,\C))$$
\end{tag}
an isomorphism of categories
\begin{tag}\label{dualhom}
$$\SMC(\A, [\B,\C]) \cong \SMC(\B, [\A,\C])$$
\end{tag}
which turns out to have a strict symmetric monoidal structure
$$\Dual_{\A,\B,\C}: [\A,[\B,\C]] \rightarrow [\B,[\A,\C]]$$
with (symmetric) monoidal inverse $\Dual_{\B,\A,\C}$. Most of the 
time we will omit the subscripts for $\Dual$.\\

To explain more precisely this result, we need to 
introduce the following terminology:\\
- We say that two functors or two natural transformations
corresponding via the isomorphism \ref{classicdual} are {\em dual}
(to each other);\\  
- We write $U$ (rather than $U_{\A,\B}$) for the forgetful
functor $\SMC(\A,\B) \rightarrow
\CAT(\A,\B)$;\\ 
- For any symmetric monoidal functor
$(F, F^0,F^2): \A \rightarrow [\B,\C]$, we let 
$F^{\ud}$ denote the dual functor 
$\B \rightarrow \CAT(\A,\C)$ of the mere functor 
$$\xymatrix{
\A 
\ar[r]^-F 
& 
\SMC(\B,\C) 
\ar[r]^-U 
& 
\CAT(\B,\C);}
$$
- For any monoidal natural transformation 
$\theta: F \rightarrow G: \A \rightarrow [\B,\C]$ between 
symmetric functors, we let         
$\theta^{\ud}$ denote the natural transformation
$F^{\ud} \rightarrow G^{\ud}: \B \rightarrow \CAT(\A,\C)$
dual of $$U * \theta: U \circ F \rightarrow U \circ 
G: \A \rightarrow \CAT(\B,\C).$$

Then we are going to exhibit 
for any symmetric monoidal functor
$F: \A \rightarrow [\B,\C]$, 
a factorisation of $F^{\ud}$ as
\begin{tag}\label{factofun}
$$\xymatrix{ 
\B
\ar[r]^-{\Fs} 
&
\SMC(\A,\C)
\ar[r]^-U
&
\CAT(\A,C)
}$$ 
\end{tag}
and define on $\Fs$ a symmetric 
monoidal structure $(\Fs, {\Fs}^0, {\Fs}^2) :\B \rightarrow [\A,\C]$.
Also for any monoidal natural transformation between symmetric functors
$\theta: (F,F^0,F^2) \rightarrow (G,G^0,G^2):\A \rightarrow [\B,\C]$,
we are going to exhibit a factorisation of 
$\theta^{\ud}: F^{\ud} \rightarrow G^{\ud}$
as
\begin{tag}\label{factonat}
$$\xymatrix{
& 
\ar@<2.5ex>@{=>}[dd]_{\theta^*} 
&
\\
\B 
\ar@/^35pt/[rr]^{\Fs}
\ar@/_35pt/[rr]_{\Gs}
&
& 
\SMC(\A,\C) 
\ar[r]^U 
&
\CAT(\A,\C)
\\
& & 
}$$
\end{tag}
such that the natural $\theta^*$ is monoidal 
$(\Fs,{\Fs}^0,{\Fs}^2) \rightarrow (\Gs,{\Gs}^0,{\Gs}^2): \B
\rightarrow [\A,\C]$.
Furthermore any symmetric monoidal functor
$(F, F^0, F^2): \A \rightarrow [\B,\C]$ is 
equal to $(\Fss, {(\Fss)}^0, (\Fss)^2)$. 
Eventually we define the isomorphism
\ref{dualhom} as the functor sending any $F$ to $F^*$ and 
any $\theta$ to $\theta^*$.\\

Let us consider any symmetric monoidal functor 
$F: \A \rightarrow [\B,\C]$ with corresponding 
$F^{\ud}: \B \rightarrow \CAT(\A,\C)$.
For any $b$, we define the symmetric monoidal structure
$( \Fs b, {(\Fs b)}^0, {(\Fs b)}^2)$ 
on $F^{\ud} b: \A \rightarrow \C$ (i.e. $\Fs b = F^{\ud}b$), 
as follows:
\begin{tag}\label{dualstr1}$\;$\\
- $(\Fs b)^0 : \un_{\C} \rightarrow \Fs (b)(\un_{\A})$
is the arrow
${(F^0)}_b: \un_{\C} \rightarrow F(\un_{\A})(b)$,
component in $b$ of the (monoidal) natural transformation
$F^0: \un \rightarrow F(\un_{\A}): \B \rightarrow \C$.\\
- For any objects $a$ and $a'$ in $\A$,
${( \Fs b)}^2_{a,a'}: \Fs b(a) \otimes \Fs b(a') \rightarrow 
\Fs b (a \otimes a')$ is the arrow
$${({(F^2)}_{a,a'})}_b: Fa(b) \otimes Fa'(b) 
\rightarrow F(a \otimes a')(b),$$ component in $b$ 
of the natural transformation 
$F^2_{a,a'}: Fa \Box Fa' \rightarrow F(a \otimes a'): \B \rightarrow
\C$.
\end{tag}

Axiom \ref{mofun3} for $(\Fs b, {(\Fs b)}^0, {(\Fs b)}^2)$ is that 
for any objects $a,a',a''$ in $\A$,
the diagram in $\C$
$$
\xymatrix{
Fa (b) \otimes 
( F(a')(b) \otimes F(a'')(b) )
\ar[r]^{\ac}
\ar[d]_{1 \otimes {{( F^2_{a',a''})}_b   }}
&
( Fa (b)  \otimes 
F(a') (b) ) \otimes F(a'')(b)
\ar[d]^{ {( F^2_{a,a'} ) }_b  \otimes 1 } 
\\
Fa(b) \otimes 
( F(a')(b) \otimes F(a'')(b) )
\ar[d]_{  {(F^2_{a, a' \otimes a''})}_b }
& 
(Fa(b)  \otimes F(a')(b)) \otimes 
F(a'')(b) 
\ar[d]^{  {(F^2_{a \otimes a', a''})}_b   } 
\\
F( a \otimes (a' \otimes a''))(b)
\ar[r]_{{  F(\ac )}_b }
&
F( (a \otimes a') \otimes a'' )(b)
}
$$
commutes.
This diagram is the evaluation in $b$ of the 
diagram in $[\B,\C]$
$$
\xymatrix{
Fa  \Box ( F(a') \Box F(a'') )
\ar[r]^{\ac}
\ar[d]_{ 1 \Box F^2_{a',a''} }
&
( Fa   \Box  F(a')  ) \Box F(a'')
\ar[d]^{ F^2_{a,a'}  \Box 1 } 
\\
Fa \Box ( F(a') \Box F(a'') )
\ar[d]_{  F^2_{a, a' \otimes a''} }
& 
(Fa  \Box F(a')) \Box F(a'')
\ar[d]^{ F^2_{a \otimes a', a''}  } 
\\
F( a \otimes (a' \otimes a''))
\ar[r]_{ F(\ac)  }
&
F( (a \otimes a') \otimes a'' )
}
$$
which commutes according to Axiom 
\ref{mofun3} for the triple $(F,F^0,F^2)$.\\

One can check in a similar way that Axioms
\ref{mofun41}, \ref{mofun42} and \ref{symofun5} for 
the triple $(\Fs b, {(\Fs b)}^0, {(\Fs b)}^2)$
are the pointwise versions 
respectively of Axioms
\ref{mofun41}, \ref{mofun42} and \ref{symofun5} for 
the triple $(F,F^0,F^2)$.\\

Similarly the naturality in $a$,$a'$ of the collection 
of arrows ${(\Fs b)}^2_{a,a'}$ can be deduced as a pointwise version
of the naturality in $a$, $a'$ of the collection 
of arrows $F^2_{a,a'}: Fa \Box Fa' \rightarrow F(a \otimes a')$
in $[\B,\C]$.\\ 


Given any arrow $f: b \rightarrow b'$ in $\B$,
we show now that the natural transformation 
$F^{\ud}f : F^{\ud}(b) \rightarrow F^{\ud}(b'): \A \rightarrow \C$
is monoidal $\Fs b \rightarrow \Fs (b')$.
We define $\Fs f$ as this last arrow of
$\SMC(\A,\C)$.\\
 
Axiom \ref{monat6} for $F^{\ud}f$ amounts to the 
commutation for any objects $a$ and $a'$ in $\A$ 
of the diagram
$$
\xymatrix{
Fa(b) \otimes F(a')(b)
\ar[r]^-{{(F^2_{a,a'})}_b}
\ar[d]_{Fa(f) \otimes F(a')(f)}
&
F(a \otimes a')(b)
\ar[d]^{F(a \otimes a')(f)}
\\
Fa(b') \otimes F(a')(b')
\ar[r]_-{{(F^2_{a,a'})}_{b'}}
&
F(a \otimes a')(b')
}
$$ in $\C$,
which does commute according to the naturality 
of $F^2_{a,a'}: Fa \Box Fa' \rightarrow F(a \otimes a')
: \B \rightarrow \C$.\\

Axiom \ref{monat7} for $F^{\ud}f$ amounts to the commutation
of 
$$
\xymatrix{
& 
F(\un_{\A})(b)
\ar[dd]^{F(\un_{\A})(f)}
\\
\un 
\ar[ru]^{{F^0}_b}
\ar[rd]_{{F^0}_{b'}}
&
\\
& 
F(\un_{\A})(b')
}
$$ in $\C$,
which holds by naturality of 
$F^0: \un \rightarrow F(\un_{\A}): \B \rightarrow \C$.\\

So far we have checked the factorisation of 
mere functors $F^{\ud} = U \circ \Fs$ as expected in \ref{factofun}. 
The monoidal structure ${(\Fs)}^0$ and ${(\Fs)}^2$
on $\Fs$ is the following.
\begin{tag}\label{dualstr2}$\;$\\
- The monoidal natural transformation
${(\Fs)}^0: \un \rightarrow \Fs( \un_{\B}): 
\A \rightarrow \C$ has component in any $a$ in $\A$, 
the arrow in $\C$ $$(Fa)^0 : \un \rightarrow F(a)(\un_{\B}).$$ 
- For any objects $b$,$b'$ in $\B$, the monoidal
natural transformation 
$${(\Fs)}^2_{b,b'}: \Fs(b) \Box \Fs(b') \rightarrow \Fs ( b \otimes b')
: \A \rightarrow \C$$ has component in any $a$ in $\A$,
the arrow in $\C$
$${(Fa)}^2_{b,b'} : Fa(b) \otimes Fa(b') \rightarrow Fa(b \otimes
b').$$
\end{tag}

The naturality $\un \rightarrow \Fs(\un_{\B}): \A \rightarrow
\C$ of ${(\Fs)}^0$ is equivalent 
to the commutation for any arrow 
$f: a \rightarrow a'$ of $\A$ of the diagram in $\C$
$$
\xymatrix{
\un
\ar@{=}[d]
\ar[r]^-{{(Fa)}^0 }
&
Fa(\un_{\B})
\ar[d]^{ {F(f)}_{\un_{\B}} }
\\
\un
\ar[r]_-{{(Fa')}^0 }
&
F(a')(\un_{\B})
}
$$
which does commute according to Axiom \ref{monat7}
for the monoidal $F(f): Fa \rightarrow Fa': \B \rightarrow \C$.\\

Axiom \ref{monat6} for ${(\Fs)}^0$ amounts to the commutation
for any objects $a, a'$ in $\A$ of the diagram in $\C$
$$ 
\xymatrix{
\un \otimes \un
\ar@{-}[r]^{\cong}
\ar[d]_{{(Fa)}^0  \otimes {(Fa')}^0}
&
\un
\ar[d]^{{( F(a \otimes a'))}^0}
\\
Fa(\un_{\B}) \otimes Fa'(\un_{\B})
\ar[r]_-{{( F^2_{a,a'})}_{\un_{\B}}   }
&
F(a \otimes a')(\un_{\B})
}
$$
which commutes according to Axiom \ref{monat7}
for the monoidal natural transformation
$F^2_{a,a'}: Fa \Box F(a') \rightarrow F(a \otimes a'): \B \rightarrow
\C$.\\

Axiom \ref{monat7} for ${(\Fs)}^0$ amounts to the equality
in $\C$ of the arrows 
${({\Fs}^0)}_{\un_{\A}}: \un \rightarrow  \Fs(\un_{\B})(\un_{\A})$ 
and ${ ( \Fs(\un_{\B}) ) }^0$. These two are respectively 
${F(\un_{\A})}^0: \un \rightarrow F(\un_{\A})(\un_{\B})$
and ${(F^0)}_{\un_{\B}}$
which are equal according to Axiom \ref{monat7} for
the monoidal natural transformation 
$F^0: \un \rightarrow F(\un_{\A}): \B \rightarrow \C$.\\

For any objects $b,b'$ in $\B$, the naturality 
$\Fs b \Box \Fs b' \rightarrow \Fs (b \otimes b'): \A \rightarrow \C$
of ${(\Fs)}^2_{b,b'}$ amounts to the commutation for any arrow 
$f: a \rightarrow a'$ of $\A$ of the diagram in $\C$
$$
\xymatrix{
Fa(b) \otimes Fa(b')
\ar[r]^-{ {(Fa)}^2_{b,b'} }
\ar[d]_{ {F(f)}_b \otimes {F(f)}_{b'} }
&
Fa ( b \otimes b')
\ar[d]^{{F(f)}_{b \otimes b'}}
\\
Fa'(b) \otimes Fa'(b')
\ar[r]_-{{(Fa')}^2_{b,b'}}
&
Fa'( b \otimes b')
}
$$
which does commute according to Axiom \ref{monat6}
for the monoidal $F(f): Fa \rightarrow Fa': \B \rightarrow \C$.\\

For any objects $b$, $b'$ in $\B$,
Axiom \ref{monat6} for ${(\Fs)}^2_{b,b'}$ amounts to 
the commutation for any objects $a$,$a'$ in $\A$ of
the diagram in $\C$
$$
\xymatrix@C=5pc{
( Fa(b) \otimes Fa(b') )
\otimes
( Fa'(b) \otimes Fa'(b'))
\ar[r]^-{  {(Fa)}^2_{b,b'} \otimes  {(Fa')}^2_{b,b'}  }
\ar@{-}[d]_{\cong}
&
Fa(b \otimes b')
\otimes 
Fa'(b \otimes b')
\ar[dd]^{{(F^2_{a,a'})}_{b \otimes b'}}
\\
( Fa(b) \otimes Fa'(b) )
\otimes 
( Fa(b') \otimes Fa'(b') )
\ar[d]_{  {(F^2_{a,a'})}_b    \otimes {(F^2_{a,a'})}_{b'}   }
&
\\
F(a \otimes a')b \otimes F(a \otimes a')b'
\ar[r]_-{  {({F(a \otimes a')}^2)}_{b,b'}    }
&
F(a \otimes a')(b \otimes b').
}
$$
This diagram commutes according to Axiom 
\ref{monat6} for the monoidal natural transformation
$F^2_{a,a'}: Fa \Box Fa' \rightarrow F(a \otimes a'): \B \rightarrow \C$.\\ 

For any objects $b,b'$ in $\B$,
Axiom \ref{monat7} for ${(\Fs)}^2_{b,b'}$ amounts to the 
equality for any objects $b$ and $b'$ in $\B$
of the two arrows in $\C$
$$
\xymatrix@C=4pc{
\un
\ar[r]^-{ {(\Fs b \Box \Fs b')}^0  } 
&
(\Fs b \Box \Fs b')(\un_{\A})
\ar[r]^-{  {( {(\Fs)}^2_{b,b'} )}_{\un_{\A}}     }
&
\Fs(b \otimes b')(\un_{\A})
}
$$
and 
${( \Fs(b \otimes b') )}^0: \un 
\rightarrow \Fs ( b \otimes b')(\un_{\A})$.\\
The first arrow rewrites
$$
\xymatrix@C=3pc{
\un_{\C} 
\ar@{-}[r]^-{\cong}
&
\un_{\C} \otimes \un_{\C}
\ar[r]^-{{F^0}_b \otimes {F^0}_{b'}}
&
F(\un_{\A})(b) \otimes F(\un_{\A})(b')
\ar[r]^-{{(F(\un_{\A}))}^2_{b,b'}  }
&
F(\un_{\A})(b \otimes b')
}
$$
and the second one is 
${F^0}_{b \otimes b'}: \un_{\C} \rightarrow F(\un_{\A})(b \otimes b')$.
These two arrows are equal according to Axiom \ref{monat6}
for the monoidal natural transformation 
$F^0: \un \rightarrow F(\un_{\A}): \B \rightarrow \C$.
(Remember that ${\un_{[\B,\C]}}^2_{b,b'}$ is the 
canonical isomorphism 
$\un_{\C} \otimes \un_{\C} \rightarrow \un_{\C}$.)\\

The naturality in $b$ and $b'$ of the collection 
of arrows
${(\Fs)}^2_{b,b'} : \Fs b \Box \Fs b' \rightarrow \Fs (b \otimes b')$
in
$[\A,\C]$ is the commutation for any arrows
$f: b_1 \rightarrow b_1'$ and $g : b_2 \rightarrow b_2'$ in $\B$
of the diagram in $[\A,\C]$
$$
\xymatrix@C=3pc{
\Fs b_1 \Box \Fs b_1'
\ar[r]^-{ {(\Fs)}^2_{b_1,b_1'} }
\ar[d]_{\Fs f \Box \Fs g}
&
\Fs (b_1 \otimes b_1')
\ar[d]^{\Fs ( f \otimes g)}
\\
 \Fs b_2 \Box \Fs b_2'
\ar[r]_-{ {(\Fs)}^2_{b_2,b_2'} }
&
\Fs (b_2 \otimes b_2')
}
$$
which is pointwise in any $a$
$$
\xymatrix@C=3pc{
Fa (b_1) \otimes Fa (b_1')
\ar[r]^-{ {(Fa)}^2_{b_1,b_1'} }
\ar[d]_{Fa(f) \otimes Fa(g)}
&
Fa (b_1 \otimes b_1')
\ar[d]^{Fa ( f \otimes g)}
\\
Fa(b_2) \otimes Fa(b_2')
\ar[r]_-{ {(Fa)}^2_{b_2,b_2'} }
&
Fa (b_2 \otimes b_2')
}
$$
that does commute according to the naturality
in $b,b'$ of the collection of arrows
$(Fa)^2_{b,b'}: Fa(b) \otimes Fa(b') \rightarrow Fa(b \otimes b')$.\\

Axioms \ref{mofun3} for the triple $( \Fs, {(\Fs)}^0, {(\Fs)}^2)$
amounts to the commutation for any objects $b,b',b''$ in $\B$ 
of the diagram in $[\A,\C]$
$$
\xymatrix{
\Fs b \Box ( \Fs b' \Box  \Fs b'') 
\ar[d]_{1 \Box {(\Fs)}^2_{b',b''}}
\ar[r]^-{\ac}
& 
(\Fs b \Box  \Fs b') \Box  \Fs b''
\ar[d]^{  {(\Fs)}^2_{b,b'} \Box 1}
\\
\Fs b \Box \Fs(b' \otimes b'')
\ar[d]_{{(\Fs)}^2_{b, b' \otimes b''}}
&
\Fs (b \otimes b') \Box \Fs b''
\ar[d]^{  {(\Fs)}^2_{b \otimes  b', b''}}
\\
\Fs ( b \otimes (b' \otimes b''))
\ar[r]_-{\Fs(\ac)}
&
\Fs ( (b \otimes b') \otimes b'').
}
$$
This diagram is pointwise in any $a$
$$
\xymatrix{
Fa(b) \otimes ( Fa(b') \otimes  Fa(b'') )
\ar[d]_{1 \otimes {(Fa)}^2_{b',b''}}
\ar[r]^-{\ac}
& 
(Fa(b) \otimes Fa(b')) \otimes  Fa(b'')
\ar[d]^{{(Fa)}^2_{b,b'} \otimes  1}
\\
Fa(b) \otimes Fa(b' \otimes b'')
\ar[d]_{{(Fa)}^2_{b, b' \otimes b''}}
&
Fa(b \otimes b') \otimes Fa( b'')
\ar[d]^{{(Fa)}^2_{b \otimes  b', b''}}
\\
Fa ( b \otimes (b' \otimes b''))
\ar[r]_-{Fa(\ac)}
&
Fa ( (b \otimes b') \otimes b'')
}.
$$
which commutes according to Axiom \ref{mofun3}
for the monoidal functor $Fa: \B \rightarrow \C$.\\

In the same way, Axioms \ref{mofun41} and \ref{mofun42}
and \ref{symofun5} for the triple $(\Fs, {(\Fs)}^0, {(\Fs)}^2)$ 
can be deduced  respectively from Axioms \ref{mofun41} 
and \ref{mofun42} and 
\ref{symofun5} for the monoidal functors $Fa$ for all objects
$a$ of $\A$.\\ 

We show now that given any monoidal natural transformation 
$\theta : F \rightarrow G: \A \rightarrow [\B,\C]$,
and any object $b$ in $\B$, the natural transformation
${\theta^{\ud}}_b: F^{\ud}b \rightarrow G^{\ud}b: \A \rightarrow \C$
is monoidal $(\Fs b, {(\Fs b)}^0, {(\Fs b)}^2 ) 
\rightarrow ( \Gs b, {(\Gs b)}^0 , {(\Gs b)}^2 )$.\\

Axiom \ref{monat6} for ${\theta^{\ud}}_b$ amounts to the 
commutation for any objects $a$, $a'$ in $\A$
of the diagram in $\C$
$$
\xymatrix{
Fa(b) \otimes Fa'(b)
\ar[d]_{{(\theta_a)}_b \otimes {(\theta_{a'})}_b }
\ar[r]^-{{( F^2_{a,a'} )}_b}
&
F(a \otimes a')(b)
\ar[d]^{ ({\theta_{a \otimes a'})}_b }
\\
 Ga(b) \otimes Ga'(b)
\ar[r]_-{{( G^2_{a,a'} )}_b}
&
G(a \otimes a')(b)
}
$$
which is the pointwise version in $b$ of 
$$
\xymatrix{
Fa  \Box Fa'
\ar[d]_{ \theta_a \Box \theta_{a'} }
\ar[r]^-{ F^2_{a,a'} }
&
F(a \otimes a')
\ar[d]^{ \theta_{a \otimes a'} }
\\
Ga \Box Ga'
\ar[r]_-{G^2_{a,a'}}
&
G(a \otimes a')
}
$$ in $[\B,\C]$,
that commutes according to Axiom \ref{monat6} for
$\theta: F \rightarrow G$.
Similarly Axiom \ref{monat7} for $\theta^{\ud}_b$
is the pointwise version in $b$ of Axiom 
\ref{monat7} for $\theta$.\\

Therefore the natural transformation $\theta^{\ud}: F^{\ud}
\rightarrow G^{\ud}$ factorises as $U * \theta^*$ as claimed 
in \ref{factonat}.
It remains to show that $\theta^*$
is monoidal $\Fs \rightarrow \Gs: \A \rightarrow [\B,\C]$.
Axiom \ref{monat6} for $\theta^*$ is that 
for any objects $b$ and $b'$ in $\B$
the diagram in $[\A,\C]$ commutes
$$
\xymatrix{
\Fs b \Box \Fs b'
\ar[r]^-{ { (\Fs) }^2_{b,b'} }
\ar[d]_{{\theta^*}_b \Box {\theta^*}_{b'}}
&
\Fs(b \otimes b')
\ar[d]^{{\theta^*}_{b \otimes b'}}
\\
\Gs b \Box \Gs b'
\ar[r]_-{ { (\Gs) }^2_{b,b'} }
&
\Gs(b \otimes b')
.}
$$
Pointwise in any $a$ in $\A$,
this diagram is
$$
\xymatrix{
Fa(b) \otimes Fa(b')
\ar[r]^-{{(Fa)}^2_{b,b'} }
\ar[d]_{{(\theta_a)}_b \Box {(\theta_a)}_{b'}}
&
Fa(b \otimes b')
\ar[d]^{{(\theta_a)}_{b \otimes b'}}
\\
Ga(b) \Box Ga( b')
\ar[r]_-{ {(Ga)}^2_{b,b'} }
&
Ga(b \otimes b')
}
$$
which commutes according to Axiom \ref{monat6} 
for the monoidal $\theta_a: Fa \rightarrow Ga: \B \rightarrow \C$.
Similarly Axiom \ref{monat7} 
for $\theta^*$ can be checked pointwise
and results from Axiom \ref{monat7} for the monoidal transformations
$\theta_a : Fa \rightarrow Ga$ for all objects $a$ of $\A$.\\

It is now straightforward to check 
that the isomorphism \ref{dualhom} admits 
a strict monoidal structure 
$\Dual_{\A,\B,\C}: [\A,[\B,\C]] \rightarrow [\B,[\A,\C]]$.
We detail briefly below the computation to establish
the equality ${(F \Box G)}^* = F^* \Box G^*: \B \rightarrow \C$
for any symmetric monoidal functors $F,G: \A \rightarrow [\B,\C]$.\\

For any object $b$ of $\B$, 
the functors $(F^* \Box G^*)(b)$
and ${(F \Box G)}^*(b)$ from $\A$ to $\C$
both  
send any arrow $f: a \rightarrow a'$ 
to the arrow 
$$
{F(f)}_b \otimes {G(f)}_b: 
F(a)(b) \otimes G(a)(b) \rightarrow F(a')(b) \otimes G(a')(b)  
$$
in $\C$.\\

For any objects $b$ in $\B$ and $a,a'$ in $\A$,
the arrows ${({(F \Box G)}^*(b))}^2_{a,a'}$
and ${((F^* \Box G^*)(b))}^2_{a,a'}$ in $\C$
are both equal to 
$$
\xymatrix{
(F(a)(b) \otimes G(a)(b)) 
\otimes 
(F(a')(b) \otimes G(a')(b))
\ar@{-}[d]^{\cong}
\\
(F(a)(b) \otimes F(a')(b))
\otimes 
(G(a)(b) \otimes G(a')(b))
\ar[d]^{{(F^2_{a,a'})}_b \otimes {(G^2_{a,a'})}_b}
\\
F(a \otimes a')(b) \otimes G(a \otimes a')(b).
}
$$

For any object $b$ in $\B$, the 
arrows ${({(F \Box G)}^*(b))}^0$
and ${((F^* \Box G^*)(b))}^0$ in $\C$
are both equal to
$$\xymatrix{
\un
\ar@{-}[r]^-{\cong}
&
\un \otimes \un
\ar[r]^-{F^0_b \otimes G^0_b}
&
F(\un_{\A})(b) \otimes G(\un_{\A})(b)
.}$$ 

For any arrow $g: b \rightarrow b'$ in $\B$,
$(F \Box G)^*$ and $F^* \Box G^*$ both send 
$g$ to the transformation between functors
$\A \rightarrow \C$ with component 
in any object $a$, the arrow
$$F(a)(g) \otimes G(a)(g):
F(a)(b) \otimes G(a)(b) 
\rightarrow F(a)(b') \otimes G(a)(b').$$ 

For any objects $b,b'$ in $\B$, the natural
transformations between functors $\A \rightarrow \C$
$${({(F \Box G)}^*)}^2_{b,b'}:
{(F \Box G)}^*(b) \Box {(F \Box G)}^*(b') 
\rightarrow 
{(F \Box G)}^*(b \otimes b')$$
and
$${(F^* \Box G^*)}^2_{b,b'}: 
(F^*(b) \Box G^*(b)) \Box (F^*(b') \Box G^*(b'))
\rightarrow 
F^*(b \otimes b')
\Box G^*(b \otimes b')$$
have both for component in any $a$ the arrow
$$
\xymatrix{
( F(a)(b) \otimes G(a)(b) )
\otimes 
( F(a)(b') \otimes G(a)(b') )
\ar@{-}[d]^{\cong}
\\
( F(a)(b) \otimes F(a)(b') )
\otimes
(G(a)(b) \otimes G(a)(b'))
\ar[d]^{ {(Fa)}^2_{b,b'} \otimes {(Ga)}^2_{b,b'} }
\\
F(a)(b \otimes b') \otimes G(a)(b \otimes b')
.}
$$ 

The natural transformations between functors $\A \rightarrow \C$
$${({(F \Box G)}^*)}^0: \un_{[\A,\C]} \rightarrow {(F \Box G)}^*(\un_{\B})$$
and
$${(F^* \Box G^*)}^0: \un_{[\A,\C]} \rightarrow (F^* \Box G^*)(\un_{\B})$$
have both for component in any object $a$ the 
arrow of $\C$
$$\xymatrix@C=4pc{
\un
\ar@{-}[r]^-{\cong}
&
\un \otimes \un
\ar[r]^-{{F(a)}^0 \otimes {G(a)}^0} 
&
F(a)(\un_{\B}) \otimes G(a)(\un_{\B}).
}$$

The computation to check that the image by the isomorphism
\ref{dualhom} 
of the unit $\un: \A \rightarrow [\B,\C]$
is $\un: \B \rightarrow [\A,\C]$ is straightforward
and left to the reader.
That the functor \ref{dualhom} preserves the canonical arrows
$\ac$, $\rc$, $\lc$ and $\syc$ is immediate since these
arrows are defined ``pointwise''.\\

Eventually to establish that $\Dual_{\A,\B,\C}$ has inverse
$\Dual_{\B,\A,\C}$ in $\SMC$, the only non immediate point to check 
is that any symmetric monoidal functor $F: \A \rightarrow [\B,\C]$
is equal to $F^{**}$. To show this
for any such $F$, it remains to check the
equality of the monoidal structures of $Fa$ and $F^{**}a$
for all objects $a$ of $\A$, and then to check the equality
of the monoidal structures of $F$ and $F^{**}$. 
Those result from the equalities
\ref{dualstr1} and \ref{dualstr2}:\\
- For any objects $a$ in $\A$ and $b,b'$ in $\B$,
${ ( F^{**}a ) }^2_{b,b'}$ $=$ ${({F^*}^2_{b,b'})}_a$ $=$
${F(a)}^2_{b,b'}$;\\
- For any objects $a$ in $\A$,
${ ( F^{**}a ) }^0$ $=$ ${({F^*}^0)}_a$ $=$ ${(Fa)}^0$;\\
- For any objects $a,a'$ in $\A$ and $b$ in $\B$,
${   ({ ( F^{**} ) }^2_{a,a'}) }_b$ $=$ 
${(F^*b)}^2_{a,a'}$ $=$ ${(F^2_{a,a'})}_b$;\\
- For any object $b$ in $\B$,
${{ ( F^{**} ) }^0}_b$ $=$ ${F^*(b)}^0$ $=$ ${F^0}_b$.
 
\begin{remark}\label{Fsbstrict}
Any symmetric monoidal functor $F: \A \rightarrow [\B, \C]$
is strict (respectively strong) if and only
if for all objects $b$ of $\B$ the functors $F^*(b): \A \rightarrow \C$
are strict (respectively strong).
\end{remark}

From now on, we also call {\em dual} the symmetric 
monoidal functors or monoidal transformations that 
correspond via the isomorphism $\Dual$.
\end{section}

\begin{section}{$F^0$ and $F^2$ are monoidal}\label{2cellsF0F2}
It is proved in this section that for any {\em symmetric} monoidal 
functor $F: \A \rightarrow \B$, the natural transformations 
$F^0: \Delta_{\un_\B} \rightarrow F \circ \Delta_{\un_{\A}}: 
\termcat \rightarrow \B$ 
and $F^2: \Ten \circ (F \times F) \rightarrow F \circ \Ten
: \A \times \A \rightarrow \B$ are
monoidal.\\

\begin{lemma}\label{lem_symofun2}
For any symmetric monoidal functor $F: \A \rightarrow \B$,
any objects $x,y,z,t$ in $\A$, the following diagram in $\B$
commutes
$$\xymatrix{
(Fx \otimes Fy) \otimes (Fz \otimes Ft)
\ar[d]_{F^2_{x,y} \otimes F^2_{z,t}}
\ar@{-}[r]^-{\cong} 
& (Fx \otimes Fz) \otimes (Fy \otimes Ft)
\ar[d]^{F^2_{x,z} \otimes F^2_{y,t}} \\
F(x \otimes y) \otimes  F(z \otimes t)
\ar[d]_{F^2_{x \otimes y, z \otimes t}}
& 
F(x \otimes z) \otimes F(y \otimes t)
\ar[d]^{F^2_{x \otimes z, y \otimes t}} \\
F((x \otimes y) \otimes (z \otimes t))
\ar@{-}[r]_-{F(\cong)} &
F((x \otimes z) \otimes (y \otimes t)). 
}$$
\end{lemma}
(Note: In the above lemma we omit to say that
the two isomorphisms $\cong$ are the canonical ones 
permuting respectively the $Fy$, $Fz$ and $y$ and $z$.)\\

\pf
Consider the following pasting of diagrams in $\B$
{\tiny
\begin{center}
$\xymatrix{
(Fx \otimes Fy) \otimes (Fz \otimes Ft)
\ar[r]^{F^2 \otimes 1}
\ar[d]_{\cong}
\ar@{}[rd]|{(a)}
&
F(x \otimes y) \otimes (Fz \otimes Ft)
\ar[r]^{1 \otimes F^2}
\ar[d]|{\cong}
\ar@{}[rrd]|{(b)}
&
F(x \otimes y) \otimes F(z \otimes t)
\ar[r]^{F^2}
&
F((x \otimes y) \otimes (z \otimes t))
\ar[d]^{\cong}
\\
((Fx \otimes Fy) \otimes Fz) \otimes Ft
\ar[d]_{\cong}
\ar@{}[rrd]|{(c)}
\ar[r]^{(F^2 \otimes 1) \otimes 1}
&
(F(x \otimes y) \otimes Fz) \otimes Ft
\ar[r]^{F^2 \otimes 1}
&
F( (x \otimes y) \otimes z ) \otimes Ft
\ar[r]^{F^2}
\ar[d]|{\cong}
\ar@{}[rd]|{(d)}
&
F( ((x \otimes y) \otimes z ) \otimes t)
\ar[d]^{\cong}
\\
(Fx \otimes (Fy \otimes Fz)) \otimes Ft
\ar[r]_{(1 \otimes F^2) \otimes 1}
&
(Fx \otimes F(y \otimes z)) \otimes Ft
\ar[r]_{F^2 \otimes 1}
&
F(x \otimes (y \otimes z)) \otimes Ft
\ar[r]_{F^2}
&
F((x \otimes (y \otimes z)) \otimes t).}
$
\end{center}
}
Here diagram $(a)$ commutes since $\B$ is 
monoidal, diagrams $(b)$ and $(c)$ commute 
according to Axiom \ref{mofun3} for $F$ and diagram
$(d)$ commutes by naturality of $F^2$.
This shows that the diagram\\
$(1)$:
{\tiny
\begin{center}
$\xymatrix{
(Fx \otimes Fy) \otimes (Fz \otimes Ft)
\ar[r]^{F^2 \otimes 1}
\ar[d]_{\cong}
&
F(x \otimes y) \otimes (Fz \otimes Ft)
\ar[r]^{1 \otimes F^2}
&
F(x \otimes y) \otimes F(z \otimes t)
\ar[r]^{F^2}
&
F((x \otimes y) \otimes (z \otimes t))
\ar[d]^{\cong}
\\
(Fx \otimes (Fy \otimes Fz)) \otimes Ft
\ar[r]_{(1 \otimes F^2) \otimes 1}
&
(Fx \otimes F(y \otimes z)) \otimes Ft
\ar[r]_{F^2 \otimes 1}
&
F(x \otimes (y \otimes z)) \otimes Ft
\ar[r]_{F^2}
&
F((x \otimes (y \otimes z)) \otimes t)
}$
\end{center}
}
commutes. 
Diagram $(1)$ with $y$ and $z$ inverted is the 
commuting diagram\\
$(2)$: 
{\tiny 
$$\xymatrix{
(Fx \otimes (Fz \otimes Fy)) \otimes Ft
\ar[r]^{(1 \otimes F^2) \otimes 1}
\ar[d]_{\cong}
&
(Fx \otimes F(z \otimes y)) \otimes Ft
\ar[r]^{F^2 \otimes 1}
&
F(x \otimes (z \otimes y)) \otimes Ft
\ar[r]^{F^2}
&
F((x \otimes (z \otimes y)) \otimes t)
\ar[d]^{\cong}
\\
(Fx \otimes Fz) \otimes (Fy \otimes Ft)
\ar[r]_{F^2 \otimes 1}
&
F(x \otimes z) \otimes (Fy \otimes Ft)
\ar[r]_{1 \otimes F^2}
&
F(x \otimes z) \otimes F(y \otimes t)
\ar[r]_{F^2}
&
F((x \otimes z) \otimes (y \otimes t))
.}$$
}
Moreover in the following pasting\\
$(3)$:\\
{\tiny
$$
\xymatrix{
(Fx \otimes (Fy \otimes Fz)) \otimes Ft
\ar[r]^{(1 \otimes F^2) \otimes 1}
\ar[d]_{\cong}
\ar@{}[rd]|{(e)}
&
(Fx \otimes F(y \otimes z)) \otimes Ft
\ar[r]^{F^2 \otimes 1}
\ar[d]|{\cong}
\ar@{}[rd]|{(f)}
&
F(x \otimes (y \otimes z)) \otimes Ft
\ar[r]^{F^2}
\ar[d]|{\cong}
\ar@{}[rd]|{(g)}
&
F( (x \otimes (y \otimes z)) \otimes t)
\ar[d]^{\cong}
\\
(Fx \otimes (Fz \otimes Fy)) \otimes Ft
\ar[r]_{(1 \otimes F^2) \otimes 1}
&
(Fx \otimes F(z \otimes y)) \otimes Ft
\ar[r]_{F^2 \otimes 1}
&
F(x \otimes (z \otimes y)) \otimes Ft
\ar[r]_{F^2}
&
F( (x \otimes (z \otimes y)) \otimes t)
}
$$
}
diagram $(e)$ commutes since $F$ is symmetric,
and diagrams $(f)$ and $(g)$ commute by naturality of
$F^2$.
A pasting of $(1)$, $(2)$ and $(3)$ above gives the 
expected result. 
\epf

\begin{lemma}\label{lem_mofun1}
For any symmetric monoidal functor $F:\A \rightarrow \B$, 
the following diagram in $\B$ commutes
\begin{center}
$\xymatrix{
\un_{\B} \ar@{-}[r]^-{\cong} \ar[dd]_{F^0} & 
\un_{\B} \otimes \un_{\B} \ar[d]^{F^0 \otimes F^0} \\
& F(\un_{\A}) \otimes F(\un_{\A})
\ar[d]^{F^2_{\un_{\A}, \un_{\A}} }\\
F(\un_{\A})
\ar@{-}[r]_-{F(\cong)}
 &
F(\un_{\A} \otimes \un_{\A}).  
}$
\end{center}
\end{lemma}
\pf
The above diagram is the following pasting
of commutative diagrams (naturality of $\rc$ in $\B$
and Axiom \ref{mofun41} for $F$)
\begin{center}
$\xymatrix{
\un_{\B} 
\ar[dd]_{F^0} 
& 
\un_{\B} \otimes \un_{\B} \ar[r]^-{F^0 \otimes F^0}
\ar[l]_{\rc} 
\ar[d]_{F^0 \otimes 1} &
F(\un_{\A}) \otimes F(\un_{\A}) 
\ar[dd]^{ F^2_{\un_{\A}, \un_{\A}} }
\\
& F(\un_{\A}) \otimes \un_{\B} \ar[ru]_{1 \otimes F^0} 
\ar[ld]_{\rc}
& \\
F(\un_{\A})
&  
& F(\un_{\A} \otimes \un_{\A})
\ar[ll]^{F(\rc)}
}$
\end{center}
\epf

\begin{lemma}\label{F2mon} For any symmetric monoidal 
$F: \A \rightarrow \B$,
the natural transformation 
$$F^2: Ten \circ (F \times F) 
\rightarrow F \circ Ten: \A \times \A
\rightarrow \B$$ is monoidal.
\end{lemma}
\pf
That $F^2$ satisfies Axiom \ref{monat6} is to say that 
the diagram in $\B$ commutes
{\small
$$
\xymatrix@C=7pc{
\Ten \circ (F \times F)(x,y) \otimes \Ten \circ (F \times F)(z,t) 
\ar[r]^-{ {(\Ten \circ (F \times F))}^2_{(x,y),(z,t)} }
\ar@{=}[d]
&
\Ten \circ (F \times F)( (x,y) \otimes (z,t))
\ar@{=}[d]
\\
(F(x) \otimes F(y)) \otimes (F(z) \otimes F(t))
\ar[d]_{F^2_{x,y} \otimes F^2_{z,t}}
& 
F(x \otimes z) \otimes F(y \otimes t)
\ar[d]^{F^2_{x \otimes z, y \otimes t}}
\\
F(x \otimes y) \otimes F(z \otimes t)
\ar@{=}[d]
&
F((x \otimes z) \otimes (y \otimes t))
\ar@{=}[d] 
\\
F \circ \Ten (x,y) \otimes F \circ \Ten (z,t)
\ar[r]_-{ {( F \circ \Ten )}^2_{(x,y),(z,t)} }
&
F \circ \Ten ( (x,y) \otimes (z,t))
}$$
}
commutes for any objects $x$, $y$, $z$, $t$ in $\A$.
Note that
${(\Ten \circ (F \times F))}^2_{(x,y),(z,t)}$ 
is the arrow
{\small
$$\xymatrix@C=3pc{ 
(F(x) \otimes F(y)) \otimes (F(z) \otimes F(t))
\ar@{-}[r]^-{\cong}
&
(F(x) \otimes F(z)) \otimes (F(y) \otimes F(t))
\ar[r]^-{F^2_{x,z} \otimes F^2_{y,t}}
&
F(x \otimes z) \otimes F(y \otimes t) 
}
$$ 
}
and that
${(F \circ \Ten)}^2_{(x,y),(z,t)}$
is the arrow 
$$
\xymatrix@C=3pc{
F(x \otimes  y) \otimes F(z  \otimes  t)
\ar[r]^-{F^2_{x \otimes y, z \otimes t}}  
&
F( (x \otimes y) \otimes (z \otimes t))
\ar@{-}[r]^-{F(\cong)}
&
F( (x \otimes z) \otimes (y \otimes t)).
}
$$
Therefore Axiom \ref{monat6} for $F^2$ is 
just Lemma \ref{lem_symofun2}.
Axiom \ref{monat7} for $F^2$ is just Lemma \ref{lem_mofun1}. 
\epf

\begin{lemma}\label{F0mon}
Given any symmetric monoidal functor $F: \A \rightarrow \B$,
the arrow $F^0: \un_{\B} \rightarrow F(\un_{\A})$
is the unique component of a monoidal natural transformation
$$\Delta_{\un_{\B}} \rightarrow F \circ \Delta_{\un_{\A}}
: \termcat \rightarrow \B.$$
\end{lemma}
\pf
Axiom \ref{monat6} for this natural transformation 
is Lemma \ref{lem_mofun1} 
and Axiom \ref{monat7} is trivial. 
\epf
\end{section}

\begin{section}{pre- and post-composition functors}\label{prepost}
In this section, $\A$, $\B$ and $\C$ stand for arbitrary
symmetric monoidal categories.
Since $\SMC$ is a 2-category, 
any 2-cell $\sigma: F \rightarrow G : \B \rightarrow \C$ in $\SMC$ 
induces a 2-cell in $\CAT$
$$\SMC(\A, \sigma): \SMC(\A,F) \rightarrow \SMC(\A,G): 
\SMC(\A,\B) \rightarrow \SMC(\A, \C),$$
and the above assignments $F \mapsto \SMC(\A,F)$ 
and $\sigma \mapsto \SMC(\A,\sigma)$ define a functor
$$\post_{\A,\B,\C}: \SMC(\B,\C) \rightarrow \CAT(\SMC(\A,\B),\SMC(\A,\C)).$$
We are going to show that the functor 
$\post_{\A,\B,\C}$ factorises as
$$\xymatrix@C=3pc{
\SMC(\B,\C) 
\ar[r]^-{\Post_{\A,\B,\C}}
&
\SMC([\A,\B],[\A,\C])
\ar[r]^-{U}
&
\CAT(\SMC(\A,\B),\SMC(\A,\C))
}$$
and that the above functor $\Post_{\A,\B,\C}$ admits a symmetric 
strict monoidal
structure $[\B,\C] \rightarrow [[\A,\B],[\A,\C]]$,
denoted $[\A,-]_{\B,\C}$.
The dual of $[\A,-]_{\B,\C}$, denoted 
$$[-,\C]_{\A,\B}: [\A,\B] \rightarrow [[\B,\C],[\A,\C]],$$
will be also briefly considered.\\

Given any symmetric monoidal functor $F: \B \rightarrow \C$,
the monoidal structure $$[\A,F] = (\SMC(\A,F), [\A,F]^0, [\A,F]^2) : 
[\A,\B] \rightarrow [\A,\C]$$
on the functor $\SMC(\A,F)$, is as follows.\\

The monoidal natural transformation ${[\A,F]}^0 : \un 
\rightarrow F \circ \un_{[\A,\B]}: 
\A \rightarrow \C$ 
is the pasting in $\SMC$ 
$$
\xymatrix{
\A \ar[rr] & & \termcat \ar[rr]^{ \Delta_{\un_{\B}} } 
\ar[rrdd]_{ \Delta_{\un_{\C}}   }
&  & \B \ar[dd]^{F}
\\
& & & \ar@{=>}[ru]^{F^0}
\\
& & & & \C
}
$$
where 2-cell $F^0$ is that
of Lemma
\ref{F0mon}.
The component of ${[\A,F]}^0$ in any $a$ in $\A$, is therefore 
$F^0: \un_{\C} \rightarrow F(\un_{\B})$.\\

For any symmetric monoidal functors $G,H: \A \rightarrow \B$,
the monoidal natural transformation 
${[\A,F]}^2_{G,H} : (F \circ G) \Box (F \circ H) \rightarrow
F \circ (G \Box H): \A \rightarrow \C$ is the pasting 
in $\SMC$
\begin{center}
$\xymatrix{
& & &
\C \otimes \C
\ar@{=>}[dd]^{F^2} 
\ar[rd]^{\Ten}
&
\\
\A \ar[r]^-{\Diag} & 
\A \times \A \ar[r]^-{G \times H} &
\B \times \B 
\ar[ru]^{F \times F}  
\ar[rd]_{\Ten} 
& 
&
\C\\
& & &
\B
\ar[ru]_{F}
&
}$
\end{center}
where the 2-cell $F^2$ is that of Lemma \ref{F2mon}.
The component of ${[\A,F]}^2_{G,H}$ in any $a$ in $\A$, is therefore
$$F^2_{G(a),H(a)}: 
FG(a) \otimes FH(a)\rightarrow F(G(a) \otimes H(a)).$$

The naturalities in $G$ and $H$ of the collection
of arrows
${[\A,F]}^2_{G,H}: (F \circ G) \Box (F \circ H) 
\rightarrow F \circ ( G \Box H)$
of $\SMC(\A,\C)$, results
actually from the fact that $\SMC$ is a 2-category.
For instance the naturality in $G$ means that for any 
monoidal transformation $\tau: G \rightarrow G': \A \rightarrow \B$,
the diagram 
$$\xymatrix{
(F \circ G) \Box (F \circ H) 
\ar[r]^-{(F * \tau) \Box 1} 
\ar[d]_{{[\A,F]}^2_{G,H}}
&
(F \circ G') \Box (F \circ H)
\ar[d]^{{[\A,F]}^2_{G',H}}\\
F \circ (G \Box H) \ar[r]_-{F * (\tau \Box 1)} &
F \circ (G'\Box H) 
}$$
in $[\A,\C]$ is commutative
and this holds by the interchange law in $\CAT$ 
for the composition of the 2-cells
\begin{center}
$\xymatrix{
& & 
\ar@{=>}[dd]^{\tau \times 1}
 & & \C \times \C 
\ar[rd]^{\Ten} 
\ar@{=>}[dd]^{F^2}
& 
\\
\A \ar[r]^-{\Diag} & 
\A \times \A
 \ar@/^37pt/[rr]^{G \times H}
 \ar@/_37pt/[rr]_{G' \times H}
 & &
\B \times \B
\ar[ru]^{F \times F}
\ar[rd]_{\Ten}
&
& \C.
\\
& & & & \B \ar[ru]_{F} & 
}$
\end{center}

Axiom \ref{mofun3} for 
the triple $(\SMC(\A,F),{[\A,F]}^2, {[\A,F]}^0)$
is the commutation for any symmetric monoidal
functors $H$,$K$,$L: \A \rightarrow \B$ of the
diagram in $[\A,\C]$
$$
\xymatrix{
(F \circ H) \Box ((F \circ K) \Box (F \circ L)) 
\ar@{-}[r]^-{\cong}
\ar[d]_{1 \Box {[\A,F]}^2_{K,L}  }
&
( (F \circ H) \Box (F \circ K) ) \Box (F \circ L)
\ar[d]^{ {[\A,F]}^2_{H,K}  \Box 1   }
\\
(F \circ H) \Box (F \circ (K \Box L))
\ar[d]_{ {[\A,F]}^2_{H, K \Box L}  }
&
(F \circ (H \Box K)) \Box (F \circ L)
\ar[d]^{ {[\A,F]}^2_{H \Box K ,L}  }
\\ 
F \circ (H \Box (K \Box L))
\ar@{-}[r]_-{F * \ac} 
&
F \circ ((H \Box K) \Box L)
.}
$$
This diagram is pointwise in any $a$
$$
\xymatrix{
FHa \otimes (FKa \otimes FLa) 
\ar@{-}[r]^-{\cong}
\ar[d]_{1 \otimes F^2_{Ka,La}  }
&
( FHa \otimes FKa ) \otimes FLa
\ar[d]^{ F^2_{Ha,Ka}  \otimes  1   }
\\
FHa \otimes F(Ka \otimes La)
\ar[d]_{F^2_{Ha, Ka \otimes La}  }
&
F(Ha \otimes Ka) \otimes FLa
\ar[d]^{ F^2_{Ha \otimes Ka ,La}  }
\\
F(Ha \otimes (Ka \otimes La))
\ar@{-}[r]_-{F(\ac)} 
&
F((Ha \otimes Ka) \otimes La)
}
$$
that commutes according to Axiom \ref{mofun3}
for $F$.\\

Similarly Axioms \ref{mofun41}, \ref{mofun42} and
\ref{symofun5} for $(\SMC(\A,F), {[\A,F]}^0, {[\A,F]}^2)$
can be checked pointwise since they hold for $F$.\\


Note the following that is immediate from the definitions.
\begin{remark}
For any symmetric monoidal functor $F:\B \rightarrow \C$,
if $F$ is strict (resp. strong) then 
$[\A,F]: [\A,\B] \rightarrow [\A,\C]$
is also strict (resp. strong). 
\end{remark}

Consider now a monoidal natural transformation
between symmetric monoidal functors 
$\sigma: F \rightarrow G: \B \rightarrow \C$.
We show that the natural transformation 
$$\SMC(\A,\sigma) : \SMC(\A,F) \rightarrow \SMC(\A,G): 
\SMC(\A,\B) \rightarrow \SMC(\A,\C)$$
is monoidal with respect to the structures 
$[\A,F]$ and $[\A,G]$.  
The resulting 2-cell $[\A,F] \rightarrow
[\A,G]: [\A,\B] \rightarrow [\A,\C]$ in $\SMC$ will be  
later denoted $[\A,\sigma]$.\\

The following pastings in $\SMC$ are equal
since Axiom \ref{monat6} for $\sigma$
states their (pointwise) equality 
\begin{tag}\label{p1monat6} 
$$
\xymatrix{
&
\ar@{=>}[d]^{\sigma \times \sigma}
\\
\B \times \B
\ar@/^38pt/[rr]^{F \times F}
\ar[rr]_{G \times G}
\ar[dd]_{\Ten}
&
&
\C \times \C
\ar[dd]^{\Ten}
\ar@{=>}[lldd]^{G^2}
\\
&&& 
\\
\B
\ar[rr]_G
&
&
\C
}
$$
\end{tag}
and
\begin{tag}\label{p2monat6}
$$
\xymatrix{
\B \times \B
\ar[rr]^{F \times F}
\ar[dd]_{\Ten}
&
&
\C \times \C
\ar[dd]^{\Ten}
\ar@{=>}[lldd]_{F^2}
\\
\\
\B 
\ar[rr]^F
\ar@/_38pt/[rr]_G
&
\ar@{=>}[d]^{\sigma}
&
\C.
\\
& 
}
$$
\end{tag}
That the natural transformation
$\SMC(\A,\sigma)$ satisfies Axiom \ref{monat6}
for the symmetric monoidal structures $[\A,F]$ and $[\A,G]$, 
is that for any symmetric monoidal functors $H, K: \A \rightarrow \B$,
the diagram below in $[\A, \C]$ commutes
\begin{center}
$\xymatrix@C=3pc{
(F \circ H) \Box (F \circ K) 
\ar[r]^-{ {[\A,F]}^2_{H,K} } 
\ar[d]_{ (\sigma * H) \Box (\sigma * K) }
&
F \circ ( H \Box K) 
\ar[d]^{\sigma * (H \Box K)}
\\
(G \circ H) \Box (G \circ K) \ar[r]_-{{[\A,G]}^2_{H,K}}
& G \circ (H \Box K)
.}$ 
\end{center}
This one does commute since the two legs of this diagram are obtained
by composing the pastings \ref{p1monat6} 
and \ref{p2monat6} by 
the arrow 
$\xymatrix{ \A \ar[r]^-{\Diag} 
& \A \times \A \ar[r]^-{H \times K}  
& \B \times \B.
}
$\\

Also the two following pastings in $\SMC$ are equal
since Axiom \ref{monat7} for $\sigma$
states their equality 
\begin{tag}\label{p1monat7}
$$
\xymatrix{
\termcat 
\ar[rr]^{\Delta_{\un_{\B}}}
\ar[rrdd]_{\Delta_{\un_{\C}}}
&&
\B
\ar@/^38pt/[dd]^G
\ar[dd]_F
\\
& \ar@{=>}[ru]^{F^0}
&
\ar@{=>}[r]^{\sigma} &
\\
& & \C
}
$$
\end{tag}
and
\begin{tag}\label{p2monat7}
$$
\xymatrix{
\termcat 
\ar[rr]^{ \Delta_{\un_{\B}} }
\ar[rrdd]_{\Delta_{\un_{\C}}}
& &
\B
\ar[dd]^G
\\
& \ar@{=>}[ru]^{G^0}
\\
& & \C
.}
$$
\end{tag}

Axiom \ref{monat7} for $\SMC(\A,\sigma)$, $[\A,F]$
and $[\A,G]$, is that 
the following diagram in $[\A,\C]$ commutes
\begin{center}
$\xymatrix{ \un_{[\A,\C]}  
\ar[r]^{ {[\A,F]}^0  }
\ar[rd]_{ {[\A,G]}^0 } 
&
F \circ \un_{[\A,\B]}
\ar[d]^{ \sigma * \un_{[\A,\B]} }\\
& G \circ \un_{[\A,\B]}
.}$
\end{center}
It does since its two legs are obtained
by composing the 2-cells \ref{p1monat7} and \ref{p2monat7} by 
the functor $\A \rightarrow \termcat$.\\


So far we have exhibited 
a factorisation of $\post_{\A,\B,\C}$ as
$$
\xymatrix{
\SMC(\B,\C) 
\ar[r]^-{\Post}
&
\SMC([\A,\B],[\A,\C])
\ar[r]^-U
& 
\CAT(\SMC(\A,\B),\SMC(\A,\C))}.$$
It is now rather straightforward to check that 
the above functor $\Post$ admits a strict symmetric
monoidal structure $[\B,\C] \rightarrow [[\A,\B],[\A,\C]]$.
This one will be denoted by $[\A,-]_{\B,\C}$, or simply
$[\A,-]$ when no ambiguity can occur. We give a few 
details of the computation below.\\

First we check that 
for any symmetric monoidal functors $F,G: \B \rightarrow \C$,
the symmetric monoidal functors 
$[\A,F \Box G]$ and $[\A,F] \Box [\A,G]:[\A,\B] \rightarrow [\A,\C]$ 
are equal.\\

Consider any 2-cell  
$\sigma: H \rightarrow K: \A \rightarrow \B$
in $\SMC$. According to the 2-naturality of the 
isomorphism \ref{2natisopro}, the 2-cells
$$\xymatrix{
&
\ar@{=>}[dd]^{\sigma}
\\
\A 
\ar@/^35pt/[rr]^{H}
\ar@/_35pt/[rr]_{K}
&
&
\B
\ar[r]^-{\Diag}
&
\B \times \B
\\
&&
}$$
and
$$\xymatrix{
&&
\ar@{=>}[dd]^{\sigma \times \sigma}
\\
\A 
\ar[r]^-{\Diag}
&
\A \times \A
\ar@/^35pt/[rr]^{H \times H}
\ar@/_35pt/[rr]_{K \times K}
&
&
\B \times \B
\\
&&
}$$
are equal in $\SMC$. By composing those with
$$\xymatrix{
\B \times \B
\ar[r]^-{F \times G}
&
\C \times \C
\ar[r]^-{\Ten}
&
\C
},$$
one obtains the equality of the images
of $\sigma$ 
by the functors $[\A, F \Box G]$ and $[\A,F] \Box [\A,G]$.\\

For any symmetric monoidal functors
$H,K: \A \rightarrow \B$, the
arrows in $[\A,\C]$ 
$${[\A, F \Box G]}^2_{H,K}: 
((F \Box G)H) \Box ((F \Box G) K) \rightarrow (F \Box G)(H \Box K)$$
and 
$${([\A, F] \Box [\A,G])}^2_{H,K}:
(FH \Box GH)
\Box
(FK \Box GK)
\rightarrow
(F(H \Box K)) \Box (G(H \Box K))
,$$
are natural transformations between functors 
$\A \rightarrow \C$ both 
with component in any $a$ the arrow
{\tiny
$$
\xymatrix@C=4pc{
(FHa \otimes GHa) \otimes (FKa \otimes GKa)
\ar@{-}[r]^-{\cong}
&
(FHa \otimes FKa) \otimes (GHa \otimes GKa)
\ar[r]^-{F^2_{Ha,Ka} \otimes G^2_{Ha,Ka}}
&
F(Ha \otimes Ka) \otimes G(Ha \otimes Ka).
}
$$
}

The arrows in $[\A,\C]$
$${ [\A,F \Box G] }^0: \un \rightarrow (F \Box G) \circ \un_{[\A,\B]}$$ 
and 
$${([\A,F] \Box [\A,G])}^0: \un \rightarrow 
( F \circ \un_{ [\A,\B] } ) \Box ( G \circ \un_{[\A,\B]} )  $$
are natural transformations 
both with component in any $a$ the arrow
$$
\xymatrix{
\un_{\C}
\ar@{-}[r]^-{\cong}
&
\un_{\C} \otimes \un_{\C}
\ar[r]^-{F^0 \otimes G^0}
&
F(\un_{\B}) \otimes F(\un_{\B}).
} 
$$

We check now that the symmetric monoidal 
functor $[\A, \un_{[\B,\C]}]$ is the unit 
of $[[\A,\B],[\A,\C]]$.\\

Consider any 2-cell  
$\sigma: H \rightarrow K: \A \rightarrow \B$
in $\SMC$. The 2-cell
$$
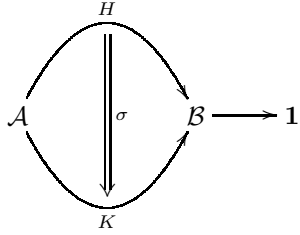
\xymatrix{
&
\ar@{=>}[dd]^{\sigma}
\\
\A
\ar@/^35pt/[rr]^H
\ar@/_35pt/[rr]_K
&
&
\B
\ar[r]
&
\termcat
\\
&&
}$$
and the identity 2-cell at the unique  
arrow $\A \rightarrow \termcat$ are equal according to the
2-naturality of the isomorphism \ref{2natisoter}.
Therefore by composing these 2-cells with the constant 
functor $\Delta_{\un_{\C}} : \termcat \rightarrow \C$,
one obtains the equality between the 2-cells
$\un_{[\B,\C]} * \sigma$, which is the image  
$\sigma$ by the functor $[\A,\un_{[\B,\C]}]$,
and the identity at $\un_{[\A,\C]}$ in $[\A,\C]$,
which is the image of $\sigma$ by the unit 
of $[[\A,\B],[\A,\C]]$.\\

For any symmetric monoidal functors 
$H,K: \A \rightarrow \B$, the 
arrows in $[\A,\C]$ 
$${[\A,\un_{[\B,\C]}]}^2_{H,K}: (\un_{[\B,\C]} \circ H)
\Box (\un_{[\B,\C]} \circ K) \rightarrow \un_{[\B,\C]} \circ (H \Box K)$$ 
and 
$${( \un_{[[\A,\B],[\A,\C]]} )}^2_{H,K}: 
\un \Box \un \rightarrow \un$$
are natural transformations both with component in
any $a$, the canonical arrow 
$\un_{\C} \otimes \un_{\C} \rightarrow \un_{\C}$.\\ 

The arrows in $[\A,\C]$
$${[\A,\un_{[\B,\C]}]}^0: \un \rightarrow \un_{[\B,\C]} \circ
\un_{[\A,\B]}$$
and $${\un_{[[\A,\B],[\A,\C]]}}^0: \un \rightarrow 
\un$$
are natural transformations both with component in any 
$a$, the identity at $\un_{\C}$.\\ 

Eventually that $[\A,-]$ preserves the canonical
arrows $\ac$, $\rc$, $\lc$ and $\syc$, results easily
from the fact that these arrows are defined pointwise in 
functor categories.
For instance the image by $[\A,-]$ of 
$\ac_{F,G,H}: F \Box (G \Box H) \rightarrow (F \Box G) \Box H$, 
for any symmetric monoidal functors $F,G,H: \B \rightarrow \C$,
is the monoidal natural transformation
$$[\A,F \Box (G \Box H)] \rightarrow [\A, (F \Box G) \Box H]: [\A,\B]
\rightarrow [\A,\C]$$
with component in any $K: \A \rightarrow \B$, the monoidal natural
transformation 
$$\ac * K:
(F \Box (G \Box H)) \circ K \rightarrow ((F \Box G) \Box H) \circ K:
\A \rightarrow \C,$$ this one has component in any $a$, the arrow
$$\ac: FKa \otimes (GKa \otimes HKa) \rightarrow (FKa \otimes GKa) \otimes
HKa.$$ Therefore the component in $K$ of $[\A,\ac_{F,G,H}]$ is 
$$\ac_{FK,GK,HK}: FK \Box (GK \Box HK) \rightarrow (FK \Box GK) \Box
HK,$$ which is also the component in $K$ of 
$$\ac_{[\A,F],[\A,G],[\A,H]}: 
[\A,F] \Box ([\A,G] \Box [\A,H]) 
\rightarrow ( [\A,F] \Box [\A,G] ) \Box [\A,H].$$
The cases for $\rc$, $\lc$ and $\syc$ are similar.\\

Since $\SMC$ is a 2-category, one has also the functor
$$\pre_{\A,\B,\C}: \SMC(\A,\B) \rightarrow \CAT(\SMC(\B,\C),\SMC(\A,\C))$$
that assigns any symmetric monoidal functor $F:\A \rightarrow \B$
to the functor 
$\SMC(F,\C) : \SMC(\B,\C) \rightarrow \SMC(\A,\C)$
and any monoidal natural transformation 
$\sigma: F \rightarrow G: \A \rightarrow \B$ to the natural transformation 
$\SMC(\sigma,\C): \SMC(F,\C) \rightarrow \SMC(G,\C)$.
This mere functor $\pre$ is the dual to $\post: \SMC(\B,\C)
\rightarrow \CAT(\SMC(\A,\B),\SMC(\A,\C))$ in the sense of \ref{classicdual}
and therefore, according to Section
\ref{inthomSMC}, it factorises as
$$
\xymatrix{
\SMC(\A,\B) 
\ar[r]^-{\Pre}
&
\SMC([\B,\C],[\A,\C]) 
\ar[r]^-U
& 
\CAT(\SMC(\B,\C),\SMC(\A,\C)).
}
$$
for a functor $\Pre$ that admits 
a symmetric monoidal structure
$[-,\C]_{\A,\B}: [\A,\B] \rightarrow  
[[\B,\C],[\A,\C]]$ dual 
in the sense of \ref{dualhom}
to 
$[\A,-]_{\B,\C}: [\B,\C] \rightarrow [[\A,\B],[\A,\C]].$\\

Note that according to Remark \ref{Fsbstrict},
for any symmetric monoidal $F: \A \rightarrow \B$,
the symmetric monoidal structure 
$[-,\C]_{\A,\B}(F):[\B,\C] \rightarrow [\A,\C]$,
which we write $[F,\C]$, is {\em strict}.
Also for any monoidal transformation 
$\sigma: F \rightarrow G: \A \rightarrow \B$,
the natural transformation 
$\SMC(\sigma,\C)$ is 
monoidal $[-,\C](\sigma):[F,\C] \rightarrow [G,\C]$ 
and we write $[\sigma,\C]$ for this last 2-cell in $\SMC$.\\

\end{section}

\begin{section}{The 2-functor $\homSMC: {\SMC}^{op} \times 
\SMC \rightarrow \SMC$}\label{funhomSMC}
In this section a ``hom'' functor
${\SMC}^{op} \times \SMC \rightarrow \SMC$ is defined.
It is denoted 
cautiously $\homSMC$ to avoid ambiguities arising from
a possible overuse of the square-brackets notation $[-,-]$.\\
  
\begin{lemma}\label{homsqcom} 
The diagram in $\SMC$
$$
\xymatrix{
& [\A,\B]
\ar[ld]_{[F,\B]}
\ar[rd]^{[\A,G]} 
&
\\
[\C,\B] 
\ar[rd]_{[\C,G]}
& 
& 
[\A,\D]
\ar[ld]^{[F,\D]}
\\
& [\C,\D]
}
$$
commutes for any $F: \C \rightarrow \A$
and $G: \B \rightarrow \D$.\\
In $\SMC$, for any 2-cells 
$\sigma: F \rightarrow F': \C \rightarrow \A$
and 
$\tau: G \rightarrow G': \B \rightarrow \D$ in $\SMC$, 
the 2-cells of $\SMC$  
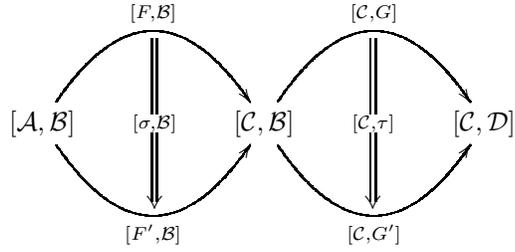
$$\xymatrix{
&
\ar@{=>}[dd]|{[\sigma,\B]}
& 
& 
\ar@{=>}[dd]|{[\C,\tau]}
\\
[\A,\B]
\ar@/^35pt/[rr]^{[F,\B]}
\ar@/_35pt/[rr]_{[F',\B]}
&
&
[\C,\B]
\ar@/^35pt/[rr]^{[\C,G]}
\ar@/_35pt/[rr]_{[\C,G']}
&
&
[\C,\D]
\\
&&&
}$$
and
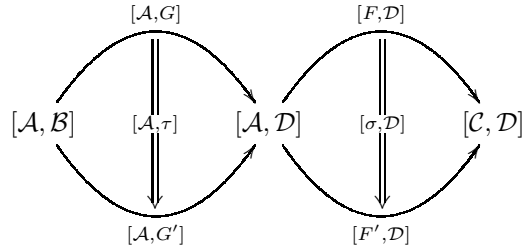
$$\xymatrix{
&
\ar@{=>}[dd]|{[\A,\tau]}
& 
& 
\ar@{=>}[dd]|{[\sigma,\D]}
\\
[\A,\B]
\ar@/^35pt/[rr]^{[\A,G]}
\ar@/_35pt/[rr]_{[\A,G']}
&
&
[\A,\D]
\ar@/^35pt/[rr]^{[F,\D]}
\ar@/_35pt/[rr]_{[F',\D]}
&
&
[\C,\D]
\\
&&&&
}$$
are equal.
\end{lemma}
\pf
To show the first assertion,
since the following diagram in $\CAT$ is commutative
$$
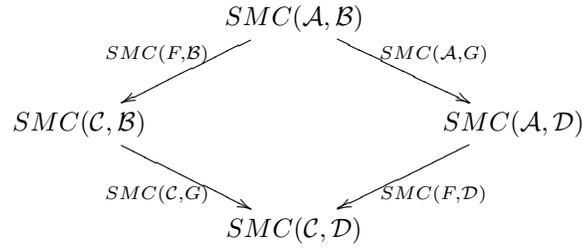
\xymatrix{
& \SMC(\A,\B)
\ar[ld]_{\SMC(F,\B)}
\ar[rd]^{\SMC(\A,G)} 
&
\\
\SMC(\C,\B) 
\ar[rd]_{\SMC(\C,G)}
& 
& 
\SMC(\A,\D)
\ar[ld]^{\SMC(F,\D)}
\\
& \SMC(\C,\D) 
}
$$
it remains to check 
that the monoidal structures of 
$[\C,G] \circ [F,\B]$ and 
$[F,\D] \circ [\A,G]$ are the same.
Straightforward computations show that 
this is the case. (The resulting monoidal structure
is described below after the definition \ref{homfundefmor}).\\

The assertion regarding the monoidal natural 
transformations $\sigma$ and $\tau$ follows from the
fact that since $\SMC$ is a 2-category, the following 
2-cells in $\CAT$ 
$$\xymatrix{
& 
\ar@{=>}[dd]|{\SMC(\sigma,1)}
&
&
\ar@{=>}[dd]|{\SMC(1,\tau)}
\\
\SMC(\A,\B)
\ar@/^35pt/[rr]^{\SMC(F,\B)}
\ar@/_35pt/[rr]_{\SMC(F',\B)}
&
&
\SMC(\C,\B)
\ar@/^35pt/[rr]^{\SMC(\C,G)}
\ar@/_35pt/[rr]_{\SMC(\C,G')}
&
&
\SMC(\C,\D)
\\
& & & & & 
}$$
and
$$\xymatrix{
& 
\ar@{=>}[dd]|{\SMC(1,\tau)}
&
&
\ar@{=>}[dd]|{\SMC(\sigma,1)}
\\
\SMC(\A,\B)
\ar@/^35pt/[rr]^{\SMC(\A,G)}
\ar@/_35pt/[rr]_{\SMC(\A,G')}
&
&
\SMC(\A,\D)
\ar@/^35pt/[rr]^{\SMC(F,\D)}
\ar@/_35pt/[rr]_{\SMC(F',\D)}
&
&
\SMC(\C,\D)
\\
& & & & & 
}$$
are equal.
\epf

As a consequence of the previous lemma, we can define
for any symmetric monoidal $F:\C \rightarrow \A$
and $G: \B \rightarrow \D$, the symmetric 
monoidal functor $[F,G]: [\A,\B] \rightarrow [\C,\D]$
as either of the composites 
\begin{tag}\label{homfundefmor}
$[\C,G] \circ [F,\B]$  
or $[F,\D] \circ [\A,G]$.
\end{tag}
Note that the monoidal structure 
of $[F,G]$ above is the following:\\
- For any symmetric monoidal functors $H,K: \A \rightarrow \B$,
the natural monoidal transformation 
$${[F,G]}^2_{H,K}: GHF \Box GKF \rightarrow G(H \Box K)F:
\C \rightarrow \D$$ has components in $c$, the arrow 
$G^2_{HFc,KFc}: 
GHFc \otimes GKFc \rightarrow G(HFc \otimes KFc)$
in $\D$.\\
- The monoidal natural transformation
$${[F,G]}^0: \un_{[\C,\D]} \rightarrow 
G \circ \un_{[\A,\B]} \circ F: \C \rightarrow \D$$
has component in any $c$ in $\C$,
the arrow $G^0: \un \rightarrow G(\un_{\B})$
in $\D$.\\    

Also, one can define for any monoidal natural transformations
$\sigma: F \rightarrow F': \C \rightarrow \A$
and 
$\tau: G \rightarrow G': \B \rightarrow \D$,
the monoidal natural transformation
$$[\sigma,\tau]: [F,G] \rightarrow [F',G']: [\A,\B] \rightarrow [\C,\D]$$
as either of the 2-cells
\begin{tag}\label{homfundef2cell}
$[\C,\tau] *  [\sigma,\B]$ or $[\sigma,\D] * [\A,\tau]$.
\end{tag}

It follows from this definition
that for any symmetric monoidal functors
$F,F',F'': \C \rightarrow \A$ and
$G,G',G'': \B \rightarrow \D$, 
and any monoidal transformations
\begin{center}
$\xymatrix{ 
F \ar[r]^{\sigma} 
&
F'
\ar[r]^{\sigma'}
&
F''
}$,
and
$\xymatrix{
G 
\ar[r]^{\tau}
& 
G'
\ar[r]^{\tau'}
&
G''
}
$,
\end{center}
the composite 
$$\xymatrix{
[F,G]
\ar[r]^-{[\sigma,\tau]}
&
[F',G']
\ar[r]^-{[\sigma',\tau']}
&
[F'',G''] : [\A,\B] \rightarrow [\C,\D]
}$$
is $[\sigma' \sigma, \tau' \tau]$.\\

Also for any symmetric monoidal 
functors $F: \C \rightarrow \A$ 
and $G: \B \rightarrow \D$,
writing $1$ for the identity natural
transformations at $F$ and $G$,
since the monoidal natural transformations
$[1,\B]: [F,\B] \rightarrow [F,\B]: [\A,\B] \rightarrow [\C,\B]$
and 
$[\C,1]: [\C,G] \rightarrow [\C,G]: [\C,\B] \rightarrow [\C,\D]$  
are identities, one has that
$[1,1]: [F,G] \rightarrow [F,G]: 
[\A,\B] \rightarrow [\C,\D]$ is 
the identity monoidal natural transformation
at $[F,G]$.\\

This is to say the following.
\begin{lemma}\label{homisfun1}
Given any symmetric monoidal categories 
$\A$,$\B$,$\C$ and $\D$ the assignments $(F,G) \mapsto [F,G]$
as in \ref{homfundefmor} and 
$(\sigma, \tau) \mapsto [\sigma,\tau]$
as in \ref{homfundef2cell} define
a functor $\SMC(\C,\A) \times \SMC(\B,\D)
\rightarrow \SMC([\A,\B],[\C,\D])$.
\end{lemma}

\begin{lemma}\label{homisfun2}
The following statements hold in $\SMC$, for any $\A$.\\
- For any $\B$, the image 
of identity 1-cell $1_{\B}$ at $\B$ by $[\A,-]$
is the identity 1-cell at $[\A,\B]$.\\
- For any 1-cells 
$F: \B \rightarrow \C$ and $G: \C \rightarrow \D$,
the 1-cells 
$[\A, G \circ F]$ and $[\A,G] \circ [\A,F]: [\A,\B] \rightarrow
[\A,\D]$ are equal.\\
- For any 1-cells
$\sigma: F \rightarrow F': \B \rightarrow \C$ 
and $\tau: G \rightarrow G': \C \rightarrow \D$,
the image by ${[\A,-]}_{\B,\D}$ of 
the composite of the 2-cell
$$
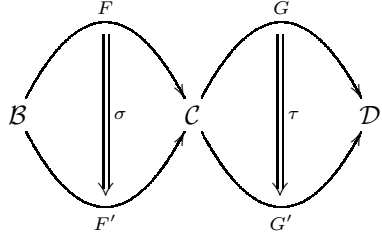
\xymatrix{
&
\ar@{=>}[dd]^{\sigma}
&
&
\ar@{=>}[dd]^{\tau}
\\
\B
\ar@/^35pt/[rr]^{F}
\ar@/_35pt/[rr]_{F'}
&
&
\C
\ar@/^35pt/[rr]^{G}
\ar@/_35pt/[rr]_{G'}
&
&
\D
\\
&&&
}
$$
is the composite of the 2-cells
$$
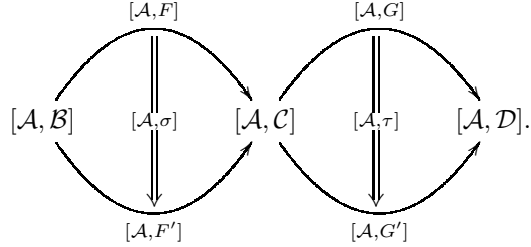
\xymatrix{
&
\ar@{=>}[dd]|{[\A,\sigma]}
& 
& 
\ar@{=>}[dd]|{[\A,\tau]}
\\
[\A,\B]
\ar@/^35pt/[rr]^{[\A,F]}
\ar@/_35pt/[rr]_{[\A,F']}
&
&
[\A,\C]
\ar@/^35pt/[rr]^{[\A,G]}
\ar@/_35pt/[rr]_{[\A,G']}
&
&
[\A,\D]
.
\\
&&&&}
$$ 
\end{lemma}
\pf
The functor $\SMC(\A,1_{\B})$ is the identity  
at the category $\SMC(\A,\B)$ and the monoidal
structure $[\A,1_\B]$ on this functor 
is strict since the identity $1_{\B}$ is. 
Therefore $[\A,1_{\B}]$ is
the identity in $\SMC$ at $[\A,\B]$.\\

For the second assertion,
since $\SMC$ is a 2-category, the
mere functors 
$\SMC(\A, G \circ F)$ and 
$\SMC(\A,G) \circ \SMC(\A,F): \SMC(\A,\B) \rightarrow
\SMC(\A,\D)$ are equal. It remains to show 
the equality of the monoidal structures
of $[\A, G \circ F]$ and $[\A,G] \circ [\A,F]$.
This is again straightforward.
For any symmetric monoidal functors
$H,K: \A \rightarrow \B$,
the monoidal natural transformations 
${([\A,G] \circ [\A,F])}^2_{H,K}$
and 
${([\A,GF])}^2_{H,K}: GFH \Box GFK \rightarrow 
GF(H \Box K): \A \rightarrow \D$ are equal,
with components in any object $a$ in $\A$, the arrow
$$\xymatrix@C=3pc{
GFHa \otimes GFKa  
\ar[r]^-{G^2_{FHa,FKa}}
&
G(FHa \otimes FKa)
\ar[r]^-{G(F^2_{Ha,Ka})}
&
G(F(Ha \otimes Ka)) 
}$$ in $\D$.\\ 
The monoidal natural transformations
${( [\A,G] \circ [\A,F] )}^0$ and 
${([\A,GF])}^0: \un_{[\A,\D]} \rightarrow GF\un_{[\A,\B]}$ are equal, 
with components in any object $a$ of $\A$, the arrow
$$
\xymatrix{
\un
\ar[r]^-{G^0}
&
G(\un_{\C})
\ar[r]^-{G(F^0)}
&
GF(\un_{\B})
}
$$
in $\D$.\\
  
The last assertion results from the fact that 
$\SMC$ is a 2-category and therefore
for any 2-cells 
$\sigma: F \rightarrow F'$ and $\tau: G \rightarrow G'$ 
in $\SMC$, the image by $\SMC(\A,-)$ of 
the composite 2-cell
$$
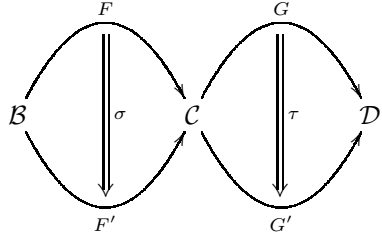
\xymatrix{
& 
\ar@{=>}[dd]^{\sigma}
&
& 
\ar@{=>}[dd]^{\tau}
\\
\B
\ar@/^35pt/[rr]^{F}
\ar@/_35pt/[rr]_{F'}
&
&
\C
\ar@/^35pt/[rr]^{G}
\ar@/_35pt/[rr]_{G'}
&
&
\D
\\
&&&&
}
$$
is the composite of the 2-cells in $\CAT$
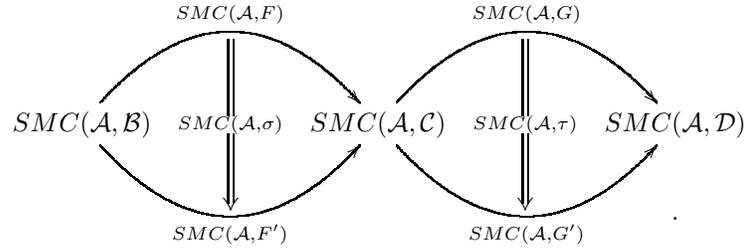
$$
\xymatrix{
& \ar@{=>}[dd]|{\SMC(\A,\sigma)} 
& & 
\ar@{=>}[dd]|{\SMC(\A,\tau)}  
\\
\SMC(\A,\B)
\ar@/^35pt/[rr]^{\SMC(\A,F)}
\ar@/_35pt/[rr]_{\SMC(\A,F')}
&
&
\SMC(\A,\C)
\ar@/^35pt/[rr]^{\SMC(\A,G)}
\ar@/_35pt/[rr]_{\SMC(\A,G')}
&
&
\SMC(\A,\D)
\\
& & & & 
.}
$$ 
\epf

One has a similar lemma for the 
assignments $[-,\A]$.
\begin{lemma}\label{homisfun3}
The following statements hold in $\SMC$ for any $\A$.\\
- For any $\B$, the 1-cell
$[1_{\B},\A]$ is the identity 1-cell at $[\B,\A]$.\\
- For any 1-cells
$F: \B \rightarrow \C$ and $G: \C \rightarrow \D$,
the 1-cells 
$[G \circ F,\A]$ and $[F,\A] \circ [G,\A]: [\D,\A] \rightarrow
[\B,\A]$ are equal.\\
- For any 2-cells 
$\sigma: F \rightarrow F': \B \rightarrow \C$ 
and $\tau: G \rightarrow G': \C \rightarrow \D$,
the image by ${[-,\A]}_{\B,\D}$ of 
the composite of the 2-cell
$$
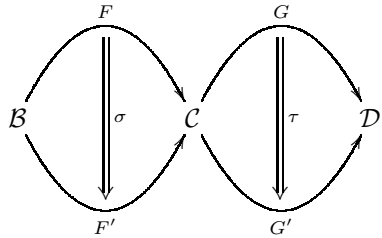
\xymatrix{
&
\ar@{=>}[dd]^{\sigma}
& & 
\ar@{=>}[dd]^{\tau}
\\
\B
\ar@/^35pt/[rr]^{F}
\ar@/_35pt/[rr]_{F'}
&
&
\C
\ar@/^35pt/[rr]^{G}
\ar@/_35pt/[rr]_{G'}
&
&
\D
\\
&&&&
}
$$
is the composite of the 2-cells
$$
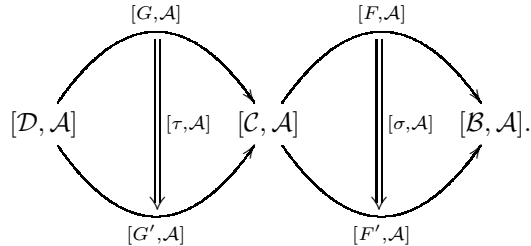
\xymatrix{
& 
\ar@{=>}[dd]^{[\tau,\A]}
& 
&
\ar@{=>}[dd]^{[\sigma,\A]}
\\
[\D,\A]
\ar@/^35pt/[rr]^{[G,\A]}
\ar@/_35pt/[rr]_{[G',\A]}
&
&
[\C,\A]
\ar@/^35pt/[rr]^{[F,\A]}
\ar@/_35pt/[rr]_{[F',\A]}
&
&
[\B,\A].
\\
&&&&
}
$$ 
\end{lemma}
The proof is similar to the one for the previous 
lemma and even simpler since all the functors
considered for the second assertion are strict.\\

\begin{proposition}\label{hom2fun}
One has a 2-functor 
$\homSMC: \SMC^{op} \times \SMC \rightarrow \SMC$
sending:\\
- any pair $(\A,\B)$ of symmetric monoidal categories
to $[\A,\B]$;\\
- any pair of symmetric monoidal functors
$(F:\C \rightarrow \A,G:\B \rightarrow \D$) 
to $[F,G]: [\A,\B] \rightarrow [\C,\D]$;\\
- any pair of 2-cells of $\SMC$
$$(\sigma: F \rightarrow F':\C \rightarrow \A, 
\tau:G \rightarrow G':\B \rightarrow \D)$$
to $[\sigma,\tau]: [F,G] \rightarrow [F',G']: [\A,\B] 
\rightarrow [\C,\D]$.
\end{proposition}
\pf
One needs to check the following three assertions.\\
1. Given any $\A$, $\B$, $\C$ and $\D$
in $\SMC$, the above assignments 
define a functor $\SMC(\C,\A) \times \SMC(\B,\D)
\rightarrow \SMC([\A,\B],[\C,\D])$.\\
2. Given any 2-cells in $\SMC$
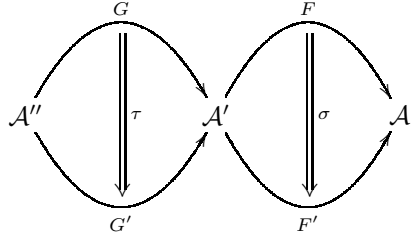
$$\xymatrix{
& \ar@{=>}[dd]^{\tau}
& 
& 
\ar@{=>}[dd]^{\sigma}
\\
\A''
\ar@/^35pt/[rr]^G
\ar@/_35pt/[rr]_{G'}
&
&   
\A'
\ar@/^35pt/[rr]^F
\ar@/_35pt/[rr]_{F'}
&
&   
\A 
\\
&&&&
}$$
and 
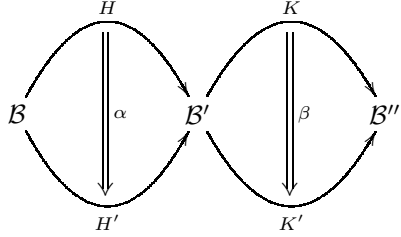
$$\xymatrix{
& 
\ar@{=>}[dd]^{\alpha}
& 
&
\ar@{=>}[dd]^{\beta}
\\
\B
\ar@/^35pt/[rr]^H
\ar@/_35pt/[rr]_{H'}
& 
&  
\B'
\ar@/^35pt/[rr]^K
\ar@/_35pt/[rr]_{K'}
&
&   
\B''
\\
&&&& 
}$$
the composite 2-cell in $\SMC$
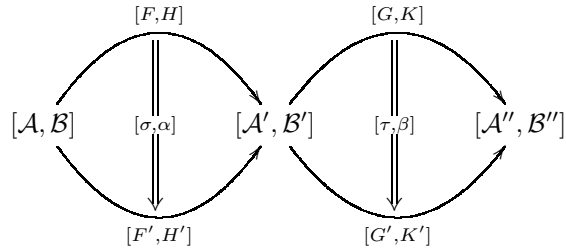
$$\xymatrix{
&
\ar@{=>}[dd]|{[\sigma, \alpha]}
& &
\ar@{=>}[dd]|{[\tau,\beta]}
\\
[\A,\B]
\ar@/^35pt/[rr]^{[F,H]}
\ar@/_35pt/[rr]_{[F',H']}
&   
&
[\A',\B']
\ar@/^35pt/[rr]^{[G,K]}
\ar@/_35pt/[rr]_{[G',K']}
&
&   
[\A'',\B''] 
\\
&&&&
}$$
is the 2-cell
$$
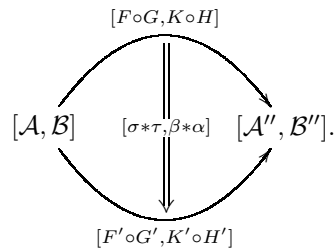
\xymatrix{
& 
\ar@{=>}@<0.9ex>[dd]|{[\sigma * \tau, \beta * \alpha]}
\\
[\A,\B]
\ar@/^35pt/[rr]^{[F \circ G, K \circ H]}
\ar@/_35pt/[rr]_{[F' \circ G',K' \circ H']}
&   
&
[\A'',\B''].
\\
& & 
} 
$$
3. Given any $\A$,
and $\B$, the 1-cell $[1_{\A},1_{\B}]$ is the identity 
at $[\A,\B]$ in $\SMC$.\\ 

\noindent Assertion 1. is Lemma \ref{homisfun1}.\\
2. is straightforward from Lemma \ref{homsqcom},
and the second and third assertions 
of Lemmas \ref{homisfun2} and \ref{homisfun3}.\\
3. is straightforward from Lemma \ref{homsqcom}
and the first assertions of Lemmas \ref{homisfun2}
and \ref{homisfun3}.
\epf

We mention here the following that is a consequence the 
2-functoriality of $\homSMC(\A,-): \SMC \rightarrow \SMC$ for any
symmetric monoidal category $\A$ and Lemma \ref{natcompfun}.
(Transpose \ref{natcompfun} to the case $\V = \CAT$ 
and $F = \homSMC(\A,-):\SMC \rightarrow \SMC$ for any symmetric 
monoidal category $\A$, to obtain the result.) 
\begin{corollary}\label{Post2natinB}
For any $\A$ in $\SMC$,
the collection of mere functors 
$$\Post_{\A,\B,\C}: 
\SMC(\B,\C) \rightarrow \SMC([\A,\B],[\A,\C])$$
for all $\B$ and $\C$, 
defines a 2-natural transformation 
between 2-functors $\SMC^{op} \times \SMC \rightarrow \CAT$.  
\end{corollary}

We can improve further this result.
\begin{corollary}\label{impronatPost2}
For any $\A$ and $\C$, 
the collection of 1-cells 
$${[\A,-]}_{\B,\C}: [\B,\C] \rightarrow [[\A,\B],[\A,\C]]$$ 
in $\SMC$, for all $\B$,
defines a 2-natural transformation between functors
$$\homSMC(-,\C) \rightarrow \homSMC(-, [\A,\C]) \circ \homSMC(\A,-)
:\SMC^{op} \rightarrow \SMC.$$ 
\end{corollary}
\pf
According to Corollary \ref{Post2natinB},
we only need to check that for any symmetric monoidal 
$F: \B' \rightarrow \B$ the following diagram
in $\SMC$ commutes
$$
\xymatrix{
[\B,\C] 
\ar[d]_{[F,\C]}
\ar[r]^-{[\A,-]}
&
[[\A,\B],[\A,\C]]
\ar[d]^{[[\A,F],1]}
\\
[\B',\C]
\ar[r]_-{[\A,-]}
&
[[\A,\B'],[\A,\C]],
}$$
but this is immediate from Corollary \ref{Post2natinB} since
all the monoidal functors involved here are strict.
\epf

\begin{lemma}\label{impronatPost}
For any $\A$ and $\B$, 
the collection of 1-cells 
${[\A,-]}_{\B,\C}: [\B,\C] \rightarrow [[\A,\B],[\A,\C]]$ in $\SMC$,
for all $\C$,
defines a 2-natural transformation between functors
$$\homSMC(\B,-) \rightarrow \homSMC([\A,\B], -) \circ \homSMC(\A,-)
:\SMC \rightarrow \SMC.$$ 
\end{lemma}
\pf 
According to Corollary \ref{Post2natinB}, 
it remains to show that for any symmetric monoidal 
functor $F: \C \rightarrow \C'$, the following diagram 
in $\SMC$ commutes
$$
\xymatrix{
[\B,\C]
\ar[r]^-{[\A,-]}
\ar[d]_{[\B,F]}
&
[[\A,\B],[\A,\C]]
\ar[d]^{[1,[\A,F]]}
\\
[\B,\C']
\ar[r]_-{[\A,-]}
&
[[\A,\B],[\A,\C']]
}
$$
and for this, it is enough to check
that the monoidal structures of 
the two legs above
are the same.\\

The monoidal natural transformations considered below, namely:\\
${( [1,[\A,F]] \circ [\A,-] )}^0$,
${([\A,-] \circ [\B,F])}^0$,
${{([1,[\A,F]] \circ [\A,-])}^2}_{H,K}$
and, 
${([\A,-] \circ [\B,F])}^2_{H,K}$,\\ are
between symmetric monoidal functors
$[\A,\B] \rightarrow [\A,\C'].$\\ 

The monoidal natural transformation 
${( [1,[\A,F]] \circ [\A,-] )}^0$ rewrites
$$  
\xymatrix@C=3pc{
\un 
\ar[r]^-{{[1,[\A,F]]}^0} 
&
[\A,F] \circ \un 
\ar@{=}[r] 
&
[\A,F] \circ [\A,\un]
}
$$ 
and has component in any $G:\A \rightarrow \B$,
the arrow in $[\A,\C']$
$$\xymatrix{
\un 
\ar[r]^-{[\A,F]^0}
&
F \circ \un
\ar@{=}[r]
&
F \circ \un \circ G
}$$
which has component in any $a$ in $\A$,
the arrow $F^0: \un \rightarrow F(\un)$ of $\C'$.\\

The monoidal natural transformation ${([\A,-] \circ [\B,F])}^0$ 
rewrites
$$
\xymatrix@C=3pc{
\un 
\ar@{=}[r]
&
[\A,\un]
\ar[r]^-{[\A,{[\B,F]}^0]}
&
[\A,F \circ \un]
}
$$
that has component in any $G: \A \rightarrow \B$,
the arrow of $[\A,\C']$
$$\xymatrix@C=3pc{
\un 
\ar@{=}[r]
&
\un \circ G
\ar[r]^-{{[\B,F]}^0 * G}
&
F \circ \un \circ G
}$$
which has component in any $a$ in $\A$,
the arrow $F^0: \un \rightarrow F(\un)$ of $\C'$.\\ 

Now consider any symmetric monoidal functors 
$H,K: \B \rightarrow \C$.\\

The monoidal natural transformation
${{([1,[\A,F]] \circ [\A,-])}^2}_{H,K}$
rewrites
{\tiny
$$\xymatrix@C=5pc{
([\A,F] \circ [\A,H]) \Box ([\A,F] \circ [\A,K])
\ar[r]^-{{[1,[\A,F]]}^2_{[\A,H],[\A,K]}}
&
[\A,F] \circ ([\A,H] \Box [\A,K])
\ar@{=}[r]
&
[\A,F] \circ [\A,H \Box K]
}$$
}
and has component in any $G: \A \rightarrow \B$,
the arrow of $[\A,\C']$
$$\xymatrix@C=4pc{
FHG \Box FKG 
\ar[r]^-{{[\A,F]}^2_{HG,KG}} 
&
F(HG \Box KG)
\ar@{=}[r]
&
F(H \Box K)G
}$$
that has component in any $a$ in $\A$, the arrow 
$$F^2_{HGa, KGa}: FHGa \otimes FKGa  \rightarrow F(HGa \otimes KGa)$$ 
of $\C'$.\\

The monoidal natural transformation 
${([\A,-] \circ [\B,F])}^2_{H,K}$ rewrites
$$\xymatrix@C=4pc{
[\A,FH] \Box [\A,FK]
\ar@{=}[r]
&
[\A,FH \Box FK]
\ar[r]^-{[\A,{[\B,F]}^2_{H,K}]}
&
[\A,F( H \Box K)]
}$$
that has component in any $G: \A \rightarrow \B$,
the arrow of $[\A,\C']$
$$
\xymatrix@C=3pc{
FHG \Box FKG
\ar@{=}[r]
&
(FH \Box FK)G
\ar[r]^-{  {[\B,F]}^2_{H,K} * G  }
&
F(H \Box K)G
}
$$ 
which has component in any $a$ in $\A$,
the arrow 
$$F^2_{HGa,KGa}:FHGa \otimes FKGa \rightarrow F(HGa \otimes KGa)$$
of $\C'$.
\epf

\begin{lemma}\label{impronatPost3}
The 2-cells in $\SMC$ 
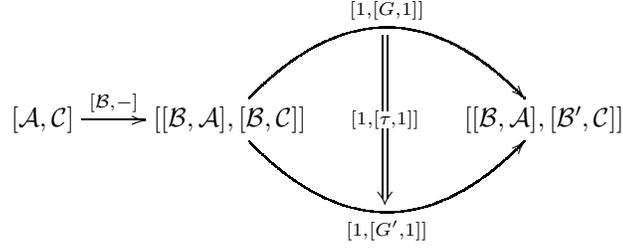
$$\xymatrix{
& & 
\ar@{=>}[dd]|{[1,[\tau,1]]}
\\
[\A,\C]
\ar[r]^-{[\B,-]}
&
[[\B,\A],[\B,\C]]
\ar@/^35pt/[rr]^{[1,[G,1]]}
\ar@/_35pt/[rr]_{[1,[G',1]]}
&
&
[[\B,\A],[\B',\C]]
\\
& & 
}$$
and
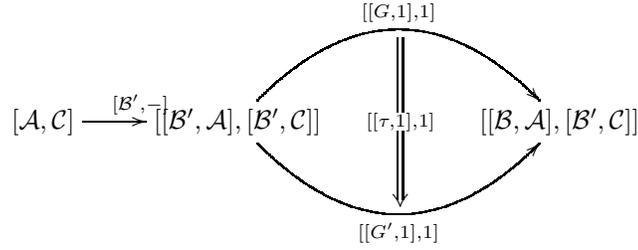
$$
\xymatrix{
& & 
\ar@{=>}[dd]|{[[\tau,1],1]}
\\
[\A,\C]
\ar[r]^{[\B',-]}
&
[[\B',\A],[\B',\C]]
\ar@/^35pt/[rr]^{[[G, 1],1]}
\ar@/_35pt/[rr]_{[[G', 1],1]}
&
&
[[\B,\A],[\B',\C]]
\\
& & 
}
$$
are equal, for any $\A$ and $\C$ and any 2-cell
$\tau: G \rightarrow G': \B' \rightarrow \B$.
\end{lemma}
\pf 
That the underlying 2-cells in $\CAT$ are equal,
is equivalent to the fact that for any 
monoidal natural transformation 
$\sigma: F \rightarrow F': \A \rightarrow \C$
the composite 2-cells in $\SMC$
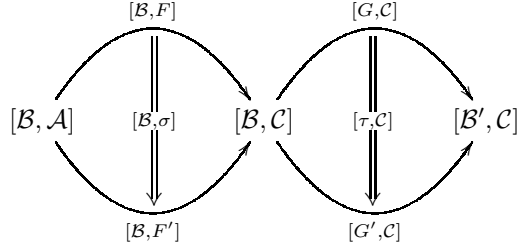
$$
\xymatrix{
&
\ar@{=>}[dd]|{[\B,\sigma]}
& 
& 
\ar@{=>}[dd]|{[\tau,\C]}
\\
[\B,\A]  
\ar@/^35pt/[rr]^{[\B,F]}
\ar@/_35pt/[rr]_{[\B,F']}
&
&
[\B,\C]
\ar@/^35pt/[rr]^{[G,\C]}
\ar@/_35pt/[rr]_{[G',\C]}
&
&
[\B',\C]
\\
&&&&
}
$$
and
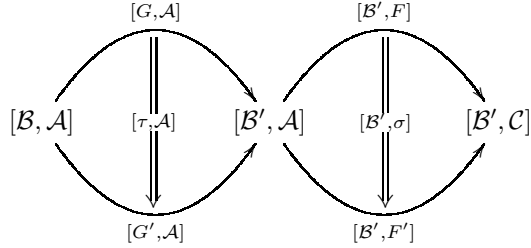
$$
\xymatrix{
& 
\ar@{=>}[dd]|{[\tau,\A]}
&&
\ar@{=>}[dd]|{[\B',\sigma]}
\\
[\B,\A]  
\ar@/^35pt/[rr]^{[G,\A]}
\ar@/_35pt/[rr]_{[G',\A]}
& &
[\B',\A]
\ar@/^35pt/[rr]^{[\B',F]}
\ar@/_35pt/[rr]_{[\B',F']}
& &
[\B',\C]
\\
&&&&
}
$$
are equal,
which is the case according to the 2-functoriality of
$\homSMC$.
To finish proving the lemma it remains to see that 
the diagram in $\SMC$ 
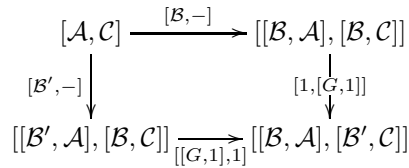
$$
\xymatrix{
[\A,\C]
\ar[r]^-{[\B,-]}
\ar[d]_{[\B',-]}
&
[[\B,\A],[\B,\C]]
\ar[d]|{[1,[G,1]]}
\\
[[\B',\A],[\B,\C]]
\ar[r]_-{[[G,1],1]}
&
[[\B,\A],[\B',\C]]
}
$$
commutes, but this is the case since all
the symmetric monoidal functors considered
here are strict. 
\epf
 
\end{section}

\begin{section}{Naturalities of the isomorphism $\Dual$}
\label{dualitynatsec}
This section tackles the question of the naturalities
of the isomorphism $\Dual$.\\

\begin{lemma}\label{dualitynatABC0}
The collection of isomorphisms defined in \ref{dualhom}
defines a 2-natural transformation
$$\SMC(\A, [\B,\C]) \cong_{\A,\B,\C} \SMC(\B, [\A, \C])$$
between 2-functors $\SMC^{op} \times \SMC^{op} \times \SMC \rightarrow
\CAT$.
\end{lemma}
\pf
To check the 2-naturality in $\C$, 
we show that given any 2-cells
$\sigma: F \rightarrow F': \A \rightarrow [\B,\C]$
and 
$\tau: h \rightarrow h': \C \rightarrow \C'$ in $\SMC$,
the 2-cells in $\SMC$ 
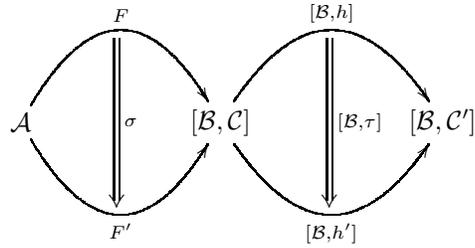
$$
\xymatrix{
&
\ar@{=>}@<0.5ex>[dd]^{\sigma}
& 
&
\ar@{=>}[dd]^{[\B,\tau]}
\\
\A
\ar@/^35pt/[rr]^{F}
\ar@/_35pt/[rr]_{F'}
&
&
[\B,\C]
\ar@/^35pt/[rr]^{[\B,h]}
\ar@/_35pt/[rr]_{[\B,h']}
&
&
[\B,\C']
\\
&&&&
}
$$
and 
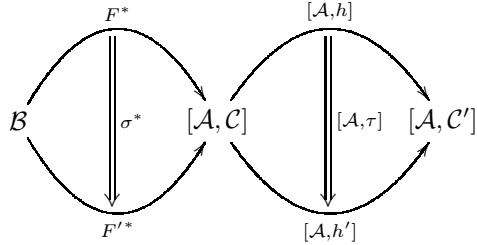
$$
\xymatrix{
&
\ar@{=>}@<0.5ex>[dd]^{\sigma^*} 
&
&
\ar@{=>}[dd]^{[\A,\tau]}
\\
\B
\ar@/^35pt/[rr]^{F^*}
\ar@/_35pt/[rr]_{{F'}^*}
&
&
[\A,\C]
\ar@/^35pt/[rr]^{[\A,h]}
\ar@/_35pt/[rr]_{[\A,h']}
&
&
[\A,\C']
\\
&&&&
}
$$
are equal.\\

Note that according to the 2-naturality in $\C$
of $U: \SMC(\B,\C) \rightarrow \CAT(\B,\C)$,
the mere natural transformation
{\small
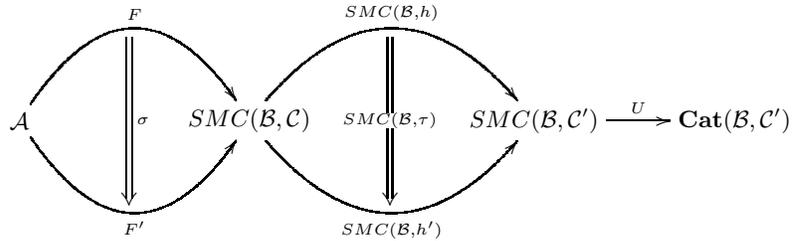
$$\xymatrix{
& 
\ar@<2ex>@{=>}[dd]^{\sigma}
& & 
\ar@{=>}[dd]|{\SMC(\B,\tau)}
\\
\A
 \ar@/^35pt/[rr]^{F}
 \ar@/_35pt/[rr]_{F'}
&
&
\SMC(\B,\C)
 \ar@/^35pt/[rr]^{\SMC(\B,h)}
 \ar@/_35pt/[rr]_{\SMC(\B,h')}
&
&
\SMC(\B,\C')
\ar[r]^-U
&
\CAT(\B,\C')
\\
& & & & &
}
$$
}
is 
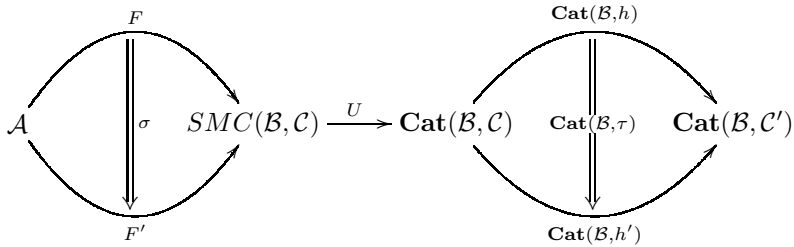
$$\xymatrix{
&
\ar@<2ex>@{=>}[dd]^{\sigma}
&
&
&
\ar@{=>}[dd]|{\CAT(\B,\tau)}
\\
\A
\ar@/^35pt/[rr]^{F}
\ar@/_35pt/[rr]_{F'}
&
&
\SMC(\B,\C)
\ar[r]^U
&
\CAT(\B,\C) 
\ar@/^35pt/[rr]^{\CAT(\B,h)}
\ar@/_35pt/[rr]_{\CAT(\B,h')}
&
&
\CAT(\B,\C')
\\
&& & & & 
}$$
which is
dual in the sense of \ref{classicdual} to
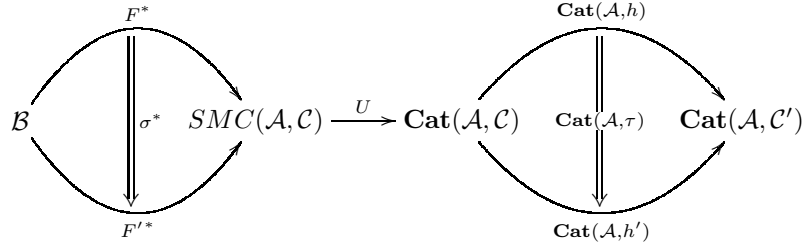
$$\xymatrix{
&
\ar@<2ex>@{=>}[dd]^{\sigma^*}
&  
& 
&
\ar@{=>}[dd]|{\CAT(\A,\tau)}
\\
\B
\ar@/^35pt/[rr]^{F^*}
\ar@/_35pt/[rr]_{{F'}^*}
&
&
\SMC(\A,\C)
\ar[r]^-U
& 
\CAT(\A,\C)
\ar@/^35pt/[rr]^{\CAT(\A,h)}
\ar@/_35pt/[rr]_{\CAT(\A,h')}
&
&
\CAT(\A,\C')
\\
& & & & & 
}$$ 
according to the 2-naturality 
of the isomorphism $\CAT(\A,\CAT(\B,\C)) \cong \CAT(\B,\CAT(\A,\C))$
in $\C$.
This last arrow is
$$
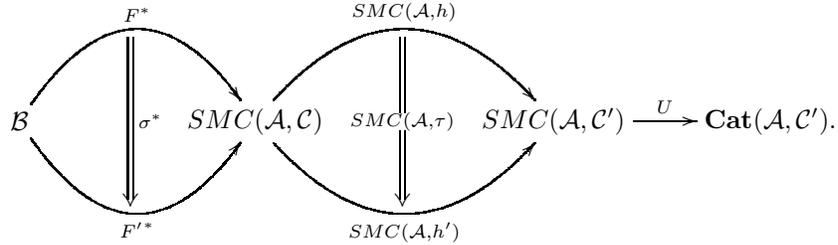
\xymatrix{
& \ar@<2ex>@{=>}[dd]^{\sigma^*} & & 
\ar@{=>}[dd]|{\SMC(\A,\tau)}
\\
\B
\ar@/^35pt/[rr]^{F^*}
\ar@/_35pt/[rr]_{{F'}^*}
&
&
\SMC(\A,\C)
\ar@/^35pt/[rr]^{\SMC(\A,h)}
\ar@/_35pt/[rr]_{\SMC(\A,h')}
&
&
\SMC(\A,\C')
\ar[r]^-U
& 
\CAT(\A,\C').\\
& & & & & 
}$$
Therefore to prove the equality of the above 2-cells $[\B, \tau] *
\sigma$ 
and $[\A,\tau] * \sigma^*$
of $\SMC$, it remains to check that $[\B,h] \circ F$ and 
$[\A,h] \circ F^*$
are dual as symmetric monoidal functors and for this it is 
enough to show the four following points.\\
$(1)$ For any objects $a$,$a'$ in $\A$ and $b$ in $\B$, the arrows
${ ( {([\B,h] \circ F)}^2_{a,a'} ) }_b$
and 
${ ( ( [\A,h] \circ F^* )(b) ) }^2_{a,a'}$ in $\C'$
are equal.\\
$(2)$ For any objects $b$,$b'$ in $\B$ and $a$ in $\A$, the arrows
${ ( {([\A,h] \circ F^*)}^2_{b,b'} ) }_a$
and 
${ (( [\B,h] \circ F )(a) ) }^2_{b,b'}$ in $\C'$
are equal.\\
$(3)$ For any object $b$ in $\B$,
the arrows
${ ([\B,h] \circ F ) }^0_b $ and ${( ([\A,h] \circ F^*)(b))}^0$
in $\C'$
are equal.\\
$(4)$ For any object $a$ in $\A$,
the arrows
${ ([\A,h] \circ F^* ) }^0_a $ and ${( ([\B,h] \circ F) (a))}^0$
in $\C'$ are equal.\\

We only check below points $(1)$ and $(3)$. The computations
for checking points $(2)$ and $(4)$ are the same
respectively as for $(1)$ and $(3)$, up to adequate
interchanges between $F$ and $F^*$ and between objects $a,a'...$
of $\A$ and objects $b, b'...$ of $\B$.\\  
 
To check $(1)$, consider any objects $a$, $a'$ in $\A$ and $b$ in $\B$. 
The arrow
$${([\B,h] \circ F)}^2_{a,a'}: 
([\B,h] \circ F)(a) \Box ([\B,h] \circ F)(a')
\rightarrow 
([\B,h] \circ F)(a \otimes a')$$ of $[\B,\C']$,
is the monoidal natural transformation  
$$\xymatrix@C=3pc{
(h \circ Fa) \Box (h \circ F(a'))
\ar[r]^-{{[\B,h]}^2_{Fa,Fa'}}
&
h \circ (Fa \Box Fa')
\ar[r]^-{h * F^2_{a,a'}}
&
h \circ F(a \otimes a'): \B \rightarrow \C'
}$$ 
which has component in $b$, 
the arrow
$$(i): \xymatrix@C=4pc{
h(Fa(b)) \otimes h(Fa'(b))
\ar[r]^-{h^2_{Fa(b),Fa'(b)}}
&
h(Fa(b) \otimes Fa'(b))
\ar[r]^-{ h ({(F^2_{a,a'})}_b) }
&
h(F(a \otimes a')(b))
}$$ in $\C'$.
On the other hand $([\A,h] \circ F^*)(b) = h \circ F^*b$ 
and the arrow $${(h \circ F^*b)}^2_{a,a'}: 
h(F^*b(a)) \otimes h(F^*b(a')) \rightarrow h(F^*b(a \otimes a'))$$
in $\C'$ is 
$$(ii): \xymatrix@C=4pc{
h(F^*b(a)) \otimes h(F^*b(a'))
\ar[r]^-{h^2_{F^*b(a), F^*b(a')}}
&
h( F^*b(a) \otimes F^*b(a') )
\ar[r]^-{h({(F^*b)}^2_{a,a'}) }
&
h(F^*b(a \otimes a')).
}$$ 
Since $F^*b(a) = Fa(b)$ and ${(F^*b)}^2_{a,a'} = {(F^2_{a,a'})}_b$,
arrows $(i)$ and $(ii)$ above are equal.\\

To check $(3)$, consider an arbitrary object $b$ in $\B$.
 The arrow
${([\B,h] \circ F)}^0$ of $[\B,\C']$
is the monoidal natural transformation
$$\xymatrix{
\un
\ar[r]^-{{[\B,h]}^0}
&
h \circ \un_{[\B,\C]}
\ar[r]^-{h * F^0}
&
h \circ (F(\un_{\A})): \B \rightarrow \C'
}$$
that has component in $b$
the arrow 
$$(i): \xymatrix{
\un_{\C'} 
\ar[r]^-{h^0}
&
h(\un_{\C})
\ar[r]^-{h(F^0_b)}
&
h(F(\un_{\A})(b))
}$$ in $\C'$. On the other hand,
$( [\A,h] \circ F^* )(b)$ is the arrow
${(h \circ F^*b)}^0 : \un \rightarrow (h \circ F^*b)(\un_{\A})$
in $\C'$, which is the composite
$$(ii):
\xymatrix@C=3pc{
\un_{\C'}
\ar[r]^{h^0}
&
h(\un_{\C})
\ar[r]^-{h( (F^*b)^0 )}
&
h(F^*b(\un_{\A})).
}
$$ 
Since ${(F^*b)}^0 = {(F^0)}_b$, arrows $(i)$ and $(ii)$ above are equal.\\

To prove the 2-naturality in $\B$, 
we show that given any 2-cells
$\sigma: F \rightarrow F': \A \rightarrow [\B,\C]$
and 
$\tau: g \rightarrow g': \B' \rightarrow \B$ in $\SMC$,
the 2-cells in $\SMC$ 
$$
\xymatrix{
& 
\ar@<0.75ex>@{=>}[dd]^{\sigma}
& 
&
\ar@{=>}[dd]^{[\tau,\C]}
\\
\A
\ar@/^35pt/[rr]^{F}
\ar@/_35pt/[rr]_{F'}
&
&
[\B,\C]
\ar@/^35pt/[rr]^{[g,\C]}
\ar@/_35pt/[rr]_{[g',\C]}
&
&
[\B',\C]
\\
& & & 
}
$$
and 
$$
\xymatrix{
& 
\ar@{=>}[dd]^{\tau}
& & 
\ar@<1ex>@{=>}[dd]^{\sigma^*}
\\
\B'
\ar@/^35pt/[rr]^{g}
\ar@/_35pt/[rr]_{g'}
&
&
\B
\ar@/^35pt/[rr]^{F^*}
\ar@/_35pt/[rr]_{F'^*}
&
&
[\A,\C].\\
& & & 
}
$$
are equal.\\

Note that the mere natural transformation
$$\xymatrix{
& 
\ar@{=>}@<2ex>[dd]^{\sigma}
&
&
\ar@{=>}[dd]|{\SMC(\tau,\C)}
\\
\A
 \ar@/^35pt/[rr]^{F}
 \ar@/_35pt/[rr]_{F'}
&
&
\SMC(\B,\C)
 \ar@/^35pt/[rr]^{\SMC(g,\C)}
 \ar@/_35pt/[rr]_{\SMC(g',\C)}
&
&
\SMC(\B',\C)
\ar[r]^-U
&
\CAT(\B',\C)
\\
& & & & 
}
$$
is 
$$\xymatrix{
& 
\ar@{=>}@<2ex>[dd]^{\sigma}
& 
& 
&
\ar@{=>}[dd]|{\CAT(\tau,\C)}
\\
\A
\ar@/^38pt/[rr]^{F}
\ar@/_35pt/[rr]_{F'}
&
&
\SMC(\B,\C)
\ar[r]^U
&
\CAT(\B,\C) 
\ar@/^35pt/[rr]^{\CAT(g,\C)}
\ar@/_35pt/[rr]_{\CAT(g',\C)}
&
&
\CAT(\B',\C)
\\
& & & & & 
}$$
which is dual in the sense of \ref{classicdual} to
$$
\xymatrix{
& 
\ar@{=>}[dd]^{\tau}
& 
&
\ar@{=>}@<2ex>[dd]^{\sigma^*}
\\
\B'
\ar@/^35pt/[rr]^{g}
\ar@/_35pt/[rr]_{g'}
&
&
\B
\ar@/^35pt/[rr]^{F^*}
\ar@/_35pt/[rr]_{{F'}^*}
&
&
\SMC(\A,\C)
\ar[r]^U
&
\CAT(\A,\C)
\\
& & & &
}
$$
according to the 2-naturality in $\B$ 
of the isomorphism $\CAT(\A,\CAT(\B,\C))  \cong \CAT(\B,\CAT(\A,\C))$.\\

Therefore to prove the equality 
of the above 2-cells $[\tau, \C] * \sigma$ and 
$\sigma^* * \tau$ of $\SMC$, it remains to check 
that $[g,\C] \circ F$ and $F^* \circ g$ are dual 
as symmetric monoidal functors 
and for this it is enough to show the    
following four points.\\
$(1)$ For any objects $a$,$a'$ in $\A$ and $x$ in $\B'$, the 
arrows
$ {( {([g,\C] \circ F)}^2_{a,a'} )}_x $ and
${((F^* \circ g)(x))}^2_{a,a'}$ in $\C$ are equal.\\
$(2)$ For any objects $a$ in $\A$ and $x$,$x'$ in $\B'$, the 
arrows
${(([g,\C] \circ F)(a))}^2_{x,x'}$
and ${( {(F^* \circ g)}^2_{x,x'} )}_a$ in $\C$ are equal.\\
$(3)$ For any object $x$ in $\B'$,
the arrows
${([g,\C] \circ F)}^0_x$ and 
${((F^* \circ g)(x))}^0$ in $\C$ are equal.\\
$(4)$ For any object $a$ in $\A$,
the arrows
${(([g,\C] \circ F)(a))}^0$ and 
${(F^* \circ g)}^0_a$ in $\C$ are equal.\\

To check point $(1)$, consider any objects $a$,$a'$ in $\A$,
and $x$ in $\B'$.
The arrow 
${([g,\C] \circ F)}^2_{a,a'}$ in $[\B',\C]$ is 
the monoidal natural transformation
$$\xymatrix{
(Fa \circ g) \Box (Fa' \circ g)
\ar@{=}[r]
&
(Fa \Box Fa') \circ g
\ar[r]^-{F^2_{a,a'} * g}
&
F(a \otimes a') \circ g : \B' \rightarrow \C
}
$$
since $[g,\C]$ is strict, and it has component in 
$x$ the arrow in $\C$
$${( F^2_{a,a'} )}_{g(x)} :  
Fa(g(x)) \otimes Fa'(g(x)) \rightarrow F(a \otimes a')(g(x))$$ 
(i.e. the component of 
$F^2_{a,a'}: Fa \Box Fa' \rightarrow F(a \otimes a'): \B \rightarrow
\C$
in $g(x)$).
This arrow is also ${(F^*(g(x)))}^2_{a,a'}$ since
$F$ and $F^*$ are dual.\\

To check point $(2)$, consider any objects $a$ in $\A$,
and $x$, $x'$ in $\B'$.
The arrow ${(( [g,\C] \circ F )(a))}^2_{x,x'}$ of $\C$ is 
$(Fa \circ g)^2_{x,x'}$ which is 
$$(i): \xymatrix@C=4pc{
Fa(g(x)) \otimes Fa(g(x')) 
\ar[r]^-{ {(Fa)}^2_{g(x),g(x')} } 
&
Fa( g(x) \otimes g(x') )
\ar[r]^-{ Fa(g^2_{x,x'}) }
&
Fa( g(x \otimes x'))
.}$$
On the other hand the arrow 
${(F^* \circ g)}^2_{x,x'}$ in $[\A,\C]$ is
the monoidal natural transformation
$$
\xymatrix@C=4pc{
F^*(g(x))
\Box 
F^*(g(x'))
\ar[r]^-{ {(F^*)}^2_{g(x),g(x')} }
&
F^*(g(x) \otimes g(x')) 
\ar[r]^-{ F^*(g^2_{x,x'})}
& 
F^*(g(x \otimes x')) : \A \rightarrow \C
}
$$
and has component in $a$
$$(ii):
\xymatrix@C=4pc{
F^*(g(x))(a)
\Box 
F^*(g(x'))(a)
\;
\ar[r]^-{ {({(F^*)}^2_{g(x),g(x')})}_a }
&
\;
F^*(g(x) \otimes g(x'))(a) 
\ar[r]^-{ {F^*(g^2_{x,x'})}_a  }
& 
F^*(g(x \otimes x'))(a)
.}
$$
Since $F$ and $F^*$ are dual, arrows $(i)$ and $(ii)$ above 
are equal.\\
 
Let us check point $(3)$. The arrow
${([g,\C] \circ F)}^0$ in $[\B',\C]$
is the monoidal natural transformation 
$$\xymatrix{
\un 
\ar@{=}[r]
&
\un_{[\B,\C]} \circ g
\ar[r]^-{F^0 * g}
&
F(\un_{\A}) \circ g : \B' \rightarrow \C
}$$ 
since $[g,\C]$ is strict.
For any object $x$ in $\B'$,
its component in $x$ is the arrow in $\C$
$${(F^0)}_{g(x)}: \un \rightarrow F(\un_{\A})(g(x))$$
that is also ${( F^*(g(x)) )}^0$ since 
$F$ and $F^*$ are dual.\\

Let us check point $(4)$. Consider any object $a$ of $\A$.
The arrow ${(([g,\C] \circ F)(a))}^0$  
of $\C$ is ${(Fa \circ g)}^0$ which is 
$$(i): \xymatrix{
\un
\ar[r]^-{{(Fa)}^0} 
&
Fa(\un_{\B})
\ar[r]^-{Fa(g^0)}
&
Fa(g(\un_{\B'}))
.}$$
On the other hand, the arrow 
${(F^* \circ g)}^0$ of $[\A,\C]$ is the monoidal natural transformation
$$
\xymatrix{
\un
\ar[r]^-{{(F^*)}^0}
&
F^*(\un_{\B})
\ar[r]^-{F^*(g^0)}
&
F^*(g(\un_{\B'})): \A \rightarrow \C
}
$$
and has component in $a$ 
$$
(ii): \xymatrix{
\un_{[\A,\C]}(a)
\ar[r]^-{{(F^*)}^0_a}
&
F^*(\un_{\B})(a)
\ar[r]^-{{F^*(g^0)}_a}
&
F^*(g(\un_{\B'}))(a)
.}
$$
Since $F$ and $F^*$ are dual, the arrows $(i)$ and $(ii)$ above are
equal.\\

The 2-naturality in $\A$ results from the 2-naturality in $\B$
since the isomorphisms \ref{dualhom} in $\CAT$,
$SMC(\A,[\B,\C]) \rightarrow \SMC(\B,[\A,\C])$ 
and $\SMC(\B,[\A,\C]) \rightarrow \SMC(\A,[\B,\C])$ are 
inverses in $\CAT$.  
\epf

Given any symmetric monoidal categories $\B$ and $\C$,
since $\Dual$ defines a 2-natural transformation in $\A$
between 2-functors $\SMC \rightarrow \CAT$
$$\SMC(\A,[\B,\C]) \rightarrow \SMC(\B,[\A,\C]),$$ 
the Yoneda Lemma (in its enriched version) indicates that 
the collection of arrows $\Dual$ in $\CAT$ above is determined 
by a unique object in $\SMC(\B, [[\B,\C],\C])$,
namely the image by $\Dual$ of the identity in $[\B,\C]$,
which we denote $\resev$. Precisely, one has the following. 
\begin{lemma}\label{chardualq}
For any symmetric monoidal categories
$\A$, $\B$ and $\C$, the diagram below 
between mere functors commute.
$$
\xymatrix{
\SMC(\A,[\B,\C]) 
\ar[rd]_{\homSMC(-,\C)}
\ar[rr]^{\Dual}
&
&
\SMC(\B,[\A,\C]).
\\
&
\SMC([[\B,\C],\C],[\A,\C])
\ar[ru]_{\SMC(q,1)}
&
}
$$
\end{lemma}
Since the collection of arrows 
$\Dual_{\A,\B,\C}: \SMC(\A,[\B,\C]) \rightarrow \SMC(\B,[\A,\C])$ 
is also 2-natural in the arguments $\B$ and $\C$, the collection
of arrows $q: \B \rightarrow [[\B,\C],\C]$ is also 2-natural
in the arguments $\B$ and $\C$. We will use this point in the next 
chapter.\\

We can improve slightly Lemma \ref{chardualq}.
\begin{lemma}\label{chardualq2}
The diagram in $\SMC$
$$
\xymatrix{
[\A,[\B,\C]] 
\ar[rd]_{[-,\C]}
\ar[rr]^{\Dual}
&
&
[\B,[\A,\C]]
\\
&
[[[\B,\C],\C],[\A,\C]]
\ar[ru]_{[q,1]}
&
}
$$ commutes for any $\A$, $\B$ and $\C$.
\end{lemma}
\pf
We check that the composite $[q,1] \circ [-,\C]$ is strict.\\

For any 1-cells $F,G: \A \rightarrow [\B,\C]$ in $\SMC$, 
the arrow 
${([q,1] \circ [-,\C])}^2_{F,G}$ in $[\B,[\A,\C]]$ is
$$\xymatrix@C=3pc{ 
([F,1] \circ q) \Box ([G,1] \circ q)
\ar@{=}[r]
&
([F,1] \Box [G,1]) \circ q
\ar[r]^-{{[-,\C]}^2_{F,G} * q}
&
[F \Box G,1] \circ q  
}$$ 
since $[q,1]$ is strict.
The above arrow 
$${[-,\C]}^2_{F,G}: [F,1] \Box [G,1] \rightarrow [F \Box G,1]$$ 
in $[[[\B,\C],\C],[\A,\C]]$
is the natural transformation with component 
in any $H: [\B,\C] \rightarrow \C$, the arrow in $[\A,\C]$
$$(H \circ F) \Box (H \circ G) 
\rightarrow H \circ (F \Box G)$$ which is the natural transformation
with  component in any $a$ of $\A$, the arrow 
$$H^2_{Fa, Ga}: HFa \otimes HGa \rightarrow H \circ (Fa \Box Ga)$$ in $\C$.
So that if the $H$ above is strict, the component in $H$ of  
${{[-,\C]}^2_{F,G}}$ is an identity.
Now for any $b$ in $\B$, the component in $b$ of 
the monoidal natural transformation 
${[-,\C]}^2_{F,G} * \resev$ is the component of ${[-,\C]}^2_{F,G}$ in 
the strict $\ev_b$ and thus it is an identity.\\

The arrow ${( [\resev,1] \circ [-,\C] )}^0$  
in $[\B,[\A,\C]]$ is 
$$
\xymatrix@C=3pc{
\un 
\ar@{=}[r]
&
\un_{[[[\B,\C],\C],[\A,\C]]} \circ \resev
\ar[r]^{ {[-,\C]}^0 * \resev }
&
[\un_{[\A,[\B,\C]]}, \C]\circ \resev
.}
$$
The above arrow ${[-,\C]}^0: \un \rightarrow [ \un_{[\A,[\B,\C]]},\C ]$
in $[  [[\B,\C],\C], [\A,\C]   ]$ is the natural transformation 
with component in any $H: [\B,\C] \rightarrow \C$, 
the arrow $\un \rightarrow H \circ \un_{[\A,[\B,\C]]}$ 
of $[\A,\C]$ which is the natural transformation 
with component in any $a$ of $\A$ the arrow 
$H^0: \un \rightarrow H(\un_{[\B,\C]})$ in $\C$.
So that if $H$ as above is strict, the component in $H$
of ${[-,\C]}^0$ is an identity.
Now for any $b$ in $\B$, the component in $b$ of 
the monoidal natural transformation 
${[-,\C]}^0 * \resev$ is the component of ${[-,\C]}^0$ in 
the strict $\ev_b$ and thus it is an identity.\\
\epf

Lemma \ref{dualitynatABC0} also implies in particular
points $(1)$,$(2)$ and $(3)$ below,
which we will often use.\\
\noindent$(1)$ Given any symmetric monoidal functors
$F: \A \rightarrow [\B,\C]$
and $G: \C \rightarrow \C'$,
the following diagram in $\SMC$ is
commutative
$$\xymatrix{
\B 
\ar[r]^{F^*}
\ar[rd]_{{([\B,G] \circ F)}^*}
&
[\A,\C]
\ar[d]^{[\A,G]}
\\
&
[\A,\C']
.}$$
\noindent$(2)$
Given any symmetric monoidal functors $F: \A \rightarrow [\B,\C]$
and $G: \B' \rightarrow \B$,
the following diagram in $\SMC$ is
commutative
$$\xymatrix{
\B' 
\ar[r]^{G}
\ar[rd]_{{([G,\C] \circ F)}^*}
&
\B
\ar[d]^{F^*}
\\
&
[\A,\C]
.}$$
\noindent$(3)$ Given any symmetric monoidal functors
$F: \A \rightarrow [\B,\C]$
and $G: \A' \rightarrow \A$, the following diagram in $\SMC$ is
commutative
$$\xymatrix{
\B 
\ar[r]^{F^*}
\ar[rd]_{{(F \circ G)}^*}
&
[\A,\C]
\ar[d]^{[G,\C]}
\\
&
[\A',\C]
.}$$ 

\begin{proposition}\label{dualitynatABC}
The collection of isomorphisms defined in \ref{dualhom}
defines a 2-natural transformation
between functors $\SMC^{op} \times \SMC^{op} \times \SMC \rightarrow
\SMC$:
$$[\A, [\B,\C]] \cong [\B, [\A, \C]].$$
\end{proposition}
\pf
For the 2-naturality in $\C$, according to Lemma
\ref{dualitynatABC0},
one just needs to check that for any arrow
$f: \C \rightarrow \C'$ the following
diagram in $\SMC$ commutes
$$
\xymatrix{
[\A,[\B,\C]] 
\ar@{-}[r]^{\cong}
\ar[d]_{[\A,[\B,f]]}
&
[\B,[\A,\C]]
\ar[d]^{[\B,[\A,f]]}
\\
[\A,[\B,\C']] 
\ar@{-}[r]_{\cong}
&
[\B,[\A,\C']]
}
$$
and for this, it just remains to show that 
the monoidal structures of 
$\Dual \circ [\A,[\B,f]]$ and $[\B,[\A,f]] \circ \Dual$
are the same.\\

Since the monoidal functor $\Dual$ is strict,
for any symmetric monoidal functors
$F,G: \A \rightarrow [\B,\C]$, 
the equality of the arrows
${(\Dual \circ [\A,[\B,f]])}^2_{F,G}$ 
and 
${( [\B,[\A,f]] \circ \Dual )}^2_{F,G}$ in $[\B,[\A,\C']]$ 
is equivalent to the fact that 
the monoidal natural transformation 
$${{[\B,[\A,f]]}^2}_{F^*,{G}^*}:
([\A,f] \circ F^*) \Box ([\A,f] \circ G^*)
\rightarrow [\A,f] \circ (F^* \Box G^*)
: \B \rightarrow [\A,\C']$$
is dual 
to
$${[\A,[\B,f]]}^2_{F,G}: 
([\B,f] \circ F) \Box ([\B,f] \circ G)
\rightarrow [\B,f] \circ (F \Box G)
: \A \rightarrow [\B,\C'].$$
Consider any objects $a$ in $\A$ and $b$ in $\B$.
The component in $a$ 
of  ${([\A,[\B,f]])}^2_{F,G}$
is the monoidal natural transformation
$${[\B,f]}^2_{Fa,Ga}:  (f \circ Fa) \Box (f \circ Ga)
\rightarrow f \circ (Fa \Box Ga): \B \rightarrow \C'$$
that has component in $b$, the arrow 
$$f^2_{Fa(b),Ga(b)}: f(Fa(b)) \otimes f(Ga(b))
\rightarrow f(Fa(b) \otimes Ga(b))$$
in $\C'$.
The component in $b$ of 
${[\B,[\A,f]]}^2_{F^*,G^*}$ is the monoidal natural 
transformation 
$${[\A,f]}^2_{F^*b,G^*b}:
(f \circ F^*b) \Box (f \circ G^*b) 
\rightarrow
f \circ (F^*b \Box G^*b): \A \rightarrow \C'
$$
that has component in $a$ the arrow 
$$f^2_{F^*b(a),G^*b(a)}:
f (F^*b(a)) \otimes f (G^*b(a)) 
\rightarrow
f (F^*b(a) \otimes G^*b(a))
$$
in $\C'$.\\

Since $\Dual$ is strict,
the equality of the arrows 
${(\Dual \circ [\A,[\B,f]] )}^0$ and 
${ ( [\B,[\A,f]] \circ \Dual )}^0$ in $[\B,[\A,\C']]$ 
amounts to the fact the monoidal natural 
transformations
$$
{[\B,[\A,f]]}^0:
\un
\rightarrow
[\A,f] \circ \un_{[\B,[\A,\C]]}:
\B \rightarrow [\A,\C']
$$
and
$$
{[\A,[\B,f]]}^0:
\un
\rightarrow [\B,f] \circ \un_{[\A,[\B,\C]]} :
\A \rightarrow [\B,\C']
$$
are dual.
Consider any objects $a$ in $\A$ and $b$ in $\B$.
The component in $a$ of
${[\A,[\B,f]]}^0$ is 
the monoidal natural transformation
$
{[\B,f]}^0:
\un
\rightarrow
f \circ \un_{[\B,\C]}:
\B \rightarrow \C'
$,
with component in $b$ the arrow
$f^0: \un_{\C'} \rightarrow f(\un_{\C})$.
The component in $b$ of 
${[\B,[\A,f]]}^0$ is the monoidal natural transformation 
${[\A,f]}^0: \un \rightarrow 
f \circ \un_{[\A,\C]}: \A \rightarrow \C'$,
with component in $a$ the arrow
$f^0: \un_{\C'} \rightarrow f(\un_{\C})$.\\

The 2-naturality in $\A$  
of $\Dual_{\A,\B,\C}: [\A,[\B,\C]] \cong [\B,[\A,\C]]$ is
actually trivial from Lemma
\ref{dualitynatABC0}. To prove it,
one just needs to check that for any 
$f: \A' \rightarrow \A$, the following
diagram in $\SMC$ commutes
$$
\xymatrix{
[\A,[\B,\C]] 
\ar@{-}[r]^{\cong}
\ar[d]_{[f,[\B,\C]]}
&
[\B,[\A,\C]]
\ar[d]^{[\B,[f,\C]]}
\\
[\A',[\B,\C]] 
\ar@{-}[r]_{\cong}
&
[\B,[\A',\C]]
}
$$
and for this, it just remains to see that 
the monoidal structures of 
$\Dual \circ [f,[\B,\C]]$ and $[\B,[f,\C]] \circ \Dual$
are the same. This is the case since all 
the symmetric monoidal functors in the two composites
above are strict.\\

The 2-naturality in $\B$ of $\Dual_{\A,\B,\C}$
results from its 2-naturality in $\A$ 
since the isomorphisms $\Dual_{\A,\B,\C}$ and $\Dual_{\B,\A,\C}$
are inverses in $\SMC$. 
\epf

The naturalities of the isomorphisms $\Dual$ are further 
investigated and the end of this section 
contains results improving on Proposition 
\ref{dualitynatABC}.\\

\begin{lemma}\label{beurk}
The two arrows of $\SMC$
$$\xymatrix{
\A 
\ar[r]^-{F}
&
[\B, \C]
\ar[r]^-{[\D,-]}
&
[[\D,\B],[\D,\C]]
}$$
and 
$$\xymatrix{
[\D,\B]
\ar[r]^-{[\D,F^{\star}]}
&
[\D, [\A, \C]]
\ar[r]^-{\Dual}
&
[\A, [\D,\C]]
}
$$
are dual (via \ref{dualhom})
for any $F: \A \rightarrow [\B,\C]$ and any $\D$.
\end{lemma}
\pf
According to Lemma \ref{chardualq2}, the above arrow $\Dual \circ [\D,F^*]$ 
rewrites
$$
\xymatrix{
[\D,\B]
\ar[r]^-{[\D,F^*]}
&
[\D,[\A,\C]]
\ar[r]^-{[-,\C]}
&
[[[\A,\C],\C],[\D,\C]]
\ar[r]^-{[q,1]}
&
[\A,[\D,\C]]
}
$$
which, according to \ref{dualitynatABC0},
has dual
$$\xymatrix@C=3pc{
\A
\ar[r]^-{q}
&
[[\A,\C],\C]
\ar[r]^-{[\D,-]}
&
[[\D,[\A,\C]], [\D,\C]]
\ar[r]^-{[[1,F^*],1]}
&
[[\D,\B],[\D,\C]]
}$$
that rewrites
\begin{tabbing}
$\xymatrix{
\A
\ar[r]^-{q}
&
[[\A,\C],\C]
\ar[r]^-{[F^*,1]}
&
[\B,\C]
\ar[r]^-{[\D,-]}
&
[[\D,\B],[\D,\C]]
}$\\
$\xymatrix{
\A
\ar[r]^-{F}
&
[\B,\C]
\ar[r]^-{[\D,-]}
&
[[\D,\B],[\D,\C]]
}$
\end{tabbing}
according to Corollary \ref{impronatPost2} and Lemma \ref{chardualq}.
\epf


\begin{lemma}
For any  
$\A$, $\B$ and $\C$, the arrow 
$[-,\C]_{\A,\B}$ in $\SMC$ is equal 
to the composite
$$\xymatrix{
[\A,\B] 
\ar[r]^-{[1,\resev]}
&
[\A,[[\B,\C],\C]]
\ar[r]^-{D}
&
[[\B,\C],[\A,\C]].      
}$$ 
\end{lemma}
\pf
$[-,\C]_{\A,\B}$ has dual 
${[\A,-]}_{\B,\C}: [\B,\C] \rightarrow [[\A,\B],[\A,\C]]$.
Since $\resev: \B \rightarrow [[\B,\C],\C]$ is the dual
of $1: [\B,\C] \rightarrow [\B,\C]$, according to Lemma \ref{beurk},
the composite $\Dual \circ [1,\resev]$ has also dual
${[\A,-]}_{\B,\C}$.
\epf


\begin{lemma}\label{tbc1}
The following diagram in $\SMC$
is commutative
$$\xymatrix{
[\D,\C] 
\ar[r]^-{[-,\B]}
\ar[d]_{[1, F^*]}
&
[[\C, \B],[\D, \B]]
\ar[d]^{[F,1]}
\\
[\D, [\A, \B]]
\ar[r]_-{\Dual}
&
[\A,[\D, \B]]
}
$$
for any $F: \A \rightarrow [\C,\B]$ and any $\D$.
\end{lemma}
\pf 
According to Lemma \ref{beurk}, the arrow
$$\xymatrix{
[\D,\C]
\ar[r]^-{[\D,F^*]}
&
[\D,[\A,\B]]
\ar[r]^-{\Dual}
&
[\A,[\D,\B]]
}$$
has dual 
$$
\xymatrix{
\A
\ar[r]^-{F}
&
[\C,\B]
\ar[r]^-{[\D,-]}
&
[[\D,\C],[\D,\B]]
.}
$$
Since the dual of 
$$[-, \B]_{\D,\C}: [\D,\C] 
\rightarrow [[\C, \B],[\D, \B]]$$
is 
$$[\D,-]_{\C, \B}:
[\C, \B] \rightarrow
[[\D,\C],[\D, \B]],$$
according to Lemma \ref{dualitynatABC0},
the dual of 
$$\xymatrix{
[\D,\C]
\ar[r]^-{[-,\B]}
&
[[\C,\B],[\D, \B]]
\ar[r]^-{[F,1]}
&
[\A,[\D,\B]]
}$$ 
is also 
$$
\xymatrix{
\A 
\ar[r]^-{F}
&
[\C,\B]
\ar[r]^-{[\D,-]}
&
[[\D,\C],[\D,\B]]
.}
$$
\epf

\begin{lemma}\label{tbc3}
The following diagram in $\SMC$ commutes
$$
\xymatrix{
& 
[\B, \D]
\ar[ld]_{[\A,-]}
\ar[rd]^{[\C,-]}
\\
[[\A,\B],[\A, \D ]]
\ar[d]_{[F^*,1]}
& 
& 
[[\C, \B],[\C, \D]]
\ar[d]^{[F,1]}
\\
[\C,[\A, \D]]
\ar[rr]_{\Dual}
&
&
[\A ,[\C, \D]]
}$$
for any $F: \A \rightarrow [\C,\B]$ and any $\D$.
\end{lemma}
\pf
According to Lemma \ref{dualitynatABC0},
for any symmetric monoidal transformation
$\sigma: G \rightarrow H: \B \rightarrow \D$,
the composite in $\SMC$
$$\xymatrix{
\C
\ar[r]^-{F^*} 
&
[\A, \B]
\ar[r]^{[\A,G]}
&
[\A,\D]
}$$
is dual to
$$\xymatrix{
\A
\ar[r]^-{F}
&
[\C,\B]
\ar[r]^{[\C,G]}
&
[\C,\D] 
}$$
and the monoidal natural transformation
$$[\A,\sigma] * F^*: 
[\A,G] \circ F^* \rightarrow [\A,H] \circ F^*
$$
is dual to
$$
[\C,\sigma] * F:
[\C,G] \circ F \rightarrow [\C,H] \circ F. 
$$
Therefore the diagram in $\CAT$ (between mere functors)
is commutative
$$
\xymatrix{
&  
\SMC(\B,\D)  
\ar[ld]_{\homSMC(\A,-)}
\ar[rd]^{\homSMC(\C,-)}
& 
\\
\SMC([\A,\B],[\A,\D])
\ar[d]_{\SMC(F^*,1)} 
& 
& 
\SMC([\C, \B],[\C,\D])
\ar[d]^{\SMC(F,1)}
\\
[\C,[\A,\D]]
\ar[rr]_{\Dual}
& 
& 
[\A,[\C,\D]].
}
$$
Now since the symmetric monoidal functors
$[\A,-]$, $[\C,-]$, $[F^*,1]$, $[F,1]$
and $\Dual$
are all strict, the result follows.
\epf

\end{section}

\begin{section}{Evaluation functors}
\label{evalfunsec}
This section gathers technical lemmas
involving the arrows $\resev$, defined in Section \ref{dualitynatsec},
and the internal hom.\\

Consider any symmetric monoidal 
categories $\A$, $\B$ and any object $a$ of $\A$.
The image of $a$ by $q: \A \rightarrow [[\A,\B],\B]$ will
be written $\ev_a$. According to Remark \ref{Fsbstrict},
since the identity symmetric monoidal functor at 
$[\A,\B]$ is strict, $\ev_a$ is strict. 
The functor $\ev_a$ sends any symmetric monoidal functor 
$F:\A \rightarrow \B$ to $Fa$ and any monoidal natural 
transformation between symmetric 
functors $\sigma: F \rightarrow G: \A \rightarrow \B$ 
to its component $\sigma_a: Fa \rightarrow Ga$ in $a$.\\

The 2-naturality of the collection of arrows
$\resev: \A \rightarrow [[\A,\B],\B]$ in the argument $\B$
is equivalent to the Lemma below.
\begin{lemma}
The following holds in $\SMC$, for any object $\A$.\\
- For any $F:\B \rightarrow \C$ the diagram 
$$\xymatrix{
\A
\ar[r]^-{\resev}
\ar[d]_{\resev}
&
[[\A,\B],\B]
\ar[d]^{[1,F]}
\\
[[\A,\C],\C]
\ar[r]_-{[[1,F],1]}
&
[[\A,\B],\C]
}
$$
commutes.\\
- For any 2-cell $\sigma: F \rightarrow G:
\B \rightarrow \C$, the 2-cells 
$$\xymatrix{
&
&
\ar@{=>}[dd]^{[1,\sigma]}
\\
\A
\ar[r]^-{\resev}
&
[[\A,\B],\B]
\ar@/^35pt/[rr]^{[1,F]}
\ar@/_35pt/[rr]_{[1,G]}
&
&
[[\A,\B],\C]
\\
&&&
}
$$
and
$$
\xymatrix{
& & 
\ar@{=>}[dd]|{[[\A,\sigma],1]} 
\\
\A
\ar[r]^-{\resev}
&
[[\A,\C],\C]
\ar@/^35pt/[rr]^{[[\A,F],1]}
\ar@/_35pt/[rr]_{[[\A,G],1]}
&
&
[[\A,\B],\C]
\\
&&
}$$
are equal.
\end{lemma}

\begin{corollary}\label{2nateva}
For any object $a$ of a symmetric monoidal category $\A$,
the following holds in $\SMC$.\\
- For any $F: \B \rightarrow \C$, the diagram 
$$
\xymatrix{
[\A,\B] 
\ar[r]^-{\ev_a} 
\ar[d]_{[\A,F]}
&
\B
\ar[d]^F
\\
[\A,\C]
\ar[r]_-{\ev_a}
&
\C
}
$$
commutes.\\
- For any 2-cell $\sigma: F \rightarrow G:
\B \rightarrow \C$, the 2-cells  
$$
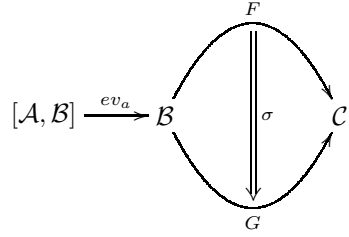
\xymatrix{
& 
\ar@{=>}[dd]^{[\A,\sigma]} 
\\
[\A,\B] 
\ar@/^35pt/[rr]^{[\A,F]}
\ar@/_35pt/[rr]_{[\A,G]}
&
&
[\A,\C]
\ar[r]^-{\ev_a}
&
\C
\\
&& 
}
$$
and
$$
\xymatrix{
& & 
\ar@{=>}[dd]^{\sigma}
\\
[\A,\B] 
\ar[r]^-{\ev_a}
&
\B 
\ar@/^35pt/[rr]^F
\ar@/_35pt/[rr]_G
& &
\C
\\
& & 
}
$$
are equal.
\end{corollary}

\begin{corollary}\label{homev3}
The diagram 
in $\SMC$ 
$$
\xymatrix{
[\B,\C] 
\ar[r]^-{[\A,-]}
\ar[rd]_{[\ev_a,1]}
&
[[\A,\B],[\A,\C]]
\ar[d]^{[1,\ev_a]}
\\
& 
[[\A,\B],\C]
}
$$
commutes for any $\A$, $\B$ and
$\C$, and any object $a$ of $\A$.
\end{corollary}
\pf
All the functors involved in the above diagram
are strict, therefore it is enough to show
that the diagram in $\CAT$ (between mere functors)
is commutative
$$
\xymatrix{
\SMC(\B,\C) 
\ar[r]^-{\Post}
\ar[rd]_{\SMC(\ev_a,1)}
&
\SMC([\A,\B],[\A,\C])
\ar[d]^{\SMC(1,\ev_a)}
\\
&
\SMC([\A,\B],\C)
.}
$$
This is exactly Corollary \ref{2nateva}.
\epf

\begin{corollary}\label{homev3dual}
The diagram in $\SMC$ 
$$
\xymatrix{
[\A,\B]
\ar[d]_{\ev_a}
\ar[r]^-{[-,\C]}
&
[[\B,\C],[\A,\C]]
\ar[d]^{[1,\ev_a]}
\\
\B
\ar[r]_-{\resev}
&
[[\B,\C],\C]
}
$$
commutes 
for any $\A$, $\B$ and
$\C$ and any object $a$ of $\A$.
\end{corollary}
\pf
According to Lemma \ref{dualitynatABC0},
$[1,\ev_a] \circ [-,\C]_{\A,\B}$ has dual
$$\xymatrix{
[\B,\C]
\ar[r]^-{[\A,-]}
&
[[\A,\B],[\A,\C]]
\ar[r]^-{[1,\ev_a]}
&
[[\A,\B],\C]
}$$ since $[-,\C]_{\A,\B}$ has dual $[\A,-]_{\B,\C}$,
whereas $\resev \circ \ev_a$ has dual 
$$\xymatrix{
[\B,\C] 
\ar[r]^-{[\ev_a,1]}
&
[[\A,\B],\C]
}$$ since $\resev$ has dual $1$.
So the result follows from Lemma \ref{homev3}.
\epf

The following result is immediate from
Lemma \ref{chardualq}.
\begin{corollary}
The diagram in $\SMC$
$$
\xymatrix{
[\A,\B] 
\ar[rr]^-{\resev} 
\ar[rd]_{[-,\C]}
&
&
[[[\A,\B], [\A,\C]],[\A,\C]]
\ar[ld]^{[[\A,-],1]}
\\
&
[[\B,\C],[\A,\C]]
&
}
$$
commutes 
for any
$\A$, $\B$ and $\C$.
\end{corollary}
\begin{corollary}\label{evhom1}
The diagram in $\SMC$ 
$$
\xymatrix{
[\B,\C]
\ar[d]_{[\A,-]}
\ar@/^60pt/[dd]^{[F,\C]}
\\
[[\A,\B],[\A,\C]]
\ar[d]_{\ev_F}
\\
[\A,\C].
}
$$ commutes
for any $F:\A \rightarrow \B$ and any $\C$.
\end{corollary}
Similarly one has from Lemma \ref{chardualq} the following.
\begin{corollary}
The diagram
in $\SMC$ 
$$
\xymatrix{
[\A,\B] \ar[rr]^-\resev
\ar[rd]_{[\C,-]} 
&
&
[[[\A,\B],[\C,\B]],[\C,\B]]
\ar[ld]^{[[-,\B],1]}\\
&
[[\C,\A], [\C,\B]]
&
}
$$
commutes 
for any  
$\A$, $\B$ and $\C$.
\end{corollary}
\begin{corollary}\label{evhom2}
The diagram in $\SMC$
$$
\xymatrix{
[\C,\A]
\ar[d]_{[-,\B]}
\ar@/^60pt/[dd]^{[\C,F]}
\\
[[\A,\B],[\C,\B]]
\ar[d]_{\ev_F}
\\
[\C,\B]
}
$$ commutes 
for any  
$F:\A \rightarrow \B$ and any $\C$. 
\end{corollary}

\begin{lemma}\label{evdual}
The diagram in $\SMC$
$$
\xymatrix{
[\A,[\B,\C]]
\ar[rr]^-{\Dual}
\ar[rd]_{\ev_a}
& 
& 
[\B,[\A,\C]]
\ar[ld]^{[1,\ev_a]}
\\
& 
[\B,\C]
& 
}
$$
commutes for any $\A$, $\B$ and $\C$ and any 
object $a$ of $\A$.
\end{lemma}
\pf
The above arrow $[1,\ev_a] \circ \Dual$
rewrites 
\begin{tabbing}
\=1. \=$\xymatrix{
[\A,[\B,\C]]
\ar[r]^-{[-,\C]}
&
[[ [\B,\C],\C ],[\A,\C]]
\ar[r]^-{[q,1]}
&
[\B,[\A,\C]]
\ar[r]^-{[1,\ev_a]}
&
[\B,\C]
}$\\\>2.
\>$\xymatrix{
[\A,[\B,\C]]
\ar[r]^-{[-,\C]}
&
[[ [\B,\C],\C ],[\A,\C]]
\ar[r]^-{[1,\ev_a]}
&
[[ [\B,\C],\C ],\C]
\ar[r]^-{[q,1]}
&
[\B,\C]
}$\\\>3.
\>$\xymatrix{
[\A,[\B,\C]]
\ar[r]^-{\ev_a}
&
[\B,\C]
\ar[r]^-{q}
&
[[ [\B,\C],\C ],\C]
\ar[r]^-{[q,1]}
&
[\B,\C]
}$\\
\>4.
\>$\xymatrix{
[\A,[\B,\C]]
\ar[r]^-{\ev_a}
&
[\B,\C]
\ar@{=}[r]
&
[\B,\C].
}$
\end{tabbing}
In the above derivation:\\
- arrow 1. is equal to 
$[1,\ev_a] \circ \Dual$ according to Lemma \ref{chardualq2};\\ 
- arrows 2. and 3. are equal due to Corollary
\ref{homev3dual};\\
- arrows 3. and 4. are equal by Lemma \ref{chardualq} 
since $\resev$ is dual to the identity.
\epf
\end{section}

\begin{section}{The tensor $\A \otimes \B$ 
for symmetric monoidal categories $\A$ and $\B$}
\label{tensorABsec}
This section is devoted to the definition
by graph and relations of the symmetric monoidal 
category {\em tensor} $\A \otimes \B$
of any two symmetric monoidal categories $\A$ and $\B$.\\ 

Let $\A$ and $\B$ stand for arbitrary symmetric monoidal 
categories.\\

We shall consider the following directed graph $\Gra_{\A,\B}$, or
simply $\Gra$ in this section. Its set of vertices $\Ver_{\A,\B}$, 
or simply $\Ver$ in this section, is the underlying set of the 
free algebra over $\Obj(\A) \times \Obj(\B)$ for the  signature consisting 
of the constant symbol $\un$ and the 2-ary symbol $\otimes$. Which is to say 
that $\Ver$ is the set of words (or more accurately trees!) of the formal 
language composed according to the following axioms:\\
- $\un$ and all pairs $(a,b)$
of $\A \times \B$, which we write $a \otimes b$,
are in $\Ver$;\\
- For any words $X$ and $Y$ in $\Ver$, 
the word $X \otimes Y$ is in $\Ver$.\\

The set of arrows of $\Gra$ is 
also a formal language defined inductively. 
$\Gra$ is the smallest graph on $\Ver$:\\
- containing the graphs 
$\Gra_{1, \A, \B}$,  
$\Gra_{2, \A, \B}$ and $\Gra_{3, \A, \B}$,
just written respectively $\Gra_1$,  
$\Gra_2$ and $\Gra_3$ in this section, 
all with set of vertices $\Ver$, and defined below;\\
- closed under the rules below 
of formation of new arrows.\\

The arrows of $\Gra_1$ are the ``canonical arrows'',
they are:\\
- for any $X,Y,Z$ in $\Ver$,
one $\ac_{X,Y,Z}: X \otimes (Y \otimes Z) 
\rightarrow (X \otimes Y) \otimes Z$
and one 
$\aci_{X,Y,Z}: (X \otimes Y) \otimes Z 
\rightarrow X \otimes (Y \otimes Z)$;\\
- for any $X$ in $\Ver$, one  
$\lc_X: \un \otimes X \rightarrow X$,
one $\lci_X: X \rightarrow \un \otimes X$,
one $\rc_X: X \otimes \un \rightarrow X$
and 
one $\rci_X: X \rightarrow X \otimes \un$;\\
- for any $X,Y$ in $\Ver$, one 
$\syc_{X,Y}: X \otimes Y \rightarrow Y \otimes X$.\\ 

$\Gra_2$ has the following set of arrows:\\ 
- for any object $b$ of $\B$, 
one $\alpha_{b}: \un \rightarrow 
\un_{\A} \otimes b$;\\
- for any object $a$ of $\A$, 
one $\beta_{a}: \un \rightarrow 
a \otimes \un_{\B}$;\\
- for any objects $a, a'$ of $\A$ and $b$ of $\B$, 
one $\gamma_{a,a',b}: 
(a \otimes b) \otimes (a' \otimes b) \rightarrow 
(a \otimes a') \otimes b$;\\
- for any objects $a$ of $\A$ and $b,b'$ of $\B$, 
one $\delta_{a,b,b'}: 
(a \otimes b) \otimes (a \otimes b') \rightarrow 
a \otimes (b \otimes b')$.\\

$\Gra_3$ consists of the following arrows:\\
- for any object $a$ of $\A$ and any arrow 
$f: b \rightarrow b'$ in $\B$, 
one $a \otimes f: a \otimes b \rightarrow a \otimes b'$;\\
- for any arrow $f: a \rightarrow a'$ of $\A$
and any object $b$ of $\B$, 
one $f \otimes b: a \otimes b \rightarrow a' \otimes b$.\\

For any $X$ in $\Ver$ and any arrow
$p: Y \rightarrow Z$, there are new arrows in $\Gra$:
$X \otimes p : X \otimes Y \rightarrow X \otimes Z$ 
and
$p \otimes X : Y \otimes X \rightarrow Z \otimes X$.\\

The set of arrows of $\Gra_1$, $\Gra_2$
and $\Gra_3$ are distinct
and arrows in $\Gra$ with distinct 
names are distinct.\\ 

Note here that since the categories $\A$ and $\B$
are small, $\Ver_{\A,\B}$ is a ``proper'' set,
and for any two $X$ and $Y$ in $\Ver_{\A,\B}$, one has a ``proper'' set
of arrows in $\Gra_{\A,\B}$ from $X$ to $Y$.\\ 

We need to consider the free category $\F_{\A, \B}$, or simply
$\F$ in this section, generated by $\Gra$. Remember that its
arrows are just concatenations
of consecutive edges of $\Gra$. This is a small
category.\\

For any $X$ of $\Ver$, one has 
two graph endomorphisms of $\Gra$, namely $X \otimes -$ and 
$- \otimes X$, sending respectively any $Y$ to $X \otimes Y$
(resp. $Y \otimes X$) and any arrow $f$ to $X \otimes f$ 
(resp. $f \otimes X$).  
The graph endomorphisms $X \otimes -$ 
and $- \otimes X$ extend both to endofunctors
of $\F$ which are still denoted respectively by
$X \otimes -$ and $- \otimes X$. Those are as follows:\\
- For any
arrow $p_1 p_2 ... p_n$ of $\F$ where 
the $p_i$'s are consecutive edges in $\Gra$,  
$X \otimes p$ is the concatenation
$(X \otimes p_1)(X \otimes p_2)... (X \otimes p_n)$,
and $p \otimes X$ is 
$(p_1 \otimes X)(p_2 \otimes X)... (p_n \otimes X)$;\\
- For any $X,Y$ in $\Ver$,
$X \otimes 1_Y = 1_{X \otimes Y} = 1_X \otimes Y$, where 
the $1_X$, $1_Y$ and $1_{X \otimes Y}$ stand for the identities
respectively at $X$, $Y$ and $X \otimes Y$ in $\F$
(i.e some empty strings).\\

For any arrow $f$ in $\A$ and any object $b$ in $\B$,
we may write $f \otimes 1$ for the arrow $f \otimes b$,
and for any object $a$ in $\A$ and any arrow $g$ 
in $\B$, we may write $1 \otimes g$ for $a \otimes g$.
Also for any arrow $p$ of $\F$ and any $X$ in $\Ver$, 
we may write $1 \otimes p$  (respectively $p \otimes 1$) 
for the arrow $X \otimes p$, (resp. $p \otimes X$).\\

We introduce now relations on arrows of $\F_{\A, \B}$
and define the {\em tensor} $\A \otimes \B$ of $\A$ and $\B$
as the category $\F_{\A,\B} / \approx$, 
quotient of $\F_{\A, \B}$ by the 
congruence $\approx$ generated by these relations. 
These relations are the $\sim$ indicated below 
from \ref{cong6} to \ref{cong12}.\\
    
For all arrows $\xymatrix{X \ar[r]^t & Y}$ and
$\xymatrix{Z \ar[r]^s & W}$ in $\Gra$,
 \begin{tag} \label{cong6}

$$\xymatrix{
&
X \otimes W
\ar[rd]^{t \otimes W}
&
\\
X \otimes Z 
\ar[ru]^{X \otimes s}
\ar[rd]_{t \otimes Z}
&
\sim
&
Y \otimes W. 
\\
&
Y \otimes Z
\ar[ru]_{Y \otimes s}
&
}$$
\end{tag}

\begin{tag}Relations stating that canonical arrows of $\A \otimes \B$
are isomorphims.\label{cong7}\end{tag}
They are the following.\\
- For all $X$, $Y$, $Z$ of $\Ver$,
$\xymatrix@C=3pc{
X \otimes (Y \otimes Z)
\ar[r]^-{\ac_{X,Y,Z}}
&
(X \otimes Y) \otimes Z
\ar[r]^-{\aci_{X,Y,Z}}
&
X \otimes (Y \otimes Z)
}
$
$\sim$ $1_{X \otimes (Y \otimes Z)}$
and
$\xymatrix@C=3pc{
(X \otimes Y) \otimes Z
\ar[r]^-{\ac_{X,Y,Z}}
&
X \otimes (Y \otimes Z)
\ar[r]^-{\aci_{X,Y,Z}}
&
(X \otimes Y) \otimes Z
}$
$\sim$ $1_{(X \otimes Y) \otimes Z}$.\\
- For all $X$ of $\Ver$, 
$\xymatrix{
\un \otimes X
\ar[r]^-{\lc_X}
&
X
\ar[r]^-{{\lci}_X }
&
\un \otimes X
}$
$\sim$ $1_{\un \otimes X}$
and
$\xymatrix{
X
\ar[r]^-{\lci_X}
&
\un \otimes X 
\ar[r]^-{\lc_X}
&
X
}$
$\sim$  
$1_X$.\\
- For all $X$ of $\Ver$,
$\xymatrix{
X \otimes \un 
\ar[r]^-{\rc_X}
&
X
\ar[r]^-{\rci_X}
&
X \otimes \un}$
$\sim$ 
$1_{X \otimes \un}$
and $\xymatrix{
X 
\ar[r]^-{ {\rci}_X } 
&
X \otimes \un
\ar[r]^-{ \rc_X } 
&
X
}$
$\sim$ 
$1_X$.\\
- For all $X, Y$ of $\Ver$,
$\xymatrix{
X \otimes Y
\ar[r]^-{s_{X,Y}}
&
Y \otimes X  
\ar[r]^-{s_{Y,X}}
&
X \otimes Y
}$
$\sim$ 
$1_{X \otimes Y}$.\\

\begin{tag}\label{cong8}
\noindent Relations giving the coherence conditions
for $\ac$, $\rc$, $\lc$ and $\syc$ in $\A \otimes \B$. 
\end{tag}
These are the following.\\ 
- For all $X$, $Y$, $Z$ and $T$ in $\Ver$,
$$\xymatrix{ X \otimes (Y \otimes (Z \otimes T)) 
\ar@{}[rrd]|{\sim}
\ar[r]^{\ac} 
\ar[d]_{1 \otimes \ac} & 
(X \otimes Y) \otimes (Z \otimes T) 
\ar[r]^{\ac}  & 
 ( ( X \otimes Y) \otimes Z) \otimes T \\
X \otimes ( ( Y \otimes Z) \otimes T) 
\ar[rr]_{\ac} & &
( X \otimes ( Y \otimes Z )) \otimes T. 
\ar[u]_{\ac \otimes 1} }$$
- For all $X, Y $ in $\Ver$,
$$\xymatrix{ X \otimes ( \un \otimes Y) \ar[rr]^{\ac} 
\ar[rd]_{1 \otimes \lc} &  \ar@{}[d]|{\sim} &  
(X \otimes \un) \otimes Y \ar[ld]^{\rc \otimes 1} \\
& X \otimes Y. &  }$$
- For all $X$ in $\Ver$,
$$\xymatrix{ X \otimes \un \ar[rr]^{\syc} \ar[rd]_{\rc} 
& \ar@{}[d]|{\sim} & \un \otimes X 
\ar[ld]^{\lc} \\ 
& X. & }$$ 
- For all $X$, $Y$ and $Z$ of $\Ver$,
$$\xymatrix{
X \otimes ( Y \otimes Z) \ar[r]^{\ac} 
\ar[d]_{1 \otimes s} & 
( X \otimes Y ) \otimes Z \ar[r]^{\syc} \ar@{}[d]|{\sim} &
Z \otimes ( X \otimes Y) \ar[d]^{\ac}\\
X \otimes ( Z \otimes Y) \ar[r]_{\ac} & 
(X \otimes Z) \otimes Y & (Z \otimes X) \otimes Y. 
  \ar[l]^{\syc \otimes 1} 
}$$

\begin{tag}\label{cong9}
Relations for the naturalities of $\ac$, $\rc$, $\lc$,
and $\syc$ in $\A \otimes \B$.\end{tag}
For instance, one has 
for any $f:X \rightarrow X'$ in $\Gra$,
and any $Y, Z$ in $\Ver$, 
$$\xymatrix@C=3pc{ 
X \otimes ( Y \otimes Z) 
\ar[r]^{\ac_{X,Y,Z}}
\ar@{}[rd]^{\sim}
\ar[d]_{f \otimes 1} &
(X \otimes Y) \otimes Z \ar[d]^{(f \otimes 1) \otimes 1}\\
X' \otimes ( Y \otimes Z) \ar[r]_{\ac_{X',Y,Z}}  &
(X' \otimes Y) \otimes Z.}
$$
We will not write here the other relations.
There are two more for the naturalities of $\ac_{X,Y,Z}$ in $Y$ and $Z$, 
one for that of $\lc_X$ in $X$,
one for that of $\rc_X$ in $X$
and two for those of $\syc_{X,Y}$
in $X$ and $Y$.\\ 

For any object $a$ in $\A$ and any arrows
$\xymatrix{b \ar[r]^{f}  & b' \ar[r]^{g} & b''}$
in $\B$, 
\begin{tag}\label{cong131}
$$\xymatrix{
a \otimes b 
\ar[rr]^{a \otimes (g \circ f)}
\ar[rd]_{a \otimes f}
&
\ar@{}[d]|{\sim}
&
a \otimes b''\\
&
a \otimes b'
\ar[ru]_{a \otimes g}
&
}$$
\end{tag}

For any object $b$ in $\B$ and any arrows 
$\xymatrix{a \ar[r]^{f}  & a' \ar[r]^{g} & a''}$
in $\A$,
\begin{tag}\label{cong133}
$$\xymatrix{
a \otimes b 
\ar[rr]^{(g \circ f) \otimes b}
\ar[rd]_{f \otimes b}
&
\ar@{}[d]|{\sim}
&
a'' \otimes b\\
&
a' \otimes b
\ar[ru]_{g \otimes b}
&
}$$

\end{tag}

For any objects $a$ in $\A$ and $b$ in $\B$,
\begin{tag}\label{cong132}
$a \otimes 1_b \sim 1_{a \otimes b}$
\end{tag}
and
\begin{tag}\label{cong134}
$ 1_a \otimes b \sim  1_{a \otimes b}$. 
\end{tag} 
where $1_b$, $1_a$ and $1_{a \otimes b}$ above
are the identities respectively at $b$ in $\B$, at $a$ in $\A$ 
and at $a \otimes b$ in $\F$.\\

For any arrows $f: a \rightarrow a'$ in $\A$ and
$g: b \rightarrow b'$ in $\B$, 
\begin{tag}\label{cong135}
$$\xymatrix{a \otimes b \ar[d]_{a \otimes g} 
\ar[r]^{f \otimes b}
\ar@{}[rd]|{\sim} & 
a' \otimes b \ar[d]^{a' \otimes g}\\
a \otimes b' \ar[r]_{f \otimes b'}& a' \otimes  b'.}$$
\end{tag}

\begin{tag}\label{cong10}
Relations for the ``naturalities'' 
of $\alpha_{b}$ in $b$ , $\beta_{a}$ in $a$,
$\gamma_{a,a',b}$ in $a$, $a'$ and $b$
and $\delta_{a, b, b'}$ in $a$, $b$ and $b'$.
\end{tag}
For instance by the relations for the ``naturality'' of $\gamma_{a,a',b}$
in $b$ it is meant that for any objects $a,a'$ in $\A$ and any 
arrow $g: b \rightarrow b'$ in $\B$, 
$$\xymatrix{ (a \otimes b) \otimes (a' \otimes b) 
\ar[d]_{ (1 \otimes g) \otimes (1 \otimes g) }
\ar[r]^-{\gamma_{a,a',b}} 
\ar@{}[rd]|{\sim}
& (a \otimes a') \otimes b 
\ar[d]^{1 \otimes g}\\
(a \otimes b')  \otimes (a' \otimes b') 
\ar[r]_-{\gamma_{a, a', b'}} 
& 
(a \otimes a') \otimes b'.
}$$  
We will not write explicitly now the seven 
other relations.\\

For any objects $a$ in $\A$ and $b$, $b'$, $b''$ in $\B$,
\begin{tag}\label{cong113}
$$\xymatrix{ 
(a \otimes b) \otimes ( (a \otimes b') \otimes (a \otimes b'') )
\ar[r]^{\ac} \ar[d]_{ 1 \otimes \delta_{a,b',b''} }
\ar@{}[rdd]|{\sim}
& 
( (a \otimes b) \otimes (a \otimes b') ) \otimes (a \otimes b'')
\ar[d]^{ \delta_{a,b,b'} \otimes 1 }
\\
(a \otimes b) \otimes (a \otimes (b' \otimes b''))
\ar[d]_{ \delta_{a,b, b' \otimes b''} } & 
(a \otimes (b \otimes b')) \otimes (a \otimes b'')
\ar[d]^{ \delta_{a, b \otimes b' ,  b''} }\\
a \otimes (b \otimes (b' \otimes b'')) 
\ar[r]_{ 1 \otimes \ac_{ b, b',  b'' } } &
a \otimes ((b \otimes b') \otimes b''). 
}$$
\end{tag}

For any objects $a$ in $\A$ and $b$ in $\B$,
\begin{tag}\label{cong1141} 
$$\xymatrix{
(a \otimes b) \otimes \un 
\ar[r]^{\rc_{a \otimes b}} 
\ar[d]_{1 \otimes \beta_a}
\ar@{}[rd]|{\sim}
& 
a \otimes b\\
(a \otimes b) \otimes (a \otimes \un_{\B}) \ar[r]_-{\delta_{a,b, \un_{\B}}} &
a \otimes (b \otimes \un_{\B}) 
\ar[u]_{1 \otimes \rc_b}
}$$
\end{tag}
and 
\begin{tag}\label{cong1142}
$$\xymatrix{
\un \otimes (a \otimes b) 
\ar[r]^{\lc_{a \otimes b}} 
\ar[d]_{\beta_a \otimes 1}
\ar@{}[rd]|{\sim}
& 
a \otimes b\\
(a \otimes \un_{\B})  \otimes (a \otimes b) \ar[r]_-{\delta_{a,\un_{\B}, b}} &
a \otimes (\un_{\B} \otimes b) 
\ar[u]_{\l \otimes \lc_b}
.}$$
\end{tag}

For any objects $a$ in $\A$ and $b$, $b'$ in $\B$,
\begin{tag}\label{cong115}
$$\xymatrix{
(a \otimes b) \otimes (a \otimes b') \ar[r]^-{ \delta_{a,b,b'} }
\ar[d]_{ \syc_{a \otimes b,a \otimes b'} } 
\ar@{}[rd]|{\sim} & 
a \otimes (b \otimes b') 
\ar[d]^{1 \otimes \syc_{b,b'}}\\
(a \otimes b') \otimes (a \otimes b) \ar[r]_-{\delta_{a,b',b}} & 
a \otimes (b' \otimes b).}$$
\end{tag}

For any objects $a$, $a'$, $a''$ in $\A$ and $b$ in $\B$,
\begin{tag}\label{cong143}
$$\xymatrix{ 
(a \otimes b) \otimes ( (a' \otimes b ) \otimes (a'' \otimes b ))
\ar[r]^{\ac} \ar[d]_{1 \otimes \gamma_{a',a'',b}}
\ar@{}[rdd]|{\sim}
& 
( (a \otimes b) \otimes (a' \otimes b) ) \otimes (a'' \otimes b)
\ar[d]^{ \gamma_{a,a',b} \otimes 1 }\\
(a \otimes b) \otimes ( (a' \otimes a'') \otimes b ))
\ar[d]_{ \gamma_{a,a' \otimes a'', b} }& 
((a \otimes a') \otimes b) \otimes (a'' \otimes b)
\ar[d]^{\gamma_{a \otimes a', a'' , b}}\\
(a \otimes (a' \otimes a'')) \otimes b 
\ar[r]_{ \ac \otimes 1 } &
( (a \otimes a') \otimes a'') \otimes b.
}$$
\end{tag}

For any objects $a$ in $\A$ and $b$ in $\B$,
\begin{tag}\label{cong1441} 
$$\xymatrix{
(a \otimes b) \otimes \un \ar[r]^{\rc_{a \otimes b}} 
\ar[d]_{1 \otimes \alpha_b}
\ar@{}[rd]|{\sim}
& 
a \otimes b\\
(a \otimes b) \otimes (\un_{\A} \otimes b) 
\ar[r]_-{ \gamma_{a,\un_{\A}, b} } &
(a \otimes \un_{\A}) \otimes b
\ar[u]_{\rc_a \otimes 1}
}$$
\end{tag}
and 
\begin{tag}\label{cong1442}
$$\xymatrix{
\un \otimes (a \otimes b) \ar[r]^{\lc_{a \otimes b}} 
\ar[d]_{\alpha_b \otimes 1} 
\ar@{}[rd]|{\sim}
& 
a \otimes b\\
(\un_{\A} \otimes b)  \otimes (a \otimes b) 
\ar[r]_-{\gamma_{ \un_{\A}, a , b } } &
(\un_{\A} \otimes a) \otimes b \ar[u]_{\lc_a \otimes 1}
.}$$
\end{tag}

For any objects $a$, $a'$ in $\A$ and $b$ in $\B$,
\begin{tag}\label{cong145}
$$\xymatrix{
(a \otimes b) \otimes (a' \otimes b) 
\ar[r]^-{\gamma_{a,a',b}} 
\ar[d]_{s_{a \otimes b, a' \otimes b}}
\ar@{}[rd]|{\sim} & 
(a \otimes a') \otimes b 
\ar[d]^{s_{a, a'} \otimes 1}\\
(a' \otimes b) \otimes (a \otimes b) 
\ar[r]_-{\gamma_{a',a,b}} & (a' \otimes a) \otimes b
.}$$  
\end{tag}

For any objects $b$,$b'$ in $\B$,
\begin{tag}\label{cong161}
$$\xymatrix{ 
\un \otimes \un 
\ar[d]_{\alpha_b \otimes \alpha_{b'}}
\ar@{-}[r]^-{\cong} 
\ar@{}[rd]|{\sim}
& \un 
\ar[d]^{ \alpha_{b \otimes b'} }\\
(\un_{\A} \otimes b) \otimes (\un_{\A} \otimes b')
\ar[r]_-{\delta_{\un_{\A}, b, b'} }& 
\un_{\A} \otimes (b \otimes b').
}$$
\end{tag}

For any objects $a$, $a'$ in $\A$,
\begin{tag}\label{cong162}
$$\xymatrix{ 
\un \otimes \un 
\ar[d]_{\beta_a \otimes \beta_{a'}}
\ar@{-}[r]^{\cong} 
\ar@{}[rd]|{\sim}
& 
\un \ar[d]^{\beta_{a \otimes a'}}\\
(a \otimes \un_{\B}) \otimes (a' \otimes \un_{\B} )
\ar[r]_-{\gamma_{a, a', \un_{\B}}} & 
(a \otimes a') \otimes \un_{\B}.
}$$
\end{tag}

\begin{tag}\label{cong17}
$\xymatrix{ \un \ar[r]^-{\beta_{\un_{\A}} } & \un_{\A} \otimes \un_{\B}  }$
$\sim$
$\xymatrix{ \un \ar[r]^-{\alpha_{\un_{\B}}} & \un_{\A} \otimes \un_{\B}}$.
\end{tag}

For any objects $a$,$a'$ in $\A$ and $b$,$b'$ in $\B$,
\begin{tag}\label{cong19}
$$\xymatrix{
((a \otimes b) \otimes (a \otimes b'))
\otimes 
((a' \otimes b) \otimes (a' \otimes b'))
\ar@{-}[r]^{\cong} 
\ar[d]_{\delta_{a,b,b'} \otimes \delta_{a',b,b'} }
\ar@{}[rd]|{\sim} 
&
((a \otimes b) \otimes (a' \otimes b))
\otimes 
((a \otimes b') \otimes (a' \otimes b'))
\ar[d]^{\gamma_{a,a',b} \otimes \gamma_{a,a',b'}}\\
(a \otimes (b \otimes b')) \otimes (a' \otimes (b \otimes b')) 
\ar[d]_{ \gamma_{a,a', b \otimes b'} } &
((a \otimes a') \otimes b) 
\otimes 
((a \otimes a') \otimes b') \ar[ld]^{\delta_{a \otimes a', b, b'}}\\
(a \otimes a') \otimes (b \otimes b'). &
 }$$
\end{tag}

\begin{tag}\label{cong12}
Expansions of all relations above by iterations of
$X \otimes -$ and $- \otimes X$ for all $X$ in $\Ver$.
\end{tag}
Which means precisely that the set of relations $\sim$ is the 
smallest set of relations on arrows of $\F$ 
containing the previous relations
(\ref{cong6} to \ref{cong19})
and satisfying the closure properties that
for any relation $f \sim g : Y \rightarrow Z$ 
and any $X$ in $\Ver$, 
one has the relations
\begin{center}
$X \otimes f \sim  X \otimes g: X \otimes Y \rightarrow X \otimes Z$
\end{center}
and 
\begin{center}
$f \otimes X \sim  g \otimes X: Y \otimes X \rightarrow Z \otimes X$.
\end{center} 

From now on, the arrows of $\A \otimes \B$,
which are $\approx$-classes, will be denoted 
with the same name as their
representatives in $\F$.\\

Now it is straightforward to check
that the category $\A \otimes \B$ admits a symmetric monoidal 
structure. It is small since $\F_{\A,\B}$ is.
According to the relations \ref{cong12}, 
the endofunctor $X \otimes -$ of $\F$
sends equivalent arrows by $\approx$ to equivalent arrows. 
Therefore there exists a unique endofunctor of the quotient 
category $\F/\approx$, still denoted 
$X \otimes -$, that makes the following diagram commute
$$
\xymatrix{ 
\F 
\ar[d]   
\ar[r]^{X \otimes -}
&  
\F 
\ar[d]\\
(\F / \approx) 
\ar[r]_{X \otimes -}
&  
(\F / \approx) 
}
$$
where the vertical arrows are the canonical ones associated
to the quotient. 
One defines similarly from the endofunctor 
$- \otimes X$ of $\F$, 
the endofunctor $- \otimes X$ of $\F / \approx$.

From the relations \ref{cong6}, an induction shows  
that for any arrows $f: X \rightarrow X'$ and
$g: Y \rightarrow Y'$ in $\F$,
the arrows 
$\xymatrix{
X \otimes Y
\ar[r]^{f \otimes Y}
&
X' \otimes Y
\ar[r]^{X' \otimes g}
&
X' \otimes Y'
}$
and 
$\xymatrix{
X \otimes Y
\ar[r]^{X \otimes g}
&
X \otimes Y'
\ar[r]^{f \otimes Y'}
&
X' \otimes Y'}$
are $\approx$-equivalent.
Therefore for any arrows $f: X \rightarrow X'$ and 
$g: Y \rightarrow Y'$ of $\F / \approx$,
$(f \otimes Y') \circ (X \otimes g) 
= (X' \otimes g) \circ (f \otimes Y)$ and
the collections of endofunctors $X \otimes -$
and $- \otimes X$ for $X$ ranging in $\Ver$,
define a bifunctor
$(\F / \approx) \times (\F / \approx) \rightarrow (\F / \approx)$.

That this bifunctor as tensor, together with the object $\un$ as unit
and the $\ac$, $\rc$, $\lc$ and 
$\syc$ as canonical arrows,
define a symmetric monoidal structure,
results from 
relations \ref{cong7} (implying that arrows 
$\ac$, $\rc$, $\lc$
and $\syc$ are isomorphisms and Axiom \ref{smcaxiom3}), 
relations \ref{cong8} (Axioms \ref{mcaxiom1}, \ref{mcaxiom2}, \ref{smcaxiom4} 
and \ref{smcaxiom5}) and 
relations \ref{cong9} (naturalities of canonical arrows).

\end{section}

\begin{section}
{The embedding $\eta: \A \rightarrow [\B, \A \otimes \B]$}
\label{embeta}
In this section a symmetric monoidal functor 
$\eta: \A \rightarrow [\B, \A \otimes \B]$  
is defined for any symmetric monoidal categories
$\A$ and $\B$. It will serve
to define the adjunction of the next section.\\

Let $\A$ and $\B$ stand here for arbitrary symmetric monoidal
categories.\\

For any object $a$ of $\A$, $\eta(a)$ denotes the following
symmetric monoidal functor $\B \rightarrow \A \otimes \B$.
It sends any object $b$ to $a \otimes b$ and any arrow
$g: b \rightarrow b'$ to 
$a \otimes g: a \otimes b \rightarrow a \otimes b'$.
These assignments define a functor 
$\B \rightarrow \A \otimes \B$, 
according to the relations
\ref{cong131} and \ref{cong132}. 
The monoidal structure on $\eta(a)$ is defined 
by ${\eta(a)}^0 = \beta_a$ and
for any objects $b,b'$ in $\B$ by
${\eta(a)}^2_{b,b'} = \delta_{a, b, b'}$. 
The relations \ref{cong10} for the naturalities
of $\delta_{a,b,b'}$ in $b$ and $b'$ give the naturalities
of ${\eta(a)}^2_{b,b'}$ in $b$ and $b'$, and the relations
\ref{cong113}, \ref{cong1141}, \ref{cong1142} and 
\ref{cong115}
give respectively Axioms \ref{mofun3}, \ref{mofun41},
\ref{mofun42} and 
 \ref{symofun5}
for the triple $(\eta(a),{\eta(a)}^0,{\eta(a)}^2)$.\\ 

For any arrow $f: a \rightarrow a'$ in $\A$, 
the collection of arrows 
$f \otimes b : a \otimes b \rightarrow a' \otimes b$, 
$b$ ranging in $\B$, forms a natural transformation 
$\eta(f):\eta(a) \rightarrow \eta(a')$ according to the 
relations \ref{cong135}. This one is moreover monoidal 
with respect to the monoidal structures of $\eta(a)$ and 
$\eta(a')$: Axioms \ref{monat6} and \ref{monat7}
for $\eta(f): \eta(a) \rightarrow \eta(a')$
result respectively from the 
 relations \ref{cong10}
for the naturality
of the $\delta_{a,b,b'}$ 
in $a$ and for the naturality of the $\beta_a$ in $a$.\\

That the assignments 
$a \mapsto \eta(a)$ and $(f: a \rightarrow a') \mapsto 
\eta(f): \eta(a) \rightarrow \eta(a')$ define a functor 
$\A \rightarrow \SMC(\B, \A \otimes \B)$ is due to the relations  
\ref{cong133}
and  \ref{cong134}.
We are going check that this functor admits the 
following symmetric 
monoidal structure $\A \rightarrow [\B,\A \otimes \B]$.
The monoidal natural transformation 
$$\eta^0: \un \rightarrow \eta(\un_{\A}): 
\B \rightarrow \A \otimes \B$$
is defined in any $b$ of $\B$ as 
$\alpha_b: \un \rightarrow \un_{\A} \otimes b$
and for any objects $a,a'$ in $\A$, the monoidal natural
transformation
$$\eta^2_{a,a'}: \eta(a) \Box \eta(a')  
\rightarrow \eta(a \otimes a'): \B \rightarrow \A \otimes \B$$ 
is defined in any $b$ of $\B$ as
$$\gamma_{a,a',b}: (a \otimes b) \otimes (a' \otimes b) \rightarrow (a \otimes a') \otimes b.$$ 

The above collection 
${\eta}^0$ is a well defined natural transformation
$\un \rightarrow \eta(\un_{\A})$
according to the relations \ref{cong10} for the naturality of the 
$\alpha_b$ in $b$. It is moreover monoidal 
with Axiom \ref{monat6} given
by the relations \ref{cong161}  
and Axiom \ref{monat7} by the relations \ref{cong17}.\\

For any objects $a,a'$ of $\A$, the above collection 
${\eta}^2_{a,a'}$
is a well defined natural transformation
$\eta(a)  \Box \eta(a')  
\rightarrow \eta(a \otimes a')$
according to the relations \ref{cong10} for 
the naturality of the 
$\gamma_{a,a',b}$ in $b$.
It is moreover monoidal 
with Axiom \ref{monat6} given by the relations 
\ref{cong19},  
and Axiom \ref{monat7} by the relations \ref{cong162}.
The collection of 
${\eta}^2_{a,a'}: \eta(a) \Box \eta(a') \rightarrow \eta(a \otimes
a')$ in $[\B,\A \otimes \B]$
is natural in $a$ and $a'$, according to the relations \ref{cong10}
for the naturalities of the $\gamma_{a,a',b}$ in $a$ and $a'$.\\

Eventually Axioms \ref{mofun3}, \ref{mofun41}, \ref{mofun42}
and \ref{symofun5} hold for the triple $(\eta, \eta^0, \eta^2)$ 
respectively due
to the relations \ref{cong143}, \ref{cong1441}, 
\ref{cong1442} and \ref{cong145}. 
\end{section}

\begin{section}{The symmetric monoidal adjunction 
$\Ext \dashv \Res : [\A \otimes \B,\C] \rightarrow [\A,[\B,\C]]$}
\label{P1}
It is established in this section the existence 
of a symmetric monoidal adjunction
$$\Ext \dashv \Res: [\A \otimes \B, \C] \rightarrow [\A,[\B,\C]]$$
such that $\Res \circ \Ext = 1$, 
for any symmetric monoidal categories $\A$, $\B$ and $\C$.\\

$\A$, $\B$ and $\C$ will stand here for arbitrary symmetric
monoidal categories.
The symmetric monoidal functor  
$\Res$ is defined as the composite in $\SMC$
$$\xymatrix{ [\A \otimes \B,\C] \ar[r]^-{[\B,-]} &
[[\B,\A \otimes \B],[\B, \C]] \ar[r]^-{[\eta, 1]} &
[\A,[\B,\C]]}.$$

We embark now for the definition of the underlying 
functor
$$\Ext: \SMC(\A,[\B,\C]) \rightarrow \SMC(\A \otimes \B, \C).$$

On objects, it is as follows.
Given any symmetric monoidal functor $F: \A \rightarrow [\B, \C]$,
it is sent by $\Ext$ to the symmetric monoidal
functor $\bar{F}: \A \otimes \B \rightarrow \C$ described below.\\

The action of $\bar{F}$ on $\Ver_{\A,\B}$ is defined by induction
according to the rules:\\
- $\bar{F}(\un) = \un$;\\
- For any objects $a$ in $\A$ and 
$b$ in $\B$, $\bar{F}(a \otimes b) = F(a)(b)$;\\
- For any objects $X,Y$ in $\A \otimes \B$,
$\bar{F}(X \otimes Y) = \bar{F}(X) \otimes \bar{F}(Y)$.\\

To define the assignment $\bar{F}$ on the arrows
of $\A \otimes \B$, one first defines the graph morphism 
denoted $\gF$,
from $\Gra_{\A,\B}$ to the underlying graph of 
$\C$. For any $X$ in $\Ver_{\A,\B}$,
$\gF(X) = \bar{F}(X)$.\\

$\gF$ on arrows of $\Gra^1_{\A,\B}$ is as follows.\\
- For any $X,Y,Z$ in $\Ver_{\A,\B}$,
$\gF$ sends 
$a_{X,Y,Z}: X \otimes (Y \otimes Z) \rightarrow (X \otimes Y) \otimes Z$
to  
{\small
$$\xymatrix@C=4pc{
\bar{F}(X \otimes (Y \otimes Z)) 
\;=\; 
\bar{F}(X) \otimes (\bar{F}(Y) \otimes \bar{F}(Z))
\;
\ar[r]^-{ \ac_{\bar{F}(X),\bar{F}(Y), \bar{F}(Z)} } 
&
\;  
(\bar{F}(X) \otimes \bar{F}(Y)) \otimes \bar{F}(Z)
\;=\;
\bar{F}( (X \otimes Y) \otimes Z) ) }$$
}
and $\gF(\bar{\ac}_{X,Y,Z})$ is 
${\ac}^{-1}_{\bar{F}(X), \bar{F}(Y), \bar{F}(Z)}$.\\
- For any $X$ in  $\Ver_{\A,\B}$,
$\gF$ sends $\lc_X: \un \otimes X \rightarrow X$ 
to
$$\xymatrix{ \bar{F}(\un \otimes X) \ar@{=}[r]
& \un \otimes \bar{F}(X) \ar[r]^-{ \lc_{\bar{F}(X)} } & 
\bar{F}(X)
}$$
and $\gF( \bar{\lc}_X )$ is ${ \lc_{\bar{F}(X)} }^{-1}$;
similarly
$\gF(\rc_X)$ is $\rc_{\bar{F}(X)}$
and $\gF(\bar{\rc}_X)$ is ${\rc}^{-1}_{\bar{F}(X)}$.\\
- For any $X,Y$ in $\Ver_{\A,\B}$,
$\gF$ sends $\syc_{X,Y}: X \otimes Y \rightarrow Y \otimes X$ 
to
$$\xymatrix@C=3pc{\bar{F}(X \otimes Y) 
\ar@{=}[r] & 
\bar{F}(X) \otimes \bar{F}(Y)
\ar[r]^-{\syc_{\bar{F}(X), \bar{F}(Y)}} &
\bar{F}(Y) \otimes \bar{F}(X) \ar@{=}[r] &
\bar{F}(Y \otimes X) 
}.$$

$\gF$ on arrows of $\Gra^2_{\A,\B}$ is as follows.\\
- For any object $b$ of $\B$,
$\gF$ sends $\alpha_b: \un \rightarrow \un_{\A} \otimes b$ 
to $$\xymatrix{ 
\bar{F}(\un) 
\ar@{=}[r] & 
\un \ar[r]^-{{F^0}_b} & 
F(\un_{\A})(b) \ar@{=}[r] &
\bar{F}(\un_{\A} \otimes b)}$$
where ${F^0}_b$ denotes the component
in $b$ of the monoidal natural transformation
$F^0: \un \rightarrow F( \un_{\A} ): \B \rightarrow \C$,
part of the monoidal structure of $F$.\\
- For any object $a$ of $\A$,
$\gF$ sends $\beta_a: \un \rightarrow a \otimes \un_{\B}$ 
to 
$$\xymatrix{ 
\bar{F}(\un) 
\ar@{=}[r] & 
\un \ar[r]^-{ {F(a)}^0 } & 
F(a)(\un_{\B}) \ar@{=}[r] &
\bar{F}( a \otimes \un_{\B} ) 
}$$
where ${F(a)}^0$
is part of the monoidal structure of $Fa: \B \rightarrow \C$.\\
- For any objects $a$, $a'$ of $\A$ and $b$ of $\B$,
$\gF$ sends $\gamma_{a,a',b} : 
(a \otimes b) \otimes (a' \otimes b) \rightarrow 
(a \otimes a') \otimes b $ to 
$$\xymatrix{
\bar{F}( (a \otimes b) \otimes (a' \otimes b) )
\ar@{=}[r] & 
F(a)(b) \otimes F(a')(b)  
\ar[r]^-{ { (F^2_{a,a'}) }_b } &
F(a \otimes a')(b) \ar@{=}[r] &
\bar{F}((a \otimes a') \otimes b)
}$$
where
${(F^2_{a,a'}) }_b$ is the component in $b$
of the monoidal transformation 
$$F^2_{a,a'}: F(a) \Box F(a') \rightarrow F (a \otimes a'): 
\B \rightarrow \C,$$
part of the monoidal structure of $F$. \\
- For any objects $a$ of $\A$ and $b$, $b'$ of $\B$,
$\gF$ sends $\delta_{a,b,b'} : 
(a \otimes b) \otimes (a \otimes b') \rightarrow 
a  \otimes (b \otimes b')$ to
$$\xymatrix{
\bar{F}( (a \otimes b) \otimes (a \otimes b') )
\ar@{=}[r] & 
F(a)(b) \otimes F(a)(b')  
\ar[r]^-{ { F(a)}^2_{b,b'} } &
F(a) (b \otimes b') \ar@{=}[r] &
\bar{F}(a \otimes (b \otimes b'))
}$$
where $F(a)^2_{b,b'}$
part of the monoidal structure of $F(a)$.\\

$\gF$ on the arrows of $\Gra^3_{\A,\B}$ is as follows.\\
- For any arrow $f: a \rightarrow a'$ in $\A$ and any object $b$ 
of $\B$, 
$\gF$ sends $f \otimes b: a \otimes b \rightarrow a' \otimes b$
to
${F(f)}_b:
F(a)(b) \rightarrow  
F(a')(b)$, the component in $b$
of the monoidal natural transformation $F(f): F(a) \rightarrow F(a')$.\\
- For any object $a$ of $\A$ and any arrow $g: b \rightarrow b'$ in $\B$, 
$\gF$ sends $a \otimes g:  a \otimes b \rightarrow a \otimes b'$
to 
$F(a)(g):
F(a)(b) \rightarrow  
F(a)(b')$.\\

$\gF$ is defined on all arrows of $\Gra_{\A,\B}$ by induction
according to the following rules.\\
- For any object $X$ of $\Ver_{\A,\B}$ and any arrow 
$p: Y \rightarrow Z$ of $\Gra_{\A,\B}$,
$\gF$ sends 
$X \otimes p: X \otimes Y \rightarrow X \otimes Z$ 
to 
$$\xymatrix{
\bar{F}(X \otimes Y) \ar@{=}[r] & \bar{F}(X) \otimes \bar{F}(Y)
\ar[r]^{ 1 \otimes \gF(p)} &
\bar{F}(X) \otimes \bar{F}(Z) \ar@{=}[r] & 
\bar{F}(X \otimes Z)
}$$
and 
$\gF$ sends
$p \otimes X : Y \otimes X \rightarrow Z \otimes X $ 
to
$$\xymatrix{
\bar{F}(Y \otimes X) \ar@{=}[r] & \bar{F}(Y) \otimes \bar{F}(X)
\ar[r]^{ \gF(p) \otimes 1} &
\bar{F}(Z) \otimes \bar{F}(X) \ar@{=}[r] & 
\bar{F}(Z \otimes X)
}.$$

The above graph morphism 
$\gF$ induces a functor $\F_{\A,\B} \rightarrow \C$,
still written $\gF$. This last functor induces a functor 
$\A \otimes \B \rightarrow \C$, namely $\bar{F}$, which we
define as the image of $F$ by $\Ext$.
To check this last point, it is sufficient to check 
that for any of the relations
$f \sim g$ on arrows of $\F_{\A,\B}$ defining $\A \otimes \B$
(from \ref{cong6} to \ref{cong12}), 
one has the equality $\gF(f) = \gF(g)$ in $\C$. Which we do below.\\

Regarding relations $\ref{cong6}$.\\
For any arrow $t: X \rightarrow Y$ and $s: Z \rightarrow W$
in $\Gra_{\A,\B}$,
$\gF$ sends  
$\xymatrix{X \otimes Z \ar[r]^{X \otimes s} &
X \otimes W \ar[r]^{t \otimes W} &  Y \otimes W}$
to
$$\xymatrix{ \bar{F}(X) \otimes \bar{F}(Z) 
\ar[r]^{ 1 \otimes \gF(s) } & 
\bar{F}(X) \otimes \bar{F}(Z)
\ar[r]^{ \gF(t) \otimes 1 } &
\bar{F}(Y) \otimes \bar{F}(W)  
}$$
and $
\xymatrix{X \otimes Z \ar[r]^{t \otimes Z} &
Y \otimes Z \ar[r]^{Y \otimes s} &  Y \otimes W }$
to
$$\xymatrix{ 
\bar{F}(X) \otimes \bar{F}(Z) 
\ar[r]^{ \gF(t) \otimes 1 } & 
\bar{F}(Y) \otimes \bar{F}(Z)
\ar[r]^{ 1 \otimes \gF(s) } &
\bar{F}(Y) \otimes \bar{F}(W)  
}$$
These two arrows are equal by the bifunctoriality of
the tensor in $\C$.\\

Relations \ref{cong7}, 
\ref{cong8} and \ref{cong9} are sent by $\gF$ to commuting diagrams
in $\C$. This results from the fact that $\C$ is symmetric
monoidal
and that $\gF$ sends:\\
- $\un$ to $\un$;\\
- objects of the form $X \otimes Y$ to $\bar{F}(X) \otimes \bar{F}(Y)$;\\
- arrows of the form $X \otimes f$ and $f \otimes X$
for any object $X$ and any arrow $f$ of $\Gra_{\A,\B}$, respectively to
$1 \otimes \gF(f)$ and to $\gF(f) \otimes 1$;\\
- the arrows $\ac$ to $\ac$,
the $\aci$ to ${\ac}^{-1}$, the $\lc$ to $\lc$,
the $\lci$ to ${\lc}^{-1}$,
the $\rc$ to $\rc$, 
the $\rci$ to ${\rc}^{-1}$ and,
the $\syc$ to $\syc$.\\

Regarding relations \ref{cong131} and \ref{cong132}.\\
For any object $a$ of $\A$ and
any arrows $\xymatrix{ b \ar[r]^f & b' \ar[r]^g & b''}$ in $\B$,
$\gF$ sends
$$\xymatrix{ 
(a \otimes b) \ar[r]^{a \otimes f} &
(a \otimes b') \ar[r]^{a \otimes g} &
(a \otimes b'')
}$$
to
$F(a)(g) \circ F(a)(f)$, that is $F(a)(g \circ f)$
since $F(a): \B \rightarrow \C$ is a functor, which is  
the image by $\gF$ of $a \otimes (g \circ f): a \otimes b \rightarrow
a \otimes b''$.\\

For any objects $a$ of $\A$ and $b$ of $\B$,
$\gF$ sends $a \otimes 1_b: a \otimes b \rightarrow a \otimes b$ 
to
$F(a)(1_b)$, which is the identity at $F(a)(b)$
since $F(a)$ is a functor, which is the identity at $\bar{F}(a \otimes b)$
and the image by $\gF$ of the identity at $a \otimes b$.\\ 

Regarding relations \ref{cong133} and \ref{cong134}.\\
For any arrows 
$\xymatrix{ a \ar[r]^f & a' \ar[r]^g & a''}$ in $\A$
and any object $b$ of $\B$,
$\gF$ sends 
$$\xymatrix{ 
(a \otimes b) \ar[r]^{f \otimes b} &
(a' \otimes b) \ar[r]^{g \otimes b} &
(a'' \otimes b)
}$$
to
$\xymatrix{ F(a)(b) \ar[r]^{{F(f)}_b} &
F(a')(b) \ar[r]^{F(g)_b} &
F(a'')(b) }$,
which is 
${F(g \circ f)}_b : F(a)(b) \rightarrow F(a'')(b)$
by functoriality of $F$,
which is the image by $\gF$ of 
$(g \circ f) \otimes b : 
(a \otimes b) \rightarrow (a'' \otimes b)$.\\

For any objects $a$ of $\A$ and $b$ of $\B$,
the image of $1_a \otimes b: a \otimes b \rightarrow a \otimes b$ 
by $\gF$ is the component in $b$ of the natural transformation
$F(1_a):F(a) \rightarrow F(a)$, which is the identity at $F(a)(b)$
since $F$ is a functor, 
which is the image by $\gF$ of the identity at $a \otimes b$.\\ 

For any arrows $f: a \rightarrow a'$ in $\A$ and
$g: b \rightarrow b'$ in $\B$,
the image by $\gF$ of Diagram \ref{cong135}
is 
$$
\xymatrix{
F(a)(b) 
\ar[r]^{{F(f)}_b}
\ar[d]_{F(a)(g)}
& 
F(a')(b)
\ar[d]^{ F(a')(g) }\\
F(a)(b')
\ar[r]_{{F(f)}_{b'}}
&
F(a')(b')
}
$$
which commutes by naturality of
$F(f):F(a) \rightarrow F(a'): \B \rightarrow \C$.\\

Regarding relations \ref{cong10}.\\

Naturality of $\alpha$.\\
For any $f: b \rightarrow b'$ in $\B$, 
the diagram in $\Gra_{\A,\B}$ 
$$\xymatrix{
\un \ar[r]^-{\alpha_b}
\ar[rd]_{ \alpha_{b'} } 
& \un_{\A} \otimes b 
\ar[d]^{1 \otimes f}\\
& \un_{\A} \otimes b'
}$$
is sent by $\gF$ to
$$\xymatrix{
\un_{\C} \ar[r]^-{F^0_b}
\ar[rd]_{ F^0_{b'} }
& F(\un_{\A})(b)
\ar[d]^{ F(\un_{\A})(f) }\\
& F(\un_{\A})(b')
}$$
which commutes by naturality of 
$F^0: \un \rightarrow F(\un_{\A}): \B \rightarrow \C$.\\

Naturality of $\beta$.\\
For any $f: a \rightarrow a'$ in $\A$, 
the diagram in $\Gra_{\A,\B}$
 $$\xymatrix{
\un \ar[r]^-{\beta_a} 
\ar[rd]_{ \beta_{a'} }
& a \otimes \un_{\B} 
\ar[d]^{f \otimes 1}\\
& a' \otimes \un_{\B}
}$$
is sent by $\gF$ to
$$\xymatrix{
\un_{\C} \ar[r]^-{{F(a)}^0} 
\ar[rd]_{ {F(a')}^0 }
& F(a)(\un_{\B})
\ar[d]^{ {F(f)}_{\un_{\B}} }\\
& F(a')(\un_{\B})
}$$
which commutes according to Axiom \ref{monat7}
for the monoidal natural transformation 
$F(f):F(a) \rightarrow F(a'): \B \rightarrow \C$.\\

Naturalities of $\gamma$.\\
For any $f:a \rightarrow c$ in $\A$ and any 
objects $a'$ of $\A$ and $b$ of $\B$, the diagram 
$$\xymatrix{
(a \otimes b) \otimes (a' \otimes b) 
\ar[r]^-{\gamma_{a,a',b}}
\ar[d]_{(f \otimes 1) \otimes 1} 
&
(a \otimes a') \otimes b 
\ar[d]^{(f \otimes 1) \otimes 1}\\
(c \otimes b) \otimes (a' \otimes b)
\ar[r]_-{\gamma_{c,a',b}} &
(c \otimes a') \otimes b
}$$
is sent by $\gF$ to
$$\xymatrix{
F(a)(b) \otimes F(a')(b) 
\ar[r]^-{ {(F^2_{a,a'})}_b }
\ar[d]_{{F(f)}_b \otimes 1} 
&
F(a \otimes a')(b) 
\ar[d]^{{F(f \otimes 1)}_b}\\
F(c)(b) \otimes F(a')(b)
\ar[r]_-{{(F^2_{c,a'})}_b} &
F(c \otimes a')(b)
}$$
which is pointwise in $b$ the diagram in $[\B,\C]$
$$\xymatrix{
F(a) \Box F(a') 
\ar[r]^-{ (F^2_{a,a'}) }
\ar[d]_{F(f) \Box 1} 
&
F(a \otimes a') 
\ar[d]^{F(f \otimes 1)}\\
F(c) \Box F(a')
\ar[r]_-{ (F^2_{c,a'}) } &
F(c \otimes a')
}$$
which commutes by naturality in $a$ 
of the collection of $F^2_{a,a'}: 
F(a) \Box F(a') \rightarrow F(a \otimes a')$
in $[\B, \C]$.\\

Similarly the images by $\gF$ of the diagram for the relations for 
the naturality of the 
$\gamma_{a,a',b}$ in $a'$ are commutative diagrams
according 
to the naturalities in $a'$ of the collection of $F^2_{a,a'}$.\\

For any objects $a,a'$ of $\A$ and any arrow $g: b \rightarrow c$ in $\B$,
the diagram
$$\xymatrix{
(a \otimes b) \otimes (a' \otimes b) 
\ar[r]^-{ \gamma_{a,a',b} }
\ar[d]_{(1 \otimes g) \otimes (1 \otimes g)} &
(a \otimes a') \otimes b \ar[d]^{ (1 \otimes g) }\\
(a \otimes c) \otimes (a' \otimes c)
\ar[r]_-{ \gamma_{a,a',c} } & 
(a \otimes a') \otimes c
}$$
is sent by $\gF$ to  
$$\xymatrix{
F(a)(b) \otimes F(a')(b) 
\ar[r]^-{ { ( F^2_{a,a'} ) }_b }
\ar[d]_{F(a)(g) \otimes F(a')(g)} &
F(a \otimes a')(b) 
\ar[d]^{F(a \otimes a')(g)}\\
F(a)(c) \otimes F(a')(c)
\ar[r]_-{{( F^2_{a,a'})}_c } & 
F(a \otimes a')(c)
}$$
which commutes since $F^2_{a,a'}$ is a 
natural transformation
$F^2_{a,a'}: F(a) \Box F(a') 
\rightarrow F(a \otimes a'): \B \rightarrow \C.$\\

Naturalities of $\delta$.\\
For any arrow $f: a \rightarrow c$ in $\A$ and any objects
$b, b'$ in $\B$, the diagram
$$
\xymatrix{
(a \otimes b) \otimes (a \otimes b') \ar[r]^-{\delta_{a,b,b'}}
\ar[d]_{ (f \otimes 1) \otimes (f \otimes 1) } &
a \otimes (b \otimes b') \ar[d]^{f \otimes 1}\\
(c \otimes b) \otimes (c \otimes b') \ar[r]_-{\delta_{c,b,b'}} &
c \otimes (b \otimes b')
}$$
is sent by $\gF$ to
$$\xymatrix{
F(a)(b) \otimes F(a)(b') \ar[r]^-{ {F(a)}^2_{b,b'} }
\ar[d]_{ {F(f)}_b \otimes {F(f)}_{b'} } &
F(a)(b \otimes b') \ar[d]^{ {F(f)}_{b \otimes b'} }\\
F(c)(b) \otimes F(c)(b') \ar[r]_-{ {F(c)}^2_{b,b'}} &
F(c)(b \otimes b')
}$$
which commutes
according to Axiom \ref{monat6}
for the monoidal natural transformation
$F(f): F(a) \rightarrow F(c): \B \rightarrow \C$.\\

For any objects $a$ in $\A$ and $b'$ in $\B$ and
any arrow $g: b \rightarrow c$ in $\B$, the diagram
$$
\xymatrix{
(a \otimes b) \otimes (a \otimes b') \ar[r]^-{\delta_{a,b,b'}}
\ar[d]_{ (1 \otimes g) \otimes 1 } &
a \otimes (b \otimes b') \ar[d]^{1 \otimes (g \otimes 1)}\\
(a \otimes c) \otimes (a \otimes b') \ar[r]_-{\delta_{a,c,b'}} &
a \otimes (c \otimes b')
}$$
is sent by $\gF$ to
$$
\xymatrix{
F(a)(b) \otimes F(a)(b') \ar[r]^-{ {F(a)}^2_{b,b'} }
\ar[d]_{ F(a)(g) \otimes 1 } &
F(a)(b \otimes b') \ar[d]^{F(a)(g \otimes 1)}\\
F(a)(c) \otimes F(a)(b') \ar[r]_-{ {F(a)}^2_{c,b'} } &
F(a)(c \otimes b')
}$$
which commutes according to the naturality of 
the collection of $F(a)^2_{b,b'}$ in $b$ since $F(a)$ is monoidal.\\

Similarly, the images by $\gF$ of diagrams for the relations 
for the naturalities of the $\delta_{a,b,b'}$ 
in $b'$ are commutative diagrams since $F(a)$ is monoidal.\\

For any objects $a$ in $\A$ 
and $b,b',b''$ in $\B$, the image of Diagram 
\ref{cong113}
by $\gF$ is 
$$
\xymatrix{
F(a)(b) \otimes ( F(a)(b') \otimes F(a)(b'') ) 
 \ar[r]^-{ \ac }
 \ar[d]_{ 1 \otimes {F(a)}^2_{b',b''} } 
&
( F(a)(b) \otimes F(a)(b') ) \otimes F(a)(b'') 
 \ar[d]^{ { F(a) }^2_{b,b'} \otimes 1 } 
\\
F(a)(b) \otimes F(a)( b' \otimes b'' ) 
 \ar[d]_{ { F(a) }^2_{b, b' \otimes b''} } 
&
F(a)(b \otimes b') \otimes F(a)(b'') 
 \ar[d]^{ {F(a)}^2_{b \otimes b',b''} } 
\\
F(a)(b \otimes (b' \otimes b'')) 
 \ar[r]_-{ F(a)(\ac_{b,b',b''}) } 
&
F(a)( ( b \otimes b' ) \otimes b'' )
}
$$
which commutes according to Axiom \ref{mofun3}
for the monoidal functor $F(a)$.\\  

For any objects $a$ in $\A$ and $b$ in $\B$,
the image by $\gF$ of Diagram 
\ref{cong1141} is 
$$
\xymatrix{
F(a)(b) \otimes \un 
\ar[r]^-{r_{F(a)(b)}}
\ar[d]_{1 \otimes {F(a)}^0} 
& F(a)(b)\\
F(a)(b) \otimes F(a)(\un_{\B})
\ar[r]_-{ {F(a)}^2_{ b,\un_{\B}  } } &
F(a)(b \otimes \un_{\B})
\ar[u]_{F(a)(\rc_b)}
}
$$
which commutes according to Axiom \ref{mofun41}
for the monoidal functor $F(a)$.\\

Similarly for any objects $a$ in $\A$ and 
$b$ in $\B$, the image by $\gF$ of
Diagram \ref{cong1142} is a commutative diagram
according to Axiom \ref{mofun42} for the monoidal 
functor $F(a)$.\\ 

Given any objects $a$ in $\A$ and $b,b'$ in $\B$, 
the image by $\gF$ of the Diagram \ref{cong115}
is 
$$
\xymatrix{
F(a)(b) \otimes F(a)(b') 
\ar[r]^-{ { F(a) }^2_{b,b'} }
\ar[d]_{s_{F(a)(b),F(a)(b')}}
&
F(a)(b \otimes b')
\ar[d]^{F(a)(s_{b,b'})}
\\
F(a)(b') \otimes F(a)(b) 
\ar[r]_-{ {F(a)}^2_{b',b} }
&
F(a)(b' \otimes b)
}
$$ 
which commutes according to Axiom \ref{symofun5}
for the symmetric functor $F(a)$.\\

For any objects $a,a',a''$ of $\A$ and $b$ of $\B$, 
the image by $\gF$ of Diagram \ref{cong143} is
$$
\xymatrix{
F(a)(b) \otimes ( (F(a')(b) \otimes F(a'')(b)  )
\ar[r]^-{ \ac }
\ar[d]_{ 1 \otimes {(F^2_{a',a''})}_b } &
( F(a)(b) \otimes F(a')(b)) \otimes  F(a'')(b)
 \ar[d]^{ {(F^2_{a,a'})}_b \otimes 1}
\\  
F(a)(b) \otimes F(a' \otimes a'')(b)
\ar[d]_{  {( F^2_{a, a' \otimes a''} )}_b } 
&
F(a \otimes a')(b) \otimes F(a'')(b) 
\ar[d]^{  {(F^2_{a \otimes a', a''})}_b }\\
F(a \otimes (a' \otimes a''))(b) 
\ar[r]_-{ {F( \ac_{a,a',a''} )}_b } &
F((a \otimes a') \otimes a'')(b) 
}
$$
which is the pointwise version
of the diagram in $[\B,\C]$
 $$\xymatrix@C=5pc{
F(a) \Box ( (F(a') \Box F(a'') )
\ar[r]^-{ \ac_{ F(a),F(a'), F(a'')}}
\ar[d]_{ 1 \Box  F^2_{a',a''} } &
( F(a) \Box F(a') ) \Box  F(a'')
 \ar[d]^{ F^2_{a,a'} \Box 1}\\  
F(a) \Box F(a' \otimes a'')
\ar[d]_{   F^2_{a, a' \otimes a''} } &
F(a \otimes a') \Box F(a'') 
\ar[d]^{  F^2_{a \otimes a', a''} }\\
F(a \otimes (a' \otimes a'')) 
\ar[r]_-{ F( \ac_{a,a',a''}) } &
F((a \otimes a') \otimes a'')  
}
$$
which commutes according to Axiom \ref{mofun3}
for the monoidal functor
$F$.\\

For any objects $a$ of $\A$ and $b$ of $\B$, the 
image by $\gF$ of Diagram \ref{cong1441}
is 
$$\xymatrix@C=3pc{
F(a)(b) \otimes \un 
\ar[r]^-{ \rc_{F(a)(b)} }
\ar[d]_{1 \otimes F^0_b}
&
F(a)(b) \\
F(a)(b) \otimes F(\un_{\A})(b) 
\ar[r]_-{ {( F^2_{a,\un_{\A}} )}_b }
&
F(a \otimes \un_{\A})(b)
\ar[u]_{ {(F(\rc_a))}_b }
}
$$
which is the pointwise version of the 
diagram in $[\B,\C]$
$$\xymatrix{
F(a) \Box \un 
\ar[r]^-{\rc_{F(a)}}
\ar[d]_{1 \Box F^0} 
&
F(a) \\
F(a) \Box F( \un_{\A} ) 
\ar[r]_-{ F^2_{ a, \un_{\A} } }
&
F(a \otimes \un_{\A})
\ar[u]_{ F(\rc_a) }
}
$$
which is commutative according 
to Axiom \ref{mofun41} 
for the monoidal functor $F$.\\
  
Similarly for any objects $a$ in $\A$ and
$b$ in $\B$, the image by $\gF$ of Diagram 
\ref{cong1442} commutes according Axiom \ref{mofun42}
for the monoidal functor $F$.\\

For any objects $a$,$a'$ of $\A$ and any $b$ of $\B$,
the image by $\gF$ of Diagram \ref{cong145}
is 
$$
\xymatrix{
F(a)(b) \otimes F(a')(b) 
\ar[r]^-{ {(F^2_{a,a'})}_b }
\ar[d]_{ \syc_{ F(a)(b), F(a)(b') } }  
&
F(a \otimes a')(b)
\ar[d]^{ {F(\syc_{a,a'})}_b} \\
F(a')(b) \otimes F(a)(b) 
\ar[r]_-{{(F^2_{a',a})}_b }
&
F(a' \otimes a)(b)
}
$$
which is the pointwise version of the diagram
in $[\B,\C]$
$$
\xymatrix{
F(a) \Box F(a') 
\ar[r]^-{ F^2_{a,a'} }
\ar[d]_{ \syc_{ F(a), F(a) } }  
&
F(a \otimes a')
\ar[d]^{F( \syc_{a,a'} ) } \\
F(a') \Box F(a) 
\ar[r]_-{ F^2_{a',a} }
&
F(a' \otimes a)
}
$$
that commutes according to Axiom \ref{symofun5}
for the symmetric monoidal functor $F$.\\

For any objects $b, b'$ in $\B$,
the image of Diagram \ref{cong161} by $\gF$
is 
$$
\xymatrix@C=3pc{
\un \otimes \un 
\ar@{-}[r]^-{\cong}
\ar[d]_{ {F^0}_b \otimes {F^0}_{b'} }
&
\un
\ar[d]^{ {F^0}_{b \otimes b'} }
\\
F(\un_{\A})(b) \otimes F(\un_{\A})(b')
\ar[r]_-{{(F(\un_{\A}))}^2_{b,b'}   }
&
F(\un_{\A})(b \otimes b')
}
$$
which commutes according to Axiom \ref{monat6}
for the monoidal natural transformation
$F^0: \un \rightarrow F(\un_{\A}): \B \rightarrow \C$.\\

For any objects $a, a'$ in $\A$,
the image of Diagram \ref{cong162} by $\gF$
is 
$$
\xymatrix{
\un \otimes \un 
\ar@{-}[r]^{\cong}
\ar[d]_{ {(F(a))}^0 \otimes {(F(a'))}^0 }
&
\un
\ar[d]^{ { (F(a \otimes a')) }^0} 
\\
F(a)(\un_{\B}) \otimes F(a')(\un_{\B})
\ar[r]_-{(F^2_{a,a'})_{\un_{\B}}}
&
F(a \otimes a')(\un_{\B})
}
$$
which commutes according to Axiom \ref{monat7} for 
the monoidal natural transformation
$$F^2_{a,a'}: F(a) \Box F(a') \rightarrow F(a \otimes a'): 
\B \rightarrow \C.$$

Regarding relations \ref{cong17}. The image by $\gF$ of 
$\beta_{\un_{\A}}: \un \rightarrow \un_{\A} \otimes \un_{\B}$
is ${(F(\un_{\A}))}^0: \un \rightarrow F(\un_{\A})(\un_{\B})$
whereas the image by $\gF$ of
$\alpha_{\un_{\B}}: \un \rightarrow \un_{\A} \otimes \un_{\B}$
is 
${F^0}_{\un_{\B}}: \un \rightarrow F(\un_{\A})(\un_{\B})$,
the component in $\un_{\B}$ of the natural transformation 
$F^0: \un \rightarrow F(\un_{\A}): \B \rightarrow \C$.
These two arrows are equal since 
Axiom \ref{monat7} for 
the natural monoidal transformation
$F^0$ states exactly their equality.\\

For any objects $a,a'$ in $\A$ and 
$b,b'$ in $\B$,
the image by $\gF$ of Diagram
\ref{cong19}
is 
$$\xymatrix{
(F(a)(b) \otimes F(a)(b'))
\otimes 
(F(a')(b) \otimes F(a')(b'))
\ar@{-}[r]^{\cong} 
\ar[d]_{{F(a)}^2_{b,b'} \otimes {F(a')}^2_{b,b'} } &
(F(a)(b) \otimes F(a')(b))
\otimes 
(F(a)(b') \otimes F(a')(b'))
\ar[d]^{ {(F^2_{a,a'})}_b \otimes {(F^2_{a,a'})}_{b'} }\\
(F(a)(b \otimes b')) \otimes (F(a')(b \otimes b')) 
\ar[d]_{ {(F^2_{a,a'})}_{b \otimes b'} } &
F(a \otimes a')(b) 
\otimes 
F(a \otimes a')(b') 
\ar[ld]^-{ { F(a \otimes a') }^2_{b, b'} }\\
F(a \otimes a')(b \otimes b') &
 }$$
which is commutative according to Axiom
\ref{monat6} for the monoidal natural transformation 
$$F^2_{a,a'}: F(a) \Box F(a') \rightarrow F(a \otimes a'): \B \rightarrow \C.$$

That the relations $\sim$ obtained by the rules of 
expansion \ref{cong12} are sent by $\gF$ to 
commuting diagrams in $\C$, can be shown by induction 
according to the functoriality of tensor in $\C$.\\
 
Note that the functor $\bar{F}: \A \otimes \B \rightarrow \C$
just defined admits a {\em symmetric strict monoidal}
structure.\\

The functor $\Ext$ is defined on arrows of $\SMC(\A,[\B,\C])$ as follows.
Given any monoidal transformation between symmetric
functors 
$\sigma: F \rightarrow G : \A \rightarrow [ \B , \C]$, it
is sent by $\Ext$ to the monoidal natural transformation
$\bar{\sigma}: \bar{F} \rightarrow \bar{G}: \A \otimes \B \rightarrow
\C$ defined by induction on the structure
of the elements of $\Ver_{\A,\B}$ according to the following
rules.\\
- $\bar{\sigma}_{\un} : \bar{F}(\un) \rightarrow \bar{G}(\un)$
is the identity at $\bar{F}(\un) = \un = \bar{G}(\un)$.\\
- For any objects $a$ of $\A$ and $b$ of $\B$,
$\bar{\sigma}_{a \otimes b}: 
\bar{F}(a \otimes b) \rightarrow \bar{G}(a \otimes b)$ is the arrow
$$\xymatrix{
\bar{F}(a \otimes b) \ar@{=}[r] 
& F(a)(b) \ar[r]^{{(\sigma_a)}_b}  
& G(a)(b) \ar@{=}[r] &
\bar{G}(a \otimes b)
}$$
where the transformation 
$\sigma_a: F(a) \rightarrow G(a): \B \rightarrow \C$ 
is the component in $a$ of $\sigma$.\\
- For any $X,Y$ in $\Ver_{\A,\B}$,
the arrow
$$\bar{\sigma}_{X \otimes Y}: \bar{F}(X \otimes Y) \rightarrow 
\bar{G}(X \otimes Y)$$
is 
$$
\xymatrix{
 \bar{F}(X \otimes Y) \ar@{=}[r]  & 
\bar{F}(X) \otimes \bar{F}(Y)
\ar[r]^{ \bar{\sigma}_X \otimes \bar{\sigma}_Y } & 
\bar{G}(X) \otimes \bar{G}(Y) \ar@{=}[r] & 
\bar{G}(X \otimes Y).
}
$$

To show that the above assignments $\bar{\sigma}$
define a natural transformation
$\bar{F} \rightarrow \bar{G}$, 
it is enough to show that 
for any $h:X \rightarrow Y$ in $\Gra_{\A,\B}$,
the diagram in $\C$
\begin{tag}\label{natbsig}
$$\xymatrix{
\bar{F}(X) 
\ar[r]^-{\gF(h)}
\ar[d]_{\bar{\sigma}_X} 
&
\bar{F}(Y)
\ar[d]^{\bar{\sigma}_Y}
\\
\bar{G}(X) 
\ar[r]_-{\gG(h)}
&
\bar{G}(Y)
}$$
\end{tag}
commutes. We shall check this now.\\

For any $X,Y,Z$ in $\Ver_{\A,\B}$, 
Diagram \ref{natbsig} for $h = \ac_{X,Y,Z}$
is 
$$
\xymatrix@C=5pc{ 
\bar{F}(X) \otimes ( \bar{F}(Y) \otimes \bar{F}(Z) )
\ar[r]^-{ \ac_{ \bar{F}(X), \bar{F}(Y), \bar{F}(Z) } }
\ar[d]_{  \bar{\sigma}_X \otimes (\bar{\sigma}_Y \otimes \bar{\sigma}_Z) }
&
(\bar{F}(X) \otimes \bar{F}(Y)) \otimes \bar{F}(Z)
\ar[d]^{ ( \bar{\sigma}_X \otimes \bar{\sigma}_Y ) \otimes \bar{\sigma}_Z }
\\
\bar{G}(X) \otimes ( \bar{G}(Y) \otimes \bar{G}(Z) )
\ar[r]_-{ \ac_{ \bar{G}(X), \bar{G}(Y), \bar{G}(Z) } }
&
(\bar{G}(X) \otimes \bar{G}(Y)) \otimes \bar{G}(Z)
}
$$
which commutes due to the naturality of $\ac$.\\

For any $X$ in $\Ver_{\A,\B}$,
Diagram \ref{natbsig} for $h = \rc_X$
is 
$$
\xymatrix{ 
\bar{F}(X) \otimes \un
\ar[r]^-{ \rc_{ \bar{F}(X)  } }
\ar[d]_{  \bar{\sigma}_X \otimes 1 }
&
\bar{F}(X)
\ar[d]^{  \bar{\sigma}_X  }
\\
\bar{G}(X) \otimes \un
\ar[r]_-{\rc_{ \bar{G}(X) } }
&
\bar{G}(X)
}
$$
which commutes due to the naturality of $\rc$.\\

Similarly, for any $X$ in $\Ver_{\A,\B}$,
Diagram \ref{natbsig} for $h = \lc_X$
is commutative due to the naturality of $\lc$.\\

For any $X,Y$ in $\Ver_{\A,\B}$,
Diagram \ref{natbsig} for $h = \syc_{X,Y}$
is 
$$
\xymatrix@C=3pc{ 
\bar{F}(X) \otimes \bar{F}(Y)
\ar[r]^-{ \syc_{ \bar{F}(X), \bar{F}(Y) } }
\ar[d]_{  \bar{\sigma}_X \otimes \bar{\sigma}_Y }
&
\bar{F}(Y) \otimes \bar{F}(X)
\ar[d]^{  \bar{\sigma}_Y \otimes \bar{\sigma}_X }
\\
\bar{G}(X) \otimes \bar{G}(Y)
\ar[r]_-{ \syc_{ \bar{G}(X), \bar{G}(Y) }  }
&
\bar{G}(Y) \otimes \bar{G}(X)
}
$$
which commutes due to the naturality of $\syc$.\\

From the above it is immediate that
Diagram \ref{natbsig} also commutes for
arrows $h$ of form $\aci_{X,Y,Z}$, $\rci_X$ and $\lci_X$.\\ 

For any object $b$ of $\B$,
Diagram \ref{natbsig} for
$h = \alpha_b: \un \rightarrow \un_{\A} \otimes b$ 
is 
$$
\xymatrix{
\un 
\ar[r]^-{F^0_b}
\ar@{=}[d] 
&
F(\un_A)(b)
\ar[d]^{{(\sigma_ {\un _{\A}})}_b}
\\
\un
\ar[r]_-{G^0_b}
&
G(\un_{\A})(b)
}
$$
which is the evaluation in $b$ of the diagram
in $[\B,\C]$
$$
\xymatrix{
\un 
\ar[r]^-{F^0}
\ar[rd]_-{G^0} 
&
F(\un_\A)
\ar[d]^{ \sigma_{\un_{\A}} }
\\
&
G(\un_{\A})
}
$$
which commutes according to Axiom \ref{monat7}
for the monoidal natural transformation 
$\sigma$.\\

For any object $a$ of $\A$,
Diagram \ref{natbsig} for 
$h = \beta_a: \un \rightarrow a \otimes \un_{\B}$
is 
$$
\xymatrix{
\un 
\ar[r]^-{{F(a)}^0}
\ar@{=}[d] 
&
F(a)(\un_{\B})
\ar[d]^{{(\sigma_a)}_{\un_{\B}}}
\\
\un
\ar[r]_-{{G(a)}^0}
&
G(a)(\un_{\B})
}
$$
which commutes according to Axiom \ref{monat7}
for the monoidal natural transformation
$\sigma_a: F(a) \rightarrow G(a): \B \rightarrow \C$.\\

For any objects $a,a'$ of $\A$
and $b$ of $\B$, Diagram
\ref{natbsig} for 
$h = \gamma_{a,a',b}: (a \otimes b) \otimes (a' \otimes b)
\rightarrow (a \otimes a') \otimes b$
is 
$$
\xymatrix{
F(a)(b) \otimes F(a')(b) 
\ar[r]^-{ {(F^2_{a,a'})}_b }
\ar[d]_{ {(\sigma_a)}_b \otimes {(\sigma_{a'})}_b  } 
&
F(a \otimes a')(b)
\ar[d]^{{(\sigma_{a \otimes a'})}_b}
\\
G(a)(b) \otimes G(a')(b) 
\ar[r]_-{ {(G^2_{a,a'})}_b }
&
G(a \otimes a')(b)
}
$$
which is the evaluation in $b$ of the diagram
in $[\B,\C]$
$$\xymatrix{
F(a) \Box F(a') 
\ar[r]^{ { F^2_{a,a'} } }
\ar[d]_{ {\sigma_a}  \Box   {\sigma_{a'} } } 
&
F(a \otimes a')
\ar[d]^{ {\sigma_{a \otimes a'} } }
\\
G(a) \Box G(a') 
\ar[r]_-{ G^2_{a,a'} }
&
G(a \otimes a')
}
$$
which commutes according to Axiom \ref{monat6}
for the monoidal natural transformation 
$\sigma$.\\

For any objects $a$ of $\A$
and $b, b'$ of $\B$,
Diagram \ref{natbsig} for 
$h = \delta_{a,b,b'}: (a \otimes b) \otimes (a \otimes b')
\rightarrow a \otimes (b \otimes b')$
is 
$$
\xymatrix{
F(a)(b) \otimes F(a)(b') 
\ar[r]^-{{F(a)}^2_{b,b'}}
\ar[d]_{ {(\sigma_a)}_b  \otimes {(\sigma_a)}_{b'}   } 
&
F(a)(b \otimes b')
\ar[d]^{{(\sigma_a)}_{b \otimes b'}}
\\
G(a)(b) \otimes G(a)(b')
\ar[r]_-{ {G(a)}^2_{b,b'}  }
&
G(a)(b \otimes b')
}
$$
which commutes according to Axiom 
\ref{monat6} for the natural transformation
$\sigma_a :F(a) \rightarrow G(a): \B \rightarrow \C$.\\

For any arrow $f: a \rightarrow a'$ in $\A$ and 
any object $b$ in $\B$, Diagram
\ref{natbsig} for 
$h = f \otimes b : a \otimes b 
\rightarrow a' \otimes b$
is 
$$
\xymatrix{
F(a)(b) 
\ar[r]^-{{F(f)}_b} 
\ar[d]_{{(\sigma_a)}_b}
& 
F(a')(b) 
\ar[d]^{{(\sigma_{a'})}_b}
\\
G(a)(b) 
\ar[r]_-{{G(f)}_b}
&
G(a')(b)
}
$$
which is the evaluation in $b$ 
of the diagram in $[\B,\C]$
$$
\xymatrix{
F(a) 
\ar[r]^-{F(f)} 
\ar[d]_{\sigma_a}
& 
F(a')
\ar[d]^{\sigma_{a'}}
\\
G(a)
\ar[r]_-{G(f)}
&
G(a')
}
$$
which commutes by naturality of 
$\sigma$.\\

For any object $a$ in $\A$ and any arrow
$g: b \rightarrow b'$ in $\B$, Diagram
\ref{natbsig} for 
$h = a \otimes g : a \otimes b 
\rightarrow a \otimes b'$
is 
$$
\xymatrix{
F(a)(b) 
\ar[r]^-{F(a)(g)} 
\ar[d]_{ {(\sigma_a)}_b }
& 
F(a)(b') 
\ar[d]^{ {(\sigma_a)}_{b'} }
\\
G(a)(b) 
\ar[r]_-{G(a)(g)}
&
G(a)(b')
}
$$
which commutes by naturality of  
$\sigma_a: F(a) \rightarrow G(a): \B \rightarrow \C$.\\

So far we have shown that Diagram \ref{natbsig}
commutes for any arrow $h$ in $\Gra_1 \cup \Gra_2 \cup \Gra_3$.
That it is also the case for all arrows $h$ of $\Gra$
is now proved by induction.\\

For any object $X$ and 
any arrow $f: Y \rightarrow Z$ in $\Gra_{\A,\B}$,
Diagram \ref{natbsig} for 
$h = X \otimes f : X \otimes Y 
\rightarrow X \otimes Z$
is 
$$
\xymatrix{
\bar{F}(X) \otimes \bar{F}(Y)
\ar[r]^-{ 1 \otimes \gF(f)    }
\ar[d]_{ \bar{\sigma}_X \otimes \bar{\sigma}_Y }
&
\bar{F}(X) \otimes \bar{F}(Z)
\ar[d]^{ \bar{\sigma}_X \otimes \bar{\sigma}_Z }
\\
\bar{G}(X) \otimes \bar{G}(Y)
\ar[r]_-{ 1 \otimes \gG(f) }
&
\bar{G}(X) \otimes \bar{G}(Z)
.}
$$
According to the functoriality of
tensor in $\C$, this diagram commutes
if the diagram
$$
\xymatrix{
\bar{F}(Y)
\ar[r]^-{ \gF(f)    }
\ar[d]_{ \bar{\sigma}_Y }
&
\bar{F}(Z)
\ar[d]^{ \bar{\sigma}_Z }
\\
\bar{G}(Y)
\ar[r]_-{ \gG(f) } &
\bar{G}(Z)
}
$$
commutes.
Similarly one shows that Diagram \ref{natbsig} for 
$h = f \otimes X : Y \otimes X 
\rightarrow Z \otimes X$
is commutative if Diagram \ref{natbsig} for $h = f$
commutes.\\

According to its inductive definition
$\bar{\sigma}: \bar{F} \rightarrow \bar{G}$ 
is trivially monoidal between
strict monoidal functors.\\

Let us give a universal characterisation of the $\bar{F}$ and $\bar{\sigma}$
that we have just been defined.\\

\begin{proposition}\label{characExt}
Given any symmetric monoidal functor $F: \A \rightarrow [\B, \C]$, 
the symmetric monoidal functor
$\bar{F}$ is the unique {\em strict} one 
$\A \otimes \B \rightarrow \C$
that renders commutative the  diagram in $\SMC$
 $$\xymatrix{
 \A 
 \ar[rd]_F
 \ar[r]^-{\eta}
 &
 [\B,\A \otimes \B]
 \ar[d]^{[\B,\bar{F}]}
\\
&
[\B,\C].
}$$
Given any 2-cell
$\sigma: F \rightarrow G: \A \rightarrow [\B,\C]$ in $\SMC$,
$\bar{\sigma}: \bar{F} \rightarrow \bar{G}$ is the unique 
2-cell in $\SMC$ such that $[\B,\bar{\sigma}] * \eta = \sigma$.
\end{proposition}
\pf
To see this, note that the commutation of the diagram 
in $\CAT$
 $$\xymatrix{
 \A 
 \ar[rd]_F
 \ar[r]^-{\eta}
 &
 \SMC( \B,\A \otimes \B )
 \ar[d]^{\SMC(\B,\bar{F})}
\\
&
\SMC(\B,\C)
}$$
is equivalent to the conjunction of the following facts:\\  
- $(1)$: For any object $a$ in $\A$, the underlying functors 
$F(a)$ and $\bar{F} \circ (\eta(a))$ are equal;\\ 
- $(2)$: For any object $a$ in $\A$, ${(\bar{F} \circ (\eta(a)))}^0 = F(a)^0$;\\
- $(3)$: For any object $a$ in $\A$, ${(\bar{F} \circ (\eta(a)))}^2 = F(a)^2$;\\ 
- $(4)$: For any arrow $f: a \rightarrow a'$
in $\A$, the natural transformations 
$F(f)$ and  $\bar{F} * (\eta(f))$ are equal.\\
Condition $(1)$ is equivalent to the conjunction of the 
following two conditions:\\
- $(1-1)$: For any objects $a$ in $\A$ and $b$ in $\B$,
$\bar{F}(a \otimes b) = F(a)(b)$;\\
and\\
- $(1-2)$: For any object $a$ in $\A$ and any arrow $g: b \rightarrow b'$
in $\B$, $\bar{F}(a \otimes g): \bar{F}(a \otimes b) \rightarrow 
\bar{F}(a \otimes b')$ is $F(a)(g): F(a)(b) \rightarrow F(a)(b')$.\\
So if $(1)$ holds, condition $(4)$ is just that:\\ 
- $(4')$ For any arrow $f: a \rightarrow a'$ in $\A$ and any object
$b$
of $\B$, 
${F(f)}_b: F(a)(b) \rightarrow F(a')(b)$ is 
equal to $\bar{F}(f \otimes b): \bar{F}(a \otimes b) 
\rightarrow \bar{F}(a' \otimes b)$.\\ 
In the case when $(1)$ holds and $\bar{F}$ is strict, condition 
$(2)$ is just equivalent to:\\
- $(2')$: For any object $a$ in $\A$, $\bar{F}$ sends
$\beta_a: \un \rightarrow a \otimes \un_{\B}$ to
${F(a)}^0: \un_{\C} \rightarrow F(a)(\un_{\B})$;\\
and condition $(3)$ is equivalent to:\\
- $(3')$: For any objects $a$ in $\A$ and $b$,$b'$ in 
$\B$, $\bar{F}$ sends the arrow 
$\delta_{a,b,b'}: (a \otimes b) \otimes (a \otimes b')
\rightarrow a \otimes ( b \otimes b')$ to ${{(Fa)}^2}_{b,b'}: 
Fa(b) \otimes Fa(b') \rightarrow Fa(b \otimes b')$.\\  

The commutation of the diagram in $\SMC$ of the proposition,
for a strict $\bar{F}$, is therefore equivalent to the conjunction of 
the above conditions
$(1-1)$, $(1-2)$, $(2')$, $(3')$ and $(4')$ and the two conditions:\\ 
- $(5)$: $([\B,\bar{F}] \circ \eta)^0 = F^0$ ;\\
and\\
- $(6)$:  $([\B,\bar{F}] \circ \eta)^2 = F^2$.\\
Condition $(5)$ is that:\\
- $(5')$: For any objects $b$ in $\B$,
$\bar{F}$ sends the arrow $\alpha_b: \un \rightarrow \un_{\A} \otimes b$
to ${F^0}_b : \un \rightarrow F(\un_{\A})(b)$.\\
Condition $(6)$ is that:\\
- $(6')$ For any objects $b$ in $\B$ and $a,a'$ in $\A$, 
$\bar{F}$ sends $\gamma_{a,a',b}: (a \otimes b) \otimes (a' \otimes b)
\rightarrow (a \otimes a') \otimes b$ to
${(F^2_{a,a'})}_b: Fa(b) \otimes Fa'(b) \rightarrow 
F(a \otimes a')(b)$.\\

Now remark that according to its inductive definition,
the monoidal functor $\bar{F}$, image of $F$ by $\Ext$, 
is the only strict one
satisfying conditions $(1-1)$, $(1-2)$, $(2')$, $(3')$,
$(4')$, $(5')$ and $(6')$ above.\\

Given monoidal transformations 
$\sigma: F \rightarrow G: \A \rightarrow [\B,\C]$,
and $\bar{\sigma}: \bar{F} \rightarrow \bar{G}$,
with $\bar{F}$, $\bar{G}$ the respective images of $F$ and $G$ by $\Ext$,
that $[\B,\bar{\sigma}] * \eta = \sigma$
just means that for any object $a$ in $\A$, the natural transformation
$\sigma_a : F(a) \rightarrow G(a): \B \rightarrow \C$
is $\bar{\sigma} * \eta(a)$, which is equivalent to the 
assertion that:\\
- $(7)$ For any objects $a$ in $\A$ and $b$ in $\B$,
the arrow $\bar{\sigma}_{a \otimes b}: \bar{F}(a \otimes b) \rightarrow
\bar{G}(a \otimes b)$ is 
${ ( \sigma_a ) }_b: Fa(b) \rightarrow Ga(b)$.\\

According to its inductive definition, 
the monoidal natural $\bar{\sigma}: \bar{F} \rightarrow \bar{G}$,
image of any $\sigma: F \rightarrow G$ by $\Ext$, is the 
unique one satisfying the condition $(7)$ above.\\
\epf

One has this alternative characterisation
from Proposition \ref{characExt} and Lemma \ref{dualitynatABC0}.
\begin{remark}\label{characExt3}
Given any symmetric monoidal $F: \A \rightarrow [\B,\C]$,
the symmetric monoidal functor 
$\bar{F}: \A \otimes \B \rightarrow \C$ is the only strict 
one that renders
commutative the diagram in $\SMC$
$$
\xymatrix{
\B 
\ar[r]^-{\eta^*}
\ar[rd]_-{F^*}
&
[\A, \A \otimes \B]
\ar[d]^{[\A,\bar{F}]}
\\
&
[\A, \C].
}
$$
\end{remark}

Eventually one has also this last characterisation 
of the $\bar{F}$ and $\bar{\sigma}$
using two diagrams in $\CAT$ rather
than one in $\SMC$.
\begin{proposition}\label{characExt2}
Given any symmetric monoidal functor $F: \A \rightarrow [\B, \C]$,
the symmetric monoidal functor 
$\bar{F}: \A \otimes \B \rightarrow \C$ is the unique 
strict one such that the following two diagrams in $\CAT$ 
commute
$$\xymatrix{
\A 
 \ar[rd]_-F
 \ar[r]^-{\eta}
 &
 \SMC(\B,\A \otimes \B)
 \ar[d]^{\SMC(\B,\bar{F})}
\\
&
\SMC(\B,\C)
}$$
and 
$$
\xymatrix{
\B
\ar[rd]_-{\Fs}
\ar[r]^-{\eta^*}
&
\SMC(\A,\A \otimes \B)
\ar[d]^{\SMC(\A,\bar{F})}
\\
&
\SMC(\A,\C).
}
$$
\end{proposition}
\pf
From the results of section \ref{inthomSMC} detailing 
the monoidal structure of $\Fs$, one checks that the commutation 
of the second diagram in $\CAT$ in the proposition implies the 
following:\\ 
- For any objects $b$ in $\B$,
$\bar{F}$ sends the arrow $\alpha_b: \un \rightarrow \un_{\A} \otimes b$
to ${F^0}_b : \un \rightarrow F(\un_{\A})(b)$;\\
- For any objects $b$ in $\B$ and $a,a'$ in $\A$, 
$\bar{F}$ sends $\gamma_{a,a',b}: (a \otimes b) \otimes (a' \otimes b)
\rightarrow (a \otimes a') \otimes b$ to
${(F^2_{a,a'})}_b: Fa(b) \otimes Fa'(b) \rightarrow 
F(a \otimes a')(b)$.
\epf

From Proposition \ref{characExt} the following becomes immediate. 
\begin{corollary}\label{coroadjten}   
The assignments 
$\bar{(-)}$ define a functor 
$\Ext: \SMC(\A,[\B,\C]) \rightarrow \SMC(\A \otimes \B, \C)$
which factorises as 
$$\xymatrix{
\SMC(\A,[\B,\C])
\ar[r]^-{\cong}
&
\SSMC(\A \otimes \B, \C)
\ar[r]
& 
\SMC(\A \otimes \B, C)},$$
where the functor on the left is an isomorphism.\\
\end{corollary}

It is also immediate from Proposition 
\ref{characExt} that the composite mere functor
$\Res \circ \Ext$ is the identity of $\SMC(\A,[ \B, \C])$.\\

Let $(-)^{\bigtriangledown}$ denote the composite functor 
$\Ext \circ \Res$. We show now the existence
of a natural transformation 
$\epsilon :  (-)^{\bigtriangledown} \rightarrow 1 : 
\SMC(\A \otimes \B,\C) \rightarrow \SMC(\A \otimes \B,\C)$.\\

For any symmetric monoidal $F: \A \otimes \B \rightarrow \C$,
the monoidal natural transformation 
$\epsilon_F: \Fb \rightarrow F: \A \otimes \B \rightarrow \C$
is defined by induction on the structure of 
the objects of $\A \otimes \B$
according to the following rules.
We shall drop the subscript $F$ and write simply $\epsilon$
when there is no ambiguity.\\
- For any objects $a$ of $\A$ and $b$ of $\B$,
$\epsilon_{a \otimes b}$ is the identity:
$\xymatrix{ \Fb (a \otimes b) \ar@{=}[r] & F(a \otimes b)}$.\\
- $\epsilon_{\un}$ is the arrow
$$\xymatrix{\Fb(\un) \ar@{=}[r] &  \un \ar[r]^{F^0} & F(\un) }.$$
- For any objects $X,Y$ of $\A \otimes \B$, 
$\epsilon_{X \otimes Y}$ is the arrow
$$\xymatrix{ \Fb(X) \otimes \Fb(Y) 
\ar[r]^-{\epsilon_X \otimes \epsilon_Y}
& 
F(X) \otimes F(Y)
\ar[r]^-{F^2_{X,Y}}  
&
F(X \otimes Y).  
}$$

To check the naturality of 
$\epsilon_F: \Fb \rightarrow F: \A \otimes \B \rightarrow \C$,
one needs to show that for any arrow $h:X \rightarrow Y$ of
$\Gra_{\A,\B}$ the following diagram in $\C$ commutes
\begin{tag}\label{epsnat}
$$\xymatrix{
\Fb(X)
\ar[d]_{\epsilon_X}
\ar[r]^-{\Fb(h)}
&
\Fb(Y)
\ar[d]^{\epsilon_Y}
\\
F(X)
\ar[r]_-{F(h)}
&
F(Y)
}
$$
\end{tag}
where $h$ also denotes the corresponding arrow 
of $\A \otimes \B$.\\

For any $X,Y,Z$ in $\Ver_{\A,\B}$, Diagram \ref{epsnat} for 
$h = \ac_{X,Y,Z}: X \otimes (Y \otimes Z) \rightarrow 
(X \otimes Y) \otimes Z$
is 
$$\xymatrix@C=4pc{
\Fb(X \otimes (Y \otimes Z))
\ar[d]_{\epsilon_{X \otimes (Y \otimes Z)}}
\ar[r]^-{\Fb(\ac_{X,Y,Z})}
&
\Fb((X \otimes Y) \otimes Z )
\ar[d]^{\epsilon_{(X \otimes Y) \otimes Z  }}
\\
F ( X \otimes (Y \otimes Z) )
\ar[r]_-{F(\ac_{X,Y,Z})}
&
F ( (X \otimes Y) \otimes Z )
}
$$
which is the external diagram in the pasting
{\small
$$\xymatrix@C=7pc{
\Fb(X) \otimes ( \Fb(Y) \otimes \Fb(Z) )
\ar[d]_{  \epsilon_X \otimes (\epsilon_Y \otimes \epsilon_Z) }
\ar[r]^{\ac_{\Fb(X),\Fb(Y),\Fb(Z)} }
&
(( \Fb(X) \otimes \Fb(Y) ) \otimes \Fb(Z) )
\ar[d]^{ ( \epsilon_X \otimes \epsilon_Y ) \otimes \epsilon_Z ) }
\\
F(X) \otimes ( F(Y) \otimes F(Z) )
\ar[d]_{1 \otimes F^2_{Y,Z}}
\ar[r]^{\ac_{F(X),F(Y),F(Z)} }
&
(( F(X) \otimes F(Y) ) \otimes F(Z) )
\ar[d]^{F^2_{X,Y} \otimes 1 }
\\
F(X) \otimes F(Y \otimes Z) 
\ar[d]_{F^2_{X,Y \otimes Z}}
&
( F(X \otimes Y) ) \otimes F(Z) 
\ar[d]^{F^2_{X \otimes Y,Z} }
\\
F( X \otimes (Y \otimes Z) )
\ar[r]_-{F(\ac_{X,Y,Z})}
&
F( (X \otimes Y) \otimes Z )
}
$$
}
where the top square commutes by naturality of
$\ac$ and the bottom one also does according to Axiom
\ref{mofun3} for $F$.\\

For any $X$ in $\Ver_{\A,\B}$, Diagram
\ref{epsnat} for 
$h = \rc_{X}: X \otimes \un \rightarrow X$
is 
$$\xymatrix{
\Fb(X \otimes \un)
\ar[d]_{\epsilon_{X \otimes \un}}
\ar[r]^-{\Fb(\rc_X)}
&
\Fb(X)
\ar[d]^{\epsilon_X}
\\
F(X \otimes \un)
\ar[r]_-{F(\rc_X)}
&
F(X)
}
$$
which is the external diagram in the pasting
$$
\xymatrix{
\Fb(X) \otimes \un 
\ar[rr]^-{\rc_{\Fb(X)}}
\ar[d]_{\epsilon_X \otimes F^0}
\ar[rd]^{\epsilon_X \otimes 1} & &  
\Fb(X) \ar[dd]^{\epsilon_X}\\
F(X) \otimes F(\un)
\ar[d]_{F^2_{X,\un}} 
& 
F(X) \otimes \un
\ar[l]^{1 \otimes F^0}
\ar[rd]^{r_{F(X)}}
&\\
F(X \otimes \un) \ar[rr]_-{F(r_X)} 
& & 
F(X)
}
$$
where the bottom left diagram commutes according to
Axiom \ref{mofun41} for $F$ 
whereas the right one commutes by naturality
of $\rc$.\\

Similarly one shows that 
for any $X$ in $\Ver_{\A,\B}$, Diagram
\ref{epsnat} for 
$h = \lc_{X}: X \otimes \un \rightarrow X$
commutes according to Axiom \ref{mofun42}
for $F$ and the naturality of $\lc$.\\

For any $X,Y$ in $\Ver_{\A,\B}$, Diagram
\ref{epsnat} for 
$h = \syc_{X,Y}: X \otimes Y \rightarrow Y \otimes X$
is 
$$\xymatrix@C=3pc{
\Fb(X \otimes Y)
\ar[d]_{\epsilon_{X \otimes Y}}
\ar[r]^-{\Fb(\syc_{X,Y}) }
&
\Fb(Y \otimes X)
\ar[d]^{\epsilon_{Y \otimes X}}
\\
F(X \otimes Y)
\ar[r]_-{F(\syc_{X,Y})}
&
F(Y \otimes X)
}
$$
which is the external diagram in the pasting
$$\xymatrix@C=4pc{
\Fb(X) \otimes \Fb(Y)
\ar[d]_{ \epsilon_X \otimes \epsilon_Y }
\ar[r]^-{ \syc_{\Fb(X),\Fb(Y)} }
&
\Fb(Y) \otimes \Fb(X)
\ar[d]^{ \epsilon_Y \otimes \epsilon_X }
\\
F(X) \otimes F(Y)
\ar[d]_{ F^2_{X,Y} }
\ar[r]^-{ \syc_{F(X),F(Y)} }
&
F(Y) \otimes F(X)
\ar[d]^{F^2_{Y,X}}
\\
F(X \otimes Y)
\ar[r]_-{F(\syc_{X,Y})}
&
F(Y \otimes X)
}
$$
where the top square commutes by naturality 
of $\syc$ and the bottom one also does according to 
Axiom \ref{symofun5} for $F$.\\ 

For any object $b$ of $\B$, Diagram
\ref{epsnat} for 
$h = \alpha_b : \un \rightarrow \un_{\A} \otimes b$
is
$$
\xymatrix{
\Fb(\un)
\ar[d]_{\epsilon_{\un}}
\ar[r]^-{ \Fb(\alpha_b) }
&
\Fb(\un_{\A} \otimes b)
\ar[d]^{ \epsilon_{(\un_{\A} \otimes b) } }
\\
F(\un)
\ar[r]_-{F(\alpha_b)}
&
F(\un_{\A} \otimes b).
}
$$
We prove below that it commutes.
The arrow $F(\alpha_b) \circ \epsilon_{\un}$ is
$\xymatrix{
\un
\ar[r]^-{F^0}
&
F(\un)
\ar[r]^-{F(\alpha_b)}
&
F(\un_{\A} \otimes b).
}$
On the other hand 
the arrow $\epsilon_{(\un_{\A} \otimes b)} \circ \Fb(\alpha_b)$
rewrites
$$\xymatrix{
\un 
\ar@{=}[r]
&
\Fb(\un)
\ar[r]^-{\Fb(\alpha_b)}
&
\Fb(\un_{\A} \otimes b)
\ar@{=}[r]
&
F(\un_{\A} \otimes b)
.}$$ 
The functor $\Fb$ sends the arrow 
$\alpha_b: \un \rightarrow \un_{\A} \otimes b$ 
of $\A \otimes \B$ 
to the arrow 
${\Res(F) }^0_b: \un \rightarrow \Res(F)(\un_{\A})(b)$ of $\C$,
component in $b$ of  
the monoidal natural transformation ${\Res(F)}^0: 
\un \rightarrow \Res(F)(\un_{\A}):\B \rightarrow \C$
which is part the monoidal structure of $\Res(F): \A \rightarrow [\B,\C]$.
$\Res(F)$ has been defined as the composite
$\xymatrix{
\A 
\ar[r]^-{\eta} 
& 
[\B, \A \otimes \B] 
\ar[r]^-{[\B,F ]}
&
[\B,\C]
}
$
thus 
${\Res(F)}^0$ is the composite arrow in $[\B,\C]$
$$\xymatrix@C=3pc{
\un
\ar[r]^-{ {[\B,F]}^0 } 
& 
[\B,F]( \un_{ [\B, \A \otimes \B]  } ) 
\ar[r]^-{ [\B,F ]( {\eta }^0) }
&
[\B,F](\eta( \un_{\A}) ).
}
$$
Now ${[\B,F]}^0: \un \rightarrow  
F \circ  \un_{[\B, \A \otimes \B] }$ 
is pointwise in $b$ the arrow 
$F^0: \un \rightarrow F(\un_{\A \otimes \B})$ in $\C$,
whereas
$
{\eta}^0: \un \rightarrow
\eta(\un_{\A}): \B \rightarrow \A \otimes \B$
takes value
in $b$ the arrow
$\alpha_b: \un \rightarrow \un_{\A} \otimes b$ in $\A \otimes \B$
and thus
$[\B,F]({\eta}^0)$ takes value in $b$ the arrow
$F(\alpha_{b}): F(\un) \rightarrow F(\un_{\A} \otimes b)$ in $\C$.\\

For any object $a$ of $\A$, Diagram
\ref{epsnat} for 
$h = \beta_a : \un \rightarrow a \otimes \un_{\B}$
is 
$$
\xymatrix{
\Fb(\un)
\ar[d]_{\epsilon_{\un}}
\ar[r]^-{\Fb(\beta_a)}
&
\Fb(a \otimes \un_{\B})
\ar[d]^{\epsilon_{(a \otimes \un_{\B})} }
\\
F(\un)
\ar[r]_-{F(\beta_a)}
&
F(a \otimes \un_{\B})
}
$$
and commutes as shown below.
The arrow $F(\beta_a) \circ \epsilon_{\un}$ rewrites
$$
\xymatrix{
\Fb(\un)
\ar@{=}[r]
&
\un
\ar[r]^-{F^0}
&
F(\un)
\ar[r]^-{F(\beta_a)}
&
F(a \otimes \un_{\B})
}
$$
whereas the arrow 
$\epsilon_{a \otimes \un_{\B}} \circ \Fb(\beta_a)$ 
rewrites
$$
\xymatrix{
\Fb(\un)
\ar[r]^-{\Fb(\beta_a)}
&
\Fb(a \otimes \un_{\B})
\ar@{=}[r]
&
F(a \otimes \un_{\B}).
}
$$
The functor $\Fb$ sends $\beta_a: \un \rightarrow a \otimes \un_{\B}$
to the arrow
${(\Res(F)(a))}^0: \un \rightarrow 
\Res(F)(a)(\un_{\B})$
in $\C$, part of the monoidal structure of $\Res(F)(a)$.
Now $\Res(F)(a)$ is the composite 
$\xymatrix{ \B 
\ar[r]^-{\eta(a)} 
& 
\A \otimes \B 
\ar[r]^-F
& \C}$ and thus the arrow 
${(\Res(F)(a))}^0$
is the composite 
$\xymatrix{
\un
\ar[r]^-{F^0}
& 
F(\un)
\ar[r]^-{F(\beta_a)} 
&
F(a \otimes \un_{\B})
}$ 
since $(\eta(a))^0 = \beta_a$.\\

For any objects $a,a'$ of $\A$ and  
$b$ of $\B$, Diagram
\ref{epsnat} for 
$h = \gamma_{a,a',b} : (a \otimes b) \otimes (a' \otimes b) 
\rightarrow (a \otimes a') \otimes b$
is 
$$
\xymatrix@C=3pc{
\Fb( (a \otimes b) \otimes (a' \otimes b) )
\ar[d]_{\epsilon_{(a \otimes b) \otimes (a' \otimes b) }}
\ar[r]^-{\Fb(\gamma_{a,a',b})}
&
\Fb((a \otimes a') \otimes b)
\ar[d]^{\epsilon_{(a \otimes a') \otimes b}}
\\
F((a\otimes b) \otimes (a' \otimes b))
\ar[r]_-{F( \gamma_{a,a',b} ) }
&
F((a \otimes a') \otimes b)
}
$$
and commutes as shown below.
The arrow $\epsilon_{ (a \otimes b) \otimes (a' \otimes b) }$
is 
$$\xymatrix@C=3pc{
\Fb((a \otimes b) \otimes (a' \otimes b))
\ar@{=}[r] 
& 
{F(a \otimes b) \otimes F(a' \otimes b)}
\ar[r]^-{F^2_{a \otimes b, a' \otimes b}} 
& 
F( (a \otimes b) \otimes (a' \otimes b)).
}$$

One has $\Fb((a \otimes a') \otimes b) = 
F((a \otimes a') \otimes b)$
and the arrow $\epsilon_{(a \otimes a') \otimes b}$
is the identity at $F((a \otimes a') \otimes b)$.
The functor $\Fb$ sends $\gamma_{a,a',b}$ to  
$${( {(\Res(F))}^2_{a,a'} ) }_b: \Res(F)(a)(b) \otimes \Res(F)(a')(b)
\rightarrow \Res(F)(a \otimes a')(b)$$ 
which is the component in $b$ of the monoidal
natural transformation
$${\Res(F)}^2_{a,a'}: \Res(F)(a) \Box \Res(F)(a') \rightarrow
\Res(F)(a \otimes a'): \B \rightarrow \C.$$
This last one is the composite
$$\xymatrix@C=4pc{
[\B,F](\eta(a))  \Box [\B,F](\eta(a'))
\ar[r]^-{ {[\B,F]}^2_{\eta(a), \eta(a')}} &
[\B,F]( \eta(a) \Box \eta(a') )
\ar[r]^-{ [\B,F]( {\eta}^2_{a,a'} ) } &
[\B,F]( \eta(a \otimes a') )  
}
$$
which component in $b$ is the arrow
$$\xymatrix@C=3pc{
F(a \otimes b) \otimes F(a' \otimes b) 
\ar[r]^-{F^2_{a \otimes b, a' \otimes b} } 
&
F( (a \otimes b) \otimes (a' \otimes b))
\ar[r]^-{ F(\gamma_{a,a',b}) }
&
F( (a \otimes a') \otimes b ). 
}
$$

For any objects $a$ of $\A$ and 
$b,b'$ of $\B$, Diagram
\ref{epsnat} for 
$h = \delta_{a,b,b'} : (a \otimes b) \otimes (a \otimes b') 
\rightarrow a \otimes (b \otimes b') $
is
$$\xymatrix@C=3pc{
\Fb((a \otimes b) \otimes (a \otimes b'))
\ar[d]_{ \epsilon_{ (a \otimes b) \otimes (a \otimes b') } }
\ar[r]^-{\Fb(\delta_{a,b,b'})}
&
\Fb(a \otimes (b \otimes b'))
\ar[d]^{\epsilon_{a \otimes (b \otimes b')}}
\\
F((a \otimes b) \otimes (a \otimes b') )
\ar[r]_-{F(\delta_{a,b,b'})}
&
F( a \otimes (b \otimes b')),
}
$$
which commutes as shown below.
$\Fb( a \otimes ( b \otimes b'))$ is 
$F(a \otimes (b \otimes b')$ and  
$\epsilon_{a \otimes ( b \otimes b')}$ 
is the identity at $F(a \otimes (b \otimes b'))$.
The arrow $\epsilon_{(a \otimes b) \otimes (a \otimes b')}$
is $$
\xymatrix@C=3pc{
\Fb((a \otimes b) \otimes (a \otimes b'))
\ar@{=}[r] & 
F(a \otimes b) \otimes F(a \otimes b')
\ar[r]^-{F^2_{a \otimes b, a \otimes b'}} &
F((a \otimes b) \otimes (a \otimes b'))
}.
$$
The functor $\Fb$ sends the arrow $\delta_{a,b,b'}$ 
to the arrow $${(\Res(F)(a))}^2_{b,b'}:
\Res(F)(a)(b)  \otimes \Res(F)(a)(b') 
\rightarrow \Res(F)(a)(b \otimes b')$$ in $\C$.
Since the functor $\Res(F)(a): \B \rightarrow \C$
is the composite
$\xymatrix{
\B 
\ar[r]^-{\eta(a)} 
& \A \otimes \B 
\ar[r]^-{F} 
& 
\C
},$  the arrow 
${(\Res(F)(a))}^2_{b,b'}$ is
$$\xymatrix{ 
\Res(F)(a)(b) \otimes \Res(F)(a)(b')
\ar@{=}[d]
\\
F(a \otimes b) \otimes F(a \otimes b')
\ar[d]^{F^2_{a \otimes b, a \otimes b'}}
\\
F((a \otimes b) \otimes (a \otimes b'))
\ar[d]^{F(\delta_{a,b,b'})} 
\\
F(a \otimes (b \otimes b'))
\ar@{=}[d]
\\ 
\Res(F)(a)(b \otimes b'). 
}$$

For any arrow $f: a \rightarrow a'$ in $\A$
and any object $b$ of $\B$,
Diagram \ref{epsnat} for 
$h = f \otimes b: a \otimes b \rightarrow a' \otimes b$
is 
$$\xymatrix{
\Fb(a \otimes b)
\ar[d]_{\epsilon_{a \otimes b}}
\ar[r]^-{\Fb(f \otimes 1)}
&
\Fb(a' \otimes b)
\ar[d]^{\epsilon_{a' \otimes b}}
\\
F(a \otimes b)
\ar[r]_-{F(f \otimes 1)}
&
F(a' \otimes b)
}
$$
and commutes since 
$\epsilon_{a \otimes b}$ (resp. $\epsilon_{a' \otimes b}$) 
is the identity at $F(a \otimes b)$ 
(resp. at $F(a' \otimes b)$)
and straightforward computation shows that the $\Fb(f \otimes b)$ 
is
$$\xymatrix@C=3pc{
\Fb(a \otimes b) 
\ar@{=}[r]
&
F(a \otimes b)
\ar[r]^-{ F(f \otimes 1)  } 
&
F(a' \otimes b)
\ar@{=}[r]
&
\Fb(a \otimes b)
}.$$
 
For any object $a$ of $\A$ and any arrow 
$g: b \rightarrow b'$ in $\B$,
Diagram \ref{epsnat} for 
$h = 1 \otimes g: a \otimes b \rightarrow a \otimes b'$ 
is 
$$\xymatrix@C=3pc{
\Fb(a \otimes b)
\ar[d]_{\epsilon_{a \otimes b}}
\ar[r]^-{\Fb(1 \otimes g)}
&
\Fb(a \otimes b')
\ar[d]^{\epsilon_{a \otimes b'}}
\\
F(a \otimes b)
\ar[r]_-{F(1 \otimes g)}
&
F(a \otimes b')
}
$$
which commutes since 
the arrows $\epsilon$ above are identities
and the image by 
$\Fb$ of $1 \otimes g: a \otimes b \rightarrow a \otimes b'$
is $\Res(F)(a)(g)$ which is $F(1 \otimes g)$.\\

So far we have shown that for any arrow $h$ in 
$\Gra_1 \cup \Gra_2 \cup \Gra_3$, Diagram \ref{epsnat}
commutes. Now we prove by induction that this is the
case for all arrows $h$ in $\Gra$.\\

For any object $X$ in $\Ver_{\A, \B}$ and 
any arrow $f: Y \rightarrow Z$ in $\Gra_{\A,\B}$,
Diagram \ref{epsnat} for 
$h = X \otimes f: X \otimes Y \rightarrow X \otimes Z$ 
is 
$$\xymatrix@C=3pc{
\Fb(X \otimes Y)
\ar[d]_{\epsilon_{X \otimes Y}}
\ar[r]^{\Fb(X \otimes f)}
&
\Fb(X \otimes Z)
\ar[d]^{\epsilon_{X \otimes Z}}
\\
F(X \otimes Y)
\ar[r]_-{F(X \otimes f)}
&
F(X \otimes Z)
}
$$
which is the external diagram in the pasting
$$\xymatrix@C=3pc{
\Fb(X) \otimes \Fb(Y)
\ar[d]_{\epsilon_X \otimes \epsilon_Y }
\ar[r]^{ 1  \otimes \Fb(f) }
&
\Fb(X) \otimes \Fb(Z)
\ar[d]^{\epsilon_X \otimes \epsilon_Z}
\\
F(X) \otimes F(Y)
\ar[d]_{F^2_{X,Y} }
\ar[r]^{ 1  \otimes F(f) }
&
F(X) \otimes F(Z)
\ar[d]^{F^2_{X,Z}}\\
F(X \otimes Y)
\ar[r]_-{F(X \otimes f)}
&
F(X \otimes Z)
}$$
By functoriality of the tensor in $\C$,
the top diagram commutes if Diagram \ref{epsnat} 
commutes for $h = f$.
The bottom one commutes
by naturality of $F^2$.\\

Similarly one shows that Diagram \ref{epsnat} 
for  $h = f \otimes X: Y \otimes X \rightarrow Z \otimes X$
commutes providing the commutation of Diagram \ref{epsnat} for $h = f$.\\

We show now that the natural transformation
$\epsilon: \Fb \rightarrow F: \A \otimes \B \rightarrow \C$ 
is monoidal.\\

It satisfies Axiom \ref{monat6}
i.e. for any $X,Y$ in $\A \otimes \B$, the diagram in $\C$ 
$$
\xymatrix{
\Fb(X) \otimes \Fb(Y) 
\ar[r]^-{ {(\Fb)}^2_{X,Y} } 
\ar[d]_{ \epsilon_X \otimes \epsilon_Y }
& \Fb(X \otimes Y)
\ar[d]^{ \epsilon_{ (X \otimes Y) }  }
\\
F(X) \otimes F(Y)
\ar[r]_-{ F^2_{X,Y} }
&
F(X \otimes Y) 
}
$$
commutes. This is the case since
$\Fb$ is strict and according to the inductive definition
of $\epsilon$.\\

The natural transformation $\epsilon$ satisfies Axiom \ref{monat7}, i.e.
the diagram in $\C$
$$
\xymatrix{
\un 
\ar[rd]_-{F^0} 
\ar[r]^-{ {(\Fb)}^0 } 
& 
\Fb(\un)
\ar[d]^{ \epsilon_{\un} }
\\
& 
F(\un)
}
$$
commutes. This holds since $\Fb$ is strict and 
by the definition of $\epsilon$ at $\un$.\\ 

Eventually, we show that the collection of monoidal 
transformations 
$\epsilon_{F} : \Fb \rightarrow F: \A \otimes \B \rightarrow \C$
for the symmetric monoidal functors $F: \A \otimes \B \rightarrow \C$,
constitutes a natural transformation 
$\Ext \circ \Res \rightarrow 1: 
\SMC(\A \otimes \B,C) \rightarrow \SMC(\A \otimes \B,\C).$ 
Given any monoidal natural transformation
$\sigma: F \rightarrow G: \A \otimes \B \rightarrow \C$ between
symmetric monoidal functors, we check by induction
on the objects of $\Ver_{\A,\B}$  
that the following diagram in $\SMC(\A \otimes \B, \C)$ commutes
\begin{tag}\label{epsnat2}
$$\xymatrix{
\Fb 
\ar[r]^{ \sigb }
\ar[d]_{ \epsilon_F }
& 
\Gb
\ar[d]^{ \epsilon_G }\\
F 
\ar[r]_{\sigma}
&
G
}$$
\end{tag}

Pointwise in $\un$, Diagram \ref{epsnat2} is 
$$\xymatrix{
\un 
\ar@{=}[r]
\ar[d]_{ F^0 } & 
\un 
\ar[d]^{ G^0 }\\
F(\un) 
\ar[r]_{ {\sigma}_{\un} } 
&
G(\un)
}$$
which commutes according to Axiom \ref{monat7} for the
monoidal $\sigma$.\\ 

For any objects $a$ in $\A$ and $b$ in $\B$,
Diagram \ref{epsnat2} pointwise in $a \otimes b$
commutes since 
$\Fb (a \otimes b) = F(a \otimes b)$,
$\Gb (a \otimes b) = G(a \otimes b)$, the components of
$\epsilon_F$ and $\epsilon_G$ in $a \otimes b$
are identities and
${\sigb}_{a \otimes b}$ is 
$$
\xymatrix{
\Fb (a \otimes b) 
\ar@{=}[r] 
&
F(a \otimes b)
\ar@{=}[r]^{ \sigma_{a \otimes b} }
&
G(a \otimes b)
\ar@{=}[r]
&
\Gb (a \otimes b)
}.$$

For any $X$ and $Y$, 
Diagram \ref{epsnat2},
pointwise in $X \otimes Y$, is the external diagram in the
pasting
$$\xymatrix{
\Fb(X) \otimes \Fb(Y) 
\ar[r]^-{ \sigb_X \otimes \sigb_Y }
\ar[d]_{ {(\epsilon_F)}_X \otimes {(\epsilon_F)}_Y }& 
\Gb(X) \otimes \Gb(Y)
\ar[d]^{ {(\epsilon_G)}_X \otimes {(\epsilon_G)}_Y }
\\
F(X) \otimes F(Y) 
\ar[r]_{\sigma_X \otimes \sigma_Y}
\ar[d]_{F^2_{X,Y}} &
G(X) \otimes G(Y)
\ar[d]^{G^2_{X,Y}}
\\
F(X \otimes Y) 
\ar[r]_-{{\sigma}_{X \otimes Y}}
&
G(X \otimes Y).
}$$
The top diagram commutes if Diagram \ref{epsnat2} 
commutes pointwise in $X$ and $Y$ and the bottom
one commutes according to Axiom \ref{monat6}
for $\sigma$.\\

\begin{remark}\label{epsilonExtis1}
For any strict (respectively strong)
symmetric monoidal functor $F: \A \otimes \B \rightarrow \C$,
$\epsilon_{F}$ is the identity (respectively an isomorphism).
\end{remark}

\begin{proposition}\label{AdjResExt}
The functor $\Res: \SMC(\A \otimes \B, \C) \rightarrow
\SMC(\A,[\B,\C])$ 
is right adjoint to $\Ext: \SMC( \A , [\B, \C]) \rightarrow 
\SMC(\A \otimes \B,\C)$. 
\end{proposition}
\pf
The unit of this adjunction is the identity 
and the counit is given in any
$F: \A \otimes \B \rightarrow \C$
by 
$\epsilon_F: (\Ext \circ \Res)(F) \rightarrow F$.
The two triangular equalities amount then 
to the facts that
$$\epsilon * \Ext : \Ext \circ \Res \circ \Ext \rightarrow \Ext$$
and
$$\Res * \epsilon : \Res \circ \Ext \circ \Res \rightarrow \Res$$
are identities.   
That $\epsilon * \Ext$ is the identity results from Remark 
\ref{epsilonExtis1}.
That $\Res * \epsilon$ is the identity
is immediate from the definition of $\epsilon$.
\epf

As the composite of two strict symmetric monoidal functors,
the functor $\Res: [\A \otimes \B,\C] \rightarrow [\A, [\B,\C]]$ 
is strict. 

\begin{proposition}\label{MAdjResExt}
One has an adjunction
$\Ext \dashv \Res: [\A \otimes \B, \C] \rightarrow [\A,[\B,\C]]$ in 
$\SMC$
where $\Res$ is strict and the underlying adjunction
in $\CAT$ is the one described in Proposition \ref{AdjResExt}.
\end{proposition}
\pf
According to Kelly's result (see \ref{monadj}), it
is enough to show the two points below:\\
- $(1)$ For any symmetric monoidal functors $F,G: \A \rightarrow [\B,\C]$
the arrow 
{\tiny
$$\xymatrix@C=5pc{
\Ext( F \Box G ) 
\ar[r]^-{\Ext(u_F \Box u_G)} 
&
\Ext( \Res \Ext (F) \Box \Res \Ext(G))
\ar[r]^-{\Ext(\Res^2_{\Ext(F),\Ext(G)})}
&
\Ext  \Res ( \Ext(F) \Box \Ext(G))
\ar[r]^-{\epsilon_{\Ext(F) \Box \Ext(G)}}
&
\Ext(F) \Box \Ext(G)
}$$ 
}
where $u$ denotes the unit of the adjunction $\Ext \dashv \Res$,
is an isomorphism;\\
- $(2)$ The arrow
$$\xymatrix@C=3pc{
\Ext(\un)
\ar[r]^-{\Ext(\Res^0)}
&
\Ext \Res (\un)
\ar[r]^-{\epsilon_{\un}}
&
\un
}
$$
is an isomorphism.\\

$(1)$ holds  since the unit $u$ of the adjunction is the identity, 
$\Res$ is strict and
$\epsilon_{\Ext(F) \Box \Ext(G)}$ is an 
isomorphism by Remark \ref{epsilonExtis1} since 
$\Ext(F)$ and $\Ext(G)$ are strict and
thus $\Ext(F) \Box \Ext(G)$ is strong.\\

$(2)$ holds since $\Res$ is strict and since 
the unit of $[\A \otimes \B, \C]$ is strong
and thus by Remark \ref{epsilonExtis1} the arrow 
$\epsilon_{\un}: \Ext \Res (\un) \rightarrow \un$ 
is an isomorphism.  
\epf

\end{section}

\begin{section}{The tensor 2-functor 
$\SMC \times \SMC \rightarrow \SMC$}\label{tens2fun}
In this section, the mapping sending any 
symmetric monoidal categories $\A$ and $\B$ to their 
tensor $\A \otimes \B$ is extended to a 2-functor
$\tenSMC: \SMC \times \SMC \rightarrow \SMC$.\\

From Proposition \ref{characExt} and its
corollary \ref{coroadjten},
one has for any 
$\A$, $\B$ and $\C$ in $\SMC$, an isomorphism in $\CAT$
$$\SMC(\A, [\B, \C]) \cong \SSMC(\A \otimes \B, \C).$$
Actually for any given $\A$ and $\B$, this isomorphism 
is 2-natural in the argument $\C$ between 2-functors 
$\SSMC \rightarrow \CAT$. 
Therefore there is a unique way of extending 
the assignments $(\A, \B) \mapsto \A \otimes \B$ 
into a 2-functor 
$\SMC \times \SMC \rightarrow \SSMC$ which makes
the collection of isomorphisms above also 2-natural in $\A$ in $\B$, i.e.
such that this collection
defines a 2-natural transformation 
between 2-functors
$\SMC \times \SMC \times \SSMC \rightarrow \CAT$.   
We shall call this extension the {\em tensor} 2-functor 
on $\SMC$. We write $\tenSMC$ for the corresponding
functor $\SMC \times \SMC \rightarrow \SMC$ and
also extend the use the binary operation symbol 
$\otimes$ to 1-cells and 2-cells to denote its images.\\  

By a simple application of the Yoneda Lemma, one 
obtains the concrete description of this tensor, as follows.
\begin{tag}
For any $F: \A \rightarrow \C$
and any $\B$ in $\SMC$, the 1-cell 
$\tenSMC(F,1_{\B}):  \A \otimes \B \rightarrow \C \otimes \B$, 
which we also write $F \otimes \B$, is the image by $\Ext$ of 
$\xymatrix{
\A
\ar[r]^-F
&
\C
\ar[r]^-{\eta}
&
[\B, \C \otimes \B].
}$
\end{tag}
\begin{tag}
For any $G: \B \rightarrow \D$
and any $\A$ in $\SMC$,
the 1-cell 
$\tenSMC(1_{\A},G): \A \otimes \B \rightarrow \A \otimes \D$,
which we also write
$\A \otimes G$,
is the image by $\Ext$
of $\xymatrix{
\A
\ar[r]^-{\eta}
&
[\D, \A \otimes \D]
\ar[r]^-{[G,1]}
&
[\B, \A \otimes \D]
}$.
\end{tag}
\begin{tag}
For any $\B$ and any 2-cell
$\sigma: F \rightarrow F': \A \rightarrow \C$
in $\SMC$, the 2-cell 
$$\tenSMC(\sigma, \B) : F \otimes \B \rightarrow F' \otimes \B: 
\A \otimes \B \rightarrow \C \otimes \B,$$
which we also write $\sigma \otimes \B$,
is the image by $\Ext$ of the 2-cell
$$\xymatrix{
& 
\ar@{=>}[dd]^{\sigma}
\\
\A
\ar@/^35pt/[rr]^{F}
\ar@/_35pt/[rr]_{F'}
&
&
\C
\ar[r]^-{\eta}
&
[\B, \C \otimes \B].
\\
&&&
}$$   
\end{tag}
\begin{tag}
For any $\A$ and any 2-cell
$\tau: G \rightarrow G': \B \rightarrow \D$
in $\SMC$,
the 2-cell $$\tenSMC(\A,\tau): \A \otimes G \rightarrow \A \otimes G'
: \A \otimes \B \rightarrow \A \otimes \D,$$
which we also write $\A \otimes \tau$,
is the image by $\Ext$ of the 2-cell
$$ 
\xymatrix{
& & 
\ar@{=>}[dd]^{[\tau,1]}
\\
\A
\ar[r]^-{\eta}
&
[\D,\A \otimes \D]
\ar@/^35pt/[rr]^{[G,1]}
\ar@/_35pt/[rr]_{[G',1]}
&
&
[\B,\A \otimes \D].
\\
&&&
}
$$
\end{tag}

Note that for any $\A$, $\B$ and $\C$ in $\SMC$, the 2-functor 
$\tenSMC(-,\B): \SMC \rightarrow \SMC$
has component in $\A$ and $\C$, a functor 
${\tenSMC(-,\B)}_{\A,\C} : \SMC(\A,\C) \rightarrow 
\SMC(\A \otimes \B, \C \otimes \B)$ 
which admits the symmetric monoidal structure 
$$\xymatrix{ 
[\A,\C]
\ar[r]^-{[1,\eta]}
&
[\A,[\B, \C \otimes \B]]
\ar[r]^{\Ext}
&
[\A \otimes \B, \C \otimes \B]
.}$$
Similarly for any 
$\A$, $\B$, $\D$ in $\SMC$, the 2-functor 
$\tenSMC(\A,-): \SMC \rightarrow \SMC$
has for component in $\B$ and $\D$, a functor 
${\tenSMC(\A,-)}_{\B,\D} : \SMC(\B,\D) \rightarrow 
\SMC(\A \otimes \B, \A \otimes \D)$ 
which admits the symmetric monoidal structure
$$\xymatrix@C=3pc{ 
[\B,\D]
\ar[r]^-{[-,\A \otimes \D]}
&
[[\D,\A \otimes \D], [\B, \A \otimes \D]]
\ar[r]^-{[\eta,1]}
&
[\A, [\B, \A \otimes \D] ]
\ar[r]^-{\Ext}
&
[\A \otimes \B, \A \otimes \D].
}$$

\end{section}

\begin{section}{Naturality issues for $\Res$ and $\Ext$}
\label{2natExtissue}
This section tackles the questions of the naturalities
of the collections of arrows  
$\Res_{\A,\B,\C}: [\A \otimes \B,\C] \rightarrow [\A,[\B,\C]]$
and 
$\Ext_{\A,\B,\C} : [\A,[\B,\C]] \rightarrow [\A \otimes \B,\C]$
in $\SMC$.\\

From the definition of the tensor in $\SMC$
the underlying functors $\SMC(\A , [\B, \C])
\rightarrow \SMC(\A \otimes \B, \C)$
of the 1-cells $\Ext_{\A,\B,\C}$ of $\SMC$ define
a 2-natural transformation between 
2-functors 
$\SMC^{op} \times \SMC^{op} \times \SSMC \rightarrow \CAT$.
Be cautious here that the domain for the third argument $\C$
of the considered 2-functors is $\SSMC$,
the 2-category with 1-cells {\em strict} functors.
This statement will be further improved by Lemmas \ref{natExtinstrictC},
\ref{2natExtinA}, \ref{2natExtinB} below as the collections of 1-cells
$\Ext_{\A,\B,\C}: [\A, [\B, \C]] \rightarrow 
[\A \otimes \B, \C]$ defines a 2-natural transformation 
between $\SMC$-valued 2-functors with domains
$\SMC^{op} \times \SMC^{op} \times \SSMC$.

\begin{lemma}\label{ResnatinC}
For any $\A$ and $\B$,
the collection of 1-cells in $\SMC$
$$\Res_{\A,\B,\C} : [\A \otimes \B,\C] \rightarrow [\A,[\B,\C]]$$
is 2-natural in $\C$.
\end{lemma}
\pf
The 1-cell $\Res_{\A,\B,\C}$ is the composite 
$$\xymatrix{
[\A \otimes \B,\C]
\ar[r]^-{[\B,-]}
&
[[\B,\A \otimes \B],[\B,\C]]
\ar[r]^-{[\eta,1]}
&
[\A,[\B,\C]]
}$$ in $\SMC$
where the collection of arrows
${[\B,-]}_{\A \otimes \B,\C}$ is 2-natural 
in $\C$ according to Lemma \ref{impronatPost}
and the collection of 
$[\eta,[\B,\C]]: [[\B,\A \otimes \B], [\B,\C]] \rightarrow 
[\A,[\B,\C]]$ is 2-natural in $\C$ 
according to the 2-functoriality of $\homSMC$.
\epf

\begin{lemma}\label{improResnatC}
The diagram in $\SMC$
$$
\xymatrix{
&
[\C,\C']
\ar[rd]^{[\B,-]}
\ar[ld]_{[\A \otimes \B,-]} 
&
\\
[[\A \otimes \B, \C],
[\A \otimes \B, \C']] 
\ar[d]_{[1,\Res]}
& 
& 
[[\B,\C],[\B,\C']]
\ar[d]^{[\A,-]}
\\
[[\A \otimes \B, \C],
[\A ,[\B, \C']]] 
& 
& 
[[\A,[\B,\C]],[\A,[\B,\C']]]
\ar[ll]^{[\Res,1]}
}
$$ 
commutes for any
$\A$, $\B$, $\C$ and $\C'$. 
\end{lemma}
\pf
Since all the functors involved in the considered 
diagram are strict, it is enough to show that the 
underlying diagram in $\CAT$
$$
\xymatrix{
&
\SMC(\C,\C')
\ar[rd]^{\homSMC(\B,-)}
\ar[ld]_{\homSMC(\A \otimes \B,-)} 
&
\\
\SMC([\A \otimes \B, \C],
[\A \otimes \B, \C']) 
\ar[d]_{\SMC(1,\Res)}
& 
& 
\SMC([\B,\C],[\B,\C'])
\ar[d]^{\homSMC(\A,-)}
\\
\SMC([\A \otimes \B, \C],
[\A ,[\B, \C']]) 
& 
& 
\SMC([\A,[\B,\C]],[\A,[\B,\C']])
\ar[ll]^{\SMC(\Res,1)}.
}
$$
commutes, which amounts to the 2-naturality in $\C$ 
of the collection 
$\Res_{\A,\B,\C}: [\A \otimes \B,\C] \rightarrow [\A,[\B,\C]]$ in
$\SMC$, that was
established in Lemma \ref{ResnatinC}.
\epf

\begin{lemma}\label{ResnatinA}
For any $\B$ and $\C$,
the collection of 1-cells in $\SMC$
$$\Res_{\A,\B,\C} : [\A \otimes \B,\C] \rightarrow [\A,[\B,\C]]$$
is 2-natural in $\A$.
\end{lemma}
\pf
Given any 1-cell $F: \A' \rightarrow \A$,
consider the pasting in $\SMC$
$$
\xymatrix{
[\A \otimes \B, \C]
\ar[d]_{[F \otimes 1,1]}
\ar[r]^-{[\B,-]}
&
[[\B, \A \otimes \B],[\B,\C]]
\ar[d]|{[[\B,F \otimes 1],1]}
\ar[r]^-{[\eta,1]}
&
[\A,[\B,\C]]
\ar[d]^{[F,1]}
\\
[\A' \otimes \B, \C]
\ar[r]_-{[\B,-]}
&
[[\B, \A' \otimes \B],[\B,\C]]
\ar[r]_-{[\eta,1]}
&
[\A',[\B,\C]]
}
$$ 
in which the left diagram commutes according 
to Corollary \ref{impronatPost2} and the right one commutes
according to the definition of $F \otimes 1$.\\

One easily adapts the previous argument to show that for any 
2-cell $\sigma: F \rightarrow F': \A' \rightarrow \A$ in $\SMC$, 
one has the equality of 2-cells 
 $[\sigma,1] * \Res = \Res * [\sigma \otimes 1,1]$ in $\SMC$.
\epf

\begin{lemma}\label{ResnatinB}
For any $\A$ and $\C$,
the collection of 1-cells in $\SMC$
$$\Res_{\A,\B,\C} : [\A \otimes \B,\C] \rightarrow [\A,[\B,\C]]$$
is 2-natural in $\B$.
\end{lemma}
\pf
Given any 1-cell $G: \B' \rightarrow \B$,
consider the pasting in $\SMC$
$$
\xymatrix{
[\A \otimes \B,\C]
\ar[r]^-{[\B,-]}
\ar[dd]_{[1 \otimes G,1]}
\ar[rd]_{[\B',-]}
&
[[\B,\A \otimes \B],[\B,C]]
\ar[rr]^{[\eta,1]}
\ar[rd]^{[1,[G,1]]}
&
&
[\A,[\B,\C]]
\ar[dd]^{[1,[G,1]]}
\\
&
[[\B',\A \otimes \B],[\B',\C]]
\ar[r]_{[[G,1],1]}
\ar[d]^{[[1,1 \otimes G],1]}
&
[[\B, \A \otimes \B],[\B',\C]]
\ar[rd]_{[\eta,1]}
\\
[\A \otimes \B',\C]
\ar[r]_-{[\B',-]}
&
[[\B',\A \otimes \B'],[\B',\C]]
\ar[rr]_{[\eta,1]}
&
&
[\A,[\B',\C]]
.}
$$
Here all the diagrams involved commute: the top 
left one commutes
according to Lemma \ref{impronatPost3}, the bottom left 
one commutes according
to Corollary \ref{impronatPost2} and the bottom right one
according to the definition 
of $1 \otimes G$.\\

One can adapt the previous argument to show for 
any 2-cell $\sigma: G \rightarrow G': \B' \rightarrow \B$, 
in $\SMC$, the equality of the 2-cells
$\Res * [[\A \otimes \sigma],\C] = [\A,[\sigma,\C]] * \Res$ in $\SMC$.
\epf

Consider any objects $\A$ and $\B$ and any 1-cell $F: \C \rightarrow \C'$ in 
$\SMC$. Lemma \ref{ResnatinC} states that there is an identity 2-cell in $\SMC$
$$
\xymatrix{
[\A \otimes \B,\C]
\ar[r]^{\Res}
\ar@{}[rd]|{=}
\ar[d]_{[\A \otimes \B,F]}
&
[\A, [\B, \C]]
\ar[d]^{[\A, [\B,F]]}
\\
[\A \otimes \B,\C']
\ar[r]_{\Res}
&
[\A, [\B, \C']]
.}
$$
Its ``mate'', that we denote
$\Xi_F$, is a 2-cell in $\SMC$
\begin{tag}\label{XiF}
$$
\xymatrix{
[\A,[\B,\C]]
\ar[r]^{\Ext}
\ar[d]_{[\A,[\B,F]]}
&
[\A \otimes \B,\C]
\ar[d]^{[\A \otimes \B,F]}
\\
[\A,[\B,\C']]
\ar[r]_{\Ext}
\ar@{=>}[ru]|{\Xi_F}
&
[\A \otimes \B,\C'].
}
$$
\end{tag}
More explicitly $\Xi_F$ is the following pasting in $\SMC$
$$
\xymatrix{
& 
\ar@{=}[d]
& 
&
\\
[\A,[\B,\C]]
\ar[r]_{\Ext}
\ar@/^40pt/[rr]^1
& 
[\A \otimes \B, \C]
\ar[r]^-{\Res}
\ar@{}[rd]|{=}
\ar[d]|{[\A \otimes \B,F]}
&
[\A,[\B,\C]]
\ar[d]|{[\A,[\B,F]]}
& 
\\
&
[\A \otimes \B, \C']
\ar@/_40pt/[rr]_{1}
\ar[r]^-{\Res}
&
[\A,[\B,\C']]
\ar[r]^{\Ext}
\ar@{=>}[d]^{\epsilon}
& 
[\A \otimes \B, \C']
\\
& 
& 
& 
}
$$
where the top identity 2-cell and $\epsilon$ are 
respectively the unit and counit of the adjunctions 
$\Ext \dashv \Res$ in the 2-category $\SMC$.\\

Considering the 2-cell above, observe that for any 
symmetric monoidal $G: \A \rightarrow [\B,\C]$, 
the component in $G$ of $\Xi_F$ is $\epsilon_{F \circ \Ext(G)}$,
which is the identity if $F$ is strict. 
Let us take note of this.
\begin{remark}
The 2-cell $\Xi_F$ is the identity for any strict $F$.
\end{remark}

Also and still according to Lemma \ref{ResnatinC},
for any 2-cell
$\sigma: F \rightarrow G: \C \rightarrow \C'$,
the two pastings in $\SMC$
$$
\xymatrix{
[\A,[\B,\C]]
\ar[rr]^{\Ext}
\ar[dd]_{[\A,[\B,F]]}
&
&
[\A \otimes \B,\C]
\ar[dd]_{[1,F]}
\ar@/^54pt/[dd]^{[1,G]}
\\
&&
\ar@{=>}[r]^-{[1,\sigma]}
&
\\
[\A, [\B,\C']]
\ar[rr]_{\Ext}
\ar@{=>}[rruu]^{\Xi_F}
&
&
[\A \otimes \B,\C']
}
$$
and
$$
\xymatrix{
&
[\A,[\B,\C]]
\ar[rr]^{\Ext}
\ar@/_54pt/[dd]_{[\A,[\B,F]]}
\ar[dd]^{[\A,[\B,G]]}
&
&
[\A \otimes \B,\C]
\ar[dd]|{[\A \otimes \B,G]}
\\
\ar@{=>}[r]^{[1,[1,\sigma]]}
&
\\
&
[\A,[\B,\C']]
\ar[rr]_{\Ext}
\ar@{=>}[rruu]_{\Xi_G}
&
&
[\A \otimes \B,\C']
}
$$
are equal.\\

As a consequence of this, one has the following result.
\begin{lemma}\label{natExtinstrictC}
For any
$\A$ and $\B$, the collection
of 1-cells $\Ext_{\A,\B,\C}: [\A,[\B,\C]]
\rightarrow [\A \otimes \B,\C]$ in $\SMC$
defines a 2-natural transformation in the argument $\C$
between the restrictions 
of the 2-functors $\homSMC(\A,-) \circ \homSMC(\B,-)$ and
$\homSMC(\A \otimes \B,-)$ to the 2-category $\SSMC$.
\end{lemma}

\begin{lemma}\label{2natExtinA}
For any $\B$ and $\C$, the collection of 1-cells 
$\Ext_{\A, \B, \C} : [\A, [\B, \C]] \rightarrow [\A \otimes \B, \C]$
of $\SMC$ defines a 2-natural transformation in the argument $\A$
$$\homSMC(-,[\B,\C]) 
\rightarrow \homSMC(-,\C) \circ \tenSMC(-,\B):\SMC^{op} \rightarrow \SMC.$$
\end{lemma}
\pf
Consider any 1-cell  $F: \A \rightarrow \A'$ in $\SMC$.
Let us see that the diagram  
$$
\xymatrix{
[\A,[\B,\C]]
\ar[r]^-{\Ext}
\ar[d]_{[F,1]}
&
[\A \otimes \B, \C]
\ar[d]^{[F \otimes 1,1]}
\\
[\A',[\B,\C]]
\ar[r]_-{\Ext}
&
[\A' \otimes \B, \C]
}
$$
in $\SMC$ commutes.
According to Lemma \ref{ResnatinA}, 
the following diagram in $\SMC$ commutes
$$
\xymatrix{
[\A \otimes \B,\C] 
\ar[r]^-{\Res}
\ar[d]_{[F \otimes 1,1]}
&
[\A,[\B,\C]]
\ar[d]^{[F ,1]}
\\
[\A' \otimes \B,\C] 
\ar[r]_-{\Res}
&
[\A',[\B,\C]]
}
$$
and the corresponding identity 
2-cell has for mate  
$$\Ext \circ [F,1] \rightarrow [F \otimes 1,1] \circ \Ext: [\A,[\B,\C]]
\rightarrow [\A' \otimes \B,\C]$$
the 2-cell $\epsilon * ( [F \otimes 1,1] \circ  \Ext)$.
In any object $H$ of $[\A,[\B,\C]]$,
this monoidal natural transformation 
has component $\epsilon_{\Ext(H) \circ (F \otimes 1)}$
which is an identity.\\

Consider now any 2-cell $\sigma: F \rightarrow G: \A' \rightarrow \A$
in $\SMC$.
According to Lemma \ref{ResnatinA}, the 2-cells
$$
\xymatrix{
&
&
\ar@{=>}[dd]^{[\sigma,1]}
\\
[\A \otimes \B, \C]
\ar[r]^{\Res}
&
[\A, [\B, \C]]
\ar@/^35pt/[rr]^{[F,1]}
\ar@/_35pt/[rr]_{[G,1]}
&
&
[\A', [\B, \C]]
\\
& & &
}
$$
and 
$$
\xymatrix{
& \ar@{=>}[dd]|{[\sigma \otimes 1,1]} & 
\\
[\A \otimes \B, \C]
\ar@/^35pt/[rr]^{[F \otimes 1,1]}
\ar@/_35pt/[rr]_{[G \otimes 1,1]}
&
&
[\A' \otimes \B, \C]
\ar[r]^{\Res}
&
[\A',[\B,\C]]
\\
& & 
}
$$
are equal.
Composing both 2-cells above on both sides by $\Ext$, one obtains 
for the first 2-cell
$$
\xymatrix{
& \ar@{=>}[dd]^{[\sigma,1]} & 
\\
[\A,[\B, \C]]
\ar@/^35pt/[rr]^{[F,1]}
\ar@/_35pt/[rr]_{[G,1]}
&
&
[\A',[\B, \C]]
\ar[r]^{\Ext}
&
[\A' \otimes \B,\C]
\\
& & 
}
$$
since $\Res \circ \Ext = 1$, whereas
for the second 2-cell, one obtains
$$
\xymatrix{
& 
&
\ar@{=>}[dd]|{[\sigma \otimes 1,1]}
\\
[\A,[\B,\C]]
\ar[r]^{\Ext}
&
[\A \otimes \B, \C]
\ar@/^35pt/[rr]^{[F \otimes 1,1]}
\ar@/_35pt/[rr]_{[G \otimes 1,1]}
& 
& 
[\A' \otimes \B, \C]
\ar[r]^{\Res}
&
[\A',[\B,\C]]
\ar[r]^{\Ext}
&
[\A' \otimes \B, \C]
\\
& & & 
}
$$
which is
$$
\xymatrix{
&
&
\ar@{=>}[dd]^{[\sigma \otimes 1,1]}
\\
[\A, [\B,\C]]
\ar[r]^{\Ext}
& 
[\A \otimes \B, \C]
\ar@/^35pt/[rr]^{[F \otimes 1,1]}
\ar@/_35pt/[rr]_{[G \otimes 1,1]}
& 
& 
[\A' \otimes \B, \C]
\\
& & 
}
$$
since the 2-cell $\epsilon * [\sigma \otimes 1, 1] * \Ext$
is an identity.
\epf

\begin{lemma}\label{2natExtinB}
For any $\A$ and $\C$, the collection of 1-cells 
$\Ext_{\A, \B, \C} : [\A, [\B, \C]] \rightarrow [\A \otimes \B, \C]$
of $\SMC$ defines a 2-natural transformation in the argument $\B$ 
$$\homSMC(\A,-) \circ \homSMC(-,\C) \rightarrow 
\homSMC(-,\C) \circ \tenSMC(\A,-):
\SMC^{op} \rightarrow \SMC.$$
\end{lemma}
\pf
Consider any $\A$ and $\C$
and any 1-cell $F: \B' \rightarrow \B$ in $\SMC$.
Let us see that the diagram 
$$
\xymatrix{
[\A,[\B,\C]]
\ar[r]^{\Ext}
\ar[d]_{[1,[F,1]]}
&
[\A \otimes \B, \C]
\ar[d]^{[1 \otimes F,1]}
\\
[\A,[\B',\C]]
\ar[r]_{\Ext}
&
[\A \otimes \B', \C]
}
$$
in $\SMC$ commutes.
The identity 2-cell corresponding to the equality of
1-cells 
$$[\A,[F,\C]] \circ \Res = \Res \circ [\A \otimes F,\C]$$  
in $\SMC$ given by Lemma \ref{ResnatinB} 
has for mate 
$$\Ext \circ [1,[F,1]] \rightarrow [1 \otimes F,1] \circ \Ext:
[\A,[\B,\C]] \rightarrow [\A \otimes \B',\C]$$
the 2-cell
$\epsilon * ([1 \otimes F,1] \circ \Ext)$ that is an identity.\\

The rest of the proof is similar to the end
the proof for Lemma \ref{2natExtinA}.
\epf
\end{section}

\begin{section}{More commuting diagrams}
\label{lasttecresults}
This section gathers technical results involving 
altogether the internal hom $\homSMC$, the tensor 
and the isomorphism $\Dual$ of $\SMC$.

\begin{lemma}\label{tbc4}
Given any symmetric monoidal categories
$\A$, $\B$, $\C$ and $\D$, the following diagram in $\CAT$
is commutative
$$
\xymatrix{
\SMC(\A,[\B,[\C,\D]]) 
\ar[rr]^{\SMC(-,\Ext)}
\ar[d]_{\Dual}
&
&
\SMC(\A,[\B \otimes \C, \D])
\ar[d]^{\Dual}
\\
\SMC(\B,[\A,[\C, \D]])
\ar[rd]_{\SMC(-,\Dual)}
&
&
\SMC(\B \otimes \C,[ \A , \D])
\ar[ld]^{\Res}
\\
&
\SMC(\B,[\C,[\A,\D]]).
&
}
$$
\end{lemma}
\pf
Considering any symmetric
monoidal functor $F: \A \rightarrow [\B,[\C,\D]]$,
we first show that the composite
$$
\xymatrix{
\B
\ar[r]^-{F^*}
&
[\A,[\C,\D]]
\ar[r]^-{\Dual}
&
[\C,[\A,\D]]
}
$$
is the image by $\Res$ of the dual of the composite 
$$\xymatrix{
\A 
\ar[r]^-{F} 
& 
[\B,[\C,\D]]
\ar[r]^-{\Ext} 
&  
[\B \otimes \C,\D]}.$$
According to Lemma \ref{chardualq},
$\Dual \circ F^*$ is the composite in $\SMC$
$$\xymatrix{
\B 
\ar[r]^-{\resev}
&
[[\B,[\C,\D]],[\C,\D]]
\ar[r]^-{[F,1]}
&
[\A,[\C,\D]]
\ar[r]^-{\Dual}
&
[\C,[\A,\D]]
}$$
whereas the dual of $\Ext \circ F$ is
$$ 
\xymatrix{
\B \otimes \C
\ar[r]^-{\resev}
&
[[[\B \otimes \C],\D],\D]
\ar[r]^-{[\Ext,1]}
&
[[\B,[\C,\D]],\D]
\ar[r]^-{[F,1]}
&
[\A,\D]
.}
$$ 
Therefore we need to show the equality of the external
legs in the pasting in $\SMC$
$$
\xymatrix{
\B 
\ar[rr]^{\eta}
\ar[dd]_{\resev}
&
&
[\C,\B \otimes \C]
\ar[d]^{[\C,\resev]}
\ar[dl]_{[-,\D]}
\\
&
[[\B \otimes \C,\D],[\C,\D]]
\ar[dl]_{[\Ext,1]}
&
[\C,[[\B \otimes \C,\D],\D]]
\ar[l]^{\Dual}
\ar[d]^{[\C,[\Ext,1]]}
\\
[[\B,[\C,\D]],[\C,\D]]
\ar[d]_{[F,1]}
\ar[rr]^-{\Dual}
&
&
[\C,[[\B,[\C,\D]],\D]]
\ar[d]^{[\C,[F,1]]} 
\\
[\A , [\C,\D] ]
\ar[rr]_-{\Dual}
& & 
[\C, [\A,\D] ]
}
$$
In this pasting all diagrams are commutative. 
According to Proposition \ref{dualitynatABC} the two bottom 
diagrams commute.
The top right triangle 
commutes according to Lemma \ref{chardualq2}.
Eventually to show the commutativity of the top left diagram in the 
pasting above, consider the dual of 
$$\xymatrix{
\B 
\ar[r]^-{\eta}
&
[\C, \B \otimes \C]
\ar[r]^-{[-,\D]}
&
[[\B \otimes \C,\D], [\C,\D]]
\ar[r]^-{[\Ext,1]}
&
[[\B,[\C,\D]],[\C,\D]]
.}$$
According to Lemma \ref{dualitynatABC0}, 
this is the arrow
$$\xymatrix{
[\B,[\C,\D]] 
\ar[r]^-{\Ext}
&
[\B \otimes \C,\D]
\ar[r]^-{[\C,-]}
&
[[\C,\B \otimes \C],[\C,\D]]
\ar[r]^-{[\eta,1]}
&
[\B,[\C,\D]]
}$$
which is $\Res \circ \Ext$, which is the identity at 
$[\B,[\C,\D]]$ and the dual of $\resev$.\\

It is rather straightforward to adapt the previous 
argument and use the 2-naturality of $\Dual$ 
to show that for any monoidal natural transformation
$\sigma$ between symmetric functors 
$F \rightarrow G: \A \rightarrow [\B,[\C,\D]]$, the monoidal 
natural transformation
$$
\xymatrix{
& \ar@{=>}@<2ex>[dd]^{\sigma^*} 
& 
& 
\\
\B
\ar@/^35pt/[rr]^{F^*}
\ar@/_35pt/[rr]_{G^*}
&
&
[\A,[\C,\D]]
\ar[r]^{\Dual}
&
[\C,[\A,\D]]
\\
& & & 
}
$$
is the image by $\Res$ of
the dual of 
$$\xymatrix{
&
\ar@{=>}@<2ex>[dd]^{\sigma} & & 
\\
\A 
\ar@/^35pt/[rr]^F
\ar@/_35pt/[rr]_G
& 
&
[\B,[\C,\D]]
\ar[r]^-{\Ext} 
&  
[\B \otimes \C,\D].
\\
& & &}$$
\epf

The image by $\Ext$ of the identity at $[\B,\C]$ is an
arrow $[\B,\C] \otimes \B \rightarrow \C$
that we denote $\Eval$.

\begin{lemma}\label{ineedit} The arrow
$[\A,-]: [\B,\C] \rightarrow [[\A,\B], [\A,\C]]$ is
equal to the composite
$$\xymatrix@C=3pc{
[\B,\C]
\ar[r]^-{[\Eval,\C]}
&
[[\A,\B] \otimes \A,\C]
\ar[r]^-{\Res}
&
[[\A,\B],[\A,\C]].
}$$
\end{lemma}
\pf
By definition of $\Eval$, the following diagram 
in $\SMC$ commutes
$$
\xymatrix{
[\A,\B]
\ar[r]^-{\eta}
\ar[rd]_{1} 
&
[\A,[\A,\B] \otimes \A]
\ar[d]^{[\A,\Eval]}
\\
&
[\A,\B]
.}
$$
Consider then the pasting of commutative diagrams
in $\SMC$
$$
\xymatrix@C=3pc{
[\B,\C]
\ar[r]^-{[\Eval,1]}
\ar[d]_{[\A,-]}
&
[[\A,\B] \otimes \A, \C]
\ar[d]_{[\A,-]}
\ar[r]^-{\Res}
&
[[\A,\B],[\A,\C]]
\\
[[\A,\B], [\A,\C]]
\ar[r]_-{[[1,\Eval],1]}
\ar@/_50pt/[rru]_{1}
&
[[\A, [\A,\B] \otimes \A], [\A,\C]]
\ar[ru]_{[\eta, 1]}
}
$$
the top left diagram commuting according to 
Corollary \ref{impronatPost2}.
\epf
\end{section}

\begin{section}
{The free symmetric monoidal category $\unc$ over the terminal
category}
\label{free1}
In this section we consider 
the free symmetric monoidal category $\unc$ over the terminal
category $\termcat$. We show that for any symmetric monoidal category $\A$,
the strict monoidal 
functor $\vv: \unc \rightarrow [\A,\A]$ sending the generator $\gen$
of $\unc$ to the identity $\A \rightarrow \A$, has 
dual $\vv^*$ which is left adjoint in $SMC$ 
to $\ev_{\gen}: [\unc,\A] \rightarrow \A$, the evaluation functor 
at $\gen$.\\
 
The set of objects of $\unc$, which we write 
$\T$, is the underlying set of
the free algebra over the one point set, with 
element denoted by $\gen$, 
for the signature consisting 
of one constant symbol denoted $\un$ and one binary
symbol denoted $\otimes$.\\ 

We write $\Grd$ 
for the graph with set of vertices 
$\T$ and set of edges defined by induction according to the 
following rules.\\
- For any $X,Y,Z$ in $\T$,
there are one edge 
$\ac_{X,Y,Z}: X \otimes (Y \otimes Z) 
\rightarrow (X \otimes Y) \otimes Z$ 
and one edge
$\aci_{X,Y,Z}: (X \otimes Y) \otimes Z 
\rightarrow X \otimes (Y \otimes Z).$\\
- For any $X$ in $\T$, there are
one edge
$\rc_X : X \otimes \un \rightarrow X$,
one edge $\rci_X: X \rightarrow X \otimes \un$,
one edge $\lc_X : \un \otimes X  \rightarrow X$ and 
one edge $\lci_X: X  \rightarrow \un \otimes X$.\\
- For any $X,Y$ in $\T$, there is
one edge $s_{X,Y}: X \otimes Y \rightarrow Y \otimes X$.\\
- For any $X$ in $\T$ and any edge $p: Y \rightarrow Z$, 
there are two new edges
$X \otimes p : X \otimes Y \rightarrow X \otimes Z$
and $p \otimes X: Y \otimes X \rightarrow Z \otimes X$.\\
- Edges in $\Grd$ with different names are different.\\

We write $\F$ for the free category 
on $\Grd$. For any object $X$ of $\Grd$, the graph morphisms
$X \otimes -$ and $- \otimes X$ extend uniquely to endofunctors of $\F$. 
The underlying category of $\unc$ is defined
as the quotient category of $\F$ by the congruence generated
by a set of relations on the arrows of $\F$
similar to those defining the tensor of Section \ref{tensorABsec}.
This time, this set is the smallest set of relations 
that is closed by the expansions of all relations 
by $X \otimes -$ and $- \otimes X$ for any object $X$,
and that contains the following.\\
- Relations expressing that the $\aci_{X,Y,Z}$ are ``inverses'' 
of the $\ac_{X,Y,Z}$, the $\rci_X$ are inverses of the $\rc_X$,
the $\lci_X$ are inverses of the $\lc_X$;\\
- The ``coherence axioms'' for symmetric monoidal categories;\\
- Relations expressing the ``naturalities'' of the collections
of arrows $\ac$, $\rc$, $\lc$ and $\syc$;\\
- Relations expressing the ``bifunctoriality 
of the tensor'', i.e.
for any edges $t: X \rightarrow Y$ and $s: Z \rightarrow W$  
in $\Grd$, the relation
$(t \otimes W)(X \otimes s) 
\sim (Y \otimes s)(t \otimes Z).$\\

For any object $X$ of $\F$, the endofunctors of $\F$ obtained
from  
$X \otimes -$ and $- \otimes X$, induce respectively 
endofunctors of $\unc$ and the collection of such 
endofunctors defines a functor
$\unc \times \unc \rightarrow \unc$. The category
$\unc$ admits a symmetric monoidal
with canonical arrows $\ac$, $\rc$, $\lc$, $\syc$
the congruence classes of edges of $\Grd$ with
the same names, with tensor the above bifunctor 
and with the object
$\un$ as unit.\\

For any symmetric monoidal category 
$\A$, by the universal property defining 
$\unc$, there exists a unique strict symmetric 
monoidal functor $v: \unc \rightarrow [\A,\A]$ 
such that $v(\gen)$ is the identity functor $\un_{\A}$
in $\A$ with its strict structure.
In particular on other objects of $\unc$, the functor $v$ is as follows:\\
- $v(\un) = \un_{[\A,\A]}$;\\
- for any $X,Y$ in $\T$, 
$v(X \otimes Y) = v(X) \Box v(Y)$.\\

Note then that since the tensor $F \Box G$
of any two strong functors $F$ and $G$ is also strong, 
one obtains by induction that
all the functors $v(X)$ above are strong. 



\begin{proposition}\label{unimonadj}
For any symmetric monoidal category $\A$,
the arrow $\vs: \A \rightarrow [\unc,\A]$ in $\SMC$, 
that is dual of $\vv: \unc \rightarrow [\A,\A]$, 
has right adjoint in $\SMC$ the evaluation at $\gen$
functor $\ev_{\gen} : [\unc, \A] \rightarrow \A$.
Moreover the composite $\ev_{\gen} \circ \vs$ is the identity
at $\A$. 
\end{proposition}
\pf
We know from Section \ref{evalfunsec} that the functor 
``evaluation at $\gen$''
$\SMC(\unc,\A) \rightarrow \A$ that 
sends any symmetric monoidal $F$ to its 
value at the generator $\gen$ and any 
monoidal $\sigma: F \rightarrow G$ to 
its component in $\gen$, admits a strict 
monoidal structure $\ev_{\gen}: [\unc,\A] \rightarrow \A$.\\

The 1-cell  $\ev_{\gen} \circ \vs$ in $\SMC$
is the identity at $\A$ since  
$\ev_{\gen} \circ \vs = \ev_{\gen}(\vv)$
by Lemma \ref{evdual}
and one has $\ev_{\gen}(\vv) = \vv(\gen) = 1_{[\A,\A]}$.\\

We show now the existence of a mere adjunction 
$\vs \dashv \ev_{\gen}: \SMC(\unc,\A) \rightarrow \A$.\\

Given any symmetric monoidal $F: \unc \rightarrow \A$, there exists
a monoidal natural transformation 
$\epsilon_F: (\vs \circ \ev_{\gen})(F) \rightarrow F$.
We will freely drop the subscript $F$ for $\epsilon$
when there is no ambiguity. It is defined by induction on the structure 
of the objects of $\unc$ according to the following rules:\\ 
- $\epsilon_{\un}: \vv(\un)(F(\gen)) \rightarrow F(\un)$
is $F^0: \un \rightarrow F(\un)$;\\
- $\epsilon_{\gen}: \vv(\gen)(F(\gen)) \rightarrow F(\gen)$
is the identity at $F(\gen)$;\\
- For any objects $X,Y$ of $\unc$, the arrow
$\epsilon_{X \otimes Y}: \vv(X \otimes Y)(F(\gen)) \rightarrow
F(X \otimes Y)$ is
$$\xymatrix{ 
\vv(X)(F(\gen)) \otimes \vv(Y)(F(\gen))
\ar[r]^-{\epsilon_X \otimes \epsilon_Y}
&
F(X) \otimes F(Y)
\ar[r]^-{F^2_{X,Y}}
&
F(X \otimes Y).
}$$

To check the mere naturality of 
$\epsilon_F: \vs (\ev_{\gen}(F)) \rightarrow F: \unc \rightarrow \A$ for a given $F$, 
it is enough to show that the following 
diagram in $\A$ commutes
\begin{tag}\label{nateps3}
$$
\xymatrix{
\vv(X)(F(\gen)) 
\ar[r]^-{{\vv(h)}_{F(\gen)}}
\ar[d]_{\epsilon_X} 
&
\vv(Y)(F(\gen))
\ar[d]^{\epsilon_{Y}}\\
F(X)
\ar[r]_{F(h)}
&
F(Y) 
}
$$
\end{tag}
for all the 
arrows $h: X \rightarrow Y$ which are equivalence classes of
edges of $\Grd$.\\

The commutation of Diagram \ref{nateps3} for 
$h = \ac_{X, Y, Z}: X \otimes (Y \otimes Z)
\rightarrow (X \otimes Y) \otimes Z$ amounts 
to the commutation of the external diagram in the following
pasting
$$
\xymatrix{
\vv(X)(F(\gen)) 
\otimes (\vv(Y)(F(\gen)) \otimes \vv(Z)(F(\gen)))
\ar[r]^{\ac}
\ar[d]_{\epsilon_{X} \otimes 
(\epsilon_{Y} \otimes \epsilon_{Z})}
&
( \vv(X)(F(\gen)) 
\otimes (\vv(Y)(F(\gen)) ) \otimes \vv(Z)(F(\gen))
\ar[d]^{
(\epsilon_{X} \otimes \epsilon_{Y}) \otimes \epsilon_{Z}
}
\\ 
F(X) \otimes (F(Y) \otimes F(Z))
\ar[d]_{1 \otimes F^2_{Y,Z}}
\ar[r]^{\ac}
&
(F(X) \otimes F(Y)) \otimes F(Z)
\ar[d]^{F^2_{X,Y} \otimes 1}
\\
F(X) \otimes F(Y \otimes Z)
\ar[d]_{F^2_{X, Y \otimes Z}}
&
F(X \otimes Y) \otimes F(Z)
\ar[d]^{F^2_{X \otimes Y} \otimes 1}
\\
F(X \otimes (Y \otimes Z))
\ar[r]_{F(\ac)}
&
F((X \otimes Y) \otimes Z)
.}
$$
In this pasting, the top diagram commutes 
according to the naturality of $\ac$ and 
the bottom one also does according to 
Axiom \ref{mofun3} for $F$.\\

Similarly one can easily check that  
Diagram \ref{nateps3} commutes for all the instances
$\rc$, $\lc$ and $\syc$ of $h$ according respectively 
to Axioms \ref{mofun41}, \ref{mofun42} and
\ref{symofun5} for $F$, and by naturality
of the canonical isomorphisms in $\A$.\\

The commutation of Diagram \ref{nateps3} for 
$h = X \otimes p$ for any object $X$ of $\T$ and 
$p: Y \rightarrow Z$ in $\Grd$
amounts, since $\vv$ is strict, to the commutation of the 
external diagram in the pasting
$$
\xymatrix@C=3pc{
\vv(X)(F(\gen)) \otimes \vv(Y)(F(\gen))
\ar[r]^{ 1 \otimes {\vv(p)}_{F(\gen)} }
\ar[d]_{\epsilon_X \otimes \epsilon_Y}
&
\vv(X)(F(\gen)) \otimes \vv(Z)(F(\gen))
\ar[d]^{\epsilon_X \otimes \epsilon_Z} 
\\
F(X) \otimes F(Y)
\ar[r]^{1 \otimes F(p)}
\ar[d]_{F^2_{X,Y}}
&
F(X) \otimes F(Z)
\ar[d]^{F^2_{X,Z}}
\\
F(X \otimes Y)
\ar[r]_{F(X \otimes p)}
&
F(X \otimes Z)
.}
$$
In this pasting the bottom diagram commutes 
by naturality of $F^2$ in its second argument and
the top diagram also does if Diagram 
\ref{nateps3} commutes for $h = p$.\\ 

Similarly one has that Diagram \ref{nateps3} for 
$h = p \otimes X$ for any object $X$ of $\T$ and 
$p: Y \rightarrow Z$ in $\Grd$ commutes 
if Diagram \ref{nateps3} already commutes
for $h = p$.\\ 

That for any symmetric monoidal 
functor $F$, the natural transformation
$\epsilon_F$ is monoidal is immediate from its definition. 
Note also that 
the transformation $\epsilon_F$ is an isomorphism for 
any strong $F$ and is the identity for any strict $F$.\\  

That the collection of arrows $\epsilon_F$ for all symmetric
monoidal $F$
defines a mere natural transformation between functors
$$\epsilon: \vs \circ \ev_{\gen} \rightarrow 1: \SMC(\unc,\A) \rightarrow \SMC(\unc,\A)$$
is that for any 2-cell 
$\sigma: F \rightarrow G: \unc \rightarrow \A$ in $\SMC$
and any $X$ in $\unc$
the following diagram in $\A$ commutes
\begin{tag}\label{eps2nat}
$$
\xymatrix{
\vv(X)(F(\gen))
\ar[d]_{ v(X)( \sigma_{\gen} )}
\ar[r]^-{ { (\epsilon_F) }_X}
&
F(X)
\ar[d]^{\sigma_X}
\\
\vv(X)(G(\gen))
\ar[r]_-{{(\epsilon_G)}_X}
&
G(X)
.}
$$
\end{tag}
We show this by induction on the structure 
of the objects $X$ of $\unc$.\\

\noindent For $X = \un$, the commutation
of Diagram \ref{eps2nat}
amounts to the commutation of 
$$
\xymatrix{
\un
\ar@{=}[d]
\ar[r]^{F^0}
&
F(\un)
\ar[d]^{\sigma_{\un}}
\\
\un \ar[r]_{G^0} 
&
G(\un)
}
$$
which is Axiom
\ref{monat7} for $\sigma$.\\

\noindent For $X = \gen$, the commutation
of \ref{eps2nat} amounts to the commutation of
$$
\xymatrix{
F(\gen)
\ar[d]_{\sigma_{\gen}}
\ar@{=}[r]
&
F(\gen)
\ar[d]^{\sigma_{\gen}}
\\
G(\gen) 
\ar@{=}[r] 
&
G(\gen)
}
$$
which is trivial.\\

\noindent For $X = Y \otimes Z$, the commutation
of \ref{eps2nat} amounts to the commutation 
of the external diagram in the pasting
$$
\xymatrix@C=4pc{
v(Y)(F\gen) \otimes v(Z)(F\gen) 
\ar[r]^-{{(\epsilon_F)}_{Y} \otimes {(\epsilon_F)}_{Z} }
\ar[d]_{v(Y)( \sigma_{\gen}) \otimes v(Z)(\sigma_{\gen}) }
& 
F(Y) \otimes F(Z)
\ar[r]^-{F^2}
\ar[d]|{\sigma_{Y} \otimes \sigma_{Z}}
&
F(Y \otimes Z)
\ar[d]^{\sigma_{Y \otimes Z}}
\\
v(Y)(G\gen) \otimes v(Z)(G\gen) 
\ar[r]_-{{(\epsilon_G)}_Y \otimes {(\epsilon_G)}_Z }
& 
G(Y) \otimes G(Z)
\ar[r]_-{G^2}
&
G(Y \otimes Z)
.}
$$
The right diagram in this pasting commutes
according to Axiom \ref{monat6} for $\sigma$
and the left one commutes if Diagram \ref{eps2nat}
commutes for $X = Y$ and for $X = Z$.\\
 
The two triangular equalities 
amount to the facts that both natural
transformations
$$\ev_{\gen} * \epsilon : \ev_{\gen} \circ \vs \circ \ev_{\gen} \rightarrow
\ev_{\gen}: \SMC(\unc,\A) \rightarrow \A$$
and\\
$$\epsilon * \vs: \vs \circ \ev_{\gen} \circ \vs \rightarrow \vs
: \A \rightarrow \SMC(\unc, \A)$$ 
are identities.\\
That $\ev_{\gen} * \epsilon$ is the identity is immediate from the definition
of $\epsilon$.
That $\epsilon * \vs$ is the identity 
results from the fact for any object $a$ in $\A$
the functor $\vs(a): \un \rightarrow \A$ is strict 
by Remark \ref{Fsbstrict} and thus 
$\epsilon_{\vs(a)}$ is an identity.\\


To check that the previous adjunction lifts
to an adjunction 
$\vs \dashv \ev_{\gen}: [\unc,\A] \rightarrow \A$
in the 2-category $\SMC$,
it remains to check that its counit  
$\epsilon: \vs \circ \ev_{\gen} \rightarrow 1$
is monoidal. This amounts to the following two points.\\
$(1)$ For any $F,G: \unc \rightarrow \A$ in $\SMC$,
the diagram in $[\unc,\A]$ 
$$
\xymatrix{
(\vs \circ \ev_{\gen}(F)) 
\Box 
(\vs \circ \ev_{\gen}(G))
\ar[rr]^-{{(\vs \circ \ev_{\gen})}^2_{F,G}}
\ar[rd]_{\epsilon_F \Box \epsilon_G}
&
&
\vs \circ \ev_{\gen}(F \Box G)
\ar[ld]^{\epsilon_{F \Box G}} 
\\
&
F \Box G
}
$$
commutes.\\
$(2)$ The diagram in $[\unc,\A]$  
$$
\xymatrix@C=3pc{
\un
\ar[r]^-{ { ( \vs \circ \ev_{\gen} ) }^0 }
\ar@{=}[rd]
&
\vs \circ \ev_{\gen} (\un_{[\unc,\A] } )
\ar[d]^{ \epsilon_{ \un_{ [\unc,\A] } } }
\\
&
\un}
$$ 
commutes.\\

Point $(1)$ above is equivalent to\\
$(1')$ For any $F,G: \unc \rightarrow \A$ and any $X$
in $\unc$, the diagram in $\A$
$$\xymatrix{
\vv(X)(F(\gen)) 
\otimes 
\vv(X)(G(\gen))
\ar[rr]^-{ {\vv(X)}^2_{F(\gen),G(\gen)} }
\ar[rd]_{ {(\epsilon_F)}_X \otimes {(\epsilon_G)}_X   }
&
&
\vv(X)(F(\gen) \otimes G(\gen))
\ar[ld]^{ {(\epsilon_{F \Box G})}_X } 
\\
&
F(X) \otimes G(X)
}$$
commutes.\\

Point $(2)$ is equivalent to\\
$(2')$ For any $X$ in $\unc$, the diagram in $\A$ 
$$\xymatrix{
\un
\ar[r]^-{{\vv(X)}^0}
\ar@{=}[rd]
&
\vv(X)(\un_{\A})
\ar[d]^{{(  \epsilon_{\un_{[\unc,\A]}} ) }_X }
\\
&
\un
}
$$
commutes.\\

Both points $(1')$ and $(2')$ can be proved by a straightforward
induction on the structure of objects of $\unc$.

\epf


\end{section}

\begin{section}{Canonical arrows and canonical diagrams}\label{canarrows}
In this section canonical arrows part of the ``lax'' symmetric monoidal
structure on the 2-category $\SMC$ are defined.
We make precise the fact that these arrows satisfy the usual coherence 
properties, are 2-natural and are isomorphisms in a lax sense.\\

Note that in this section we shall sometimes use for convenience
the notation $\A\B$ to denote the tensor $\A \otimes \B$ of any 
symmetric monoidal categories $\A$ and $\B$.\\

Given any symmetric monoidal categories 
$\A$, $\B$ and $\C$ one has the composite functor
$$
\gamma^1_{\A,\B,\C}:
\xymatrix{
\SSMC( \A \otimes \B, \C)
\ar@{-}[r]^-{\cong}
&
\SMC(\A,[\B,\C])
\ar[r]^-{\Dual}
&
\SMC(\B,[\A,\C])
\ar@{-}[r]^-{\cong}
&
\SSMC(\B \otimes \A,\C)
}
$$
By the definition of the tensor in $\SMC$ 
and the 2-naturality of $\Dual$ (Lemma \ref{dualitynatABC0}),
the collection of $\gamma^1_{\A,\B, \C}$ defines a 2-natural 
transformation between $\CAT$-valued functors 
with domain $\SMC^{op} \times \SMC^{op} \times \SSMC$.
By Yoneda, the 2-natural $\gamma^1$ corresponds
to a 2-natural transformation, 
namely $S_{\B,\A}: \B \otimes \A \rightarrow \A \otimes \B$
between 2-functors $\SMC \times \SMC \rightarrow \SSMC$.
For any $\A$ and $\B$, 
the 1-cell $S_{\A,\B}: \A \otimes \B \rightarrow \B \otimes \A$ 
is the image by 
$\Ext$ of the 1-cell
$\eta^*: \A \rightarrow [\B, \B \otimes \A]$ 
dual of $\eta: \B \rightarrow [\A, \B \otimes \A]$.\\

Given any symmetric monoidal categories $\A$, $\B$ and $\C$
and $\D$, one defines $\gamma^2_{\A, \B, \C, \D}$ as the composite 
functor
$$\xymatrix{
\SSMC(( \A \otimes \B) \otimes \C, \D) 
\ar@{-}[d]^{\cong}
\\
\SMC(\A \otimes \B, [ \C, \D])
\ar[d]^{\Res}
\\
\SMC(\A, [ \B , [ \C, \D ]] )
\ar[d]^{\SMC( \A,\Ext )}
\\
\SMC(\A, [ \B \otimes \C, \D])
\ar@{-}[d]^{\cong}
\\
\SSMC(\A \otimes ( \B \otimes \C) , \D)
}$$
and one defines $\gamma^3_{\A,\B,\C,\D}$ as the composite functor
$$\xymatrix{
\SSMC( \A \otimes (\B \otimes \C), \D ) 
\ar@{-}[d]^{\cong}
\\
\SMC(\A, [ \B \otimes \C, \D ])
\ar[d]^{\SMC(\A, \Res)}
\\
\SMC(\A, [ \B , [ \C, \D ]] )
\ar[d]^{\Ext}
\\
\SMC(\A \otimes \B, [ \C, \D ])
\ar@{-}[d]^{\cong}
\\
\SSMC(( \A \otimes \B ) \otimes \C, \D).
}$$

According to Lemmas \ref{ResnatinA}, \ref{ResnatinB},
\ref{ResnatinC}, \ref{2natExtinA}, \ref{2natExtinB}
and \ref{natExtinstrictC}, both collections 
$\gamma^2_{\A, \B, \C,\D}$ 
and 
$\gamma^3_{\A,\B,\C,\D}$ define 2-natural transformations
between 2-functors 
$$\SMC^{op} \times \SMC^{op} \times \SMC^{op} \times \SSMC \rightarrow
\CAT.$$
Therefore by Yoneda, the collections $\gamma^2_{\A, \B, \C, \D}$ 
and $\gamma^3_{\A, \B, \C, \D}$ correspond 
respectively to 2-natural transformations, 
namely $$A_{\A,\B,\C}: \A \otimes (\B \otimes \C)
\rightarrow (\A \otimes \B) \otimes \C$$ 
and 
$${A'}_{\A,\B,\C}: (\A \otimes \B) \otimes \C
\rightarrow \A \otimes (\B \otimes \C),$$
both between 
$\SSMC$-valued 2-functors with domains
$\SMC \times \SMC \times \SMC$.\\

For any $\A$,$\B$,$\C$ in $\SMC$, $A_{\A,\B,\C}$
is the image of the identity at $(\A \otimes \B) \otimes \C$
by $\gamma^2_{\A,\B,\C,(\A \otimes \B) \otimes \C}$.
This is to say that it is the image by $\Ext$ of
\begin{center}
$\xymatrix{
\A
\ar[d]^{\eta}
\\
[ \B, \A \otimes \B]
\ar[d]^{[\B,\eta]}
\\
[ \B , [\C, (\A \otimes \B) \otimes \C]]
\ar[d]^{\Ext}
\\
[ \B \otimes \C, (\A \otimes \B) \otimes \C].
}$
\end{center}
The 1-cell above is $\Res(A_{\A,\B,\C})$ and its dual
is strict
according to Remark \ref{Fsbstrict} since for all objects $a$ of $\A$,
the symmetric monoidal functor 
$\Res(A)(a): \B \otimes \C \rightarrow (\A \otimes \B) \otimes \C$
is strict as an image by $\Ext$.\\

For any $\A$, $\B$, $\C$ in $\SMC$, the 1-cell
${A'}_{\A, \B, \C} : (\A \otimes \B) \otimes \C 
\rightarrow \A \otimes ( \B \otimes \C )$ is the image of the 
identity at $\A  \otimes (\B \otimes \C)$ by
$\gamma^3_{\A,\B,\C, \A \otimes (\B \otimes \C)}$.
This is to say that 
$A'_{\A,\B,\C}: (\A \otimes \B) \otimes \C \rightarrow \A \otimes (\B \otimes
\C)$
is the image by $\Ext \circ \Ext$ of the arrow
$$\xymatrix{
\A
\ar[r]^-{\eta}
&
[\B \otimes \C, \A \otimes (\B \otimes \C)]
\ar[r]^-{\Res}
&
[\B,[\C, \A \otimes (\B \otimes \C)]].
 }$$\\

Given any symmetric monoidal categories $\A$ and $\B$, 
one defines $\gamma^4_{\A,\B}$ as the composite
functor
$$
\xymatrix@C=4pc{
\SSMC(\A \otimes \unc, \B)
\ar@{-}[r]^-{\cong}
&
\SMC(\A,[\unc,\B])
\ar[r]^-{\SMC(\A, \ev_{\gen})}
&
\SMC(\A,\B)
}
$$
and 
$\gamma^5_{\A,\B}$ as
the composite functor
$$
\xymatrix{
\SSMC(\unc \otimes \A,\B)
\ar[r]^-{\gamma^1_{\unc, \A, \B}}
&
\SSMC(\A \otimes \unc, \B)
\ar[r]^-{\gamma^4_{\A,\B}}
&
\SMC(\A, \B)
.}
$$
According to Lemma \ref{2nateva}, both collections $\gamma^4_{\A,\B}$
and $\gamma^5_{\A,\B}$ are 2-natural transformations 
between 2-functors $\SMC \times \SSMC \rightarrow \CAT$.
Therefore by Yoneda $\gamma^4_{\A,\B}$ and 
$\gamma^5_{\A,\B}$ correspond respectively
to the 2-natural transformations 
between 2-functors $\SMC \rightarrow \SMC$
$$R'_{\A}:
\xymatrix{
\A 
\ar[r]^-{\eta} 
&
[\unc, \A \otimes \unc]
\ar[r]^-{\ev_{\gen}}
&  
\A \otimes \unc}$$
and 
$$L'_{\A}: 
\xymatrix{\A 
\ar[r]^-{\eta^*} 
&
[\unc, \unc \otimes \A]
\ar[r]^-{\ev_{\gen}}
&
\unc \otimes \A
}$$
that satisfies moreover the following.
\begin{lemma}\label{SRpLp}
The diagram in $\SMC$
$$\xymatrix{\A \otimes \unc  \ar[rr]^{S} & & 
\unc \otimes \A  \\  
& \A 
\ar[lu]^{R'}
\ar[ru]_{L'}
& }$$
commutes for any $\A$.
\end{lemma}

Given any symmetric monoidal category
$\A$, the symmetric monoidal functor $L_{\A}: \unc \otimes \A \rightarrow
\A$ is defined as the image by $\Ext$ of
the symmetric monoidal functor $\vv: \unc \rightarrow [\A,\A]$   
defined in Section \ref{free1}.
The symmetric monoidal functor $R_{\A}: \A \otimes \unc \rightarrow \A$ 
is the image by $\Ext$ of the symmetric monoidal functor 
$\vs: \A \rightarrow [\unc,\A]$,
dual of $\vv$ via the isomorphism \ref{dualhom}.\\

When no ambiguity can occur, we shall
omit the subscripts $\A$, $\B$, ..., for the ``canonical'' 
arrows $A$, $A'$, $R$, $R'$, $L$, $L'$ and $S$.\\

For any $\A$, $\B$, $\C$ in $\SMC$,
the composite functor $\Dual_{\B,\A,\C} \circ \Dual_{\A,\B, \C}$
is the identity and thus the functor
$\gamma^1_{\B, \A, \C} \circ \gamma^1_{\A, \B, \C}$ is the identity
at $\SSMC(\A \otimes \B, \C)$. Therefore by Yoneda, one has the following.
\begin{lemma}\label{waxiom3}
For any $\A$ and $\B$, the composite $S_{\B, \A} \circ S_{\A,\B}$
is the identity at $\A \otimes \B$ in $\SMC$. 
\end{lemma}

\begin{lemma}\label{ResSD}
The diagram in 
$\SMC$ 
$$
\xymatrix{
[\A \otimes \B,\C]
\ar[r]^{[S,1]}
\ar[d]_{\Res}
&
[\B \otimes \A,\C]
\ar[d]^{\Res}
\\
[\A,[\B,\C]]
\ar[r]_{\Dual}
&
[\B,[\A,\C]]
}
$$
commutes for any $\A$, $\B$ and $\C$.
\end{lemma}
\pf
Consider the pasting of diagrams
in $\SMC$
$$\xymatrix{
[\A \otimes \B,\C] 
\ar[rr]^{[S,1]}
\ar[d]_{[\B,-]}
\ar[rd]^{[\A,-]}
& 
&
[\B \otimes \A,\C]
\ar[d]^{[\A,-]}
\\
[[\B, \A \otimes \B],[\B,\C]]
\ar[d]_{[\eta,1]}
&
[[\A,\A \otimes \B],[\A,\C]]
\ar[rd]_{[\eta^*,1]}
\ar[r]^{[[1,S],1]}
&
[[\A,\B \otimes \A],[\A,\C]]
\ar[d]^{[\eta,1]}
\\
[\A,[\B,\C]]
\ar[rr]_{\Dual}
&
&
[\B,[\A,\C]]
.}
$$
All the diagrams above are commutative, the top
one according to Corollary \ref{impronatPost2}
and the left-bottom one 
according to Lemma \ref{tbc3}.
\epf

According to the 2-naturality of $S$, one can define
for any symmetric monoidal categories $\A$, $\B$ and $\C$ 
the 2-cell $T_{\A,\B,\C}$ in $\SMC$ 
as either one of the two composites
$$
\xymatrix{
(\A \otimes \B) \otimes \C
\ar[r]^{S}
&
\C \otimes (\A \otimes \B)
\ar[r]^{1 \otimes S}
&
\C \otimes (\B \otimes \A)
}
$$
or
$$
\xymatrix{
(\A \otimes \B) \otimes \C
\ar[r]^{S \otimes 1}
&
(\B \otimes \A) \otimes \C
\ar[r]^{S}
&
\C \otimes (\B \otimes \A).
}
$$

\begin{lemma}
For any $\A$, $\B$ and $\C$,
the image by $\Res \circ \Dual \circ \Res$
of the 1-cell 
$A: \A \otimes (\B \otimes \C) \rightarrow 
(\A \otimes \B) \otimes \C$ in $\SMC$
equals to
$$
\xymatrix{
\B
\ar[r]^-{\eta^*}
&
[\A,\A \otimes \B]
\ar[r]^-{[1,\eta]}
&
[\A,[\C,(\A \otimes \B) \otimes \C]]
\ar[r]^-{\Dual}
&
[\C,[\A,(\A \otimes \B) \otimes \C]].
}
$$
\end{lemma}
\pf
The image by $\Res$ of 
$A_{\A,\B,\C}$ is the 1-cell
$$\xymatrix{
\A 
\ar[r]^-{\eta}
&
[\B,\A \otimes \B]
\ar[r]^-{[1,\eta]}
&
[\B,[\C,(\A \otimes \B) \otimes \C]]
\ar[r]^-{\Ext}
&
[\B \otimes \C, (\A \otimes \B) \otimes \C]
.}
$$
According to Lemma \ref{tbc4}, the image by $\Res$ of its dual 
is 
$$
\xymatrix{
\B
\ar[r]^-{F^*}
&
[\A,[\C,(\A \otimes \B) \otimes \C]]
\ar[r]^-{\Dual}
&
[\C,[\A,(\A \otimes \B) \otimes \C]]
}
$$
where $F$ is the 1-cell
$\xymatrix{
\A 
\ar[r]^-{\eta}
&
[\B,\A \otimes \B]
\ar[r]^-{[1,\eta]}
&
[\B,[\C,(\A \otimes \B) \otimes \C]]
}.$
According to Lemma \ref{dualitynatABC0} this $F$ has
dual 
$$
\xymatrix{
\B
\ar[r]^-{\eta^*}
&
[\A, \A \otimes \B]
\ar[r]^-{[1,\eta]}
&
[\A,[\C,(\A \otimes \B) \otimes \C]]
} 
$$ 
\epf
According to the previous lemma and 
Lemmas \ref{beurk} and \ref{dualitynatABC0}, one obtains
\begin{corollary}\label{RdRd}
For any $\A$, $\B$ and $\C$, the
image by 
$\Dual \circ \Res \circ \Dual \circ \Res$ 
of the 1-cell $A: \A \otimes (\B \otimes \C)
\rightarrow (\A \otimes \B) \otimes \C$ 
in $\SMC$
equals to
$$
\xymatrix{
\C
\ar[r]^-{\eta^*}
&
[\A \otimes \B,(\A \otimes \B) \otimes \C]
\ar[r]^-{[\A,-]}
&
[[\A,\A \otimes \B], [\A, (\A \otimes \B) \otimes \C]]
\ar[r]^-{[\eta^*,1]}
&
[\B,[\A,(\A \otimes \B) \otimes \C]].
}
$$
\end{corollary}

We can now relate the 1-cells $A$ and $A'$.
\begin{lemma}\label{AandAprel}
The diagram in $\SMC$
$$\xymatrix{
(\A \otimes \B) \otimes \C
\ar[r]^{A'}
\ar[d]_{T_{\A,\B,\C}}
&
\A \otimes (\B \otimes \C)
\\
\C \otimes (\B \otimes \A)
\ar[r]_{A}
&
(\C \otimes \B) \otimes \A
\ar[u]_{T_{\C,\B,\A}}
}
$$
commutes for any $\A$, $\B$ and $\C$.
\end{lemma}
\pf
Both legs of the diagram are strict functors.
We show that their images by $\Res$ are strict
and that their images by $\Res \circ \Res$ are equal.\\
 
The 1-cell
$$
\xymatrix{
(\A \otimes \B) \otimes \C
\ar[r]^S
& 
\C \otimes (\A \otimes \B)
\ar[r]^{1 \otimes S}
&
\C \otimes (\B \otimes \A)
\ar[r]^A
&
(\C \otimes \B) \otimes \A
}
$$
has image by $\Res$
$$
\xymatrix{
\A \otimes \B
\ar[r]^-{\eta^*}
&
[\C, \C \otimes (\A \otimes \B)]
\ar[r]^-{[1, 1 \otimes S]}
&
[\C,\C \otimes (\B \otimes \A)]
\ar[r]^-{[1,A]}
&
[\C,(\C \otimes \B) \otimes \A]
}
$$
which is
$$
\xymatrix{
\A \otimes \B
\ar[r]^-{S}
&
\B \otimes \A
\ar[r]^-{{\eta}^*}
&
[\C,\C \otimes (\B \otimes \A)]
\ar[r]^-{[1,A]}
&
[\C,(\C \otimes \B) \otimes \A]
}
$$
or 
$$
\xymatrix{
\A \otimes \B
\ar[r]^-{S}
&
\B \otimes \A
\ar[r]^-{{\Res(A)}^*}
&
[\C,(\C \otimes \B) \otimes \A].
}
$$
Observe that this 1-cell is strict since both $S$  
and ${\Res(A)}^*$ are. 
The image by $\Res$ of the 1-cell 
${(\Res(A))}^* \circ S$ above 
is 
$$
\xymatrix@C=3pc{
\A 
\ar[r]^-{\eta^*}
&
[\B,\B \otimes \A]
\ar[r]^-{[1,{\Res(A)}^*]}
&
[\B,[\C,(\C \otimes \B) \otimes \A]]
}
$$
which is the dual of $\Res \circ \Dual \circ \Res(A_{\C,\B,\A})$,
that is 
$$
\xymatrix{
\A
\ar[r]^-{\eta^*}
&
[\C \otimes \B, (\C \otimes \B) \otimes \A]
\ar[r]^-{[\C,-]}
&
[[\C,\C \otimes \B],[\C, (\C \otimes \B) \otimes \A]]
\ar[r]^-{[\eta^*,1]}
&
[\B,[\C,(\C \otimes \B) \otimes \A]]
}
$$
according to Corollary \ref{RdRd}.\\

On the other hand the image 
by $\Res$ of the 1-cell
$$\xymatrix{
(\A \otimes \B) \otimes \C
\ar[r]^-{A'}
&
\A \otimes (\B \otimes \C)
\ar[r]^-{1 \otimes S}
&
\A \otimes (\C \otimes \B)
\ar[r]^-S
&
(\C \otimes \B) \otimes \A
}$$
is the strict 1-cell
$[1,S] \circ [1,1 \otimes S] \circ \Res(A'): \A \otimes \B 
\rightarrow [\C,(\C \otimes \B) \otimes \A]$
which image by $\Res$ is
{\small
$$
\xymatrix@C=3pc{
\A
\ar[r]^-{\eta}
&
[\B \C, \A  (\B  \C)]
\ar[r]^-{\Res}
&
[\B,[\C, \A (\B  \C)]]
\ar[r]^-{[1,[1,1 \otimes S]]}
&
[\B,[\C, \A  (\C  \B)]]
\ar[r]^-{[1,[1,S]]}
&
[\B,[\C, (\C  \B) \A]]
}.
$$
}
This last 1-cell rewrites successively as
{\tiny
\begin{tabbing}
\=1.
$
\xymatrix{
\A
\ar[r]^-{\eta}
&
[\B \C, \A (\B \C)]
\ar[r]^-{[1, 1 \otimes S]}
&
[\B \C, \A (\C \B)]
\ar[r]^-{[1,S]}
&
[\B \C, (\C  \B) \A]
\ar[r]^-{\Res}
&
[\B,[\C, (\C  \B)  \A]]
}
$\\
\>2.
$
\xymatrix{
\A
\ar[r]^-{\eta}
&
[\C  \B, \A  (\C  \B)]
\ar[r]^-{[S, 1]}
&
[\B  \C, \A  (\C  \B)]
\ar[r]^-{[1,S]}
&
[\B  \C, (\C  \B) \A]
\ar[r]^-{\Res}
&
[\B,[\C, (\C \B) \A]]
}
$\\
\>3.
$
\xymatrix{
\A
\ar[r]^-{\eta}
&
[\C \B, \A  (\C  \B)]
\ar[r]^-{[1, S]}
&
[\C  \B, (\C  \B) \A]
\ar[r]^-{[S,1]}
&
[\B  \C, (\C  \B)  \A]
\ar[r]^-{\Res}
&
[\B,[\C, (\C  \B) \A]]
}
$\\
\>4.
$
\xymatrix{
\A
\ar[r]^-{\eta^*}
&
[\C  \B, (\C \B) \A]
\ar[r]^-{[S,1]}
&
[\B  \C, (\C  \B)  \A]
\ar[r]^-{\Res}
&
[\B,[\C, (\C  \B) \A]]
}
$\\
\>5.
$
\xymatrix{
\A
\ar[r]^-{\eta^*}
&
[\C \B, (\C  \B) \A]
\ar[r]^-{\Res}
&
[\C, [\B, (\C  \B)  \A]]
\ar[r]^-{\Dual}
&
[\B,[\C, (\C  \B) \A]]
}
$\\
\>6.
$
\xymatrix{
\A
\ar[r]^-{\eta^*}
&
[\C \B, (\C  \B) \A]
\ar[r]^{[\B,-]}
&
[[\B,\C \B], [\B,(\C  \B) \A]]
\ar[r]^-{[\eta,1]}
&
[\C, [\B, (\C  \B)  \A]]
\ar[r]^-{\Dual}
&
[\B,[\C, (\C  \B) \A]]
}
$\\
\>7.
$
\xymatrix{
\A
\ar[r]^-{\eta^*}
&
[\C \B, (\C  \B) \A]
\ar[r]^-{[\C,-]}
&
[[\C, \C  \B]], [\C, (\C  \B)  \A]]
\ar[r]^-{[\eta^*,1]}
&
[\B,[\C, (\C  \B)  \A]]}.$
\end{tabbing}
}
\noindent In the above derivation, arrows 4. and 5. are equal
due to Lemma \ref{ResSD} and arrows 6. and 7. are equal due to
Lemma \ref{tbc3}.
\epf

\begin{lemma}\label{Ainverse2}
There exists a 2-cell
$A \circ A' \rightarrow 1: (\A \otimes \B) \otimes \C
\rightarrow (\A \otimes \B) \otimes \C$ in $\SMC$, for any 
$\A$, $\B$ and $\C$.
\end{lemma}
\pf
The image by $\Res \circ \Res$ of the composite $A \circ A'$
is 
$$\xymatrix@C=4pc{
\A 
\ar[r]^-{\Res(\Res(A'))}
&
[\B,[\C,\A \otimes (\B \otimes \C)]]
\ar[r]^-{[1,[1,A]]}
&
[\B,[\C,(\A \otimes \B) \otimes \C]]
}$$
which rewrites successively as 
{\small
\begin{tabbing}
\=$\xymatrix{
\A 
\ar[r]^-{\eta}
&
[\B \C, \A (\B  \C)]
\ar[r]^-{\Res}
&
[\B,[\C,\A  (\B  \C)]]
\ar[r]^-{[1,[1,A]]}
&
[\B,[\C,(\A  \B) \C]].
}
$\\
\>$\xymatrix{
\A 
\ar[r]^-{\eta}
&
[\B  \C, \A  (\B  \C)]
\ar[r]^-{[1,A]}
&
[\B  \C, (\A  \B)  \C]
\ar[r]^-{\Res}
&
[\B,[\C, (\A  \B) \C]]
}
$\\
\>$\xymatrix{
\A 
\ar[r]^-{\Res(A)}
&
[\B \C, (\A  \B)  \C]
\ar[r]^-{\Res}
&
[\B,[\C, (\A  \B) \C]]
}
$\\
\>$\xymatrix{
\A 
\ar[r]^-{\eta}
&
[\B, \A  \B]
\ar[r]^-{[1,\eta]}
&
[\B,[\C,(\A  \B) \C]]
\ar[r]^-{\Ext}
&
[\B  \C, (\A  \B)  \C]
\ar[r]^-{\Res}
&
[\B,[\C,(\A  \B)  \C]]
}
$\\
\>$\xymatrix{
\A 
\ar[r]^-{\eta}
&
[\B, \A  \B]
\ar[r]^-{[1,\eta]}
&
[\B,[\C,(\A  \B) \C]]
}
$
\end{tabbing}
}
Therefore
$\Res( \Res ( A \circ A'))$ is the 
image by $\Res$ of 
$\eta: \A \otimes \B \rightarrow [\C, (\A \otimes \B) \otimes \C]$.
The 1-cell
$\Res(A \circ A')$, that is $[\C,A] \circ \Res(A')$, is strict
since the 1-cells $A$ and $\Res(A')$ are both strict.
Thus there exists a 2-cell $\Res(A \circ A') \rightarrow \eta$
and since $A \circ A'$ is strict, there exists a 2-cell 
$A \circ A' \rightarrow 1$. 
\epf

\begin{corollary}\label{Ainverse}
There exists a 2-cell
$A' \circ A \rightarrow 1: \A \otimes (\B \otimes \C) \rightarrow
\A \otimes (\B \otimes \C)$ in $\SMC$, for any $\A$, $\B$ and $\C$.
\end{corollary}
\pf 
Such a 2-cell is obtained as 
the pasting 
$$
\xymatrix{
\A \otimes (\B \otimes \C)
\ar@{}[rd]|{=}
\ar[r]^A
\ar[d]_{ T^{-1}_{\C,\B,\A} }
&
(\A \otimes \B) \otimes \C
\ar@{}[rd]|{=}
\ar[r]^{A'}
\ar[d]|{T_{\A,\B,\C}}
&
\A \otimes (\B \otimes \C)
\ar[d]^{T^{-1}_{\A,\B,\C}}
\\
(\C \otimes \B) \otimes \A
\ar[r]_{A'}
\ar@/_40pt/[rr]_1
&
\C \otimes (\B \otimes \A)
\ar[r]_A
\ar@{=>}[d]
&
(\C \otimes \B) \otimes \A
\\
&
&
}
$$
where the top diagrams commute according to 
Lemma \ref{AandAprel} and 
the bottom 2-cell comes from Lemma \ref{Ainverse2}.
\epf

\begin{lemma}\label{waxiom12}
The diagram in $\SMC$
$$\xymatrix{
((\A \otimes  \B) \otimes \C) \otimes \D 
\ar[r]^{A'}
\ar[d]_{A' \otimes 1}
& 
(\A \otimes  \B) \otimes ( \C \otimes \D ) 
\ar[r]^{A'}
& 
\A \otimes (\B \otimes (\C \otimes \D))\\
(\A \otimes (\B \otimes \C)) \otimes \D \ar[rr]_{A'} && 
\A \otimes ((\B \otimes \C)  \otimes \D) \ar[u]_{1 \otimes A'}
}$$
commutes for any $\A, \B, \C$ and $\D$. 
\end{lemma}
\pf
Both legs of the diagram are strict.
We show that their images by $\Res$ are strict,
that their images by $\Res \circ \Res$ are strict
and that their images by $\Res \circ \Res \circ \Res$ 
are equal.\\
 
The image by $\Res$ of the 1-cell
$$\xymatrix{   
((\A \otimes \B) \otimes \C) \otimes \D
\ar[r]^{A'}
&
(\A \otimes \B) \otimes (\C \otimes \D)
\ar[r]^{A'}
&
\A \otimes (\B \otimes (\C \otimes \D))
}$$
is the strict 1-cell 
$$
\xymatrix{
(\A \otimes \B) \otimes \C
\ar[r]^-{\Res(A')}
&
[\D, (\A \otimes \B) \otimes (\C \otimes \D)]
\ar[r]^-{[1,A']}
&
[\D, \A \otimes ( \B \otimes (\C \otimes \D))].
}
$$
This last one has image by $\Res$ 
$$
\xymatrix{
\A \B
\ar[r]^-{\eta}
&
[\C  \D, (\A  \B) (\C \D)]
\ar[r]^-{\Res}
&
[\C,[\D, (\A  \B)  (\C  \D)]]
\ar[r]^-{[1,[1,A']]}
&
[\C,[\D, \A (\B  (\C  \D))]]
}
$$
which rewrites successively as
\begin{tabbing}
\=$\xymatrix{
\A \B
\ar[r]^-{\eta}
&
[\C  \D, (\A  \B) (\C \D)]
\ar[r]^-{[1,A']}
&
[\C  \D, \A (\B  (\C  \D))]
\ar[r]^-{\Res}
&
[\C,[\D, \A  (\B  (\C \D))]]
}
$\\
\>$\xymatrix{
\A \otimes  \B
\ar[r]^-{\Res(A')}
&
[\C \otimes  \D, \A \otimes (\B \otimes (\C \otimes \D))]
\ar[r]^-{\Res}
&
[\C,[\D, \A \otimes (\B \otimes (\C \otimes \D))]].
}
$
\end{tabbing}
This last 1-cell is strict and has image 
by $\Res$ 
$$
\xymatrix@C=4pc{
\A
\ar[r]^-{\Res(\Res(A'))}
&
[\B,[\C \D, \A (\B (\C \D))]]
\ar[r]^-{[1,\Res]}
&
[\B, [\C, [\D, \A ( \B  (\C \D))]]]
}
$$
which rewrites successively as
{\tiny 
\begin{tabbing}
\=1.
$\xymatrix@C=1pc{
\A
\ar[r]^-{\eta}
&
[\B(\C \D), \A (\B (\C \D))]\;
\ar[r]^-{[\C\D,-]}
&
\;
[[\C\D,\B(\C \D)], [\C\D,\A (\B (\C \D))]]
\ar[r]^-{[\eta,1]}
&
[\B,[\C \D, \A (\B (\C \D))]]
\ar[r]^-{[1,\Res]}
&
[\B, [\C, [\D, \A ( \B  (\C \D))]]]
}
$\\
\>2. 
$\xymatrix{
\A
\ar[r]^-{\eta}
&
[\B(\C \D), \A (\B (\C \D))]
\ar[r]^-{[\C\D,-]}
&
[[\C\D,\B(\C \D)], [\C\D,\A (\B (\C \D))]]
\ar[r]^-{[1,\Res]}
&
[[\C\D,\B(\C \D)], [\C,[\D,\A (\B (\C \D))]]   ]
\; ...}$
\\
$\xymatrix{...
\ar[r]^-{[\eta,1]}
&
[\B, [\C, [\D, \A ( \B  (\C \D))]]]
}
$\\
\>3. 
$\xymatrix{
\A
\ar[r]^-{\eta}
&
[\B(\C \D), \A (\B (\C \D))]
\ar[r]^-{[\D,-]}
&
[[\D,\B(\C \D)], [\D,\A (\B (\C \D))]]
\ar[r]^-{[\C,-]}
&
[[\C,[\D,\B(\C \D)]], [\C,[\D,\A (\B (\C \D))]]]
\ar[r]^-{[\Res,1]}
&
...
}$\\
$\xymatrix{ ...\;
[[\C\D,\B(\C \D)], [\C,[\D,\A (\B (\C \D))]]]
\ar[r]^-{[\eta,1]}
&
[\B, [\C, [\D, \A ( \B  (\C \D))]]]
}
$
\end{tabbing}
}
\noindent the above arrows 2. and 3. being equal according to Lemma
\ref{improResnatC}.\\

On the other hand the 1-cell
$$
\xymatrix{
((\A \B) \C) \D
\ar[r]^-{A' \otimes 1}
&
(\A (\B \C)) \D
\ar[r]^-{A'}
&
\A ( (\B \C) \D  )
\ar[r]^-{1 \otimes A'}
&
\A ( \B (\C \D))
}
$$
has image by $\Res$
$$\xymatrix{
(\A \B) \C
\ar[r]^-{A'}
&
\A (\B \C)
\ar[r]^-{\eta}
&
[\D, ((\A  (\B  \C)) \D]
\ar[r]^-{[1,A']}
&
[\D, \A  ( (\B \C) \D)]
\ar[r]^-{[1,1 \otimes A']}
&
[\D, \A ( \B (\C \D))]
}
$$
which is 
$$
\xymatrix{
(\A \B) \C
\ar[r]^-{A'}
&
\A (\B \C)
\ar[r]^-{\Res(A')}
&
[\D, \A  ( (\B \C) \D)]
\ar[r]^-{[1,1 \otimes A']}
&
[\D, \A ( \B (\C \D))].
}
$$
Note that this last 1-cell is strict.
It has image by $\Res$ the 1-cell
$$
\xymatrix@C=3pc{
\A \B
\ar[r]^-{\Res(A')}
&
[\C, \A (\B \C)]
\ar[r]^-{[1,\Res(A')]}
&
[\C,[\D, \A ( (\B \C) \D )]]
\ar[r]^-{[1,[1, 1 \otimes A']]}
&
[\C,[\D, \A ( \B (\C \D))]]
}
$$
which is also strict and has image by $\Res$
{\small
$$
\xymatrix@C=4pc{
\A
\ar[r]^-{\Res ( \Res (A'))}
&
[\B,[\C, \A (\B \C)]]
\ar[r]^-{[1,[1,\Res(A')]]}
&
[\B,[\C,[\D, \A ((\B \C) \D )]]]
\ar[r]^-{[1,[1,[1,1 \otimes A']]]}
&
[\B,[\C,[\D, \A (\B (\C \D) )]]]
}.
$$
}
This last arrow rewrites successively
{\tiny
\begin{tabbing}
\=1.
$\xymatrix@C=3pc{
\A
\ar[r]^-{\eta}
&
[\B \C, \A (\B \C)]
\ar[r]^-{\Res}
&
[\B,[\C, \A (\B \C)]]
\ar[r]^-{[1,[1,\Res(A')]]}
&
[\B,[\C,[\D, \A ((\B \C) \D )]]]
\ar[r]^-{[1,[1,[1,1 \otimes A']]]}
&
[\B,[\C,[\D, \A (\B (\C \D) )]]]
}
$\\
\>2.
$\xymatrix@C=3pc{
\A
\ar[r]^-{\eta}
&
[\B \C, \A (\B \C)]
\ar[r]^-{[1,\Res(A')]}
&
[\B\C, [\D, \A( (\B \C) \D) ] ]
\ar[r]^-{[1,[1, 1 \otimes A']]}
&
[\B \C , [\D, \A ( \B (\C  \D)) ] ]    
\ar[r]^-{\Res}
&
[\B,[\C , [\D, \A ( \B (\C  \D)) ] ]] 
}$
\\
\>3.
$\xymatrix{
\A
\ar[r]^-{\eta}
&
[(\B \C)\D, \A ((\B \C)\D)]
\ar[r]^-{\Res}
&
[\B\C, [\D, \A( (\B \C) \D) ] ]\;\;
\ar[r]^-{[1,[1, 1 \otimes A']]}
&
\;\;
[\B \C , [\D, \A ( \B (\C  \D)) ] ]    
\ar[r]^-{\Res}
&
[\B,[\C , [\D, \A ( \B (\C  \D)) ] ]] 
}$
\\
\>4.$\xymatrix{
\A
\ar[r]^-{\eta}
&
[(\B \C)\D, \A ((\B \C)\D)]
\ar[r]^-{[1,1 \otimes A']}
&
[(\B \C)\D, \A ( \B (\C  \D)) ]     
\ar[r]^-{\Res}
&
[\B \C , [\D, \A ( \B (\C  \D)) ]]
\ar[r]^-{\Res}
&
[\B,[\C , [\D, \A ( \B (\C  \D)) ] ]] 
}
$\\
\>5.$\xymatrix{
\A
\ar[r]^-{\eta}
&
[\B (\C\D), \A (\B (\C\D))]
\ar[r]^-{[A',1]}
&
[(\B \C)\D, \A ( \B (\C  \D)) ]    
\ar[r]^-{\Res}
&
[\B \C , [\D, \A ( \B (\C  \D)) ]]
\ar[r]^-{\Res}
&
[\B,[\C , [\D, \A ( \B (\C  \D)) ] ]] 
}
$\\
\>6.
$
\xymatrix{
\A
\ar[r]^-{\eta}
&
[\B(\C\D), \A (\B (\C  \D))]
\ar[r]^-{[A' ,1 ]}
&
[(\B\C)\D, \A ( \B (\C  \D) ) ]
\ar[r]^-{[\D,-]}
&
[[\D, (\B\C) \D],
[\D,\A ( \B (\C \D)) ]]
\ar[r]^-{[\eta,1]}
&
...
}$\\
\>$\xymatrix{ ...\;
[\B\C,[\D,\A ( \B (\C \D))]]
\ar[r]^-{\Res}
&
[
\B,
[\C, [\D,\A ( \B (\C \D)) ]]
]
}
$\\
\>7.
$
\xymatrix{
\A
\ar[r]^-{\eta}
&
[\B(\C\D), \A (\B (\C  \D))]
\ar[r]^-{[\D , - ]}
&
[[\D,\B(\C\D)],[\D, \A (\B (\C  \D))]]
\ar[r]^-{[[1,A'],1]}
&
[[\D, (\B\C) \D],
[\D,\A ( \B (\C \D)) ]]
\ar[r]^-{[\eta,1]}
&
...
}$\\
\>$\xymatrix{ ...\;
[\B\C,[\D,\A ( \B (\C \D))]]
\ar[r]^-{\Res}
&
[
\B,
[\C, [\D,\A ( \B (\C \D)) ]]
]
}
$\\
\>8.
$
\xymatrix@C=3pc{
\A
\ar[r]^-{\eta}
&
[\B(\C\D), \A (\B (\C  \D))]
\ar[r]^-{[\D , - ]}
&
[[\D,\B(\C\D)],[\D, \A (\B (\C  \D))]]
\ar[r]^-{[\Res(A'),1]}
&
[\B\C,[\D,\A ( \B (\C \D))]]
\ar[r]^-{[\C,-]}
&
...
}$\\
\>$\xymatrix{ ...\;
[[\C,\B\C], [\C, [\D,\A ( \B (\C \D)) ]]]
\ar[r]^-{[\eta,1]}
&
[
\B,
[\C, [\D,\A ( \B (\C \D)) ]]
]
}
$\\
\>9.$\xymatrix{
\A
\ar[r]^-{\eta}
&
[\B(\C\D), \A (\B (\C  \D))]
\ar[r]^-{[\D , - ]}
&
[[\D,\B(\C\D)],[\D, \A (\B (\C  \D))]]
\ar[r]^-{[\C,-]}
&
[[\C,[\D,\B(\C\D)]],[\C,[\D, \A (\B (\C  \D))]]]
\; ...
}$\\
\>$\xymatrix@C=4pc{... 
\ar[r]^-{[\Res (\Res(A')),1]} 
&
[
\B,
[\C, [\D,\A ( \B (\C \D)) ]]
]
}
$\\
\>10. $\xymatrix{
\A
\ar[r]^-{\eta}
&
[\B(\C\D), \A (\B (\C  \D))]
\ar[r]^-{[\D , - ]}
&
[[\D,\B(\C\D)],[\D, \A (\B (\C  \D))]]
\ar[r]^-{[\C,-]}
&
[[\C,[\D,\B(\C\D)]],[\C,[\D, \A (\B (\C  \D))]]]
\; ...
}$\\
\>$\xymatrix{
...
\ar[r]^-{[\Res,1]}
&
[ [\C \D,\B( \C\D) ],[\C,[\D, \A (\B (\C  \D))]]] 
\ar[r]^-{[\eta,1]}
&
[
\B,
[\C, [\D,\A ( \B (\C \D)) ]]
]
.}
$
\end{tabbing}
}
In the above derivation the equalities between
arrows 6. and 7. and between arrows 8. and 9. hold due 
to Corollary \ref{impronatPost2}.
\epf

\begin{corollary}\label{waxiom1}
The diagram in $\SMC$
$$\xymatrix{
\A \otimes ( \B \otimes ( \C \otimes \D ))) 
\ar[d]_{1 \otimes A}
\ar[r]^A
& 
(\A \otimes  \B) \otimes ( \C \otimes \D ) 
\ar[r]^A
& 
(( \A \otimes \B) \otimes \C) \otimes \D\\
\A \otimes ((\B \otimes \C) \otimes \D) \ar[rr]_A && 
( \A \otimes (\B \otimes \C) ) \otimes \D \ar[u]_{A \otimes 1}
}$$
commutes for any $\A, \B, \C$ and $\D$.
\end{corollary}
\pf
Consider the two pastings below in $\SMC$
where all the diagrams involved commute according to
Lemma \ref{AandAprel} and the naturalities of $S$ and $A'$.
$$ 
\xymatrix{
\A(\B(\C\D))
\ar[r]^{1 \otimes A}
\ar[d]_{1 \otimes T^{-1}}
&
\A((\B\C)\D)
\ar[d]_{1 \otimes T}
\ar[rr]^A
\ar[rd]^{T^{-1}}
&
&
(\A(\B\C))\D
\ar[rr]^{A \otimes 1}
\ar[d]_T
\ar[rd]^{T^{-1} \otimes 1}
&
&
((\A\B)\C)\D
\ar[d]^{T \otimes 1}
\\
\A((\D\C)\B)
\ar[r]^{1 \otimes A'}
\ar[d]_{S}
&
\A(\D(\B\C))
\ar[d]_{S}
&
(\D(\B\C))\A
\ar[ld]^{(1 \otimes S) \otimes 1}
\ar[r]^{A'}
&
\D((\B\C)\A)
\ar[d]_{1 \otimes (S \otimes 1)}
&
((\C\B)\A)\D
\ar[ld]^S
\ar[r]^{A' \otimes 1}
&
(\C(\B\A))\D
\ar[d]^S
\\
((\D\C)\B)\A
\ar[r]_{A' \otimes 1}
&
(\D(\C\B))\A
\ar[rr]_{A'}
&
&
\D((\C\B)\A)
\ar[rr]_{1 \otimes A'}
&
&
\D(\C(\B\A))
}$$
$$
\xymatrix{
&
\A(\B(\C\D))
\ar[r]^A
\ar[d]^{T^{-1}}
\ar[ld]_{1 \otimes T^{-1}}
&
(\A\B)(\C\D)
\ar[d]_{T}
\ar[rd]^{T^{-1}}
\ar[rr]^A
&
&
((\A\B)\C)\D
\ar[d]_T
\ar[rd]^{T \otimes 1}
\\
\A((\D \C) \B)
\ar[rd]_S
&
((\C\D)\B)\A
\ar[r]^{A'}
\ar[d]^{(S \otimes 1) \otimes 1}
&
(\C\D)(\B\A)
\ar[rd]_{S \otimes 1}
&
(\D\C)(\A\B)
\ar[d]^{1 \otimes S}
\ar[r]^{A'}
&
\D(\C(\A\B))
\ar[d]_{1 \otimes (1 \otimes S)}
&
(\C(\B \A))\D
\ar[ld]^S
\\
&
((\D\C) \B) \A
\ar[rr]_{A'}
&
&
(\D\C)(\B\A)
\ar[r]_{A'}
&
\D(\C(\B\A))
.}
$$
\epf

\begin{lemma}\label{waxiom42}
The diagram
in $\SMC$
$$\xymatrix{ (\A \otimes \B) \otimes \C \ar[r]^{A'} \ar[d]_{S \otimes 1} &
\A \otimes (\B \otimes \C) \ar[r]^{S} & 
(\B \otimes \C) \otimes \A
\ar[d]^{A'}\\
(\B \otimes \A) \otimes \C \ar[r]_{A'} &
\B \otimes (\A \otimes \C) 
& 
\B \otimes (\C \otimes \A)
\ar[l]^{1 \otimes S} 
}$$
commutes for any $\A$, $\B$ and $\C$.
\end{lemma}
\pf
Both legs of this diagram are strict.
We show that their images by $\Res$ are strict
and that their images by $\Dual \circ \Res \circ \Res$
are equal.\\ 

The 1-cell
$$\xymatrix{  
(\A \otimes \B) \otimes \C
\ar[r]^{S \otimes 1}
&
(\B \otimes \A) \otimes \C
\ar[r]^{A'}
&
\B \otimes (\A \otimes \C)
}$$
has image by $\Res$ 
$$
\xymatrix{
\A \otimes \B
\ar[r]^-{S}
&
\B \otimes \A
\ar[r]^-{\eta}
&
[\C, (\B \otimes \A) \otimes \C]
\ar[r]^-{[1,A']}
&
[\C, \B \otimes (\A \otimes \C)]
}
$$
which is
$$\xymatrix{
\A \otimes \B
\ar[r]^-{S}
&
\B \otimes \A
\ar[r]^-{\Res(A')}
&
[\C, \B \otimes (\A \otimes \C)]
}.$$
This last 1-cell is strict and has
image by $\Res$ the 1-cell
$$
\xymatrix@C=3pc{
\A 
\ar[r]^-{\eta^*}
&
[\B,\B \otimes \A]
\ar[r]^-{[1,\Res(A')]}
&
[\B,[\C,\B \otimes (\A \otimes \C)]]
}
$$
which, by Lemma \ref{dualitynatABC0}, has dual 
$$
\xymatrix@C=3pc{
\B 
\ar[r]^-{\eta}
&
[\A, \B \otimes \A]
\ar[r]^-{[1,\Res(A')]}
&
[\A,[\C,\B \otimes (\A \otimes \C)]]
}
$$
that is $\Res(\Res(A'))$.
On the other hand
the 1-cell
$$
\xymatrix{
(\A \B) \C
\ar[r]^-{A'}
&
\A (\B \C)
\ar[r]^-S
&
(\B \C) \A
\ar[r]^-{A'}
&
\B (\C \A)
\ar[r]^-{1 \otimes S}
&
\B (\A \C)
}
$$
has image by $\Res$
$$
\xymatrix{
\A \B
\ar[r]^-{\Res(A')}
&
[\C, \A (\B \C)]
\ar[r]^-{[1,S]}
&
[\C, (\B \C) \A]
\ar[r]^-{[1,A']}
&
[\C, \B (\C \A)]
\ar[r]^-{[1, 1 \otimes S]}
&
[\C, \B (\A \C)]
}
$$
which is strict. This last 1-cell has image by $\Res$
{\tiny
$$
\xymatrix{
\A
\ar[r]^-{\eta}
&
[\B \C, \A (\B \C)]
\ar[r]^-{\Res}
&
[\B, [\C, \A (\B \C)]]
\ar[r]^-{[1,[1,S]]}
&
[\B, [\C, (\B \C) \A]]
\ar[r]^-{[1,[1,A']]}
&
[\B,[\C, \B (\C \A)]]
\;
\ar[r]^-{[1,[1,1 \otimes S]]}
&
\;
[\B,[\C, \B (\A \C)]]
}
$$
}
which rewrites
{\tiny
\begin{tabbing}
\=$\xymatrix{
\A
\ar[r]^-{\eta}
&
[\B \C , \A (\B \C)]
\ar[r]^-{[1,S]}
&
[\B \C, (\B \C) \A]
\ar[r]^-{[1,A']}
&
[\B \C, \B (\C \A)]
\ar[r]^-{[1,1 \otimes S]}
&
[\B \C, \B (\A \C)]
\ar[r]^-{\Res}
&
[\B,[ \C, \B (\A \C)]]
}$\\
\>$\xymatrix{
\A
\ar[r]^-{\eta^*}
&
[\B \C, (\B \C) \A]
\ar[r]^-{[1,A']}
&
[\B \C, \B (\C \A)]
\ar[r]^-{[1, 1 \otimes S]}
&
[\B \C, \B (\A \C)]
\ar[r]^-{[\C,-]}
&
[[\C, \B \C], [\C, \B (\A \C)]]
\ar[r]^-{[\eta,1]}
&
[\B, [\C, \B (\A \C)]]  
}$
\end{tabbing}
}
According to Lemma \ref{dualitynatABC0} 
and Lemma \ref{beurk}
this last 1-cell has dual
$$
\gamma: \xymatrix{
\B
\ar[r]^-{\eta}
&
[\C , \B \C]
\ar[r]^-{[1,F]}
&
[\C, [\A, \B (\A \C)]]
\ar[r]^-{\Dual}
&
[\A, [\C, \B (\A \C)]]
}
$$ 
where $F: \B \otimes \C \rightarrow [\A, \B \otimes (\A \otimes \C)]$ 
stands for the 1-cell
dual 
of 
$$
\xymatrix{
\A
\ar[r]^-{\eta^*}
&
[\B \C, (\B \C) \A]
\ar[r]^-{[1,A']}
&
[\B \C, \B( \C \A)]
\ar[r]^-{[1,1 \otimes S]}
&
[\B \C, \B (\A \C)].
}
$$
This $F$ is
$$\xymatrix{
\B \C
\ar[r]^-{\eta}
&
[\A, (\B \C) \A]
\ar[r]^-{[1,A']}
&
[\A, \B (\C \A)]
\ar[r]^-{[1,1 \otimes S]}
&
[\A, \B (\A \C)]
}$$ 
and
therefore the above 1-cell $\gamma$ is
{\small
$$
\xymatrix{
\B
\ar[r]^-{\eta}
&
[\C , \B \C]
\ar[r]^-{[1,\eta]}
&
[\C, [\A, (\B \C) \A]]
\ar[r]^-{[1,[1,A']]}
&
[\C, [\A, \B (\C \A)]]\;
\ar[r]^-{[1,[1,1 \otimes S]]}
&
\;
[\C, [\A, \B (\A \C)]]
\ar[r]^-{\Dual}
&
[\A, [\C, \B (\A \C)]]
}
$$
}
which rewrites successively as
\begin{tabbing}
\=1. \xymatrix@C=4pc{
\B
\ar[r]^-{\Res(\Res(A'))}
&
[\C, [\A, \B (\C \A)]]
\ar[r]^-{[1,[1,1 \otimes S]]}
&
[\C, [\A, \B (\A \C)]]
\ar[r]^-{\Dual}
&
[\A, [\C, \B (\A \C)]]
}\\
\>2. \xymatrix{
\B
\ar[r]^-{\eta}
&
[\C \A, \B (\C \A)]
\ar[r]^-{\Res}
&
[\C, [\A, \B (\C \A)]]\;\;
\ar[r]^-{[1,[1,1 \otimes S]]}
&
\;\;
[\C, [\A, \B (\A \C)]]
\ar[r]^-{\Dual}
&
[\A, [\C, \B (\A \C)]]
}\\
\>3. \xymatrix{
\B
\ar[r]^-{\eta}
&
[\C \A, \B (\C \A)]
\ar[r]^-{[1,1 \otimes S]}
&
[\C \A, \B (\A \C)]
\ar[r]^-{\Res}
&
[\C, [\A, \B (\A \C)]]
\ar[r]^-{D}
&
[\A, [\C, \B (\A \C)]]
}\\
\>4. \xymatrix{
\B
\ar[r]^-{\eta}
&
[\A \C, \B (\A \C)]
\ar[r]^-{[S,1]}
&
[\C \A, \B (\A \C)]
\ar[r]^-{\Res}
&
[\C, [\A, \B (\A \C)]]
\ar[r]^-{D}
&
[\A, [\C, \B (\A \C)]]
}
\\
\>5. \xymatrix{
\B
\ar[r]^-{\eta}
&
[\A \C, \B (\A \C)]
\ar[r]^-{\Res}
&
[\A, [\C, \B (\A \C)]].
}
\end{tabbing}
In the above derivation
arrows 4. and 5. are equal due to Lemma \ref{ResSD}.
\epf

\begin{corollary}\label{waxiom4}
The diagram
in $\SMC$ 
$$\xymatrix{ \A \otimes (\B \otimes \C) \ar[r]^A \ar[d]_{1 \otimes S} &
(\A \otimes \B) \otimes \C \ar[r]^{S} & 
\C \otimes ( \A \otimes \B )
\ar[d]^{A}\\
\A \otimes (\C \otimes \B) \ar[r]_A &
(\A \otimes \C) \otimes \B 
& 
(\C \otimes \A) \otimes \B
\ar[l]^{S \otimes 1} 
}$$
commutes
for any $\A$, $\B$ and $\C$.
\end{corollary}
\pf
Consider the pasting of commutative diagrams in $\SMC$
$$
\xymatrix{
\A(\B\C)
\ar[rr]^A
\ar[ddd]_{1 \otimes S}
\ar[rd]_{T^{-1}}
& 
& 
(\A\B)\C
\ar[d]^{T}
\ar[rr]^S
&
&
\C(\A\B)
\ar[ddd]^A
\ar[ld]^{T^{-1}}
\\
&
(\C\B)\A 
\ar[d]_{S \otimes 1}
\ar[r]^{A'}
&
\C(\B \A)
\ar[r]^S
&
(\B\A)\C
\ar[d]^{A'}
&
\\
&
(\B\C)\A
\ar[r]_{A'}
&
\B(\C\A)
&
\B(\A\C)
\ar[l]^{1 \otimes S}
&
\\
\A(\C\B)
\ar[rr]_{A}
\ar[ru]^{T^{-1}}
& & 
(\A \C) \B
\ar[u]^T
& & 
(\C \A) \B.
\ar[ll]^{S \otimes 1}
\ar[lu]_T
}
$$
\epf

\begin{lemma}\label{waxiom5}
The diagram in $\SMC$
$$\xymatrix{\A \otimes \unc \ar[rd]_{R} \ar[rr]^{S} & & 
\unc \otimes \A \ar[ld]^{L}  \\  
& \A & }$$
commutes for any $\A$.
\end{lemma}
\pf
By definition of $L$ and $S$, the functor $L \circ S$ is strict 
and the two diagrams involved in the 
pasting in $\SMC$ below commute
\begin{center}
$\xymatrix{
\A 
\ar[r]^{\eta}
\ar@{=}[d]
&
[\unc, \A \otimes \unc]
\ar[d]^{[\unc,S]}
\\
\A
\ar[r]^{\eta^*}
\ar[rd]_{\vs }
&
[\unc, \unc \otimes \A]
\ar[d]^{[\unc,L]}
\\
&
[\unc,\A]
.}
$
\end{center}
This shows that $L \circ S$ is the image by $\Ext$ of $\vs$.
By definition, $R$ is also that image.
\epf

\begin{lemma}\label{waxiom22}
The diagram in $\SMC$
$$
\xymatrix{
(\A \otimes \unc) \otimes \C
\ar[rr]^{A'}
& 
& 
\A \otimes (\unc \otimes \C)
\\
&
\A \otimes \C
\ar[lu]^{R' \otimes 1}
\ar[ru]_{1 \otimes L'} 
&
}
$$
commutes for all $\A$ and $\C$.
\end{lemma}
\pf
Both legs of the diagram are strict, we show 
that their images by $\Res$ are equal.\\

The image by $\Res$ of $A' \circ (R' \otimes 1)$ 
is
{\tiny
\begin{tabbing}
\=1. 
\=$ 
\xymatrix{
\A 
\ar[r]^-{R'}
&
\A \otimes \unc
\ar[r]^-{\eta} 
&
[\C, (\A \otimes \unc) \otimes \C]
\ar[r]^-{[1,A']}
&
[\C, \A \otimes (\unc \otimes \C)]
}
$
\\
\>2.
\>
$
\xymatrix{
\A 
\ar[r]^-{\eta}
&
[\unc, \A \otimes \unc]
\ar[r]^-{\ev_{\gen}} 
&
\A \otimes \unc
\ar[r]^-{\Res(A')}
&
[\C, \A \otimes (\unc \otimes \C)]
}
$
\\
\>3.
\>$
\xymatrix@C=3pc{
\A 
\ar[r]^-{\eta}
&
[\unc, \A \otimes \unc]
\ar[r]^-{[1,\Res(A')]} 
&
[\unc,[\C,\A \otimes (\unc \otimes \C)]]
\ar[r]^-{\ev_{\gen}}
&
[\C, \A \otimes (\unc \otimes \C)]
}
$
\\
\>4.
\>
$
\xymatrix{
\A 
\ar[r]^-{\eta}
&
[\unc \otimes \C, \A \otimes (\unc \otimes \C)]
\ar[r]^-{\Res} 
&
[\unc,[\C,\A \otimes (\unc \otimes \C)]]
\ar[r]^-{\ev_{\gen}}
&
[\C, \A \otimes (\unc \otimes \C)]
}
$
\\
\>5.
\>
$
\xymatrix{
\A 
\ar[r]^-{\eta}
&
[\unc \otimes \C, \A \otimes (\unc \otimes \C)]
\ar[r]^-{[\C,-]} 
&
[[\C,\unc \otimes \C],[\C,\A \otimes (\unc \otimes \C)]]
\ar[r]^-{[\eta,1]}
&
[\unc,[\C,\A \otimes (\unc \otimes \C)]]
\ar[r]^-{\Dual}
&
...
}$\\
\>
\>
$
\xymatrix{
...
\;
[\C,[\unc,\A \otimes (\unc \otimes \C)]]
\ar[r]^-{[1,\ev_{\gen}]}
&
[\C, \A \otimes (\unc \otimes \C)]
}
$
\\
\>6.
\>
$
\xymatrix{
\A 
\ar[r]^-{\eta}
&
[\unc \otimes \C, \A \otimes (\unc \otimes \C)]
\ar[r]^-{[\unc,-]} 
&
[[\unc, \unc \otimes \C],[\unc, \A \otimes (\unc \otimes \C)]]
\ar[r]^-{[\eta^*,1]}
&
[\C,[\unc,\A \otimes (\unc \otimes \C)]]
\ar[r]^-{[1,\ev_{\gen}]}
&
[\C, \A \otimes (\unc \otimes \C)]
}
$
\\
\>7.
\>
$
\xymatrix{
\A 
\ar[r]^-{\eta}
&
[\unc \otimes \C, \A \otimes (\unc \otimes \C)]
\ar[r]^-{[\unc,-]} 
&
[[\unc, \unc \otimes \C],[\unc, \A \otimes (\unc \otimes \C)]]
\ar[r]^-{[1,\ev_{\gen}]}
&
[ [\unc, \unc \otimes \C], \A \otimes (\unc \otimes \C)]
\ar[r]^-{[\eta^*,1]}
&
[\C, \A \otimes (\unc \otimes \C)]
}
$
\\
\>8.
\>
$
\xymatrix{
\A 
\ar[r]^-{\eta}
&
[\unc \otimes \C, \A \otimes (\unc \otimes \C)]
\ar[r]^-{[\ev_{\gen},1]} 
&
[ [\unc, \unc \otimes \C], \A \otimes (\unc \otimes \C)]
\ar[r]^-{[[1,S],1]}
&
[[\unc, \C \otimes \unc], \A \otimes (\unc \otimes \C)]
\ar[r]^-{[\eta,1]}
&
[\C , \A \otimes (\unc \otimes \C)]
}
$
\\
\>9.
\>
$
\xymatrix{
\A 
\ar[r]^-{\eta}
&
[\unc \otimes \C, \A \otimes (\unc \otimes \C)]
\ar[r]^-{[S,1]} 
&
[ \C \otimes \unc, \A \otimes (\unc \otimes \C)]
\ar[r]^-{[\ev_{\gen},1]}
&
[[\unc, \C \otimes \unc], \A \otimes (\unc \otimes \C)]
\ar[r]^-{[\eta,1]}
&
[\C , \A \otimes (\unc \otimes \C)]
}
$
\\
\>10.
\>
$
\xymatrix{
\A 
\ar[r]^-{\eta}
&
[\unc \otimes \C, \A \otimes (\unc \otimes \C)]
\ar[r]^-{[S,1]} 
&
[ \C \otimes \unc, \A \otimes (\unc \otimes \C)]
\ar[r]^-{[R',1]}
&
[\C, \A \otimes (\unc \otimes \C)]
}
$
\\
\>11.
\>
$
\xymatrix{
\A 
\ar[r]^-{\eta}
&
[\unc \otimes \C, \A \otimes (\unc \otimes \C)]
\ar[r]^-{[L',1]}
&
[\C, \A \otimes (\unc \otimes \C)]
}
$
\end{tabbing}
}
and this last arrow is the image of $L'$ by $\Res$.\\

In the above derivation:\\
- arrows 2 and 3 are equal by Corollary \ref{2nateva},\\
- arrows 4 and 5 are equal by Lemma \ref{evdual},\\
- arrows 5 and 6 are equal by Lemma \ref{tbc3},\\
- arrows 7 and 8 are equal because $\Res(S) = \eta^*$
 and by Corollary \ref{homev3},\\
- arrows 8 and 9 are equal by Corollary \ref{2nateva}.
\epf

\begin{corollary}
The diagram in $\SMC$
$$
\xymatrix{
\A \otimes (\unc \otimes \C)
\ar[rr]^{A}
& 
& 
(\A \otimes \unc) \otimes \C
\\
&
\A \otimes \C
\ar[lu]^{1 \otimes L'}
\ar[ru]_{R' \otimes 1} 
&
}
$$
commutes for all $\A$ and $\C$.
\end{corollary}
\pf
Consider the pasting in $\SMC$ where all
diagrams commute.
$$
\xymatrix{
\A(\unc \C) 
\ar[rrrr]^A
\ar[rd]^{T^{-1}}
\ar[rdd]|{1 \otimes S}
& 
& 
& 
&
(\A \unc) \C
\ar[ld]_T
\ar[ldd]|{S \otimes 1}
\\
&
(\C \unc) \A
\ar[rr]^{A'}
\ar[d]|S
& 
&
\C (\unc \A)
\ar[d]|S\\
& 
\A (\C \unc)
& 
\C \A
\ar[d]|S
\ar[lu]|{R' \otimes 1}
\ar[ru]|{1 \otimes L'}
&
(\unc \A) \C
\\
& & \A \C
\ar[lu]|{1 \otimes R'}
\ar[ru]|{L' \otimes 1}
\ar@/^40pt/[lluuu]^{1 \otimes L'}
\ar@/_40pt/[rruuu]_{R' \otimes 1}
}
$$
\epf

For any  
$F: \A \rightarrow \B$ the following diagram
in $\SMC$ 
$$
\xymatrix{
[\unc,\A] 
\ar[r]^{[1,F]}
\ar[d]_{\ev_{\gen}}
&
[\unc,\B]
\ar[d]^{\ev_{\gen}}
\\
\A
\ar[r]_F
&
\B
}
$$
commutes according to Corollary \ref{2nateva},
so that by considering the mate of this 
identity 2-cell one obtains a 2-cell, namely
$\delta_{F}$ as follows
\begin{tag}\label{2cdelta}
$$
\xymatrix{
\A 
\ar[r]^{F}
\ar[d]_{\vs}
&
\B
\ar[d]^{\vs}
\ar@{=>}[ld]
\\
[\unc,\A]
\ar[r]_{[1,F]}
&
[\unc,\B].
}
$$
\end{tag}

\begin{lemma}\label{Rinverse}
In $\SMC$, one has the following for any $\A$.\\
- The composite 
$\xymatrix{
\A \ar[r]^-{R'}
&
\A \otimes \unc
\ar[r]^-{R}
&
\A
}$ is the identity 
at $\A$.\\
- There exists a 1-cell 
$F : \A \otimes \unc \rightarrow \A \otimes \unc$ with 
two 2-cells
$F \rightarrow 1: \A \otimes \unc \rightarrow \A \otimes \unc$ 
and    
$F \rightarrow R' \circ R$.
\end{lemma}
\pf
To check that the composite $R \circ R'$ is the identity
at $\A$ in $\SMC$, consider the pasting in $\SMC$ of 
commutative diagrams
$$
\xymatrix{
\A
\ar@/^35pt/[rr]^{\vs}
\ar[rd]_{R'}
\ar[r]^{\eta}
&
[\unc, \A \otimes \unc]
\ar[d]^{\ev_{\gen}}
\ar[r]_{[1,R]}
&
[\unc,\A]
\ar[d]^{\ev_{\gen}}
\\
&
\A \otimes \unc
\ar[r]_R
&
\A
}
$$
the bottom right diagram here commuting 
according to Corollary \ref{2nateva}. 
One concludes since $\ev_{\gen} \circ \vs = 1$
as shown in Proposition \ref{unimonadj}.\\

We prove now the existence of
a 1-cell $F: \A \otimes \unc \rightarrow \A \otimes \unc$ in $\SMC$
with two 2-cells
$F \rightarrow 1$ and $F \rightarrow  R' \circ R$.
For this, we exhibit a 1-cell
$G: \A \rightarrow [\unc, \A \otimes \unc]$ 
with two 2-cells
one from $G$ to $\eta$ and the other
one from $G$ to the composite
$$\xymatrix{
\A
\ar[r]^-{\vs}
&
[\unc,\A]
\ar[r]^-{[1,R']}
&
[\unc, \A \otimes \unc]
}$$
which is actually the image by $\Res$ of $R' \circ R$.
To see that this is sufficient, suppose that such 2-cells exist.
Then let  $F:\A \otimes \unc \rightarrow \A \otimes \unc$
be the image by $\Ext$ of $G$.
Since the image by $\Ext$ of $\eta$ is the identity 
at $\A \otimes \unc$ one obtains 
a 2-cell $F \rightarrow 1_{\A \otimes \unc}$. 
The 2-cell $G \rightarrow \Res(R' \circ R)$  
corresponds to a 2-cell $F \rightarrow R' \circ R$
via the adjunction $\Ext \dashv \Res$.\\

The 1-cell $G$ in question is 
$$\xymatrix{
\A 
\ar[r]^-{R'}
&
\A \otimes \unc
\ar[r]^-{\vs}
&
[\unc, \A \otimes \unc]
}
$$
and the 2-cells are
$$
\xymatrix{ 
&
&
&
\\
\A
\ar[r]^{\eta}
\ar@/_35pt/[rr]_{R'}
&
[\unc,\A \otimes \unc]
\ar[r]^{\ev_{\gen}}
\ar@{}[d]|{=}
\ar@/^35pt/[rr]^1
&
\A \otimes \unc
\ar[r]^{\vs}
\ar@{=>}[u]_{\epsilon}
&
[\unc, \A \otimes \unc]
\\
& & & 
}
$$ 
where $\epsilon$ is the co-unit of the adjunction
of Proposition \ref{unimonadj}
and the $\delta_{R'}$
$$
\xymatrix{
\A 
\ar[d]_{\vs}
\ar[r]^{R'}
&
\A \otimes \unc
\ar[d]^{\vs}
\ar@{=>}[ld]
\\
[\unc, \A]
\ar[r]_{[1,R']}
&
[\unc, \A \otimes \unc]
}
$$
defined in \ref{2cdelta}.
\epf

\end{section}

\begin{section}{A symmetric monoidal closed structure on $\SMC_{ / \sim}$}
\label{closedquo}
This section contains a proof of the following result:
the category $\SMC_{/\sim}$,
quotient of $\SMC$ by the congruence generated by its 2-cells
admits a symmetric monoidal closed structure.\\

Given any small category $\C$, consider the equivalence
on its set of objects generated by the relation consisting of 
the pairs $(x,y)$ such that there exists an arrow $x \rightarrow y$ 
in $\C$.
Its classes are the so-called {\em connected components} 
of $\C$ and we write $[x]$ for the connected component 
of the object $x$. The set of connected components of $\C$
is denoted $\pi(\C)$.
Given any functor $F: \C \rightarrow \C'$, one obtains
a map $\pi(F) : \pi(\C) \rightarrow \pi(\C')$
sending any connected component $[x]$ to $[F(x)]$.
Note then that:\\
- For any functors  $F:\A \rightarrow \B$ and $G: \B \rightarrow \C$
one has
$\pi (G \circ F) = \pi(G) \circ \pi(F)$;\\
- For any category $\A$, the map $\pi(1_{\A})$ for the identity
functor $1_{\A}$ at $\A$ is the identity map at $\pi(\A)$;\\
- For any natural transformation $\sigma: F \rightarrow G: \A
\rightarrow \B$, one has
$\pi(F) = \pi(G): \pi(\A) \rightarrow \pi(\B)$.\\

For convenience we shall consider further 
mere categories as locally discrete 2-categories
and mere functors as 2-functors.
With this convention, according to the 
remark above, the assignments $\pi$ above define 
a 2-functor $\pi: \CAT \rightarrow \Set$.\\

Given a 2-category $\bA$, consider
the equivalence $\sim$ on its 1-cells generated
by the relation consisting of the pairs $(f, g)$ 
with same domains and codomains
and such that there exists a 2-cell $f \rightarrow g$. 
This is equivalence is compatible with the
composition of 1-cells. We write $\bA_{/\sim}$ for 
the category with the same objects as $\bA$, with arrows
$x \rightarrow y$, the equivalence classes $f^{\sim}$
by $\sim$ of 1-cells $f:x \rightarrow y$ in $\bA$,
and with identities and composition induced
by those of $\bA$:\\
- For any object $x$ of $\bA$, the identity at $x$ in $\bA / \sim$
is the $\sim$-class of the identity at $x$ in $\bA$;\\
- For any arrows $\xymatrix{\ar[r]^f & \ar[r]^g &}$ in $\bA$,
the composite $g^{\sim} \circ f^{\sim}$ in $\bA/\sim$ is
the $\sim$-class of $g \circ f$.\\  
One has a 
2-functor $p_{\bA}: \bA \rightarrow \bA_{/\sim}$,
that is the identity on the set of objects, sends
any 1-cell $f:x \rightarrow y$ to its equivalence class $f^{\sim}$
by $\sim$ and any 2-cell to an identity.
For any 2-functor $F: \bA \rightarrow \bB$, there exists
a unique functor $F_{/\sim}: \bA_{/\sim} \rightarrow \bB_{/\sim}$
that renders commutative the diagram of 2-functors
$$
\xymatrix{
\bA
\ar[r]^F
\ar[d]_{p_{\bA}}
&
\bB
\ar[d]^{p_{\bB}}
\\
\bA_{/\sim}
\ar[r]_{F_{/\sim}}
&
\bB_{/\sim}.
}
$$

\begin{remark}\label{piquo}
For any 2-category $\bA$, the following assertions hold.\\
(i) For any objects $x,y$ of $\bA$,
there is an isomorphism between sets
$$\bA_{/\sim}(x,y)  \cong \pi( \bA (x,y) )$$
(sending $f^{\sim}$ to $[f]$, for any 1-cell $f:x \rightarrow y$).\\   
(ii) For any 1-cell $f:x' \rightarrow x$ in $\bA$,
the following diagram of maps commutes
$$
\xymatrix@C=3pc{
\pi(\bA(x,y))
\ar[r]^{\pi(\bA(f,y))}
\ar@{-}[d]_{\cong}
&
\pi(\bA(x',y))
\ar@{-}[d]^{\cong}
\\
\bA_{/\sim}(x,y)
\ar[r]_{\bA_{/\sim}(f^{\sim},y)}
&
\bA_{/\sim}(x',y).
}
$$
(iii) For any 1-cell $g:y \rightarrow y'$ in $\bA$,
the following diagram of maps commutes
$$
\xymatrix@C=3pc{
\pi(\bA(x,y))
\ar[r]^{\pi(\bA(x,g))}
\ar@{-}[d]_{\cong}
&
\pi(\bA(x,y'))
\ar@{-}[d]^{\cong}
\\
\bA_{/\sim}(x,y)
\ar[r]_{\bA_{/\sim}(x,g^{\sim})}
&
\bA_{/\sim}(x,y').
}
$$
\end{remark}

Note that for any two arbitrary 2-categories 
$A$ and $B$, the category 
 $(\bA_{/\sim}) \times (\bB_{/\sim})$
is isomorphic  to ${(\bA \times \bB)}_{/ \sim}$.
Therefore from any 2-functor 
$H: \bA \times \bB \rightarrow \bC$ one obtains
a functor  
$$\cH: 
\xymatrix{ 
(\bA_{/\sim}) \times (\bB_{/\sim})
\ar@{-}[r]^-{\cong}
&
{(\bA \times \bB)}_{/ \sim}
\ar[r]^-{H_{/\sim}}
&
\bC_{/\sim}}.$$
which is more concretely the following:\\
- For any objects $x$ of $\bA$ and $y$ of $\bB$,
$\cH(x,y)$ is $H(x,y)$;\\
- For any 1-cells $f: x \rightarrow x'$ in $\bA$
and $g: y \rightarrow y'$ in $\bB$, $\cH(f^{\sim},g^{\sim})$ 
is $H(f,g)^{\sim}$.\\
In particular for the 2-functor 
$\tenSMC : \SMC \times \SMC \rightarrow \SMC$
defined in Section \ref{tens2fun}, one obtains
the functor 
$\check{\tenSMC}: \SMC_{/\sim} \times 
\SMC_{/\sim} \rightarrow \SMC_{/\sim}$ which we denote 
by the 2-ary symbol $\cotimes$.\\ 

We can now formulate our result.
\begin{theorem}\label{mainresult}
The category
$\SMC_{/\sim}$ admits the symmetric monoidal closed structure
$$(\SMC_{/\sim}, \cotimes,
\unc, A^{\sim}, R^{\sim}, L^{\sim}, S^{\sim}).$$
For any symmetric monoidal category $\B$,
the right adjoint to 
$- \cotimes \B: \SMC_{/\sim} \rightarrow \SMC_{/\sim}$ 
sends any symmetric monoidal category $\C$ to 
$[\B,\C]$ and
for any symmetric monoidal functor 
$F: \C \rightarrow \C'$, sends $F^{\sim}$
to ${[\B,F]}^{\sim}:[\B,\C] \rightarrow [\B,\C']$. 
\end{theorem} 
\pf
According to the 2-functoriality of $p_{\SMC}$, one 
obtains the following results.
The arrows $A^{\sim}$ are 
isomorphisms with inverses
the ${A'}^{\sim}$
according to Lemma \ref{Ainverse2} 
and Corollary \ref{Ainverse}.
The $R^{\sim}$ are isomorphisms with inverses the ${R'}^{\sim}$
according to Lemma \ref{Rinverse}.
The collection of arrows ${(A_{\A,\B,\C})}^{\sim}$ (respectively 
${({A'}_{\A,\B,\C})}^{\sim}$) is 
natural in $\A$, $\B$ and $\C$ since the collection
of $A_{\A,\B,\C}$ (respectively ${A'}_{\A,\B,\C}$) is natural
in $\A$, $\B$ and $\C$.
The collections 
${(R'_{\A})}^{\sim}$ and  
${(R_{\A})}^{\sim}$  
are natural in $\A$
according to the naturality of the collection 
of $R'_{\A}$.   
The collection ${(L'_{\A})}^{\sim}$
is natural in $\A$
according to the naturality of the collection 
of $L'_{\A}$. 
The arrows $L^{\sim}$ have inverses the $L'^{\sim}$
since the $R^{\sim}$ and ${R'}^{\sim}$ are inverses and
according to Lemmas \ref{SRpLp}, \ref{waxiom5} and \ref{waxiom3}.
Therefore the collection 
${(L_{\A})}^{\sim}$  
is also natural in $\A$.
The collection of ${(S_{\A,\B})}^{\sim}$ is natural
in $\A$ and $\B$ since the collection $S_{\A,\B}$ is.
Then the coherence axioms 
\ref{mcaxiom1}, \ref{mcaxiom2}, \ref{smcaxiom3}, 
\ref{smcaxiom4} and \ref{smcaxiom5} hold for
$\ac = A^{\sim}$, $\rc = R^{\sim}$, $\lc = L^{\sim}$
and $\syc = S^{\sim}$ according respectively to 
the points \ref{waxiom1}, \ref{waxiom22}, 
\ref{waxiom3}, \ref{waxiom4} and \ref{waxiom5}.
\\

Let us consider any symmetric monoidal categories
$\A$, $\B$ and $\C$.
According to Proposition \ref{MAdjResExt}, 
one has 2-cells in $\CAT$
$$1 \rightarrow \Res \circ \Ext:  
\SMC(\A,[\B,\C]) \rightarrow \SMC(\A,[\B,\C])$$
and 
$$\Ext \circ \Res \rightarrow 1: \SMC(\A \otimes \B, \C) 
\rightarrow \SMC(\A \otimes \B,\C).$$
Their images by $\pi: \CAT \rightarrow \Set$
being identities, one obtains
according to Remark \ref{piquo}(i) the isomorphism in $\Set$
\begin{tag}\label{closedadj}
$$
\xymatrix{
\SMC_{/\sim}(\A @ \B,\C) 
\ar@{}[r]|{\cong}
&
\pi ( \SMC(\A \otimes \B,\C) )
\ar[r]^{\pi(\Res)}
&
\pi(\SMC(\A,[\B,\C]))
\ar@{}[r]|{\cong}
&
\SMC_{/\sim}(\A,[\B,\C])
.}
$$
\end{tag}
 
This collection of isomorphisms is actually natural in $\A$.
To see this, consider
any symmetric monoidal functor $F: \A' \rightarrow \A$,
and any symmetric monoidal categories $\B$ and $\C$.
According to Lemma \ref{ResnatinA}, the following 
diagram in $\CAT$ is commutative
$$
\xymatrix{
\SMC(\A \otimes \B,\C)
\ar[r]^-{\Res}
\ar[d]_{\SMC(F \otimes 1, 1)}
&
\SMC(\A, [\B,\C])
\ar[d]^{\SMC(F,1)}
\\
\SMC(\A' \otimes \B, \C)
\ar[r]_-{\Res}
&
\SMC(\A',[\B,\C]).
}
$$
By definition of $\cotimes$ and according to Remark
\ref{piquo} (i) and (ii), applying $\pi$ to the
above yields the diagram of maps
$$
\xymatrix{
\SMC_{/\sim}( \A @ \B, \C) 
\ar@{-}[r]^{\cong}
\ar[d]_{\SMC_{/\sim} (F^{\sim} \cotimes 1, 1 ) }
&
\SMC_{/\sim}( \A,[\B,\C]) 
\ar[d]^{ \SMC_{/\sim}( F^{\sim},1 ) }
\\
\SMC_{/\sim}( \A' @ \B, \C) 
\ar@{-}[r]_{\cong}
&
\SMC_{/\sim}( \A',[\B,\C]). 
}
$$

For any symmetric monoidal category $\B$ there exists a unique
functor $\SMC_{/\sim} \rightarrow \SMC_{/\sim}$ taking 
values $[\B,\C]$ for any symmetric monoidal category $\C$ 
and that renders 
the collection of isomorphisms \ref{closedadj} also natural in $\C$.  
For any arrow $F: \C \rightarrow \C'$, 
this functor sends the class $F^{\sim}$ to 
the image by $\pi(\Res)$ of
$$\xymatrix{
[\B,\C] \otimes \B
\ar[r]^-{\Eval^{\sim}}
&
\C
\ar[r]^-{F^{\sim}}
&
\C'
}
$$
and this image is ${[\B,F]}^{\sim}$ according to Lemma 
\ref{ineedit}.
\epf

\end{section}

\end{document}